\documentclass[twoside, 10pt]{book}

\usepackage{amsmath}
\usepackage{amssymb}
\usepackage{a4}

\input epsf

\makeatletter
\newcounter{mijnlijst}
\newenvironment{lijst}%
	{\begin{list}
		{
	         \hfill
		 (\roman{mijnlijst})
		 }
		{
                \parskip0cm
		\topsep0cm
                \leftmargin0.0cm
		\rightmargin0cm
                \itemindent1.0cm
                \labelwidth1.0cm
                \itemsep0cm
                \parsep0cm
		\usecounter{mijnlijst}
		}
	}%
	{
	\end{list}}
\makeatother
\newcounter{llijst}
\newenvironment{looplijst}%
        {\begin{list}
                {
                 \arabic{llijst}.
                 \hfill
                 }
                {
                \topsep0cm
                \leftmargin.8cm
                \rightmargin0cm
                \itemsep0cm
                \parsep0cm
                 \partopsep0cm
                \usecounter{llijst}
                }
        }%
        {\end{list}}
\newcommand{\fr}[1]{\mbox{$\frac{1}{ #1}$}}  

\newcommand{\beq}{\begin{eqnarray*}}
\newcommand{\eeq}{\end{eqnarray*}}
\newcommand{\bmap}{\begin{displaymath}\begin{array}{rccl}}
\newcommand{\emap}{\end{array}\end{displaymath}}
\newcommand{\bib}[1]{\bibitem[#1]{#1}}
\newcounter{proof}
\newcommand{\eis}{\mbox{$\quad = \quad$}\bigskip}    
\newcommand{\N}{\mathbb N}			
\newcommand{\Z}{\mathbb Z}			
\newcommand{\R}{\mathbb R}			
\newcommand{\C}{\mathbb C}			
\renewcommand{\P}{\mathbb P}			
\newcommand{\U}{\mathcal U}			
\newcommand{\desda}{\mbox{$\Leftrightarrow$}}   
\newcommand{\return}{{\bf return}}

\newcommand{\step}{\refstepcounter{proof}{\bf Step \theproof.~}}

\usepackage{amsthm}
\usepackage{makeidx}
\usepackage[noend]{algorithmic}
\usepackage{epic}

\newcommand{\elabel}[1]{\label{#1}}

\newcommand{\conf}{{\mbox {\sl {\bf conf}}_n}}
\newcommand{\confp}{{\mbox {\sl {\bf conf}}_{n'}}}
\newcommand{\conft}{{\mbox {\sl {\bf conf}}_2}}
\newcommand{\confd}{{\mbox {\sl {\bf conf}}_3}}
\newcommand{\clconf}{{\mbox {\sl {\bf clconf}}_n}}
\newcommand{\cda}{{\mbox {\sl CDA}_n}}
\newcommand{\cdat}{{\mbox {\sl CDA}_3}}
\newcommand{\cuat}{{\mbox {\sl CUA}_3}}
\newcommand{\cuaf}{{\mbox {\sl CUA}_4}}
\newcommand{\cdau}{{\mbox {\sl CDA}_n^{\U}}}
\newcommand{\cua}{{\mbox {\sl CUA}_n}}
\newcommand{\FMd}{{\mbox {\sl FM}_d}}
\newcommand{\FMt}{{\mbox {\sl FM}_2}}
\newcommand{\CONF}{{\mbox {\sl CONF}_n}}
\newcommand{\CLCONF}{{\mbox {\sl CLCONF}_n}}
\newcommand{\CONFt}{{\mbox {\sl CONF}_3}}

\newcommand{\XAH}{{\mbox {\sl {\bf XAH}}[n]}}
\renewcommand{\XAH}{{\mbox {\sl XAH}_n}}
\newcommand{\xah}{{\mbox {\sl XAH}_n}}
\newcommand{\xaht}{{\mbox {\sl XAH}_3}}
\newcommand{\xedah}{{\mbox {\sl XEDAH}_n}}

\newcommand{\DA}{{\mbox {\sl DA}_n}}
\newcommand{\DAt}{{\mbox {\sl DA}_3}}
\newcommand{\UA}{{\mbox {\sl UA}_n}}
\newcommand{\pda}{\psi_{{\text{DA}}_n}}
\newcommand{\pdat}{\psi_{{\text{DA}}_3}}
\newcommand{\fpdat}{\pdat(\CONFt)/\sim_{\{p,r\}}}
\newcommand{\pua}{\psi_{{\text{UA}}_n}}

\newcommand{\draw}{{\mbox {\sl {\bf draw}}}}
\newcommand{\readoff}{{\mbox {\sl {\bf read}}}}
\newcommand{\scr}{{\mbox {\sl {\bf Scr}}_n}}
\newcommand{\dom}{{\mbox {\sl {\bf Dom}}_n}}
\newcommand{\domt}{{\mbox {\sl {\bf Dom}}_3}}
\newcommand{\rng}{{\mbox {\sl {\bf Rng}}_n}}
\newcommand{\rot}{{\mbox {\sl {\bf rot}}}}
\newcommand{\ah}{{\mbox {\sl AH}_n}}
\newcommand{\edah}{{\mbox {\sl EDAH}_n}}

\makeindex

\newtheoremstyle{break}
  {9pt}
  {9pt}
  {\itshape}
  {}
  {\bfseries}
  {}
  {\newline}
  {}

\newtheoremstyle{breakrm}
  {9pt}
  {9pt}
  {}
  {}
  {\bfseries}
  {}
  {\newline}
  {}

\makeatletter
\if@twoside
  \def\ps@headings{%
      \let\@oddfoot\@empty\let\@evenfoot\@empty
      \def\@evenhead{\thepage\hfil\slshape\leftmark}%
      \def\@oddhead{{\slshape\rightmark}\hfil\thepage}%
      \let\@mkboth\markboth
    \def\chaptermark##1{%
      \markboth {%
        \ifnum \c@secnumdepth >\m@ne
            \@chapapp\ \thechapter. \ %
        \fi
        ##1}{}}%
    \def\sectionmark##1{%
      \markright {%
        \ifnum \c@secnumdepth >\z@
          \thesection. \ %
        \fi
        ##1}}}
\else
  \def\ps@headings{%
    \let\@oddfoot\@empty
    \def\@oddhead{{\slshape\rightmark}\hfil\thepage}%
    \let\@mkboth\markboth
    \def\chaptermark##1{%
      \markright {\MakeUppercase{%
        \ifnum \c@secnumdepth >\m@ne
            \@chapapp\ \thechapter. \ %
        \fi
        ##1}}}}
\fi
\makeatother

\theoremstyle{break}
\newtheorem{claim}{Claim}[chapter]

\newtheorem{definition}[claim]{Definition}
\newtheorem{lemma}[claim]{Lemma}

\newtheorem{property}[claim]{Property}
\newtheorem{proposition}[claim]{Proposition}
\newtheorem{theorem}[claim]{Theorem}

\newtheorem{corollary}[claim]{Corollary}

\theoremstyle{breakrm}
\newtheorem{remark}[claim]{Remark}
\newtheorem{eexample}[claim]{Example}
\newtheorem{notation}[claim]{Notation}
\newtheorem{alg}[claim]{Algorithm}

\newcommand{\construction}{\setcounter{proof}{0}\refstepcounter{claim}%
	 {\bf Construction \theclaim}}

\newcommand{\bfindex}[1]{{\bf #1}\index{#1}}

\newcommand{\be}[2]{\beta^{\mathit{#1}\,}_{\mathit{#2}\,}}
\newcommand{\abe}[2]{|\beta^{\mathit{#1}\,}_{\mathit{#2}\,}|}
\newcommand{\al}[2]{\alpha^{\mathit{#1}\,}_{\mathit{#2}\,}}
\newcommand{\ho}[2]{h^{\mathit{#1}\,}_{\mathit{#2}\,}}
\newcommand{\K}[2]{K^{\mathit{#1}\,}_{\mathit{#2}\,}}

\renewcommand{\be}[2]{\beta^{#1\,}_{#2\,}}
\renewcommand{\abe}[2]{|\beta^{#1\,}_{#2\,}|}
\renewcommand{\al}[2]{\alpha^{#1\,}_{#2\,}}
\renewcommand{\ho}[2]{h^{#1\,}_{#2\,}}
\renewcommand{\K}[2]{K^{#1\,}_{#2\,}}

\begin{document}

\begin{titlepage}
\begin{center}
{\huge \bf Limits of Voronoi Diagrams}
\end{center}

\vfill

\begin{center}
{\large Proefschrift}
\end{center}

\begin{center}
\parbox{11cm}{ter verkrijging van de graad van doctor aan de Universiteit
Utrecht op gezag van de Rector Magnificus, Prof.~dr.~W.H.~Gispen,
ingevolge het besluit van het College voor Promoties in het
openbaar te verdedigen op donderdag 10 oktober 2002 des ochtends te 10.30 uur}
\end{center}
\vspace{4cm}
\begin{center}
door
\end{center}

\begin{center}
{\large Roderik Cornelis Lindenbergh}
\end{center}

\begin{center}
geboren op 12 juni 1970, te Aerdenhout
\end{center}
\vspace{2cm}
\newpage
\thispagestyle{empty}
\noindent
\begin{tabular}{ll}
promotor: &Prof.~Dr.~D.~Siersma\\
copromotor: &Dr.~W.~L.~J.~van~der~Kallen\\
\end{tabular}

\medskip
\noindent
\begin{tabular}{ll}
 Faculteit der Wiskunde en  Informatica\\
  Universiteit Utrecht\\
\end{tabular}

\vfill
\noindent
Dit proefschrift werd mede mogelijk gemaakt met financi\"ele
steun van de\\
 Nederlandse Organisatie voor Wetenschappelijk Onderzoek.

51M20, 52C35, 57N80, 68U05.

\medskip
\noindent
ISBN 90-393-3137-5

\end{titlepage}

\setcounter{tocdepth}{1}

\tableofcontents

\chapter{Introduction.}

In this thesis we study sets of points in the plane
and their Voronoi diagrams, in particular
when the points coincide.   
We bring together two ways of studying point sets that 
have received a lot of attention 
in recent years: Voronoi diagrams and compactifications
of configuration spaces. 
We study moving and colliding points and this enables us to introduce
`limit Voronoi diagrams'. 
We define several compactifications 
by considering geometric properties of pairs and triples
of points. In this way we are able to define a smooth, real
version of the Fulton-MacPherson compactification. We show
how to define Voronoi diagrams on elements of these
compactifications  and describe the connection
with the limit Voronoi diagrams.

\subsection*{Voronoi diagrams and supermarkets.}

Consider all supermarkets in a city.
We divide the city in sectors by considering
the closest supermarket:
all people in the sector of some supermarket are closer  
to this supermarket than to any other supermarket.
The {\it Voronoi diagram} of the supermarkets
is this subdivision of the city into sectors. The {\it Voronoi
cell} of one supermarket is just the sector of the
supermarket.
An example of a Voronoi diagram of six supermarkets,
numbered $1$ to $6$, 
is given in Figure \ref{fsupervor}.

 \begin{figure}[!ht]
 \begin{center}
 \setlength{\unitlength}{.75cm}
 \begin{picture}(5,4)
 \put(2.5,2){\makebox(0,0)[cc]{
         \leavevmode\epsfxsize=5\unitlength\epsfbox{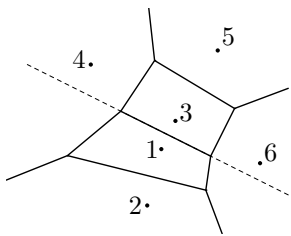}}}
 \put(2.55,1.5){\makebox(0,0)[l]{1}}
 \put(2.25,0.5){\makebox(0,0)[l]{2}}
 \put(3.15,2.15){\makebox(0,0)[l]{3}}
 \put( 1.25,3.1){\makebox(0,0)[l]{4}}
 \put(3.9,3.5){\makebox(0,0)[l]{5}}
 \put(4.65,1.45){\makebox(0,0)[l]{6}}
 \end{picture}
 \caption{\elabel{fsupervor}
 The Voronoi diagram of six supermarkets, labeled $1, 2, 3, 4, 5, 6$.}
 \end{center}
 \end{figure}

In the figure we see some boundaries, for example the boundary 
between supermarket $1$ and supermarket $3$. People living
on this boundary are at equal distance between those supermarkets.
The boundary is part of
a line, the dotted line in the picture, 
which is called the {\it bisector} of
$1$ and $3$, because on one side of the line  people
are closer to supermarket $1$ and on the other side
people are closer to supermarket $3$. So, the bisector 
divides the plane into two sectors or {\it half-planes}. The sector 
containing supermarket $1$ is called
{\it Voronoi half-plane} $vh(1,3)$ and the sector
on the other side of the bisector is Voronoi
half-plane  $vh(3,1)$. We can express the
sector $V(1)$ of supermarket $1$ in terms of these 
Voronoi half-planes:
\begin{equation}
  V(1)  ~=~  vh(1,2) \cap vh(1,3) \cap 
vh(1,4) \cap vh(1,5) \cap vh(1,6)
\elabel{ev1}
\end{equation}
Of course this formula just states that you are in the sector
of supermarket $1$ if you are closer to supermarket $1$
than to $2$, closer to supermarket $1$ than to $3$, 
etcetera. If we want to code the bisector of 
$1$ and $3$, we need only two ingredients: one
point on the bisector and the {\it angle} that the bisector makes
with the horizontal line. This will be important later on.

\subsection*{Driving supermarkets: changing Voronoi diagrams.}

In the Dutch countryside food is supplied by
little supermarket lorries, or supervans.
These vans drive through the countryside looking for customers.
In the part of the countryside we are considering there are
five supervans. If we know their positions at a certain moment,
we can determine the Voronoi diagram of the five vans
at that moment. But, if the  vans drive continuously the 
Voronoi diagram of the vans changes continuously
as well. In Chapter~\ref{chlimit} we model
these driving supervans and their changing Voronoi diagrams. 
For every supervan there is
a curve that gives at any time $t$  the position of the supervan.
If we want to know the Voronoi diagram at time $t$, we specify $t$
in the curves describing the position of the vans.  Using 
positions given by the curves we compute the Voronoi diagram.

\begin{figure}[!ht]
\begin{center}
\setlength{\unitlength}{1.02cm}
\begin{picture}(12,4)
\put(1.5,2.3){\makebox(0,0)[cc]{
        \leavevmode\epsfysize=3\unitlength\epsfbox{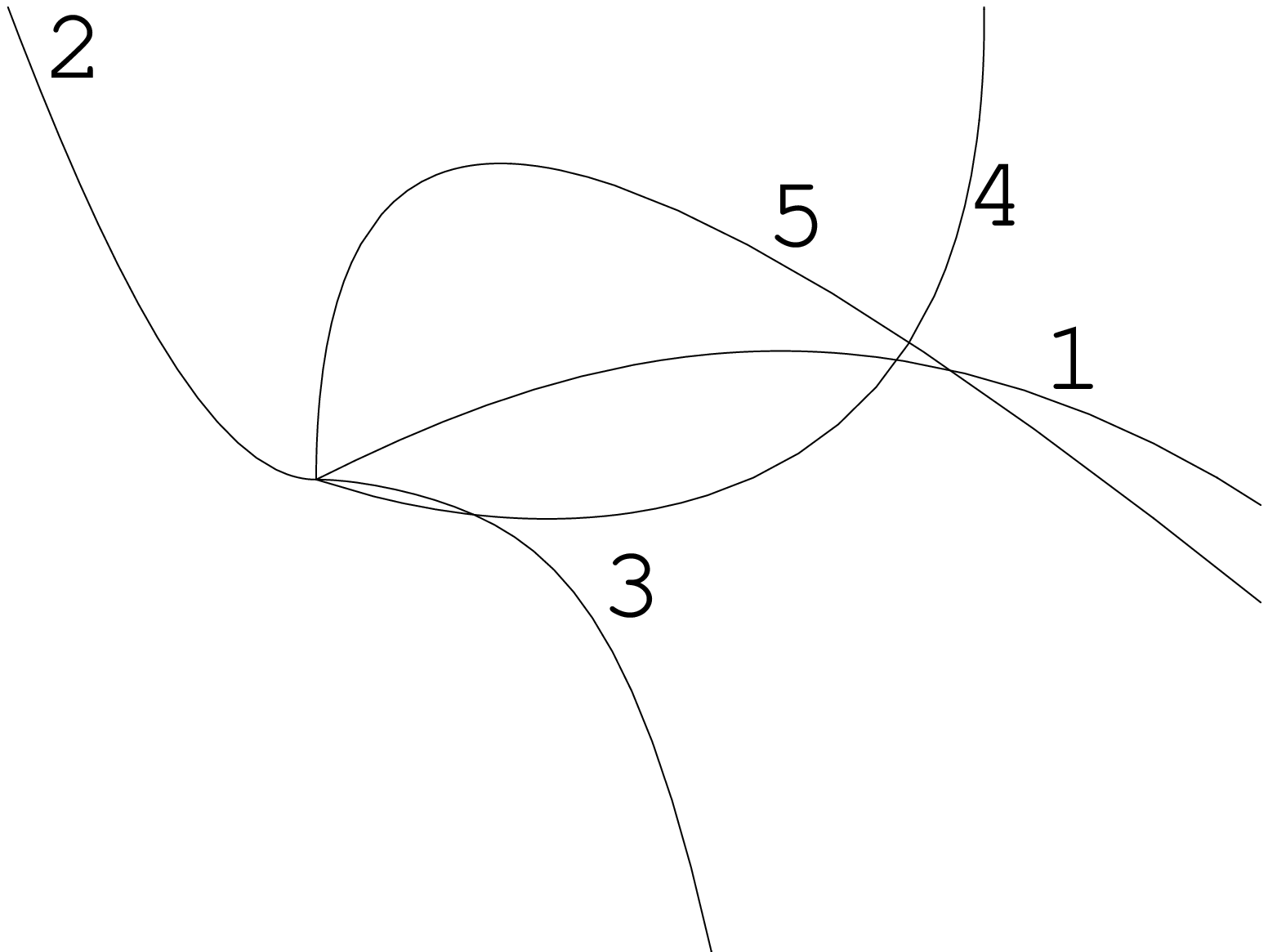}}}
\put(6,2.3){\makebox(0,0)[cc]{
        \leavevmode\epsfysize=3\unitlength\epsfbox{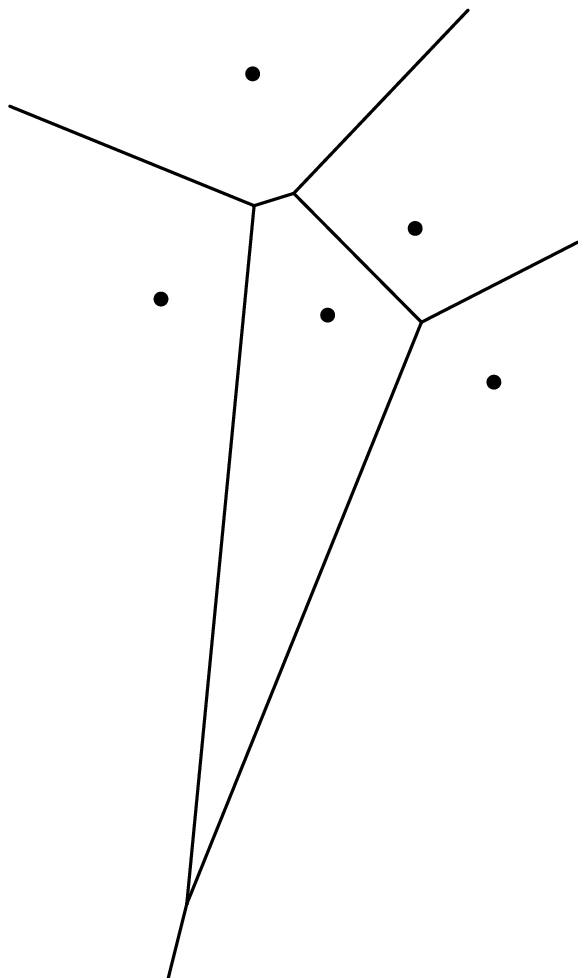}}}
\put(10.4,2.3){\makebox(0,0)[cc]{
        \leavevmode\epsfysize=3\unitlength\epsfbox{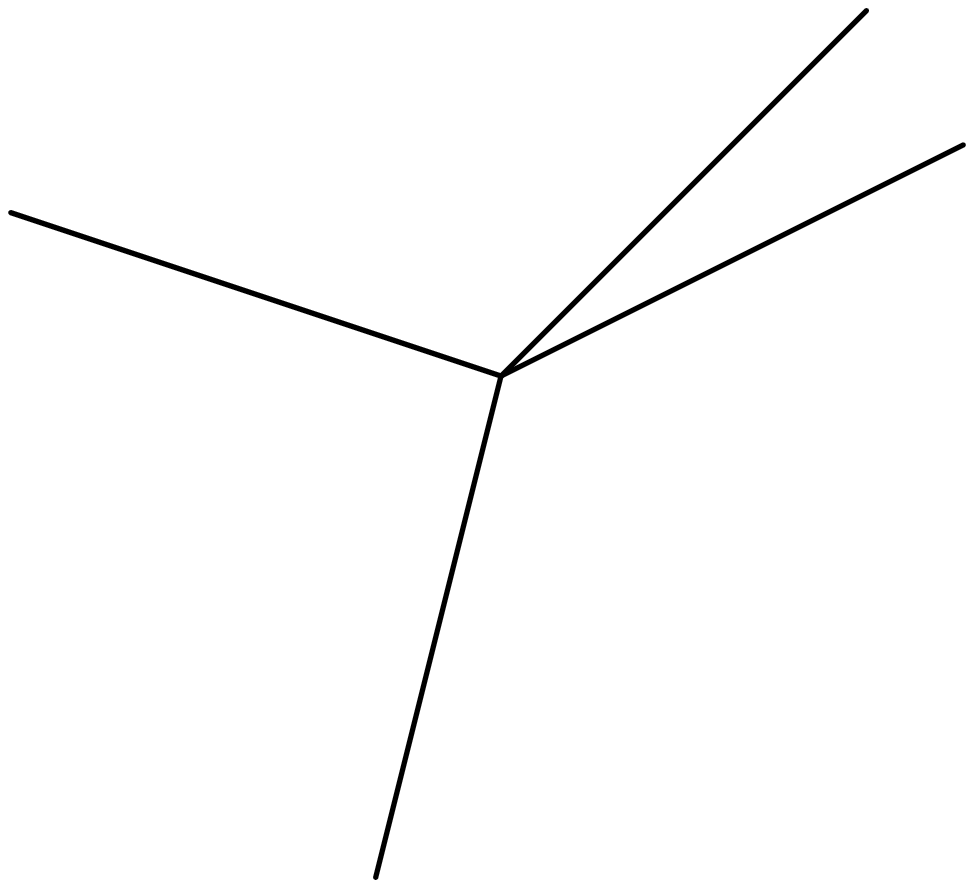}}}
\put(9.75,1.05){\makebox(0,0)[l]{{\tt 2}}}
\put(8.8,2.75){\makebox(0,0)[l]{{\tt 2}}}
\put(8.95,3.35){\makebox(0,0)[l]{{\tt 5}}}
\put(11.3,3.85){\makebox(0,0)[l]{{\tt 5}}}
\put(11.7,3.45){\makebox(0,0)[l]{{\tt 1}}}
\put(11.95,2.88){\makebox(0,0)[l]{{\tt 4}}}
\put(10.3,1.0){\makebox(0,0)[l]{{\tt 4}}}
 \put(1.5,0.1){\makebox(0,0)[l]{a}}
 \put(6,0.1){\makebox(0,0)[l]{b}}
 \put(10.4,0.1){\makebox(0,0)[l]{c}}
\put(5.4,2.75){\makebox(0,0)[l]{{\tt 2}}}
\put(5.8,3.7){\makebox(0,0)[l]{{\tt 5}}}
\put(6.5,3.2){\makebox(0,0)[l]{{\tt 1}}}
\put(6,2.6){\makebox(0,0)[l]{{\tt 3}}}
\put(6.7,2.4){\makebox(0,0)[l]{{\tt 4}}}
\end{picture}
\caption{a.\ The positions of five driving supervans before collision.
b.\ The Voronoi diagram of the vans,  just before collision,
and c. at collision.
\elabel{fvorlimit}}
\end{center}
\end{figure}

There is only one problem.
At one day a very unfortunate accident happens: all vans run into
each other at $t=0$. How should we define the Voronoi diagram of
the vans at $t=0$? In this particular example
we decide to analyze the situation just before the accident 
happens. That is, we let time run backwards and investigate what
is going on for small negative $t$.
 The route covered by the five vans before 
the accident is shown in Figure \ref{fvorlimit}.a. A Voronoi
diagram of the positions of the vans just before the accident
is given in Figure \ref{fvorlimit}.b. 

Our strategy is to
define a Voronoi diagram at $t=0$ that is consistent with the situation
just before $t=0$. We call such a diagram a `limit Voronoi 
diagram' and  an example of this that resembles
the Voronoi diagram in Figure \ref{fvorlimit}.b is
shown in Figure \ref{fvorlimit}.c. Compare the two diagrams.
Note that the directions of the outgoing edges in the
two diagrams are very similar. The Voronoi cell
of van $3$, in the  middle  of the diagram of
Figure Figure \ref{fvorlimit}.b seems to have completely 
disappeared in  Figure \ref{fvorlimit}.c. These two remarks
address questions that we answer in Chapter \ref{chlimit}
and the chapters after that: 
what information is needed to to create a diagram as shown
in figure \ref{fvorlimit}.c? And, 
can we still say something about the Voronoi cells that
seem to disappear in the limit, like the cell of
van $3$? 

\subsection*{Voronoi diagrams for coinciding points.}

\begin{figure}[!ht]
\begin{center}
\setlength{\unitlength}{0.8cm}
\begin{picture}(15,5)
\put(1.5,2.8){\makebox(0,0)[cc]{
        \leavevmode\epsfxsize=2.0\unitlength\epsfbox{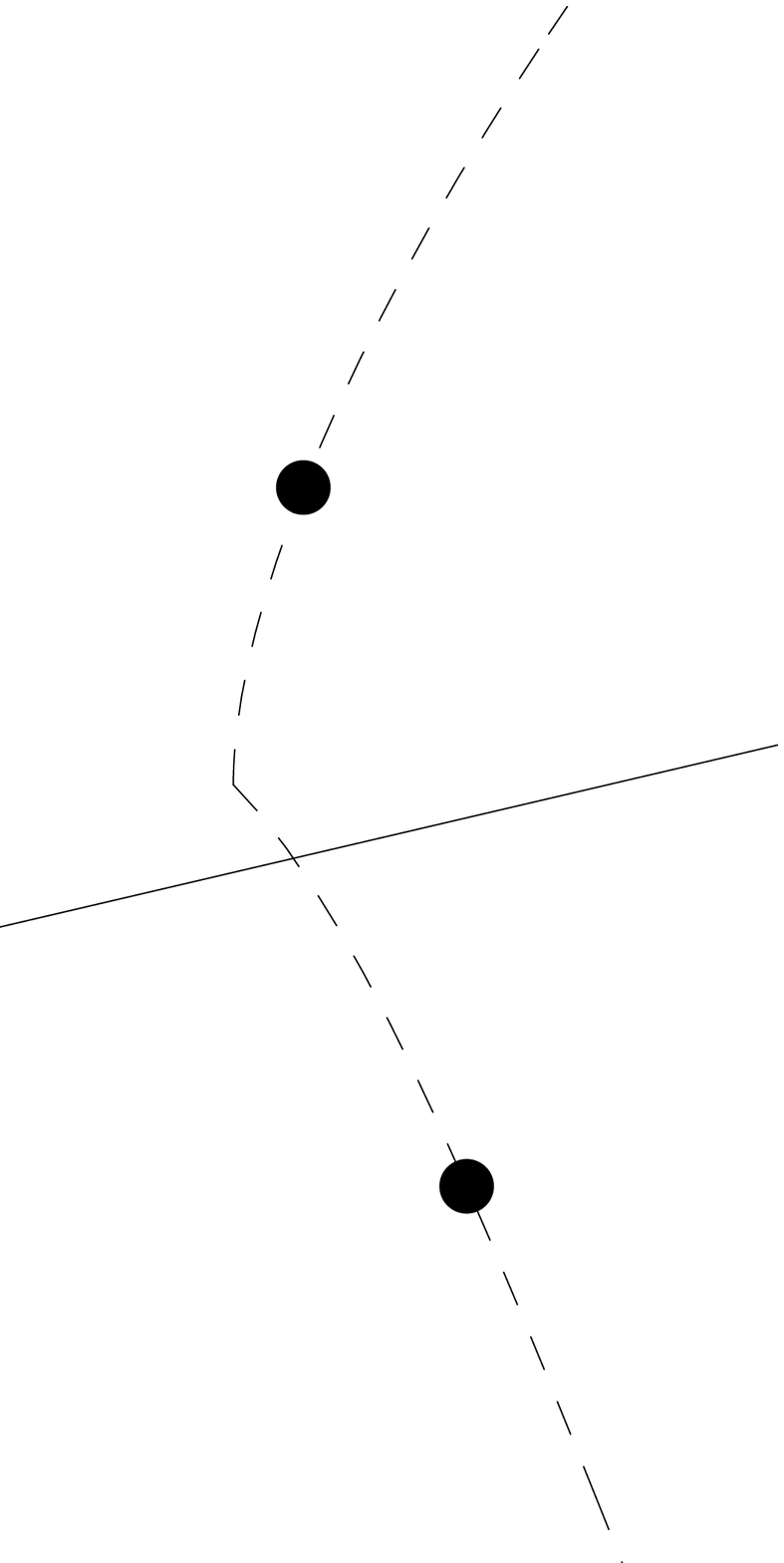}}}
\put(7.5,2.8){\makebox(0,0)[cc]{
        \leavevmode\epsfxsize=2.0\unitlength\epsfbox{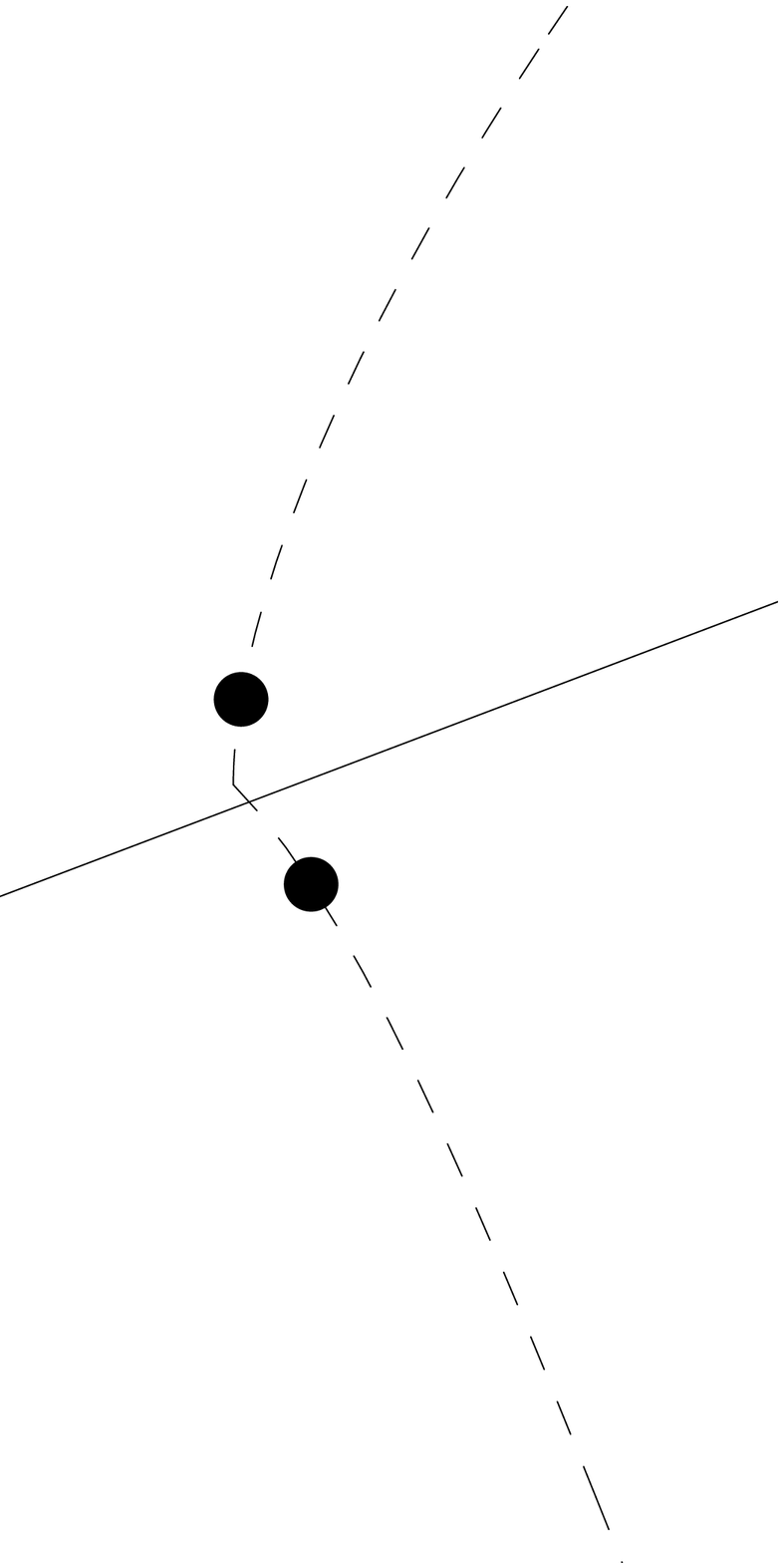}}}
\put(13.5,2.8){\makebox(0,0)[cc]{
        \leavevmode\epsfxsize=2.0\unitlength\epsfbox{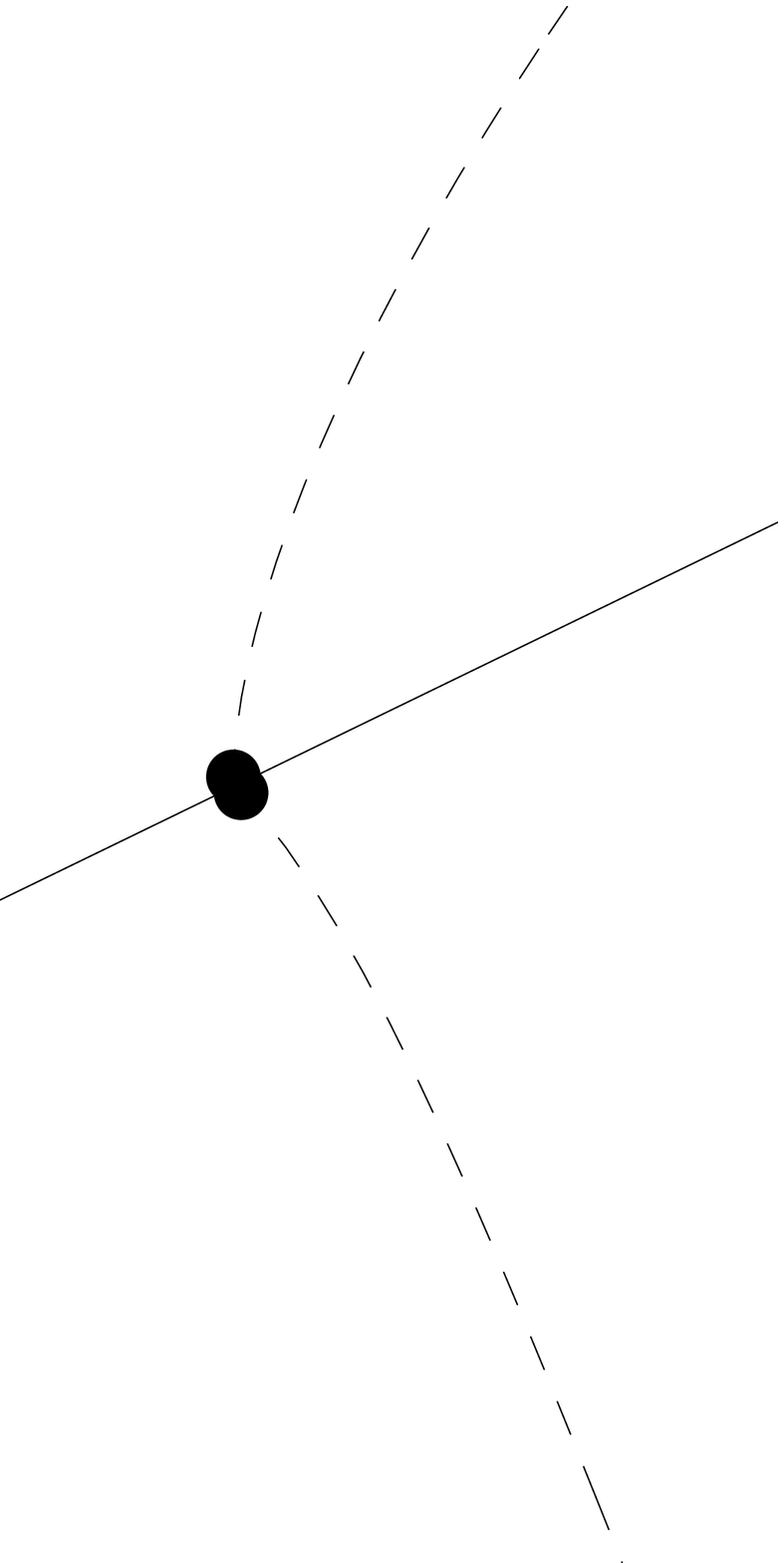}}}
\put(0.8,3.6){\makebox(0,0)[l]{ 1}}
\put(1.2,1.7){\makebox(0,0)[l]{ 2}}
\put(6.7,3.1){\makebox(0,0)[l]{ 1}}
\put(6.9,2.4){\makebox(0,0)[l]{ 2}}
\put(12.7,3.0){\makebox(0,0)[l]{ 1}}
\put(13.0,2.5){\makebox(0,0)[l]{ 2}}
 \put(1.5,0.0){\makebox(0,0)[l]{a}}
 \put(7.5,0.0){\makebox(0,0)[l]{b}}
 \put(13.5,0.0){\makebox(0,0)[l]{c}}
\end{picture}
\caption{The Voronoi diagram of two coinciding points.
 \elabel{fthreebisectors}}
\end{center}
\end{figure}
We restrict ourselves for a moment to the case of two points, $p_1$ 
and $p_2$.
The position of each point at time $t$ is given by a curve.
Suppose that the two points meet at $t=0$.
Such situation is shown in Figure \ref{fthreebisectors}.a-c. 
In Figures \ref{fthreebisectors}.a and b, the two
points are still distinct.  Therefore we can draw the line $l_{12}$
that passes through both points. This line makes
some angle $\alpha_{12}$ with the horizontal axis.  So this angle
$\alpha_{12}$ is in some sense the {\it angle} of the points $p_1$ and $p_2$. 
As the positions of the points depend on time $t$, the angle
$\alpha_{12}$ also depends on time, that is $\alpha_{12} = \alpha_{12}(t)$.
The Voronoi diagram of points $p_1$ and
$p_2$ is determined by the line at equal distance from the points $p_1$
and $p_2$. That line is exactly the bisector of $p_1$ and $p_2$: 
the line perpendicular to line $l_{12}$ passing through 
the middle of the line segment $12$. If $p_1$ and $p_2$ coincide, we define
the middle of the line segment $12$ as the point $p_1=p_2$ itself. And
we define the angle $\alpha_{12}(0)$  as the limit for small negative
$t$ of $\alpha_{12}(t)$. Now we define the bisector at $t=0$ in terms 
of this angle
and this middle point. That is, the bisector of $p_1=p_2$ is the line
passing through $p_1=p_2$ perpendicular to the direction $\alpha_{12}(0)$.
But this implies that we have created a Voronoi diagram 
for the two coinciding points $p_1$ and $p_2$!

We drop this particular example but conclude the following:
we can define a limit Voronoi diagram for two 
coinciding points $p_1$ and $p_2$ if we know the following information:\\
 \phantom{kip} -- the position $p_1=p_2$ of the coinciding points.\\
 \phantom{kip} -- an angle $\alpha_{12}$ mod $2 \pi$.\\
In this way we can define a limit Voronoi diagram for 
an arbitrary number of coinciding points as well: as long as we have
for every pair of coinciding points $p_i$ and $p_j$
a position $p_i=p_j$ and  an 
angle $\alpha_{ij}$, we can define the bisector of 
$p_i$ and $p_j$. And, using the
bisector we can determine the two Voronoi half-planes
$vh(p_i,p_j)$ and $vh(p_j,p_i)$. Now we are done,
as any Voronoi cell can be expressed as an intersection
of half-planes, as we saw in Equation \ref{ev1}.

We work along these lines in Chapter \ref{chlimit}. 
For example, the curves
describing the positions of the two
points $p_1=p_1(t)$ and $p_2=p_2(t)$ in Figure \ref{fthreebisectors}
 are given by
$p_1(t) = (t^2, t + t^2 -  .3 t^3)$ and $p_2(t) = (t,-t-3 t^2 + 2 t^3)$.
We only allow curves given by pairs of polynomials in $t$.
We call points described by such curves {\it polynomial sites}. 
After having defined Voronoi diagrams for polynomial sites we
show how to determine the Voronoi diagrams without having to compute
all bisectors.

First we assume that we have a set of $n$ polynomial sites
that all coincide at $t=0$. So, this is the situation in 
Figure \ref{fvorlimit}.a,  where five polynomial sites cluster together.
Leaving out polynomial site
$p_3(t)$ in the beginning will not change the shape of the  limit diagram as we see
in figure \ref{fvorlimit}.c. We show in Section \ref{sposcell} which sites exactly
can be left out without changing the boundary of the diagram. 
So,  we want to characterize those sites that are somehow
in the interior of the cluster. Of course we could as well  
characterize the sites that are on the exterior of the cluster.
This is done in terms of the convex hull of the sites in the cluster
for small enough positive~$t$. 

\begin{figure}[!ht]
\begin{center}
\setlength{\unitlength}{0.9cm}
\begin{picture}(11,3.5)
\put(1.5,2.0){\makebox(0,0)[cc]{
        \leavevmode\epsfxsize=3\unitlength\epsfbox{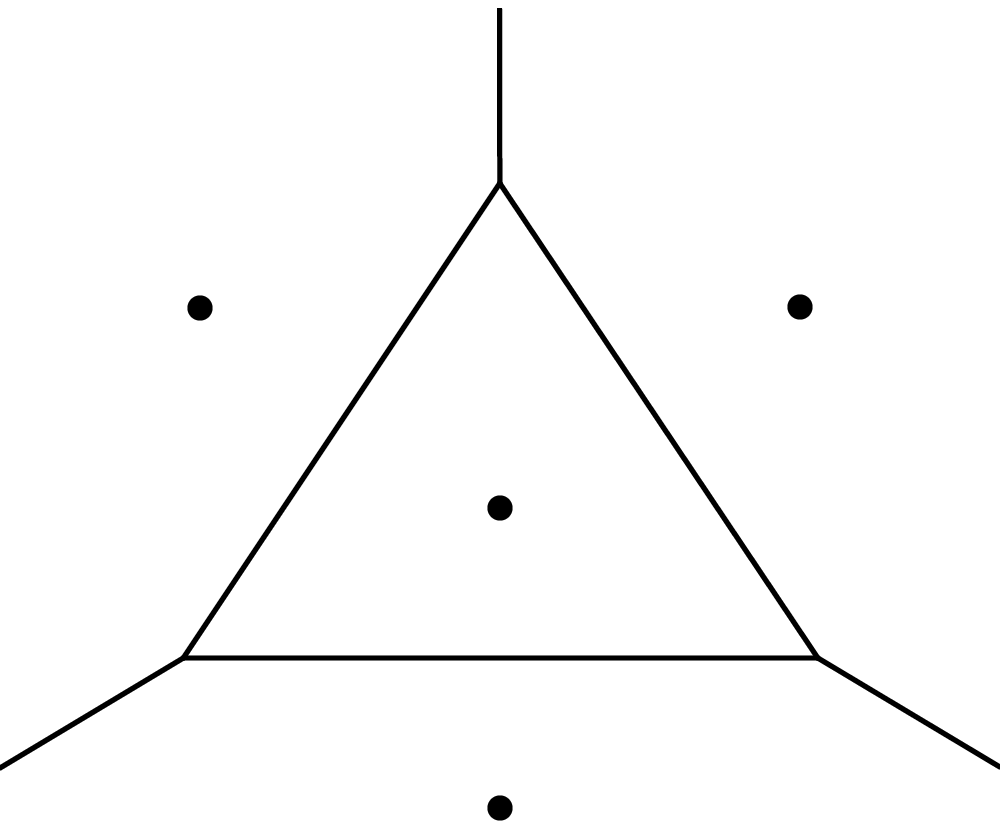}}}
\put(5.5,1.79){\makebox(0,0)[cc]{
        \leavevmode\epsfysize=2.5\unitlength\epsfbox{ilimitright.eps}}}
\put(9.5,1.98){\makebox(0,0)[cc]{
        \leavevmode\epsfxsize=3\unitlength\epsfbox{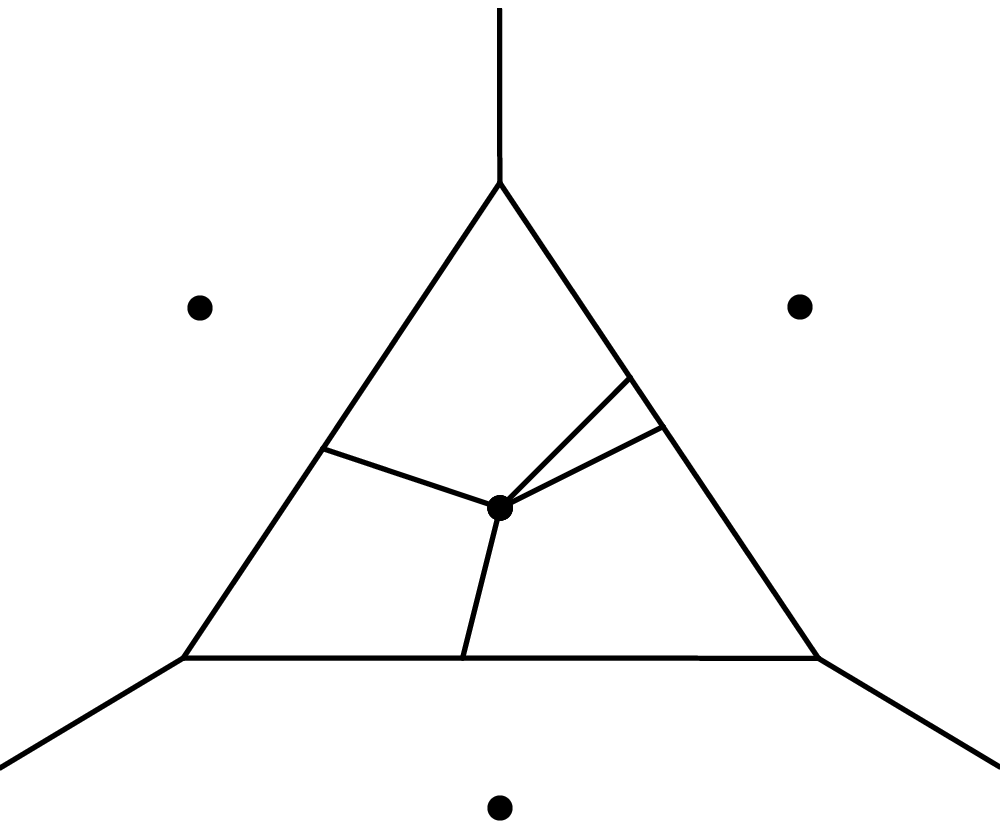}}}
\put(3.5,1.85){\makebox(0,0)[l]{ +}}
\put(7.2,1.85){\makebox(0,0)[l]{ =}}
\put(1.5,0.1){\makebox(0,0)[l]{ a}}
\put(5.5,0.1){\makebox(0,0)[l]{ b}}
\put(9.5,0.1){\makebox(0,0)[l]{ c}}
\put(0.2,2.6){\makebox(0,0)[l]{ 6}}
\put(2.5,2.6){\makebox(0,0)[l]{ 8}}
\put(1.15,1.1){\makebox(0,0)[l]{ 7}}
\put(8.2,2.6){\makebox(0,0)[l]{ 6}}
\put(10.5,2.6){\makebox(0,0)[l]{ 8}}
\put(9.15,1.1){\makebox(0,0)[l]{ 7}}
\end{picture}
\caption{Plugging  Voronoi diagrams.
 \elabel{foneplug}}
\end{center}
\end{figure}

For a general set $S(t)$ of $n$ polynomial sites, the
boundary of the Voronoi diagram at $t=0$ is determined by {\it plugging}.
We demonstrate plugging in Figure \ref{foneplug} where 
we compute the Voronoi diagram of 
$S(t) = \{p_1,\dots, p_8\}$, for $p_6 = (-3,2)$, $p_7 = (0,-3)$
and $p_8=(3,2)$. The other polynomial sites form a cluster in the
origin. In fact we used this cluster already for the supervan example
in Figure \ref{fvorlimit}.
The cluster locations of $S(t)$ at $t=0$ are the distinct positions 
of the sites at $t=0$. So, we have four cluster locations.
The Voronoi diagram of the cluster locations is given in 
Figure \ref{foneplug}.a. 
Next, the cell of every cluster location is filled,
if necessary. Three cells in our example correspond to a single point,
so no filling is needed there.  
The limit Voronoi diagram of the points clustering in the origin,
is shown in Figure \ref{foneplug}.b. We simply plug this
diagram in the appropriate cell in order to obtain
the limit Voronoi diagram of $S(t)$. The result is shown
in Figure \ref{foneplug}.c.

The limit Voronoi diagrams are really new diagrams: in general 
they can not be realized as classic Voronoi diagrams, for example
Figure \ref{foneplug}.c. In classic Voronoi
diagrams, every Voronoi cell has positive area, for example. This is not
true for limit Voronoi diagrams.

\subsection*{Compactifications of configuration spaces.}

A collection of $n$ points in $\R^2$ is often called a {\it configuration}.
The {\it configuration space} of $n$ distinct points in $\R^2$
is just the space that consists of all possible configurations
of $n$ distinct, ordered points. Suppose for example that $n=3$. Any element 
of the configuration space consists of three lebeled 
points $p_1,p_2, p_3 \in \R^2$ 
such that:  $p_1$ is distinct from $p_2$; $p_1$ is distinct from $p_3$;
and $p_2$ is distinct from $p_3$.  A natural way of describing 
the configuration $(p_1,p_2,p_3)$ is  by listing  the coordinates 
of $p_1$, $p_2$, and $p_3$. But if we list all six coordinates in
single file, we obtain an element $c=(p_1,p_2,p_3)$ in $\R^6$.
That is, the configuration space of three distinct point in the
plane is part of a  six dimensional space. In fact it
is six dimensional, as almost all elements of $\R^6$ can be seen
as some configuration $(p_1,p_2,p_3)$ of three distinct points.

We are interested in the Voronoi diagram of $n$ points in the plane.
That is, we want to define a Voronoi diagram for every configuration
in a configuration space of $n$ points. We are especially interested
in possible Voronoi diagrams for point sets that contain coinciding
points. We call such point sets {\it degenerate} configurations.
The idea of a compactification of a configuration space is as follows:
we want to construct a space that encodes
all possible configurations of $n$ distinct points, both
non-degenerate and degenerate.
Degenerate configurations
should be on the
boundary of this space. By adding this boundary
the degenerate configurations are incorporated. For a bounded space,
adding the boundary is the same as compactifying, which explains
the name `compactification'.
One important reason to compactify is the hope to be able to extend
some definition, in our case the definition of Voronoi diagram,
to the degenerate configurations. This extension would 
give access to limit objects which in our case are limits of Voronoi
diagrams.

We denote the set of ordered $n$ tuples of all pairwise distinct 
points in the plane
by $\CONF$.  Although this is not very useful for us,
an easy example is the compactification $(\R^2)^n$ of 
the configuration space $\CONF$ itself. 
The problem with this compactification is that it gives very
little information on degenerate configurations: for two coinciding
points $p_1=p_2$,  the only point that is added by compactifying
is the point $(p_1,p_2)$. 
Consider the two points $p_1(t)$ and $p_2(t)$ where 
$p_2(t) = p_1(t) + t(\cos \alpha, \sin \alpha)$. 
If $t$ goes to zero, $p_2(t)$ will coincide with $p_1(t)$. Describe the
Voronoi diagram of $p_1(0)$ and $p_2(0)$ as in Figure 
\ref{fthreebisectors}.
Then every distinct value of $\alpha$ corresponds to a 
distinct direction of the bisector of the points $p_1(0)$ 
and $p_2(0)$. But this means that there is not one 
Voronoi diagram corresponding to $p_1(0)=p_2(0)$ but a
complete collection of diagrams, parameterized  by $\alpha$~mod $2 \pi$.

\subsection*{Compactifying using angles.}

By now, the attentive reader should be convinced that
it is not a strange idea to use angles for
a suitable compactification.
In Chapter~\ref{changles} we may consider all angles between
$n$ points:  we start with $n$ distinct, labeled points
in the plane. For every two points $p_i$ and $p_j$ 
with $i< j$ we
consider the directed line $l_{ij}$ passing through the two points.
This line makes some angle $\alpha_{ij}$ mod $2 \pi$ with 
the horizontal axis. The {\it angle map} $\pda$ maps a configuration
of $n$ distinct points to the $\binom{n}{2}$ directed angles 
$\alpha_{ij}$ with $i<j$. 
Let $T$ be a triangle with vertices $1$, $2$ and $3$.
If we know the angles $\alpha_{12}$, $\alpha_{13}$
and $\alpha_{23}$, we know in fact the shape of the triangle $T$. 
So, $\pdat(\CONFt)$ describes in a way all distinct triangles
on three points.
We describe explicitly this image $\pdat(\CONFt)$ in Chapter~\ref{changles}.
We explain the consequences  of the action of the symmetric group $S_3$ 
on the labels of the three vertices. Finally we give geometric interpretations 
of boundary points
of $\pdat(\CONFt)$ in terms of degenerate configurations.

The first compactification that we encounter is $\cda$, the 
{\it compactification} of {\it directed angles}.
This space is defined as the closure of the graph of the angle map
$\pda$, so a point $\gamma_n \in \cda$ consists of $n$  points
and $\binom{n}{2}$ angles between those points. If $\gamma_n$
is a boundary point of $\cda$, then not all points need
to be distinct. In any case we know for every two points
$p_i,p_j \in \gamma_n$ an angle $\alpha_{ij}$ and the position(s) 
$p_i \neq p_j$ or $p_i=p_j$. This implies that for 
$\gamma_n$ a Voronoi diagram $V(\gamma_n)$ is defined.

In Chapter~\ref{chcont} we prove a continuity theorem
for Voronoi diagrams of data sets from $\cda$.
It states essentially that
two data sets $\gamma_n, \eta_n \in \cda$ that are
close in the Euclidean metric,
have Voronoi diagrams whose boundaries
are close in the Hausdorff metric. The Hausdorff
metric is very suitable to
compare images. Two sets $A$ and $B$ in the
plane are within Hausdorff distance $r$ iff
$r$ is the smallest number such that any point
of $A$ is within distance $r$ from some point in $B$
and vice versa.

\subsection*{An algebraic description:  from angles to slopes.}

The line through the points $(x_1,y_1)$ and $(x_2,y_2)$
is given algebraically by the equation
\begin{equation}
	y - y_1 ~ = ~ \frac{y_2-y_1}{x_2-x_1}(x-x_1) \elabel{eslope}
\end{equation}
The quantity $\frac{y_2-y_1}{x_2-x_1}$ is of course the {\it slope} of 
the line.
Instead of writing down the angle $\alpha_{ij}$ for every
two points $p_i$ and $p_j$, we consider all slopes $a_{ij}$.
It turns out that there exists an algebraic relation between the six 
possible slopes for four distinct points  $p_0$, $p_1$, $p_2$ and $p_3$. 
This relation is called the {\it six-slopes} formula,
compare Figure \ref{fsixslopes}. The formula is given
by $\Delta_{0123} = 0$, where
$$
\Delta_{0123} =(a_{01}-a_{12})(a_{02}-a_{23})(a_{03}-a_{13})
                -(a_{01}-a_{13})(a_{02}-a_{12})(a_{03}-a_{23}).
$$
We prove in Chapter \ref{changles} that there also exists a
{\it triangle} formula. This is a relation
$t_{ijk}=0$ between slopes and $x$-coordinates of the points 
$p_i$, $p_j$, and $p_k$.
So by now we have two formulas involving points and slopes:
one holds for any three points and the other holds for any 
four points.  
Suppose we have some configuration $c$
of points and slopes such that all triangle formulas and all
six-slopes formulas hold.
A question that we consider in Chapter \ref{changles}
is the following:  is such configuration $c$ always a limit
of non-degenerate configurations. 
To that end we introduce a variety that consists of exactly those
configurations such that all triangle formulas and all
six-slopes formulas hold. 
That is: the algebraic variety $T_n$ is just the set of zeroes common to
all triangle formulas and all six-slopes formulas.
It serves as an algebraic counterpart of the compactification
$\cua$ of {\it undirected angles}.   We prove
that $T_n$ is not a smooth variety and give a geometric
interpretation for the non-singular configurations.

 \begin{figure}[!ht]
 \begin{center}
 \setlength{\unitlength}{1cm}
 \begin{picture}(6,2.4)
 \put(3,1.2){\makebox(0,0)[cc]{
         \leavevmode\epsfxsize=6\unitlength\epsfbox{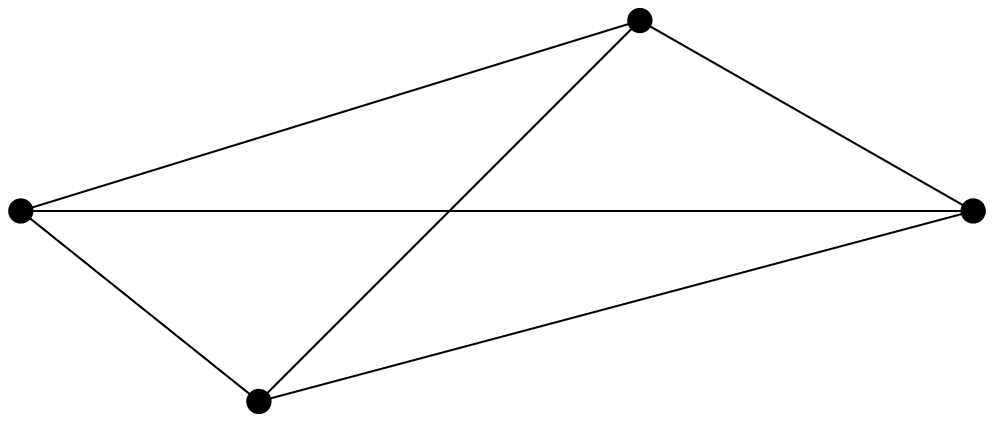}}}
 \put(0.5,0.5){\makebox(0,0)[l]{$a_{01}$}}
 \put(1.8,1.4){\makebox(0,0)[l]{$a_{02}$}}
 \put(1.85,2.05){\makebox(0,0)[l]{$a_{03}$}}
 \put(3.9,0.45){\makebox(0,0)[l]{$a_{12}$}}
 \put(3.43,1.7){\makebox(0,0)[l]{$a_{13}$}}
 \put(4.75,2.0){\makebox(0,0)[l]{$a_{23}$}}
 \end{picture}
 \caption{\elabel{fsixslopes}
 For four distinct points the six-slopes formula holds.}
 \end{center}
 \end{figure}

\subsection*{A smooth and clickable compactification.}

Look at Figure \ref{fvorlimit} again. On the right, a limit Voronoi
diagram of the points $p_1, \dots, p_5$ is displayed. For every
point $p_i$ a Voronoi cell is visible, except for $p_3$.
One could say that the Voronoi cell $V(p_3)$ is so small that
we cannot see it. A solution in this case could be  to rescale or magnify 
the picture somehow until the cell $V(p_3)$ becomes visible. 
After magnifying enough we would get a picture that is very similar
to Figure \ref{fvorlimit}.b
In more complicated
configurations we might need to rescale or zoom in several times
at several positions in order to distinguish every cell. 

\begin{figure}[!ht]
\begin{center}
\setlength{\unitlength}{1cm}
\begin{picture}(6,6)
\put(3,5){\makebox(0,0)[cc]{
        \leavevmode\epsfxsize=2\unitlength\epsfbox{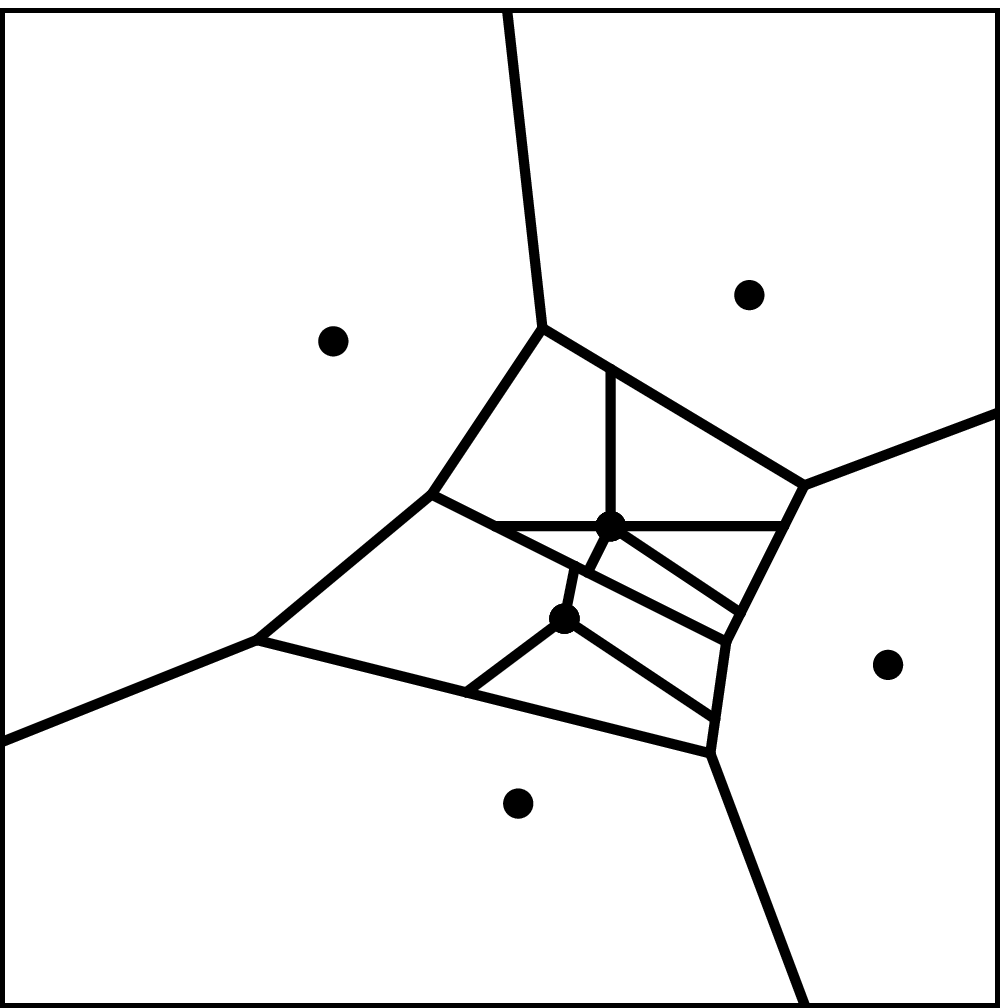}}}
\put(1.5,1){\makebox(0,0)[cc]{
        \leavevmode\epsfxsize=2\unitlength\epsfbox{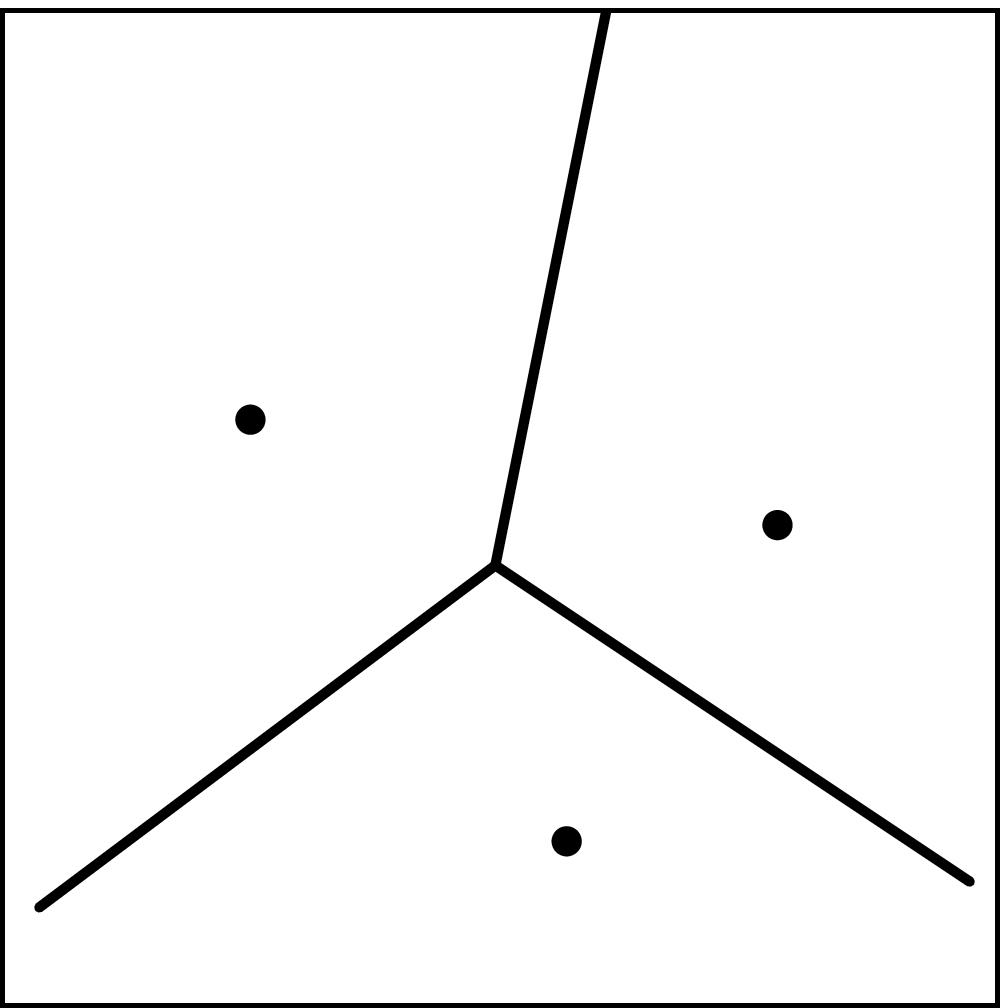}}}
\put(4.5,1){\makebox(0,0)[cc]{
        \leavevmode\epsfxsize=2\unitlength\epsfbox{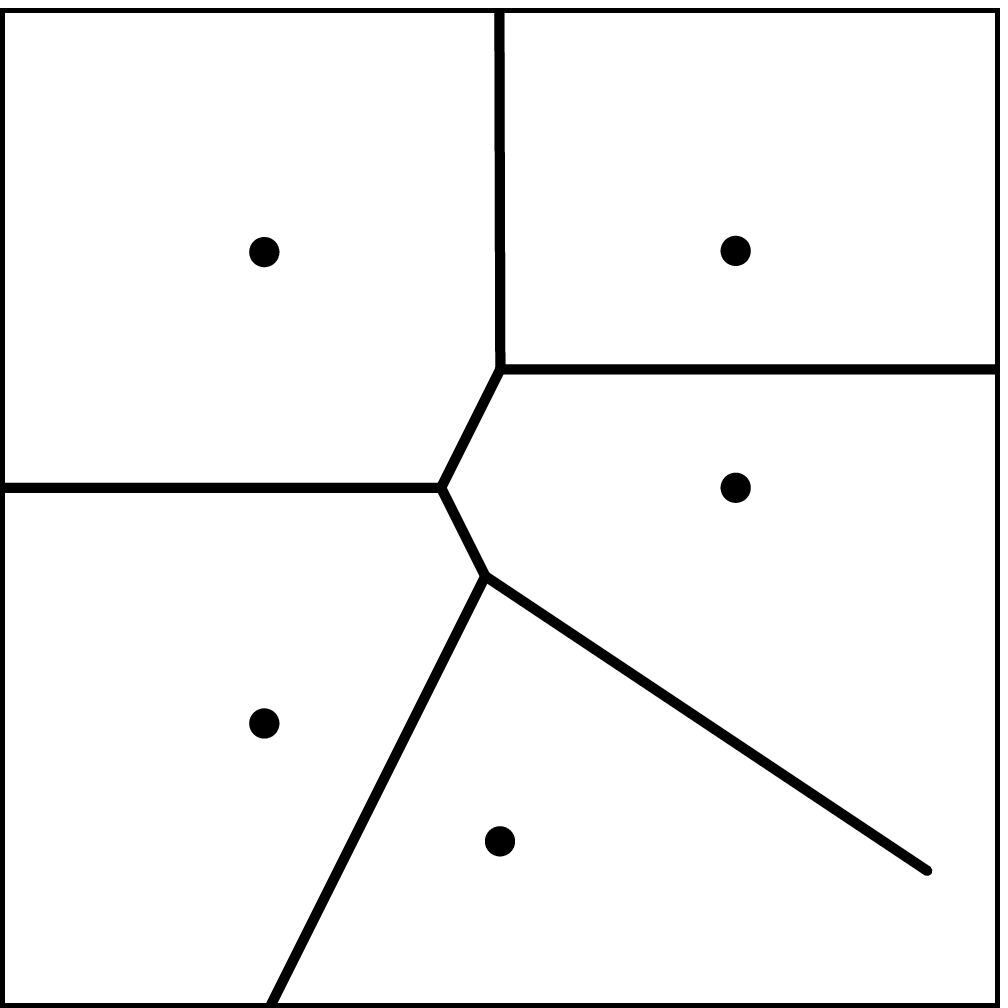}}}
%
%
\put(2.6,5.5){\makebox(0,0)[c,c]{1}}
\put(3.1,4.2){\makebox(0,0)[c,c]{2}}
\put(3.6,5.7){\makebox(0,0)[c,c]{5}}
\put(3.8,4.4){\makebox(0,0)[c,c]{6}}
\put(1.45,.4){\makebox(0,0)[c,c]{3}}
\put(0.9,1.3){\makebox(0,0)[c,c]{7}}
\put(2.1,1.3){\makebox(0,0)[c,c]{12}}
\put(4.8,.3){\makebox(0,0)[c,c]{4}}
\put(3.8,.6){\makebox(0,0)[c,c]{8}}
\put(3.9,1.7){\makebox(0,0)[c,c]{9}}
\put(5.1,.8){\makebox(0,0)[c,c]{10}}
\put(5,1.8){\makebox(0,0)[c,c]{11}}
\unitlength=1mm
\linethickness{.3mm}
\dottedline{1}(17,20)(32,48)
\dottedline{1}(43,20)(33,49)
\end{picture}
 \caption{\elabel{fclick}
A clickable Voronoi diagram in a clickable configuration.}
\end{center}
\end{figure}

This shows the need of a `clickable' compactification on which we can define
`clickable' Voronoi diagrams. An example is given in Figure \ref{fclick}.
In the top `screen', six points are visible. But the two points in the middle
are in fact clusters of points. One cluster consists of the points $p_3$, $p_7$
and $p_{12}$. This can be `seen' by clicking on the cluster:
The bottom left screen appears, displaying these points separately. The bottom
right screen appears after clicking on the other cluster and displays the
points $p_4$, $p_8$, $p_9$, $p_{10}$ and $p_{11}$. 
 
The Fulton-MacPherson compactification is a well-known compactification
of the configuration space of $n$ distinct labeled points
that has such a `clickable' structure. 
To analyze Voronoi diagrams we do not need the full 
generality of this compactification. That is why we construct a
real version of this compactification, incorporating some
ideas of Kontsevich-Soibelman. Let $c \in \CONF$ be some configuration
of $n$ distinct points. 
Kontsevich and Soibelman write down both the angles $\alpha_{ij}$ 
between any two points, and, for every ordered triple of points $(p_i,p_j,p_k)$,
the ratio $\be{ik}{ij}=\frac{|p_i-p_k|}{|p_i-p_j|}$. 
They take the closure
of the image of $\CONF$ under the map that assigns all angles
and ratios to a set of distinct points. This results in a manifold
with corners $\FMt(n)$.

\begin{figure}[!ht]
\begin{center}
\setlength{\unitlength}{1cm}
\begin{picture}(5,3.2)
\put(2.5,1.6){\makebox(0,0)[cc]{
        \leavevmode\epsfxsize=4\unitlength\epsfbox{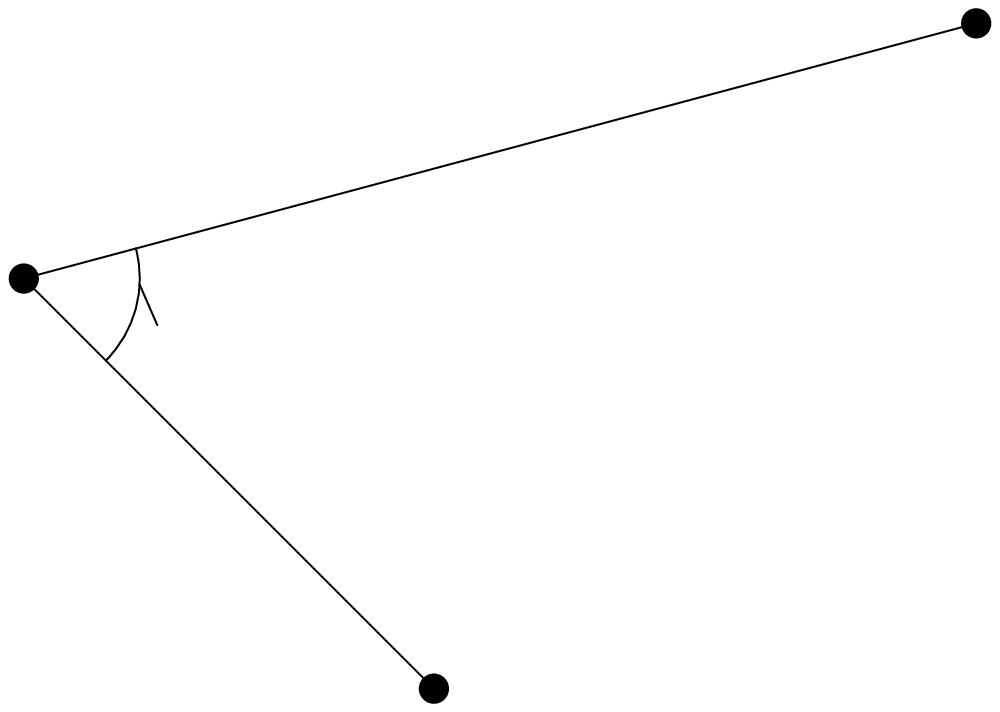}}}
\put(1.35,1.8){\makebox(0,0)[l]{$\al{ik}{ij}$}}
\put(2.4,2.8){\makebox(0,0)[l]{$\be{ik}{ij}$}}
\put(0.2,2){\makebox(0,0)[l]{$p_i$}}
\put(2.5,0.5){\makebox(0,0)[l]{$p_j$}}
\put(4.7,2.9){\makebox(0,0)[l]{$p_k$}}
\end{picture}
 \caption{\elabel{fhookintro}
The hook $\ho{ik}{ij}=(\be{ik}{ij},\al{ik}{ij})$ hinged at $p_i$ from $p_j$ to $p_k$.
The little arrow indicates the positive direction.}
\end{center}
\end{figure}

In Chapter~\ref{chmodel} we show how to adapt their approach so that
the resulting compactification is a smooth manifold. 
Instead of just a ratio, we write down a {\it hook} for every
triple of points.  A hook $\ho{ik}{ij}$ consists of this ratio 
$\be{ik}{ij}$ together with an angle $\al{ik}{ij}$ 
between the legs of the hook, cf.\ Figure 
\ref{fhookintro}.  
We define a space $\xah$ as the closure of the image of $\CONF$ under the
map
\begin{displaymath}
\begin{array}{rccl}
        \psi_{\ah}:& \CONF & \rightarrow &
                 (\R / \pi \Z)^{\binom{n}{2}} \times
                        ( (-\infty, \infty] \times \R/ 2 \pi \Z) / \sim_k)
                                ^{6\binom{6}{3}},\\
	& (p_1, \dots, p_n) & \mapsto &
		((\alpha_{ij})_{1 \leq i < j \leq n},(\be{ik}{ij},\al{ik}{ij})),
	\quad i, j, k~ \text{pairwise distinct} 	
\end{array}
\end{displaymath}
where $\sim_k$ denotes the identification of 
$(\be{ik}{ij}, \al{ik}{ij})$ with $(-\be{ik}{ij}, \al{ik}{ij}+\pi)$.
 
Any point $x \in \xah$ that is added to $\CONF$ by taking the closure has some
ratio $\be{ik}{ij}(x)$ equal to zero. We think of these points
as corresponding to degenerate configurations, that is,  configurations
that have some coinciding points.  We exploit these zero ratios
to define a series of {\it screens} corresponding to $x$.
A screen is just a copy of the plane $\R^2$. In Figure
\ref{fclick} for example,  there are three screens. We `fill'
these screens with the degenerate configuration. We can do it so 
that any two
points in the configuration are separated in at least one screen.
Moreover, the process is such that it shows how to write 
$\xah$ locally as the graph of a function. This proves that
$\xah$ is indeed a smooth manifold. 

The close connection between $\xah$ and the manifold
with corners $\FMt$ has two
immediate applications. First of all there is Theorem \ref{tfiber}
that reveals the corner structure of $\FMt$: the corners appear
automatically if we describe the natural map from $\FMt$ to
$\xah$.   

Moreover, we can apply this same method of defining and filling screens
to the manifold with corners $\FMt(n)$. That is, 
we add a family of screens to a degenerate configurations $x \in \FMt(n)$.
But $\FMt(n)$ contains all angles mod $2 \pi$ between two points.
Therefore we can, without ambiguities, define a clickable Voronoi diagram as
shown in Figure \ref{fclick} for every $x\in \FMt(n)$.
We conclude Chapter \ref{chmodel} by showing how to jump between several
descriptions and models in this thesis. 
Suppose we start with some $x \in \xah$, for example.
We can always construct a set $S(t)$ of polynomial sites 
$\{ p_1(t), \dots, p_n(t) \}$ with the following 
property: if we write down all angles and ratio's
for the sites in $S(t)$, we exactly obtain $x$.

\subsection*{Comparing prices: higher order Voronoi diagrams.}

In what preceded we have extended  Voronoi diagrams by limit 
Voronoi diagrams.  Another extension or generalization of the classic
Voronoi diagrams are the {\it higher order} Voronoi diagrams.
We discuss some results on higher order Voronoi diagrams in 
Chapter~\ref{chposet}.

People in the Netherlands like to compare prices. That is,
they have no fixed supermarket, but frequent the, say, five
closest supermarkets to pick the bargains. For those
people we have to divide the city into different sectors that
we call $5$-sectors.
In one $5$-sector, the five closest supermarkets are the same.
So, if people know which $5$-sector they live in, they know the five closest
supermarkets.  An example is given in Figure \ref{fvoronoi1}.
In Figure \ref{fvoronoi1}.a the 6 supermarkets of above are
displayed again. In Figure \ref{fvoronoi1}.b the four different
$5$-sectors are drawn. So, in sector $13456$, the five closets
supermarkets are supermarkets $1$, $3$, $4$, $5$ and $6$.
Or, if we put it the other way around: in $5$-sector
$13456$, supermarket $2$ is the most far away supermarket.
This diagram is the {\it fifth order} Voronoi diagram of
the six supermarkets.

\begin{figure}[!ht]
\begin{center}
\setlength{\unitlength}{1cm}
\begin{picture}(10,4.2)
\put(2,2.5){\makebox(0,0)[cc]{
        \leavevmode\epsfxsize=4\unitlength\epsfbox{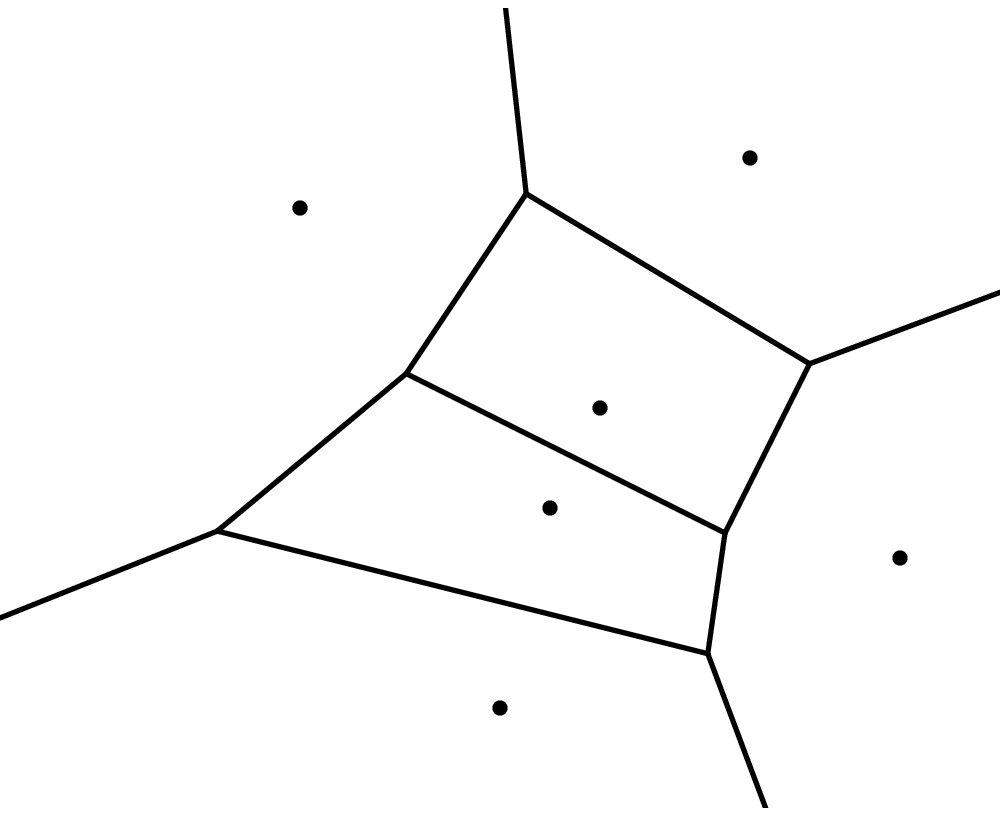}}}
\put(2.3,1.9){\makebox(0,0)[l]{1}}
\put(1.8,1.3){\makebox(0,0)[l]{2}}
\put(2.6,2.6){\makebox(0,0)[l]{3}}
\put( .98,3.4){\makebox(0,0)[l]{4}}
\put(3.1,3.7){\makebox(0,0)[l]{5}}
\put(3.7,2.1){\makebox(0,0)[l]{6}}
\put(8,2.5){\makebox(0,0)[cc]{
        \leavevmode\epsfxsize=4\unitlength\epsfbox{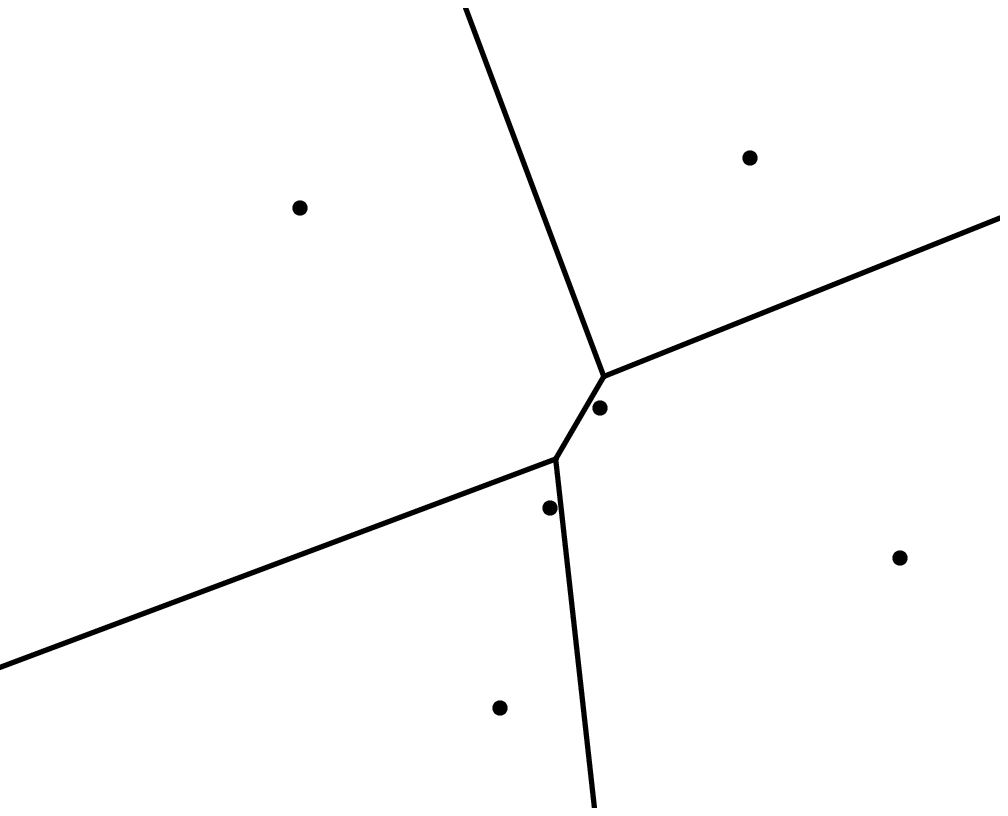}}}
\put(6.1,2.85){\makebox(0,0)[l]{12345}}
\put(9.2,4.0){\makebox(0,0)[l]{13456}}
\put(8.8,2.3){\makebox(0,0)[l]{12356}}
\put( 6.5,1.4){\makebox(0,0)[l]{12346}}
\put(2,0){\makebox(0,0)[l]{a}}
\put( 8,0){\makebox(0,0)[l]{b}}
\end{picture}
\caption{\elabel{fvoronoi1}
a.\ The Voronoi diagram of six supermarkets $1, \dots, 6$, and
b.\ the fifth order Voronoi diagram of the same supermarkets.}
\end{center}
\end{figure}

Note that in Figure \ref{fvoronoi1}.b only four distinct
$5$-sectors are present. This implies that two supermarkets
(which ones?) are nowhere the most far away supermarket.
We can change however the position of the supermarkets
in such a way that every supermarket is somewhere
the most far away supermarket. (How to do this?
Put the supermarkets  on the vertices of a convex
$6$-gon.)  This shows that if we change the position
of the supermarkets, the set of $5$-sectors may change as
well. In Chapter \ref{chlimit} we write down
all $5$-sectors. Moreover, we write down all
$4$-sectors, all $3$-sectors and all $2$-sectors. And we
add all $1$-sectors and all $6$-sectors, although they
are trivial (why?). If we collect everything 
we get a collection of subsets of $\{1,2,3,4,5,6\}$.
Or more general: we start with a set $S$ of
$n$ distinct points in the plane. The {\it Voronoi poset} $\Pi(S)$
consists exactly of those subsets of $\{1, \dots, n\}$
that have a non-empty Voronoi cell. 

 Although
the $5$-sectors may change if we change the position
of the supermarkets, there are certain invariants.
An invariant is just a number that is independent of the position
of the supermarkets.
In Theorem \ref{theuler} we prove for example the following:
suppose the number of supermarkets is odd. Then
the number of odd sectors equals the number of
even sectors. Or more formal:
for almost every configuration $S$ the total number of cells
in all even order Voronoi diagrams of $S$
equals the total  number of cells in all odd order Voronoi diagrams of $S$.

\chapter{Voronoi diagrams.}

In this chapter we give a short introduction
to Voronoi diagrams.
Most material covered in this chapter can be found in good books
on computational and discrete geometry. A very readable general introduction
to computational geometry is the book by De Berg, Kreveld,
Overmars and Schwarzkopf, see \cite{BKOS}. The book
by Edelsbrunner, see \cite{Ed}, is more in-depth, but less suited
as an introduction. Okabe, Boots and Sugihara wrote
a monograph on Voronoi diagrams for a broad audience, see \cite{OBS}. 
Most recent is the chapter on Voronoi diagrams written by
Aurenhammer and Klein  in the Handbook of Computational
Geometry (\cite{AK}). Some concepts that are used in the rest of this thesis are
shortly discussed in an overview paper on computational topology, cf.\ \cite{DEG}.
 
Apparently, Dirichlet in 1850 and  Voronoi in 1908, cf. \cite{Vo},
were the first that used a notion of Voronoi diagram. They
considered Voronoi diagrams of regular point sets, associated
to quadratic forms. Consult \cite{CF} for a recent and recommendable  book on Voronoi
diagrams of quadratic forms. 

\section{Convex hull.}
\elabel{sconvexhull}

A subset $A$ of the plane is \bfindex{convex} if for 
any two points $p, q \in A$ the line segment $pq$ is
contained in $A$ as well.  The \bfindex{convex hull}
$CH(A)$ of a set $A$ is the smallest convex set 
containing $A$. Any two non-coinciding points
$p=(p_x,p_y)$ and $q=(q_x,q_y)$ define two \bfindex{hull half-planes}
$hh_{pq}$ and $hh_{qp}$ where: 
\beq
	hh_{pq} & := & \{ r=(r_x,r_y) \in \R^2 \,\mid\,  
		\text{Det}(p,q,r) \,\geq\, 0 \}.
\eeq
The determinant $\text{Det}(p,q,r)$ is given by: 
\beq
\text{Det}(p,q,r) & =  &  \left| \begin{array}{ccc}
                1 & p_x & p_y \\
                1 & q_x & q_y \\
                1 & r_x & r_y
        \end{array}
        \right|.
\eeq
Let $S$ be a finite set of distinct points in $\R^2$.
We can write $CH(S)$ as an  intersection of hull half-planes.
More precisely,   $CH(S)$ equals  the intersection of those half-planes
defined by points in $S$ that contain all points of $S$
\beq
	CH(S) & = & \bigcap_{p_i, p_j \in S} 
		\{ h_{p_ip_j} ~:~ S \subset h_{p_ip_j} \}.
\eeq
As a non-empty intersection of half-planes, $CH(S)$ is a convex polygon. 
We can represent this polygon by listing its consecutive vertices
in clockwise order. In this way  $CH(S)$ corresponds to a cyclically ordered 
list of points from $S$. 
\begin{figure}[!ht]
\begin{center}
\setlength{\unitlength}{1em}
\begin{picture}(16,8)
\put(8,4){\makebox(0,0)[cc]{
        \leavevmode\epsfxsize=14\unitlength\epsfbox{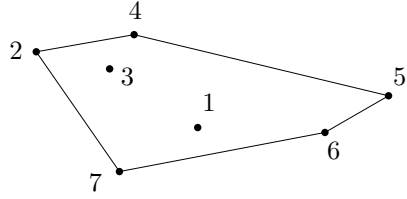}}}
\put(7.8,4){\makebox(0,0)[l]{1}}
\put(.5,6){\makebox(0,0)[l]{2}}
\put(4.7,5){\makebox(0,0)[l]{3}}
\put(5,7.5){\makebox(0,0)[l]{4}}
\put(15,5){\makebox(0,0)[l]{5}}
\put(12.5,2.2){\makebox(0,0)[l]{6}}
\put(3.5,1){\makebox(0,0)[l]{7}}
\end{picture}
\caption{\elabel{fhull}The convex hull of the points $1, \dots, 7$.}
\end{center}
\end{figure}
\begin{eexample} 
In Figure \ref{fhull}, $CH(S) \leftrightarrow \{2,4,5,6,7\}$.
\end{eexample}
Suppose the points in $S$ start moving around in the plane. Assume that
at the start of the motion, no three points in $S$ are collinear. 
A small enough disturbance of the points in $S$ does not change
the ordered list of points of $CH(S)$. Only when a point $p_k$ 
passes a convex hull edge $p_ip_j$, the list changes (simultaneously
$\text{Det}(p_i, p_j, p_k) = 0$). We call this change in $CH(S)$
a {\bf convex hull event}. \index{convex hull!event}

\section{Voronoi diagram.}
\elabel{svordia}

Start again with a set $S = \{ p_1, \dots, p_n \}$ of $n$ distinct points in the plane.
The \bfindex{Voronoi cell} $V(p_i)$ of a point $p_i \in S$ is defined as
\beq
	V(p_i) := \{ q \in \R^2: d(p_i, q) ~\leq d(p_j, q), ~i \neq j \}.
\eeq
Here, $d(p,q)$ denotes the ordinary Euclidean distance between $p$ and $q$.
Note that we define a Voronoi cell as a closed subset of
$\R^2$, in contrast to the choice made in \cite{Ed,BKOS}.
The \bfindex{Voronoi diagram} $V(S)$ of $S$ is the family of subsets of $\R^2$
consisting of the Voronoi cells and all of their intersections. 
The boundary of a Voronoi
cell consists of \bfindex{Voronoi edges} and 
\index{Voronoi vertex} {\bf Voronoi vertices}.
A point $q \in \R^2$ is on the Voronoi edge $e(p_i,p_j)$ if
$d(q,p_i) = d(q,p_j)$ and $d(q,p_k) \geq d(q,p_i)$ if $k \neq i,j$. 
A point $q \in \R^2$ is a Voronoi vertex if it is present on at least
two Voronoi edges.
A circle $C$ is an {\bf empty circle} \index{circle!empty} 
with respect to $S$ if there are no points of $S$
inside the circle.
For any three points $p_i$, $p_j$, and $p_k$, that are not 
collinear, there exists a unique circle $C_{ijk}$ passing through
$p_i$, $p_j$, and $p_k$. A circle $C_{ijk}$ is
a \bfindex{Voronoi circle} if it is an empty circle. 

\begin{figure}[!ht]
\begin{center}
\setlength{\unitlength}{1em}
\begin{picture}(16,10)
\put(8,5){\makebox(0,0)[cc]{
        \leavevmode\epsfxsize=14\unitlength\epsfbox{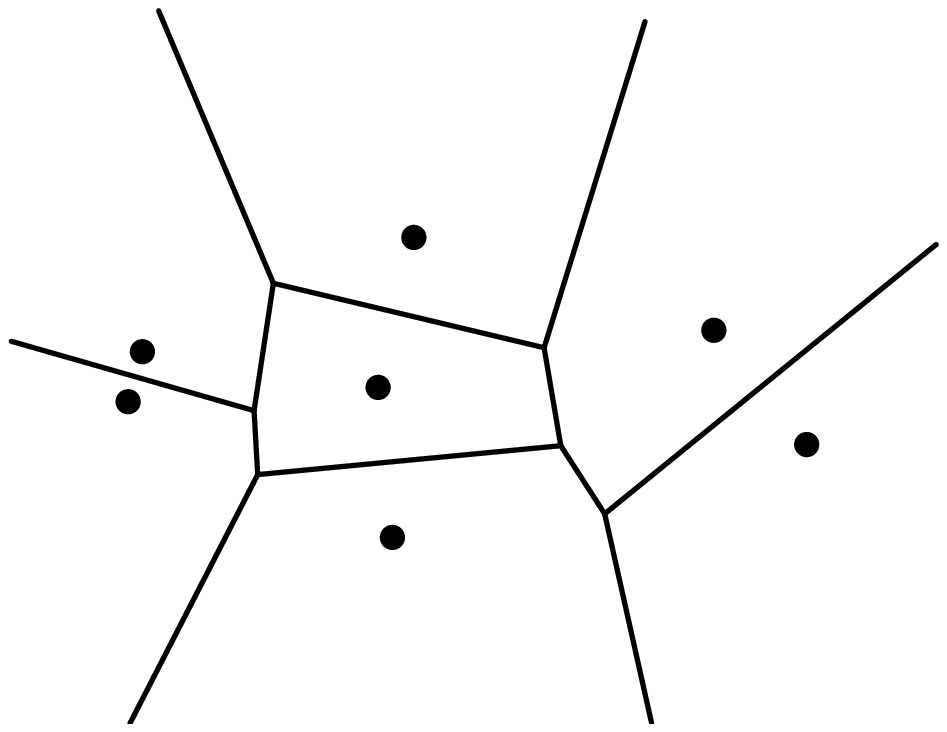}}}
\put(7,2){\makebox(0,0)[l]{1}}
\put(7.3,5){\makebox(0,0)[l]{2}}
\put(13.5,4){\makebox(0,0)[l]{3}}
\put(12,6){\makebox(0,0)[l]{4}}
\put(2.5,4){\makebox(0,0)[l]{5}}
\put(3,6){\makebox(0,0)[l]{6}}
\put(7,7.8){\makebox(0,0)[l]{7}}
\end{picture}
\caption{\elabel{fvoronoi}The Voronoi diagram of some points $1, \dots, 7$.}
\end{center}
\end{figure}

\begin{lemma}
$q \in \R^2$ is a Voronoi vertex if and only if $q$ is the center of
a Voronoi circle.
\end{lemma}

\begin{proof}
\cite{BKOS}, Theorem 7.4.
\end{proof}

The \bfindex{bisector} $B(p_i,p_j)$ of two points $p_i$ and $p_j$ is the line
equi-distant to $p_i$ and $p_j$. It is 
perpendicular to the line segment $p_ip_j$ passing
through $\frac{1}{2}(p_i + p_j)$.
A point $q$ is in the  \bfindex{Voronoi half-plane} $vh(p_i,p_j)$  if it is not closer to
$p_j$ than to $p_i$:
\beq
	vh(p_i, p_j) & = & \{q \in \R^2 \,\mid\, d(q,p_i) \leq d(q, p_j)\}.
\eeq
As a consequence,  the Voronoi half-plane $vh(p_i,p_j)$ is bounded 
by the bisector
$B(p_i,p_j)$.  Any Voronoi cell is an intersection of half-planes
\beq
	V(p_i) & = & \bigcap_{j \neq i} vh(p_i, p_j).
\eeq
Therefore, any Voronoi cell $V(p_i)$ is convex and is 
either bounded or unbounded.

\begin{lemma}
\elabel{lunbound}
$V(p_i)$ is unbounded if and only if $p_i \in \delta CH(S)$, where
$\delta CH(S)$ denotes the boundary of the convex hull $CH(S)$.
\end{lemma}

\begin{proof}
\cite{OBS}, Property V2.
\end{proof}

\begin{lemma}
\elabel{lbicenter}
Let $p_i$, $p_j$ and $p_k$ be three distinct points that are not collinear.
Let $c$ be the center of the circle through $p_i$, $p_j$ and $p_k$.
Then: 
\beq
	c & = & B(p_i,p_j) \cap B(p_i,p_k).
\eeq
\end{lemma}

\begin{proof}
The point $c$ is equi-distant to $p_i, p_j$ and $p_k$. 
\end{proof}

\section{Topological changes.}
\elabel{stopchange}

\begin{figure}[ht]
 \begin{center}
\setlength{\unitlength}{1cm}
\begin{picture}(5,4.2)
\put(2.5,2.25){\makebox(0,0)[cc]{
        \leavevmode\epsfxsize=5\unitlength\epsfbox{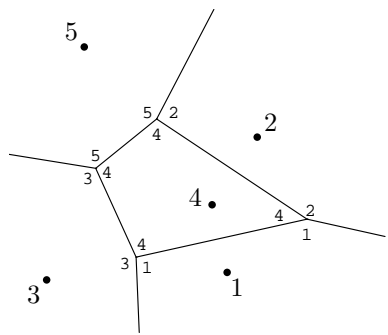}}}
\put(3.0,0.7){\makebox(0,0)[l]{1}}
\put(3.45,2.9){\makebox(0,0)[l]{2}}
\put(.3,.6){\makebox(0,0)[l]{3}}
\put(2.45,1.9){\makebox(0,0)[l]{4}}
\put(0.8,4.1){\makebox(0,0)[l]{5}}
\end{picture}
 \end{center}
 \caption{\label{typeeps}A diagram of type $\{134, 142, 245, 354\}$.}
\end{figure}

\begin{definition}
\elabel{dtype}
Let $V(S)$ be a Voronoi diagram. Represent every vertex $x$ of $V(S)$ 
as an ordered list of labels
of the points of $S$ on the Voronoi circle $C(x)$ of $x$.
The order of the labels corresponds to the cyclic, clockwise order
of the points on $C(x)$.
The set of all these lists is the \bfindex{type} of the Voronoi diagram $V(S)$.
\end{definition}

\begin{eexample}
\elabel{extype}
A diagram of type $\{134, 142, 245, 354\}$ is depicted in
Figure \ref{typeeps}.
\end{eexample}

A point set $S$ is in \bfindex{general position} iff no three points
are collinear and no four points are cocircular. If $S$ is in general
position, then every Voronoi circle has exactly three points
on its boundary. When the points in $S$ start to move around in the
plane, the type of $V(S)$ changes exactly when the configuration
of empty circles changes. This can happen generically in two ways.

\begin{figure}[!ht]
\begin{center}
\setlength{\unitlength}{1cm}
\begin{picture}(8,2.3)
\put(4,1.4){\makebox(0,0)[cc]{
        \leavevmode\epsfxsize=8\unitlength\epsfbox{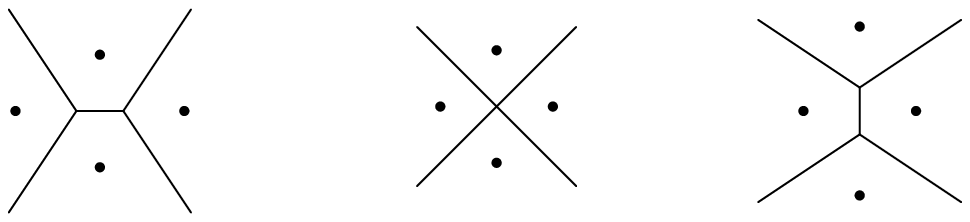}}}

\put(0,1.4){\makebox(0,0)[l]{a}}
\put(0.9,2.2){\makebox(0,0)[l]{b}}
\put(1.8,1.4){\makebox(0,0)[l]{c}}
\put(0.9,.6){\makebox(0,0)[l]{d}}
\put(3.3,1.4){\makebox(0,0)[l]{a}}
\put(4.1,2.2){\makebox(0,0)[l]{b}}
\put(4.9,1.4){\makebox(0,0)[l]{c}}
\put(4.1,.6){\makebox(0,0)[l]{d}}
\put(6.25,1.4){\makebox(0,0)[l]{a}}
\put(7,2.4){\makebox(0,0)[l]{b}}
\put(7.7,1.4){\makebox(0,0)[l]{c}}
\put(7,.4){\makebox(0,0)[l]{d}}
\end{picture}
 \caption{\elabel{fcircleevent}A circle event.}
\end{center}
\end{figure}

The first is when two empty circles $C_{abd}$ and $C_{bcd}$ coincide. This
      is a \index{circle!event} {\bf circle event}, see Figure 
      \ref{fcircleevent}. Before the event, 
      $a$, $b$, $c$, and $d$ define two empty circles $C_{abd}$
      and $C_{bcd}$. If  $c$ moves to the left in
      the leftmost figure, $a$, $b$, $c$ and $d$ become 
      cocircular. If $c$ continues moving left, one arrives
      in the situation of the rightmost figure, where
      $C_{abc}$ and $C_{acd}$ are the empty circles. The corresponding change in
      type is given by:
$$ \{\dots, abd, bcd, \dots \} \quad \rightarrow \quad \{\dots, abcd, \dots\} \quad \rightarrow \quad
	\{\dots, abc, acd, \dots \}. $$

\begin{figure}[!ht]
\begin{center}
\setlength{\unitlength}{1cm}
\begin{picture}(10,2.4)
\put(5,1.4){\makebox(0,0)[cc]{
        \leavevmode\epsfxsize=10\unitlength\epsfbox{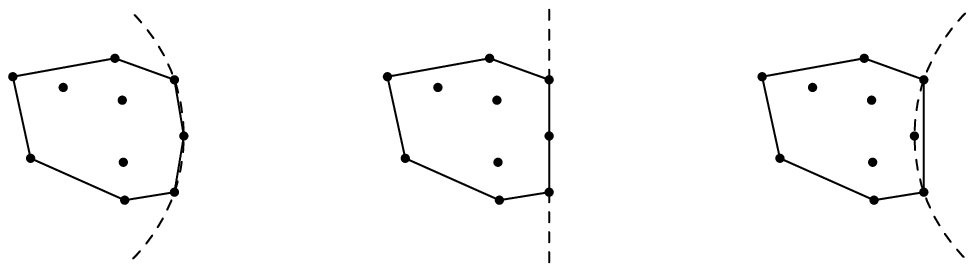}}}
\put(2.0,2.1){\makebox(0,0)[l]{a}}
\put(2.1,1.5){\makebox(0,0)[l]{b}}
\put(2.0,.8){\makebox(0,0)[l]{c}}
\put(5.8,2.1){\makebox(0,0)[l]{a}}
\put(5.8,1.5){\makebox(0,0)[l]{b}}
\put(5.8,.8){\makebox(0,0)[l]{c}}
\put(9.6,2.1){\makebox(0,0)[l]{a}}
\put(9.0,1.5){\makebox(0,0)[l]{b}}
\put(9.6,.8){\makebox(0,0)[l]{c}}
\end{picture}
 \caption{\elabel{fhullevent}A convex hull event.}
\end{center}
\end{figure}

The other way by which the type of a Voronoi diagram can change
is by means of a convex hull event, see Figure \ref{fhullevent}.
Consider the circle defined by the points $a,b$ and $c$.
In the figure on the left, $b \in \delta CH(S)$. The circle
$C_{abc}$ contains all other points of $S$ in its interior.
Suppose $b$ moves to the left. At some stage, $b$
passes through the line segment $ab$. At this moment,
the circle $C_{abc}$ swaps over, and becomes empty, as in the picture
on the right. This means that in  the type of $V(S)$, an extra term
$acb$ appears. 

\section{Delaunay triangulation.}

\begin{figure}[!ht]
\begin{center}
\setlength{\unitlength}{1em}
\begin{picture}(16,10)
\put(8,5){\makebox(0,0)[cc]{
        \leavevmode\epsfxsize=14\unitlength\epsfbox{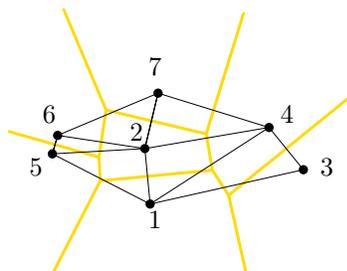}}}
\put(7,2){\makebox(0,0)[l]{1}}
\put(6.3,5.3){\makebox(0,0)[l]{2}}
\put(13.5,4){\makebox(0,0)[l]{3}}
\put(12,6){\makebox(0,0)[l]{4}}
\put(2.5,4){\makebox(0,0)[l]{5}}
\put(3,6){\makebox(0,0)[l]{6}}
\put(7,7.8){\makebox(0,0)[l]{7}}
\end{picture}
\caption{\elabel{fdelaunay}The Delaunay triangulation of the points
         introduced in Figure \ref{fvoronoi}.}
\end{center}
\end{figure}

Starting from  a Voronoi diagram $V(S)$ we define the \bfindex{Delaunay graph},
$DG(S)$, of $S$. The vertices of $DG(S)$ are exactly the points
in $S$. Two vertices $p_i$ and $p_j$  are connected by an edge
in $DG(S)$ exactly if there exists an  edge $e(p_i,p_j)$ of
positive length in the Voronoi diagram $V(S)$.
Let $D(S)$ be the straight line embedding of the Delaunay graph $DG(S)$.
For point sets  $S$ that are in general position, $D(S)$ 
triangulates the convex hull $CH(S)$. For this reason,
the straight line embedding $D(S)$ of the Delaunay graph $DG(S)$ is
called the \bfindex{Delaunay triangulation} of $S$. This terminology
is also used for point sets $S$ that are not in general position.
In the latter case $D(S)$ is not necessarily a triangulation.
An example of a Delaunay triangulation is presented in Figure \ref{fdelaunay}.

The Delaunay triangulation $D(S)$ is \bfindex{dual} to the
Voronoi diagram $V(S)$ in the following sense: vertices in the
Voronoi diagram correspond to faces in the Delaunay triangulation,
while Voronoi cells correspond with vertices of $D(S)$. 
As a consequence, the effect of a convex hull event  and
a circle event on the Delaunay triangulation $D(S)$ is easily found. 
At the  convex hull event
depicted  in Figure \ref{fhullevent} a triangle $T_{abc}$ appears
in the Delaunay triangulation $D(S)$.
Moreover, for any triple of points $p_i$, $p_j$, and $p_k$ 
in a point set $S$ it 
follows that $p_i$, $p_j$, and $p_k$ gives a triangle in
$D(S)$ if and only if the circle $C_{ijk}$ is a Voronoi circle.
Indeed, in the circle event of Figure \ref{fcircleevent}
an edge $bd$ of the Delaunay triangulation flips over in an edge $ac$.  
This proves the following lemma.

\begin{lemma}
\elabel{ltype}  Let $S$ be a set of distinct points in $\R^2$ 
in general position. The type of the Voronoi diagram $V(S)$ is
the list of triangles in the Delaunay triangulation $D(S)$.
\end{lemma}

The following characterization of Delaunay triangulations explains why
Delaunay triangulations are often used in the generation of altitude maps
and contour plots. The \bfindex{minimal angle} of a triangulation
is the smallest angle that occurs in any of the triangles of the
triangulation. In Figure \ref{fdelaunay}, the minimal angle is the
angle $\angle 526$.

\begin{lemma}
Let $S$ be a point set in general position.
The Delaunay triangulation $D(S)$ is that triangulation that maximizes
the minimal angle  over all triangulations of $S$.
\end{lemma}

\begin{proof}
\cite{BKOS}, Theorem 9.9.
\end{proof}

\section{Geometric transformations.}

\subsection{The lifting transformation.}
\elabel{lifting}

\begin{figure}[!ht]
\begin{center}
\setlength{\unitlength}{1em}
\begin{picture}(12,12)
\put(6,6){\makebox(0,0)[cc]{
        \leavevmode\epsfxsize=14\unitlength\epsfbox{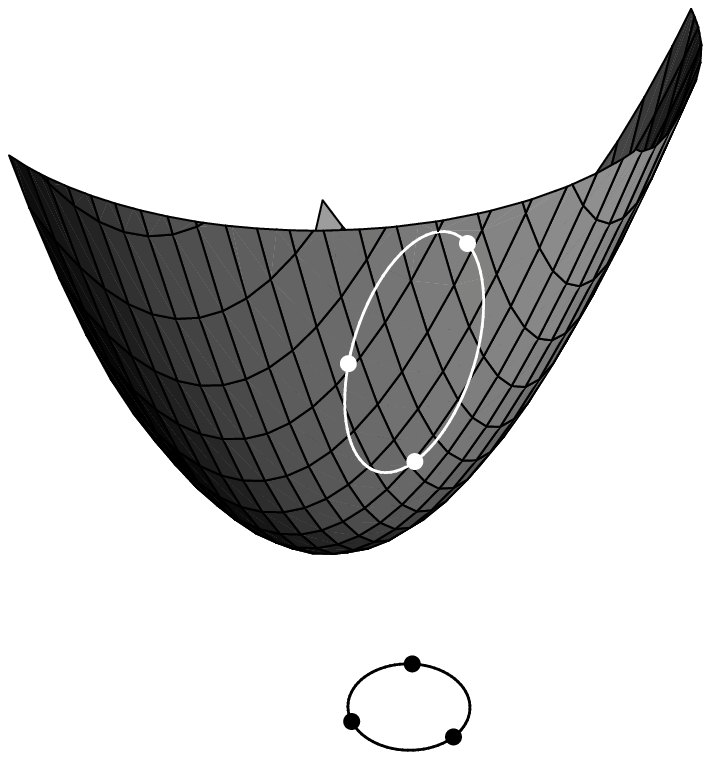}}}
\end{picture}
\caption{\elabel{fparaboloid}Lifting $(x,y)$ to $(x,y,x^2+y^2)$.}
\end{center}
\end{figure}

There is an easy correspondence between Delaunay triangulations of point
sets in $\R^2$
and convex hulls in $\R^3$: the \bfindex{lifting transformation}. 
This transformation is established by the map
\bmap
        \psi: & \R^2 & \rightarrow & \R^3, \\
        & (x,y) & \mapsto &  (x,y,x^2+y^2).
\emap
It maps a point in the plane onto a point of the paraboloid $P$ 
defined by $z=x^2+y^2$,
see Figure \ref{fparaboloid}.  The map $\psi$ has two properties that are important to us.

\begin{lemma}
\elabel{lparaboloid}
Let $C$ be a circle in the $xy$-plane.
\begin{lijst}
\item $\psi(C) = P \cap H_C$, where $H_C$ is a plane in $\R^3$.
\item A point $q$ is inside $C$ if and only if $\psi(q)$ below $H_C$.
\end{lijst}
\end{lemma}

\begin{proof}
The proof follows directly from combining the equation of a circle in the plane and the
equation of the paraboloid $P$.
\end{proof}

Let $S$ be a set of points in the $xy$-plane and $S'$ the set of 
images of the points in $S$ on
the paraboloid $P$. By Lemma \ref{lparaboloid}, $\psi$ maps empty circles 
defined by $S$ onto faces of the lower convex hull of $S'$. Therefore, the 
projection of the lower hull faces of $CH(S')$ onto the $xy$-planes 
gives exactly the Delaunay triangulation $D(S)$.

Let $S$ be a set of $n$ points in general position. 
A subset $A$ of $k$
points of $S$  is a $\mathbf{k}${\bf -set} \index{k-set@$k$-set} 
if it can be separated from the complementary set $B = S\setminus A$ 
of $n-k$ points by a plane $V_A$. Here we say that $V_A$ 
\bfindex{separates} $A$ from $B$ if $V_A$ can be oriented so that
all points in $A$ are on the positive side of $A$, while all points of $B$
are on the negative side.
  Lemma \ref{lparaboloid}
shows that
any set of $k$ points contained in a circle in $\R^2$ can be mapped
to a $k$-set in $\R^3$. 

\subsection{Mapping points to planes.}

A set of $n$ planes $\mathcal{V}$ defines a subdivision of $\R^3$
into connected pieces of dimension $0,1,2$ or $3$. This subdivision 
is the \bfindex{arrangement} $\mathcal{A}(\mathcal{V})$ 
of $\mathcal{V}$. Arrangements form a topic of their
own, cf.\ \cite{OT}.  
We assume that any point in any plane $V_i \in \mathcal{V}$ for
$i=1, \dots, n$  can be written as $(x,y, f_{{\mathcal{V}_i}}(x,y))$,
for some linear function $f_{\mathcal{V}_i}: \R^2 \rightarrow \R$.
That is, $\mathcal{V}$ consists of non-vertical planes only. 
We say that $p=(p_x,p_y,p_z) \in \R^3$ is \bfindex{above} plane $\mathcal{V}_i$ iff
$p_z > f_{\mathcal{V}_i}(p_x,p_y)$; similarly for \bfindex{below}. 
The $\mathbf{k}${\bf -level} \index{k-level@$k$-level} of the arrangement
$\mathcal{A}(\mathcal{V})$ consists of 
  those points in $\R^3$ above or in $k-1$ planes and 
below or in $n-k$ planes of $\mathcal{V}$. 
The $n$-level of an arrangement $\mathcal{A}(\mathcal{V})$ is also called
the \bfindex{upper envelope} of $\mathcal{A}(\mathcal{V})$.

\begin{figure}[!ht]
\begin{center}
\setlength{\unitlength}{1em}
\begin{picture}(12,10.4)
\put(6,5.3){\makebox(0,0)[cc]{
        \leavevmode\epsfxsize=12\unitlength\epsfbox{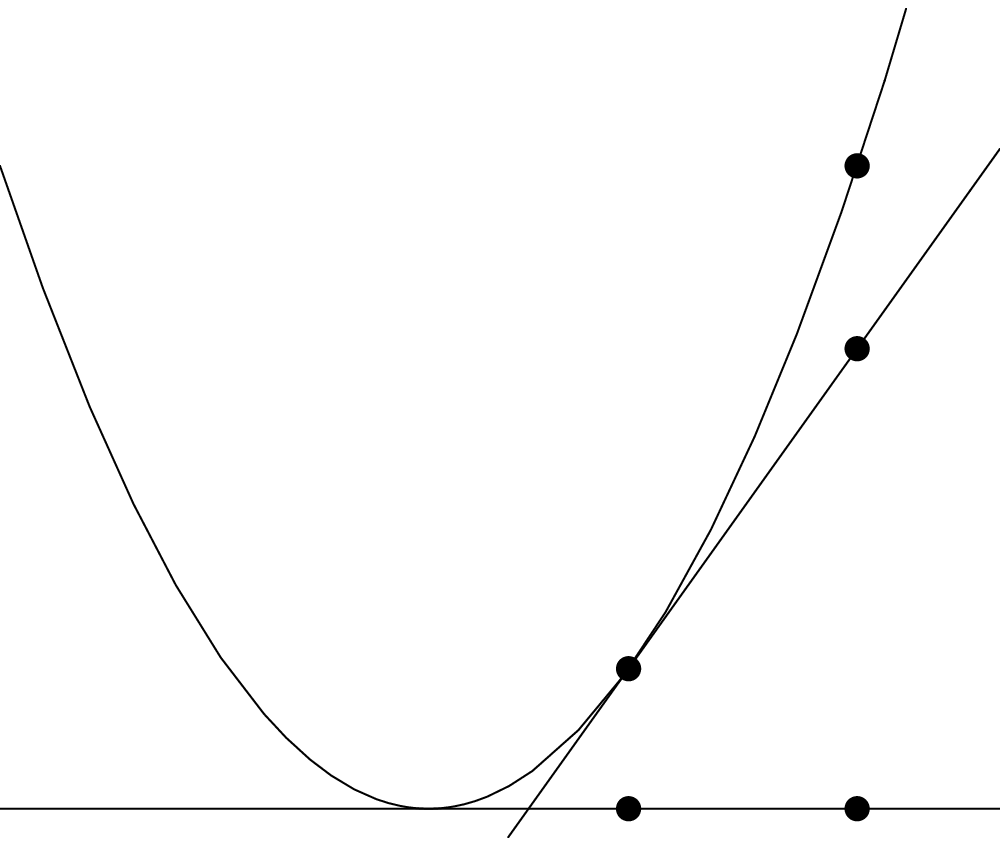}}}
 \put(6.8,0){\makebox(0,0)[l]{$p$}}
 \put(11,0){\makebox(0,0)[l]{$q$}}
\put(5.5,3){\makebox(0,0)[l]{$\psi(p)$}}
\put(11,6.5){\makebox(0,0)[l]{$h_p(q)$}}
\put(8,9){\makebox(0,0)[l]{$\psi(q)$}}
\end{picture}
\caption{\elabel{fd2para}Lifting points to planes.}
\end{center}
\end{figure}

Using the map $\psi$ again, we show that we can determine the Voronoi diagram
of a set $S$ of points  in $\R^2$ by computing the boundary  of an upper 
envelope in $\R^3$.  Consider 
Figure \ref{fd2para} that demonstrates the principle one dimension
lower.  To any point $p$ in a fixed point set $S$ we associate 
a plane $h_p$ in $\R^3$. That is, we map $p$ to the unique plane $h_p$
tangent to the paraboloid $P$ at the point $\psi(p)$. 
Lemma \ref{llift} shows that the set of planes $\{h_p\,|\,p \in S\}$ completely encodes
the relative distances of points $q \in \R^2$ to points in $S$.

\begin{lemma}
\elabel{llift}
 Let $q \in \R^2$ and let $h_p(q)$ be the intersection
of the vertical  line through $q$ and the plane $h_p$. Then
\beq
d(q,p)^2 & = & \psi(q) - h_p(q).
\eeq
\end{lemma} 

\begin{proof}
An easy computation or \cite{Ed}, Observation 13.3.
\end{proof}

Applying this result, we compute the Voronoi diagram $V(S)$ as follows.
Let $q \in \R^2$. Suppose that $q$ is in the interior $V(p_i)^{\circ}$ of
the Voronoi cell $V(p_i)$. 
By Lemma \ref{llift} this is equivalent to $h_{p_i}$ being the
first hyperplane that one encounters from the set
$\{h_{p_j}\,\mid\, p_j \in S\}$ if one goes downwards in negative vertical 
direction, starting from $\psi(q)$.  But this means that 
$V(S)$ is exactly the projection of the boundary of the upper envelope
of the arrangement $\{h_{p_j}\,\mid p_j \in S\}$ on $\R^2$.

We can connect the Delaunay triangulation $D(S)$
of a point set $S$ with the lower convex hull of a set $S'$.
Meanwhile the Voronoi diagram $V(S)$ can be obtained from the
upper envelope of an arrangement.
Remember that the Delaunay triangulation $D(S)$  is dual to  
the Voronoi diagram $V(S)$.
By now it may not come as  a surprise, that there exists a duality transfer 
from  lower convex hulls to upper envelopes and vice versa. For more
information consult \cite{Ed, BKOS}.

\chapter{A Voronoi poset.}
\elabel{chposet}

\hyphenation{pa-ra-me-ter mo-ving cor-res-pond}

Given a set $S$ of $n$ points in general position, we consider
all $k$-th order Voronoi diagrams on $S$, for $k=1,\dots,n$,
simultaneously. We recall symmetry relations for
the number of cells, number of vertices and number of
circles of certain orders.  We introduce a poset $\Pi(S)$  that 
consists of the $k$-th order Voronoi cells for all $k=1,\dots,n$,
that occur for some set $S$.  We prove that there exists a rank function
on $\Pi(S)$ and moreover that the number of elements of odd rank equals
the number of elements of even rank of $\Pi(S)$, provided that
$n$ is odd. 

\section{Introduction.}

The dynamics of Voronoi diagrams in the plane is well understood. When
$n-1$ points are fixed and one point is moving continuously  somewhere
inside the convex hull, combinatorial changes of the Voronoi diagram correspond
to changes in the configuration of empty circles, see Chapter~1 and 
\cite{AGMR} for example. Changes in the configuration of non-empty circles
correspond to combinatorial changes of higher order Voronoi 
diagrams. Here the $k$-th order Voronoi diagram associates to
each subset of size $k$ of generating sites that region in the plane that
consist of points closest to these $k$ sites.

We consider all $k$-th order Voronoi diagrams simultaneously for 
$k$ between $1$ and $n$. We do so by introducing the  
Voronoi poset $\Pi(S)$ of a set $S$ of $n$ distinct sites in the plane.
The poset consists of all sets of labels that correspond to 
a subset of sites that defines some non-empty Voronoi cell
in some $k$-th order Voronoi diagram.

Higher order Voronoi diagrams have been investigated by numerous 
people. Many results are published in an article by D.T.~Lee,
see \cite{Le}. A survey is given in Edelsbrunners book on
algorithms in combinatorial geometry, see \cite{Ed}. In
Chapter~1 we have introduced $k$-sets and we have shown that there
exists a map $\psi$ that changes the point inside circle
relation in $\R^2$ into a point below hyperplane relation in $\R^3$.
It turns out that these circles containing points serve
as a `building block' for higher order Voronoi diagrams as we 
discuss in full detail in Section \ref{skvoronoi}. As a consequence, 
formulas counting $k$-sets in $\R^3$ can be applied in the counting
of vertices, edges and cells of higher order Voronoi diagrams.
Instead of considering circles that contain a fixed number of, say, $k$ points,
one can also consider circles that contain at most $k$ points.
This is done in \cite{GHK}.

Let $T$ be a set of $n$ points
in $\R^3$ in general position that are the vertices of a convex polytope.
Sharir, \cite{Sh}, Lemma 4.4 and Clarkson and Shor, \cite{CS}, 
Theorem~3.5
prove that the number of $k$-sets of $T$ is given by
$ 2(k+1)(n-k-2)$.
They prove this formula using probabilistic methods that we do not discuss
here. 

 \begin{figure}[!ht]
 \begin{center}

 \leavevmode\epsfxsize=4.0cm\epsfbox{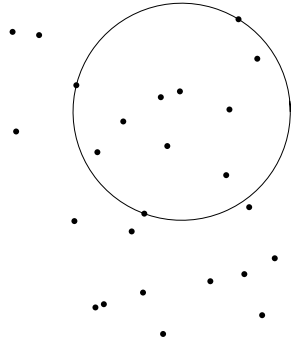}

 \end{center}
 \caption{An invariant for circle configurations. \elabel{cieps}}
 \end{figure}

This formula can also be derived in the context of $k$-th order
Voronoi diagrams from Lee's results as has been observed by
several people,  
see again Clarkson and Shor, \cite{CS} or Andrzejak et al.\ , \cite{AAHSW}.
We give this derivation explicitly and state  in Theorem \ref{ciplus}
that 
\beq
	c_i + c_{n-i-3} & =  &2(i+1)(n-i-2),
\eeq
where $c_i$ denotes the number of circles defined by a set $S$ of
$n$ points in general position in the plane, containing
exactly $i$ points of $S$. For an illustration, see Figure
\ref{cieps}.  Moreover we explicitly derive 
similar formulas for the number of cells $f_k$ in the $k$-th order
Voronoi diagram $V_k(S)$, see Lemma \ref{fiplus}, and the number 
of vertices $v_k$, see Lemma \ref{tildev}, in $V_k(S)$.
\beq
 f_k + f_{n-k+1} & = &  2k(n-k+1) +1 - n,	  \\
 v_k + v_{n-k} & =  & 4k(n-k).
\eeq

These `symmetry relations' are independent of the
particular position of the sites in $S$, provided $S$
is in general position: while the number 
of cells in some $k$-th order Voronoi diagram may change,
depending on the configuration, the sum of the   number of cells in the 
$k$-th order diagram and the number of cells in the $(n-k+1)$-th
diagram remains constant.

In Section \ref{svoronoiposet} we introduce the Voronoi poset
mentioned above and prove that $\Pi(S)$ has a rank function.
As an application  of the symmetry relations we prove in Theorem \ref{theuler}
that the number of elements of odd rank in $\Pi(S)$ equals 
the number of even rank, provided that $n$ is odd.

The Voronoi poset of a set $S$ of $n$ moving points seems a natural object 
to study as changes of the poset 
occur exactly at those moments where $S$ is not in general position. 
As there are tight connections between higher order  Voronoi diagrams,
$k$-levels in certain arrangements in $\R^3$ and certain $k$-sets in $\R^3$,
the study of the Voronoi poset may have applications in these areas
as well.

\section{Higher order Voronoi diagrams.}
\elabel{skvoronoi}

\subsection{Definition of $k$-th order Voronoi diagram.}

Let $S=\{p_1, \dots, p_n\}$ be a set of $n$ distinct points in $\R^2$ 
in general position.
Let $0 \leq k \leq n$. For
every point $p$ in the plane we ask for the $k$ nearest points from 
$S$. That is, we look for a subset $A \subset S$, such that
\begin{displaymath}
|A| = k, \qquad \forall x \in A, \quad  \forall y \in S-A : \quad d(p,x) 
	\quad \leq \quad d(p,y).
\end{displaymath}
For two points in $\R^2$, we define a half-plane
$$ h(x,y) ~:=~ \{ p \in \R^2 \mid d(x,p) \leq d(y,p) \}.$$
The \bfindex{Voronoi cell}
\index{Voronoi cell! order k} of $A \subset S$ of \bfindex{order} $|A|$
 is the intersection of half-planes
\beq
V(A) & := & \bigcap_{x \in A,~y\in S-A} h(x,y),
\eeq
whenever this intersection is not empty. As an intersection of
half-planes, $V(A)$ is a convex polygon. 

\begin{remark}
It is left as an exercise to the reader to show that assuming
general position implies that a Voronoi cell is not  a line
segment or a single point.
\end{remark}

The {\bf $\mathbf{k}$-th order Voronoi diagram}
\index{Voronoi diagram! order k}
 is the   subdivision of $\R^2$, induced
by the set of Voronoi cells of order $k$. For later purposes, we identify 
the $k$-th order Voronoi diagram with the set of non empty $k$-th order
Voronoi cells. 
\beq
	V_k(S) & := & \{V(A) \mid A \subset S, 
		~|A| = k, ~V(A) \neq \emptyset \}. 
\eeq
\begin{figure}[!ht]
\begin{center}
\setlength{\unitlength}{1cm}
\begin{picture}(12,3.4)
\put(6,2.0){\makebox(0,0)[cc]{
        \leavevmode\epsfxsize=12\unitlength\epsfbox{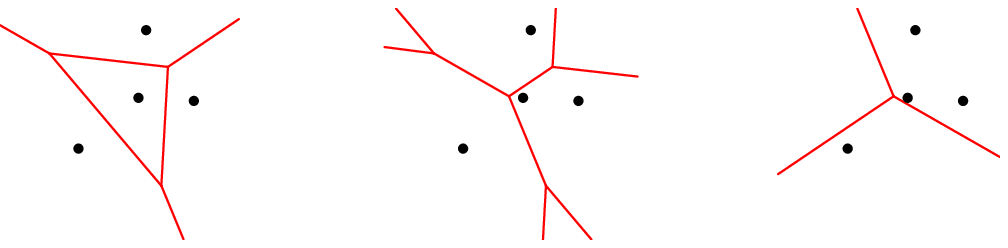}}}
\put(1.55,3.3){\makebox(0,0)[l]{$1$}}
\put(2.5,2.5){\makebox(0,0)[l]{$2$}}
\put(1.45,2.2){\makebox(0,0)[l]{$3$}}
\put(0.75,1.6){\makebox(0,0)[l]{$4$}}
\put(5.6,3.1){\makebox(0,0)[l]{$13$}}
\put(7.2,3.1){\makebox(0,0)[l]{$12$}}
\put(4.6,3.1){\makebox(0,0)[l]{$14$}}
\put(6.9,1.8){\makebox(0,0)[l]{$23$}}
\put(5.0,1.6){\makebox(0,0)[l]{$34$}}
\put(6.6,0.8){\makebox(0,0)[l]{$24$}}
\put(11.5,2.8){\makebox(0,0)[l]{$123$}}
\put(9.6,2.7){\makebox(0,0)[l]{$134$}}
\put(10.6,1.6){\makebox(0,0)[l]{$234$}}
\put(1.5,0){\makebox(0,0)[l]{a.}}
\put(6,0){\makebox(0,0)[l]{b.}}
\put(10.5,0){\makebox(0,0)[l]{c.}}
\end{picture}
 \caption{\elabel{fthreeinone}A first, second and third order
Voronoi diagram.}
\end{center}
\end{figure}

\begin{eexample}
Let $S=\{(45,86),(76,40),(40,42),(1,9)\}$. Figure \ref{fthreeinone}.a.\ 
shows the first order Voronoi diagram. The
second order diagram is shown in Figure \ref{fthreeinone}.b.\ and  the third order 
diagram in Figure \ref{fthreeinone}.c.\
All non-empty Voronoi cells are indicated by their
generators. Note that not all possible triples occur.
\end{eexample}

\begin{remark}
A planar graph that represents a point-face dual of the $k$-th order 
Voronoi diagram can be constructed as follows, cf. \cite{AS}. 
Write down for every $A \subset S$ with $|A|=k$ and $V_k(A) \neq \emptyset$
its \bfindex{centroid} $c(A)$, defined by $c(A)=(1/k)\sum_{p \in A} p$. 
Two centroids $C(A)$ and $c(B)$ are connected by an edge exactly iff
$V_k(A)$ and $V_k(B)$ share an edge. 

\end{remark}

\subsection{Circles and higher order Voronoi diagrams.}

In this section, we state some elementary properties of 
higher order Voronoi diagrams.
Every edge in $V_k(S)$ is part of some bisector $B(a,b)$, 
with $a,b \in S$. The Voronoi
vertices are exactly those points that are in the centers of the circles 
determined by
three points from $S$. Therefore, under our general position assumption,
every Voronoi vertex has valency three. The following theorem describes 
the local situation around a Voronoi vertex. The symbol
$\bigodot_{a,b,c}$ denotes the circle passing through the points
$a$, $b$, and $c$.

\begin{theorem}
\elabel{thcenter}
Let $x$ be the center of $\bigodot_{a,b,c}$, for $a,b,c \in S$, let 
\begin{eqnarray*} H & = & \{~z \in S~|~d(x,z) < d(x,a)~\},
\end{eqnarray*}
and let $k=|H|$. Then $x$ is a Voronoi vertex of $V_{k+1}(S)$
and $V_{k+2}(S)$. The Voronoi edges and cells that contain $x$ are
given in Figure \ref{abc}. Moreover, all Voronoi vertices are
of this form.

 \begin{figure}[!ht]
 \begin{center}

 \leavevmode\epsfxsize=9.6cm\epsfbox{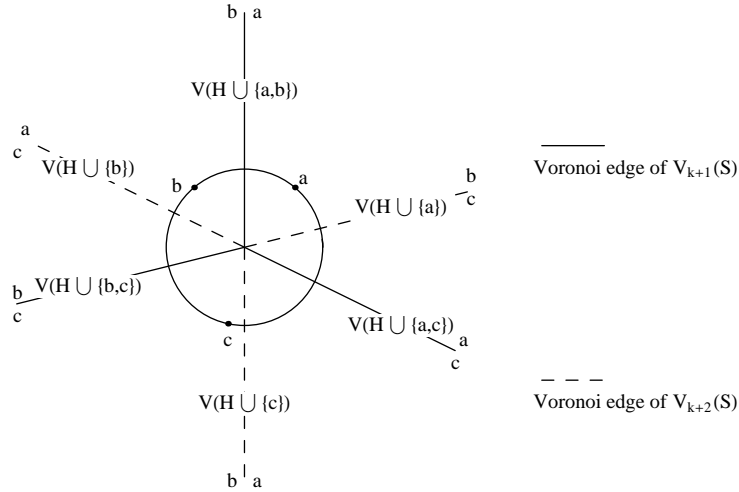}

 \end{center}
 \caption{The Voronoi diagram around $x$. \elabel{abc}}
 \end{figure}

\end{theorem}

\begin{proof} \cite{De}, Theorem 1 and Theorem 2. \end{proof}

Let $a,b,c$ and $H$ be as defined in Theorem \ref{thcenter}. We
define the \bfindex{order} \index{circle!order} of a circle  $\bigodot_{a,b,c}$ as $|H|$.
Notation: $|\bigodot_{a,b,c}| ~ :=~  |H|$. An {\bf order k Voronoi
circle} 
\index{Voronoi circle!order $k$}
$\bigodot_{a,b,c}$ is a circle through three points
$a,b$ and $c$ from $S$ that contains exactly $k$ points from $S-\{a,b,c\}$.
In fact, from all $\binom{n}{3}$ Voronoi circles $\bigodot_{a,b,c}$ 
and all sets $H_{a,b,c}$, compare  Theorem \ref{thcenter},  almost enough
information is provided to construct all $k$-th order Voronoi diagrams
$V_k(S)$ for $k=1, \dots,n-1$. 

\begin{alg}Voronoi diagrams of all orders.
\elabel{algkvoronoi}
\end{alg}

\begin{algorithmic}[1]
\REQUIRE set $S$ of $n$ points in general position.
\ENSURE  all $k$-th order Voronoi diagrams $V_k(S)$ for $k =1, \dots, n-1$.
\STATE Compute all circles $\bigodot_{a,b,c}$ defined by $S$. 
\STATE Compute all sets $H_{a,b,c}$ defined by $S$.
\STATE Take all circles of order $k-1$ and order $k-2$. The centers
	of these circles are exactly the vertices of $V_k(S)$.
\STATE Theorem \ref{thcenter} gives for every vertex the three incident
      edges and the three incident cells.
\STATE Two vertices are connected by an edge iff the two vertices
      have two incident cells in common.  Skip the edge if it
      is used.
\STATE Edges that are not skipped are unbounded edges. 
      Their direction and orientation still have to be computed.
      The direction is simply the direction of the bisector containing the edge.
      The orientation follows from Figure \ref{abc}. 
\end{algorithmic}

\begin{remark}
Denote the number of circles of order
$k$ by $c_k$ and the number of vertices in a $k$-th order Voronoi diagram by $v_k$. As a consequence of Theorem \ref{thcenter} we get
\begin{eqnarray}
\label{cv}
	v_k & = & c_{k-1} ~+~ c_{k-2}.
\end{eqnarray}
\end{remark}

\subsection{Counting vertices, edges and cells.}

The following theorem shows that the  total number of vertices,
edges and Voronoi cells does not depend
on the positions of the points in $S$, assuming general position.

\begin{theorem}
Let $v_k$, $e_k$, and $f_k$ denote the number of vertices, edges and
cells in $V_k(S)$ for some set $S$ of size $n$ in general 
position. The total number of vertices, edges and cells
in the Voronoi diagram of all orders are as follows.
\begin{lijst}
\item $\sum_{k=1}^n v_k$ $=$ $\frac{1}{3} n(n-1)(n-2)$.\vspace{.1cm}
\item $\sum_{k=1}^n e_k$ $=$ $\frac{1}{2} n(n-1)^2$.\vspace{.1cm}
\item $\sum_{k=1}^n f_k$ $=$ $\frac{1}{6} n(n^2+5)$.
\end{lijst}
\end{theorem}

\begin{proof} We prove the three claims.
\begin{lijst}
\item Every circle center defined by three distinct sites from
$S$ is a Voronoi vertex in some $k$-th and $(k+1)$-th order Voronoi
diagram. As there are $\binom{n}{3}$ distinct circles, the first 
claim follows.
\item Consider the arrangement of bisectors $\mathcal{A}(S)$.
Fix one bisector $B(a,b)$.
As $S$ is in general position, we may assume that the
bisector $B(a,b)$ is divided into $n-1$ line segments by the Voronoi circle
centers $ abx_3, abx_4, \dots, abx_n$,
where we write $S=\{a,b,x_3,\dots,x_n\}$. Every line segment is an edge in some
$k$-th order Voronoi diagram. As there are $\binom{n}{2}$ distinct
bisectors, claim {\it (ii)} follows.
\item The Euler  formula, $v_k -e_k +f_k = 1$, holds for every order.
Therefore
\beq
\sum_{k=1}^n f_k & = & n + \sum_{k=1}^n e_k - \sum_{k=1}^n v_k,
\eeq
which completes the proof.
\qedhere
\end{lijst}
\end{proof}

The number of vertices, edges and cells in $V_k(S)$ depends
on the configuration of $S$ as the ordinary Voronoi diagram shows.
The following theorem gives expressions for these numbers.
Let $f_k^{\infty}$ denote the
number of unbounded cells in the $k$-th order Voronoi diagram. By definition
$f_0^{\infty} := 0$.

\begin{theorem}
\label{euler}
\label{fk}
Let $S$ be in general position. Then the number of vertices, edges and cells
in the $k$-th order Voronoi diagram can be expressed as follows.
\begin{lijst}
\item  $v_k$  $=$  $2(f_k - 1) - f_k^{\infty}.$ 
\item  $e_k$  $=$  $3(f_k - 1) - f_k^{\infty}.$
\item  $f_k$  $=$  $(2 k-1)n - (k^2-1) - \sum_{i=1}^{k} f_{i-1}^{\infty}.$
\end{lijst}
\end{theorem}

\begin{proof}
\cite{Ed, Le}.
\end{proof}

Note that $f_n=1$. 
Substituting $k=n$ in  the expression for $f_k$ in 
Theorem \ref{fk} yields the following equation for the total
number of unbounded  cells: 
\begin{eqnarray}
\label{totalunbound}
	\sum_{i=1}^{n} f_{i-1}^{\infty} & = & n(n-1).
\end{eqnarray} 

The unbounded cells in the $k$-th order Voronoi diagram can be 
characterized as follows: let $\overline{pq}$ denote the
line segment with endpoints $p$ and $q$  and $l_{pq}$ the
line through $p$ and $q$.

\begin{property}
\elabel{unbounded}
 A cell $V(A)$ of the $k$-th order Voronoi diagram
$V_k(S)$ is unbounded if and only if one of the following two
conditions holds.
\begin{lijst}
	\item There exists a line $l$ that separates $A$ from $S-A$.
	\item There exist two consecutive points $p$ and $q$,
              with $p,q \in S-A$, on $\delta CH(S-A)$ such that the
              points in $A-\overline{pq}$ are in the open half plane 
	      defined by $l_{pq}$ opposite to $CH(S-A)$.
\end{lijst}
\end{property}

\begin{proof} \cite{OBS}, Property OK4.
\end{proof}

Under the general position assumption, we only need to
consider condition {\it (i)} in Property \ref{unbounded}. It is clear that
in this case the following symmetry holds:
\begin{equation}
\label{dual}
	f_k^{\infty} \eis f_{n-k}^{\infty}.
\end{equation}

\subsection{Circle events and hull events.}

In Section \ref{stopchange} we have discussed the type of a Voronoi diagram in
connection with circle events and convex hull events. 
We generalize these notions to $k$-th order Voronoi diagrams.
Let $S$ be a set of points in general position and fix some $k \in 1,\dots,n$.
As in \cite{Le}  we call a Voronoi vertex of $V_k(S)$ \bfindex{old}
if its corresponding circle has order $k-2$ and 
\bfindex{new} if its corresponding circle has order $k-1$.
We will also use the words old and new to indicate order $k-2$ and
order $k-1$ circles with respect to $V_k(S)$.

\begin{eexample}
In the classic or  first order diagram, vertices correspond with empty circles, 
so all vertices are new.
\end{eexample}

Represent every vertex $x$ of $V_k(S)$  as an ordered list of labels
of the points of $S$ on the $k-1$ or $k-2$ order Voronoi circle
$C(x)$ corresponding to $x$. Order the labels as in Definition \ref{dtype}.
The \bfindex{type} of a $k$-th order Voronoi diagram $V_k(S)$
consists of two lists: a list of ordered labels of old vertices
and a list of ordered labels of new vertices.

\begin{eexample}
 Consider the configuration in Figure \ref{n=4}.a. 
The circles $124$ and $143$ are empty, while the circles $123$ and $243$
both contain one point. Therefore, the type of the second order
Voronoi diagram is given by: $\{\{124,143\},\{123,243\}\}$.
\end{eexample}

Suppose that the points in $S$ start moving.
$S$ is in general position, as long a no four points are cocircular and
no three points are collinear.
If four points become cocircular, then two circles $C_{abd}$ and $C_{bcd}$, 
both surrounding  a subset
$T \subset S$ of size $|T|=k$, coincide and change into two circles
$C_{abc}$ and $C_{acd}$, see Figure \ref{fcircleevent}. 
This corresponds to a change of the Voronoi vertices $abd$ and $bcd$
into Voronoi vertices $abc$ and $acd$ in both $V_{k+1}(S)$ and $V_{k+2}(S)$. 
This is the generalization 
of a circle event and is called a {\bf $\mathbf{k}$-th order circle event}.
\index{circle!event!order k@order $k$} So, a zero order
circle event is just a circle event.

\begin{figure}[!ht]
\begin{center}
\setlength{\unitlength}{1cm}
\begin{picture}(11.5,2.7)
\put(5,1.5){\makebox(0,0)[cc]{
        \leavevmode\epsfxsize=10\unitlength\epsfbox{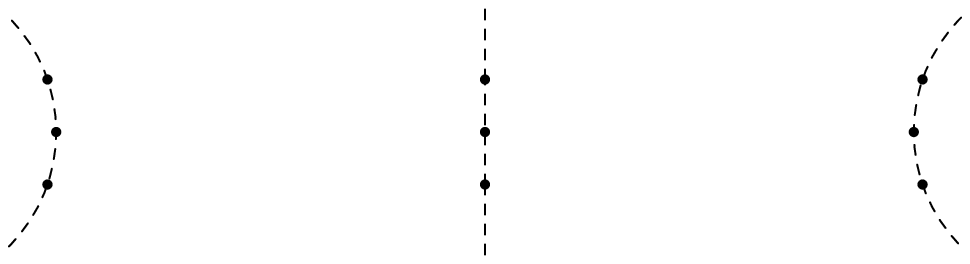}}}
\put(0,1.5){\makebox(0,0)[l]{$\mathbf{k}$}}
\put(0.8,2.1){\makebox(0,0)[l]{a}}
\put(0.9,1.6){\makebox(0,0)[l]{b}}
\put(0.8,1){\makebox(0,0)[l]{c}}
\put(1.55,1.5){\makebox(0,0)[l]{$\mathbf{n-k-3}$}}
\put(4.5,1.5){\makebox(0,0)[l]{$\mathbf{k}$}}
\put(5.2,2.1){\makebox(0,0)[l]{a}}
\put(5.2,1.6){\makebox(0,0)[l]{b}}
\put(5.2,1){\makebox(0,0)[l]{c}}
\put(5.8,1.5){\makebox(0,0)[l]{$\mathbf{n-k-3}$}}
\put(8.6,1.5){\makebox(0,0)[l]{$\mathbf{k}$}}
\put(9.6,2.1){\makebox(0,0)[l]{a}}
\put(9.0,1.6){\makebox(0,0)[l]{b}}
\put(9.6,1){\makebox(0,0)[l]{c}}
\put(9.7,1.5){\makebox(0,0)[l]{$\mathbf{n-k-3}$}}
\end{picture}
 \caption{\elabel{fhigherhull}A $k$-th order convex hull event.}
\end{center}
\end{figure}

If three points become collinear, an order $k$ circle $abc$ changes into
an order $n-k-3$ circle $abc$, see Figure \ref{fhigherhull}.  This implies
that the vertex $abc$ disappears in $V_{k+1}(S)$ and $V_{k+2}(S)$ 
and appears in $V_{n-k-1}(S)$ and $V_{n-k-2}(S)$. This is the generalization
of a hull event and is called a {\bf $\mathbf{k}$-th order hull event}.
\index{hull event!order k@order $k$} So, a zero order
hull event is just a hull event.
At a $k$-th order circle event the type of $V_{k+1}(S)$ and $V_{k+2}(S)$
changes. At a $k$-th order hull event, the type of
$V_{k+1}(S), V_{k+2}(S), V_{n-k-1}(S)$ and $V_{n-k-2}(S)$  changes.

\section{The Voronoi poset.}
\elabel{svoronoiposet}

\subsection{Definition and examples.}

Fix a labeling of the sites in $S$  and identify
a  set of sites $A \subset S$ that defines a non-empty  Voronoi cell $V(A)$
with the set of labels $L(A) \subset [n]$ of the sites  in $A$.
A subset $L$ of $[n]$ may or may not correspond to some
Voronoi cell $V(A_L)$.
For $k=1$ we retain the ordinary Voronoi diagram, which implies
the correspondence
\begin{displaymath}
 V_1(S)  ~\leftrightarrow~ \{ \{1\}, \{2\}, \dots, \{n\} \}.
\end{displaymath}

We define $V_0(S) = \{ \emptyset \}$. The set $\{ \{1, \dots, n\} \}$
corresponds to $V_n(S)$.
We consider the set of all Voronoi cells that appear for
a  given set $S$ of points and  call the set of corresponding labels the
\bfindex{Voronoi poset} $\Pi(S)$ of $S$:
\begin{displaymath}
        \Pi(S) ~:=~ \bigcup_k~\{~L(A)~|~V(A) \in V_k(S)~\}.
\end{displaymath}

This definition also makes sense when we drop the general position
assumption.

We order the
elements in the poset by set inclusion of the sets $L(A)$. This yields a
partially ordered set. For more on partially ordered sets consult \cite{Zi}.
 The poset is bounded since we have the empty set as
$\hat{0}$, the unique minimal element, and the set $[n]$ as $\hat{1}$, the
unique maximal element. In general, a poset is called \bfindex{graded}  if it is bounded
and if every maximal chain has equal length. We show that $\Pi(S)$ is graded.
Below we give an example showing that $\Pi(S)$  is in general  not
a lattice.

\begin{property}
$\Pi(S)$ is graded.
\end{property}

\begin{proof}
We show that $r(L(A)) = |L(A)|$ is a \bfindex{rank function}
for $\Pi(S)$. A rank function maps an element $x$ from a poset to a unique
level in such a way that the level corresponds to the length of any maximal 
chain from $x$ to $\hat{0}$.
 Let $L(A) \in \Pi(S)$, with $|L(A)| = k$. Every
point $x \in V(A)$ has all $k$ points from $A$
as its $k$ nearest neighbors. Order those points with respect to
their distance to $x$.  As we assumed general position it is always possible
to change the choice of $x$ in such a way that this order is strict.
By removing at each step  the furthest point still available,
we obtain a chain of length $k$ that descends to $\hat{0}$. 
\end{proof}

We analyze the two smallest cases, assuming general position.

\begin{eexample}
For $n=3$ there is one poset, the full poset on $[3]$. That is,
$$\Pi_3(S)=\{\emptyset, 1,2,3,12,13,23,123\}.$$
\end{eexample}

\begin{eexample}
\elabel{extwofour}
For $n=4$ there are two essentially distinct posets as is evident from
the circles defined by four points. Since $n=4$, Voronoi
circles are of order one or two. Eulers formula rules out
four circles of order one. At the same time less than two circles
would not yield enough cells in the first order diagram.
\end{eexample}

\begin{figure}[!ht]
\begin{center}
\setlength{\unitlength}{1cm}
\begin{picture}(10.5,3.8)
\put(2,2){\makebox(0,0)[cc]{
        \leavevmode\epsfxsize=3.5\unitlength\epsfbox{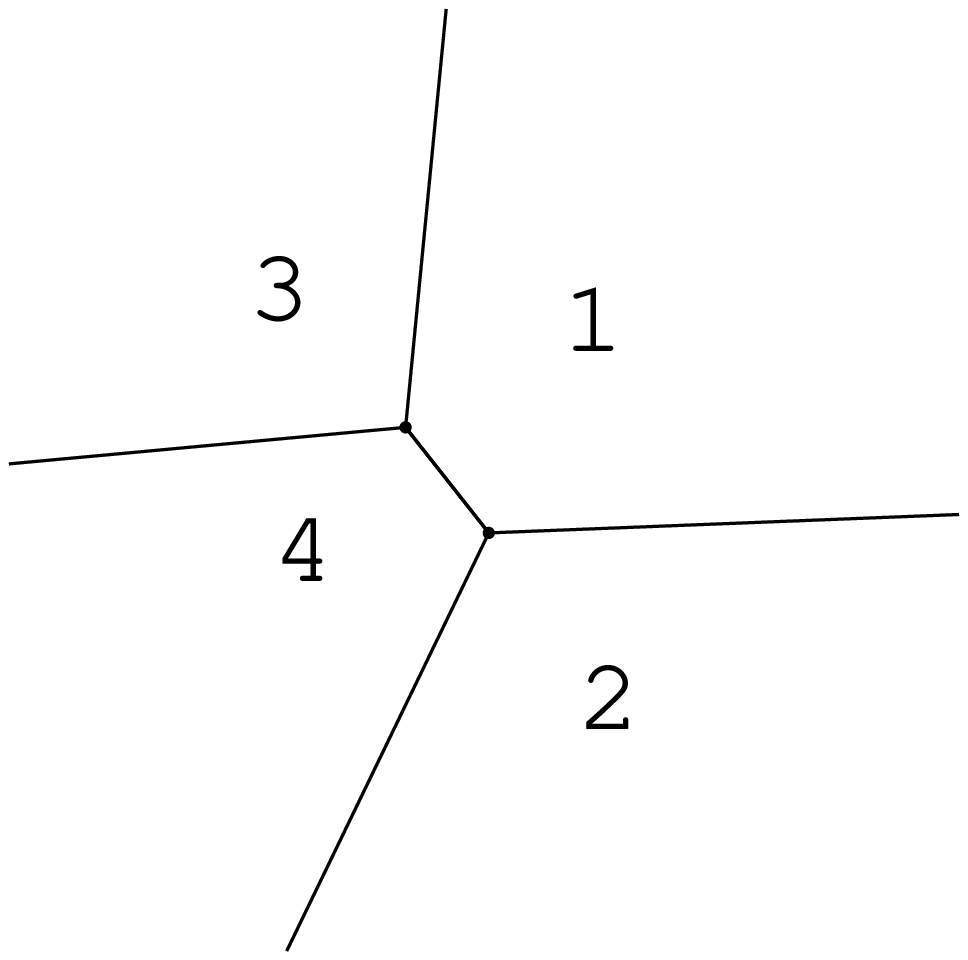}}}
\put(8.5,2){\makebox(0,0)[cc]{
        \leavevmode\epsfxsize=3.5\unitlength\epsfbox{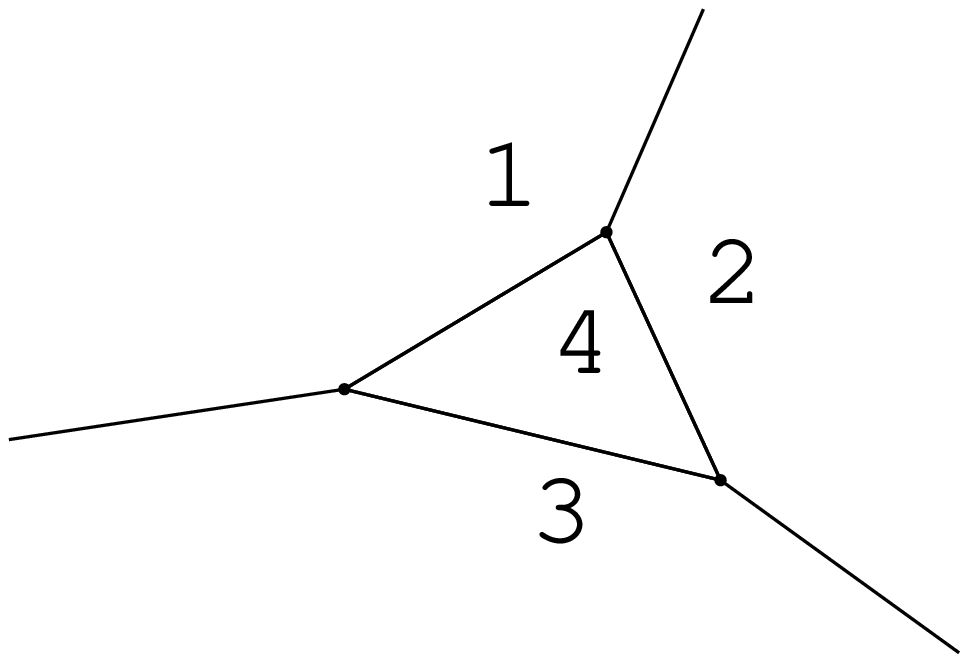}}}
\put(2.0,0){\makebox(0,0)[l]{a.}}
\put(8.5,0){\makebox(0,0)[l]{b.}}
\end{picture}
\caption{\elabel{n=4}The two distinct first order Voronoi diagrams.}
\end{center}
\end{figure}

Figure \ref{n=4}.a has two empty Voronoi circles. Clearly, subset
$23$ is the only subset of $[4]$ missing, leaving us with :

$$\Pi_{4}(S_1)=\{\emptyset, 1,2,3,4,12,13,14,24,34,123,124,134,234,1234\}.$$

This example shows that the Voronoi poset is in general not a lattice.
A lattice requires that every two elements of the poset have a
unique minimal upper bound. In this example, the elements $2$ and $3$ have two 
minimal upper bounds, namely $123$ and $234$.

In Figure \ref{n=4}.b.\ there are three empty Voronoi circles. 
The cell $123$ cannot appear
in the third order diagram, but all other subsets of $[4]$ do appear, thus:

$$\Pi_{4}(S_2)=\{\emptyset, 1,2,3,4,12,13,14,23,24,34,124,134,234,1234\}.$$

\subsection{The order complex of the Voronoi poset.}

The standard way to associate a topological space to a finite poset $(P,\leq)$
is by means of the \bfindex{order complex} $\Delta(P)$ of the poset, see
\cite{Bj,Wa}.
The order complex is the simplicial complex of all nonempty chains 
of $P$. A \bfindex{chain} of $P$ of \bfindex{length} $k$
 is a totally ordered subset
$$ x_0 ~<~ x_1 ~<~ x_2 ~<~ \dots ~<~ x_k, $$ 
of elements $x \in P$. The well-known \bfindex{geometric realization}
associates a topological space with a simplicial complex.

As a  Voronoi poset on a set of  $n$ points $p_1, \dots, p_n$ always has a 
unique maximal element $\{1, \dots,n \}$, the geometric realization
of the order complex is a cone and therefore contractible. This
shows that the topological space that we have associated with
$S$ is homotopy equivalent with a point, and therefore not
very interesting.

More promising is to consider the complement, that is
the \bfindex{anti Voronoi poset} $aP(S)$, consisting of those 
subsets of $\{1, \dots, n\}$ that are not
in the Voronoi poset. Another possibility is to consider
the arrangement of bisectors.

\section{Symmetry relations.}

Given a set $S$ of sites, we count for every order $k$ the number
of vertices $v_k$, the number of edges $e_k$ and the number of
non empty Voronoi cells $f_k$. The \index{f-vector@$f$-vector}
 $\mathbf{f}${\bf -vector} of $\Pi(S)$ is
the vector $\{f_1, f_2, \dots, f_n\}$. The $c$- and 
$e$-vector are defined analogously.  Note that the $f$-vector
of $\Pi(S)$ may change if the position of the sites in $S$ changes.

\begin{eexample}
Consider the two configurations $S_1$ and $S_2$ on four points, 
presented in Example \ref{extwofour}. The $f$-vector of $\Pi(S_1)$
equals $(1,4,5,4,1)$, while the $f$-vector of $\Pi(S_2)$
equals $(1,4,6,3,1)$.
\end{eexample}

\subsection{Symmetry in the number of cells.}

It turns out that a symmetry exists in the $f$-vectors.

\begin{lemma}
\label{fiplus}
Consider the $f$-vector of $\Pi(S)$, where $|S|=n$.
Then $f_{k}+f_{n-k+1}$ is a constant independent of the position
of the points in $S$.  More precisely,
\begin{eqnarray}
\label{si}
	f_{k}+f_{n-k+1} & = & 2k(n-k+1) + 1 - n.
\end{eqnarray}
\end{lemma}

\begin{proof} We apply Theorem \ref{fk} to $f_{k}$ and $f_{n-k+1}$:
\begin{eqnarray*}
f_k + f_{n-k+1} & = & 
	(2k-1)n - k^2 + 1 - \sum_{i=1}^{k} f_{i-1}^{\infty} \\
&& + \quad (2(n-k+1)-1)n - (n-k+1)^2 +1  - 
	\sum_{i=1}^{n-k+1} f_{i-1}^{\infty}, \\
& = & 2kn - 2k^2 + 2k + 1 - n + n(n-1) - (~\sum_{i=1}^{k} f_{i-1}^{\infty}
		+ \sum_{i=1}^{n-k+1} f_{i-1}^{\infty}~).
\end{eqnarray*}
We join the two sums by applying Symmetry Equation \ref{dual}.\
and evaluate the result by using Equation \ref{totalunbound}.\
\beq
\sum_{i=1}^{k} f_{i-1}^{\infty}~+~\sum_{i=1}^{n-k+1} f_{i-1}^{\infty} 
\eis  \sum_{i=1}^{n} f_{i-1}^{\infty} \eis n(n-1).
\eeq
The lemma follows from combining the two equations above. \end{proof}

\subsection{Symmetry in the number of vertices.}

A similar equation holds for the number of vertices of a collection
of Voronoi diagrams $V_k(S)$, for  $k=1,\dots,n-1$.

\begin{lemma}
\label{tildev} Let $S$ be a  set of points in general position with $|S|=n$.
Let $v_k$ denote the number of vertices in the $k$-th order
Voronoi diagram. Then:
\begin{eqnarray}
	v_k + v_{n-k} & = & 4k(n-k) - 2n.
\end{eqnarray}
\end{lemma}

\begin{proof}
Using Theorem \ref{euler} we write $v_k+v_{n-k}$ in terms of
numbers of cells. Next we regroup and apply Symmetry Equation \ref{dual}.
After applying Theorem \ref{fk} we combine using symmetry again. Finally,
using $\sum_{i=1}^n f_{i-1}^{\infty}=n(n-1)$  completes the proof.
\begin{align*}
v_k + v_{n-k} &=  2(f_k-1) - f_k^{\infty} + 2(f_{n-k}-1) - f_{n-k}^{\infty},  \\
              &=  2 ( f_k + f_{n-k} -2 - f_k^{\infty} ), \\
              &=  2 ( n^2 -2 n + 2 kn - 2 k^2 -
                      ( \sum_{i=1}^k f_{i-1}^{\infty} + \sum_{i=1}^{n-k} f_{i-1}^{\infty} 
                      + f_k^{\infty} ) ), \\
	      &=  2 ( n^2 - 2n + 2k n -2 k^2 - \sum_{i=1}^n f_{i-1}^{\infty} ), \\
              &=  -2 n + 4 kn - 4 k^ 2, \\
              &=  4k(n-k)-2n. \qedhere
\end{align*}
\end{proof}

\subsection{Symmetry in the number of Voronoi circles.}

Recall that the order of a Voronoi circle equals the number
of points of $S$ contained in its interior. We define
the  $\mathbf{c}${\bf -vector}
\index{c-vector@$c$-vector} of $S$ as the vector
\beq
 c(S) & =  & \{ c_0, c_1, \dots, c_{n-3} \}, 
\eeq
where $c_i$ denotes the number of circles of order $i$. The following
theorem states that for $n$ arbitrary points in general position, the
number of circles containing exactly  $i$ points on their inside plus
the number of circles containing exactly $i$ points on their outside is
 constant.  We prove this by applying the above results.

\begin{theorem}
\label{ciplus}
Consider the $c$-vector of $\Pi(S)$, where $|S|=n$.
Then $c_i+c_{n-i-3}$ is a constant independent of the position
of the points in $S$. 
More precisely,
\begin{eqnarray}
\label{ci}
        c_i+c_{n-i-3} & = & 2(i+1)(n-2-i), \\\nonumber
                      & = & 2i(n-i-3) + 2(n-2).
\end{eqnarray}
\end{theorem}

\begin{proof} We prove the theorem by induction.

\noindent{\bf [i=0]. }
\hyphenation{di-men-sio-nal}
We use the \bfindex{lifting transformation}. This transformation 
changes the point-inside-circle
relation in 2-dimensional space in a point-below-plane relation
in 3-dimensional space. See also Section \ref{lifting}.
The lifting transformation map $\phi$
is defined by
\bmap
	\phi: &  \R^2 &  \rightarrow &  \R^3, \\
	& (x,y)  & \mapsto & (x,y,x^2+y^2).
\emap
It lifts points in the plane to the unit paraboloid in three-space.
As every circle defined by $S$ in the plane contains 
only three points from $S$, every hyperplane defined by $\phi(S)$
contains only three points from $\phi(S)$ as well.
The number $c_0$  of empty circles of $S$ in the plane equals the number
of facets of the lower hull of $\phi(S)$ in three dimensions.
At the same time, the number $c_{n-3}$ of
circles that contain all other points of $S$ equals the number of facets 
of the upper hull of $\phi(S)$. All images of points in  $S$
under $\phi$ are part of the convex hull of $\phi(S)$.  Since the convex hull
of a point set of $n$ points consists of $2n-4$ facets, if every
facet is a triangle,  see \cite{BKOS}, Theorem 11.1, the claim follows.

\noindent{\bf [induction step]. }
We deduce the expression for $c_k + c_{n-k-3}$ by applying
Equation \ref{cv}, followed by combining  Lemma \ref{tildev} 
and the induction hypothesis:
\begin{align*}
c_k + c_{n-k-3} &= c_{k-1} + c_k + c_{n-k-3} + 
	c_{n-k-2} - (c_{k-1} + c_{n-k-2}),\\
                &= v_{k+1} + v_{n-(k+1)} - (c_{k-1} + c_{n-k-2}), \\
		&= 2(2(k+1)-1)(n-(k+1))-2(k+1)\\
                &  \qquad -(2(k-1+1)(n-2-(k-1))),\\
		&= 2(k+1)(n-2-k).
\qedhere
\end{align*} 
\end{proof}

Let $\tilde{f_k} := f_{k}+f_{n-k+1}$ and $\tilde{c_i} := c_i+c_{n-i-3}$.
By the {\bf  reduced $\mathbf{f}$-vector},
\index{f-vector@$f$-vector!reduced} denoted $\tilde{f}$, we mean the vector
 of $\tilde{f_k}$'s for all distinct $k$.
That is
\begin{eqnarray*}
       \tilde{f} \quad := \quad  \{ \tilde{f_0}, \tilde{f_1}, \dots,
             \tilde{f}_{\lfloor\frac{n-1}{2}\rfloor} \}.
\end{eqnarray*}
$\tilde{c}$ is defined similarly.  As a consequence of Lemma \ref{fiplus}
and  Theorem \ref{ciplus},
$\tilde{f}$ and  $\tilde{c}$ are only dependent on $n$.
 
\begin{eexample}
As an example we present the reduced $f$- and $c$-vectors for $n\in\{3,\dots,12\}$.
\end{eexample}

\begin{center}
\begin{tabular}{rll}\medskip
n & $\tilde{f}$ & $\tilde{c}$\\
3 & (4, 6)      &   (2)   \\
4 & (5, 9)      &   (4)  \\
5 & (6, 12, 14)   & (6, 8)  \\
6 & (7, 15, 19)         & (8, 12) \\
7 & (8, 18, 24, 26)     & (10, 16, 18)  \\
8 & (9, 21, 29, 33)     & (12, 20, 24)  \\
9 & (10, 24, 34, 40, 42)        &  (14, 24, 30, 32) \\
10 & (11, 27, 39, 47, 51)       &   (16, 28, 36, 40) \\
11 & (12, 30, 44, 54, 60, 62)   &   (18, 32, 42, 48, 50) \\
12 & (13, 33, 49, 61, 69, 73)   &   (20, 36, 48, 56, 60)
\end{tabular}
\end{center}

\begin{remark}
Computer calculations did not suggest any similar symmetry relation
for the number of edges.
\end{remark}

\subsection{Relations between cells and circles.}

\begin{corollary}
$\tilde{f_i}$= $\tilde{f_0} + \tilde{c}_{i-1}$ = $\tilde{c}_{i-1} + n +1$.
\end{corollary}

\begin{proof}
This follows directly from Lemma \ref{fiplus} and Theorem \ref{ciplus}.
\end{proof}

\begin{property}
Let $f_i^{\infty}$ denote the number of unbounded cells in the $i$-th order diagram
and let $c_i$ denote the number of circles of order $i$:
\begin{eqnarray}
\label{fiinfci}
	f_i^{\infty} + (c_{i-1}-c_{i-2}) & = & 2(n-i).
\end{eqnarray}
\end{property}

\begin{proof}
We prove the property by induction.

\noindent{\bf [i=1] }
$c_{-1}$ is zero by definition.
The number of vertices $v_1$ in the first order Voronoi diagram equals
the number of circles of order zero, $c_0$. The claim follows
from applying Theorem \ref{fk}:
\begin{align*}
	f_1^{\infty} + ( c_0 - c_{-1} ) &=  f_1^{\infty} + v_1, \\
	                   &=  f_1^{\infty} + 2(f_1-1) - f_1^{\infty}, \\				   &=  2(n-1).
\end{align*}

\noindent{\bf [induction step] }
Assume we have proved that
\begin{eqnarray*}
	f_i^{\infty} + (c_{i-1} - c_{i-2}) & = & 2(n-i). 
\end{eqnarray*}
We rewrite this, by using induction, as
\begin{eqnarray}
\label{ci-}
	c_{i-1} & = & 2ni-i(i+1)-\sum_{k=1}^{i+1} f_{k-1}.
\end{eqnarray}
Evaluate $c_i-c_{i-1}$:
\begin{align}
\label{cici-}
c_i-c_{i-1} &= ( c_i + c_{i-1} ) -2 c_{i-1}, \\\nonumber
            &= v_{i+1} - 2 c_{i-1}, \\\nonumber
	    &= 2( f_{i+1} - 1) - f_{i+1}^{\infty} -2 c_{i-1}.
\end{align}
Substituting this expression for $c_i-c_{i-1}$ and applying Theorem 
\ref{fk} and Equation \ref{ci-} proves the claim:
\begin{align*}
 f_{i+1}^{\infty} + (c_i - c_{i-1}) & ~=~  2( f_{i+1}-1-c_{i-1}), \\
                                    & ~=~  2(n-i-1).
\qedhere
\end{align*}
\end{proof}

\begin{corollary}
\label{cdetf}
The $c$-vector totally determines the $f$-vector. The correspondence 
is given by
\begin{eqnarray*}
	f_k & = & n - k + 1 + c_{k-2}.
\end{eqnarray*}
\end{corollary}

\begin{proof}
Applying Equation \ref{fiinfci} we get
\begin{eqnarray*}
	\sum_{i=1}^{k} f_{i-1}^{\infty} & = & (k-1)(2n-k) - c_{k-2}.
\end{eqnarray*}
The claim follows from evaluating Theorem \ref{fk} using the expression
above.
\end{proof}

\section{Even versus odd order cells.}

Given a grading on a set of objects, it is common to consider the 
\bfindex{Poincar\'e polynomial} $P(t)$ of the grading. The $i$-th 
coefficient of this polynomial equals the number of objects
of grade $i$. In our case, the objects are the elements of
the Voronoi poset $\Pi(S)$, while the grading is given by the 
rank function on the poset. Recall that the rank of an element $x$
in $\Pi(S)$ is just the order $k$ of the Voronoi diagram in
which $x$ occurs as a cell. The $i$-th coefficient of the
Poincar\'e polynomial $P(t)$ is given by $f_i$,
as $f_i$ gives the number of cells in the $i$-th order diagram $V_i(S)$.
So, the Poincar\'e polynomial $P(t)$ of $\Pi(S)$ with respect to 
our rank function is given by
\beq
	P(t) & = & f_{0} + f_1 t + f_2 t^2 + \dots + f_n t^n.
\eeq

As an application of the symmetry relations we compare the number of
cells in the even order Voronoi diagrams with the number of cells 
in the odd order diagrams. In terms of the Poincar\'e polynomial
$P(t)$ of above, the following result can also be formulated as
\beq
	P(-1) & = & 0.
\eeq

\begin{theorem}
\elabel{theuler}  
Let $S$ be a set of points in general position with $|S|=n \geq 3$.
Assume $n$ is odd.
In this case, the number of cells in the even order Voronoi diagrams equals
the number of cells in the odd order Voronoi diagrams.
\end{theorem}

\begin{proof}
Write $\tilde{f}_i ~=~ f_{i} + f_{n-i+1}$. We show that $A=0$, where:
\beq
	A & = & -f_0 + f_1 - f_2 + \dots - f_{n-1} + f_n.
\eeq
So $A$ is the number of cells in the odd order diagrams minus the number of cells
in the even order diagrams:
\beq
	A & = & -f_0 + \tilde{f}_1 + \frac{1}{2} 
		\tilde{f}_{\frac{n+1}{2}} + t_n,
\eeq
where
\beq
 t_n &  := & \sum_{i=2}^{\frac{n-1}{2}} (-1)^{i+1} \tilde{f}_i. 
\eeq
Clearly, $f_0 = 1$, as $f_0$ counts the empty set. $\tilde{f}_1$ is the number
of points in $S$ plus the number of cells in $V_n(S)$, so $\tilde{f}_1 = n+1$. 
Applying Equation \ref{si} gives:
\beq
  \tilde{f}_{\frac{n+1}{2}} & =  & -(-1)^{\frac{n+1}{2}} \frac{n^2+3}{4} .
\eeq
Straightforward calculations show that:
\beq
  t_n & = &  (-1)^{\frac{n+1}{2}} \frac{n^2+3}{4} - n,
\eeq
from which it follows that:
\begin{equation*}
A  \eis  -1+n+1-(-1)^{\frac{n+1}{2}} \frac{n^2+3}{4}
    +(-1)^{\frac{n+1}{2}} \frac{n^2+3}{4}-n \eis 0. \qedhere 
\end{equation*}
 \end{proof}

The claim of Theorem \ref{theuler} does not hold when $n$ is even. However,
the following result does hold.

\begin{lemma}
Let $S$ be a set of points in general position, with $|S|=n \geq 3$.
Assume $n$ is even. Let $A(S)$ denote the number of cells in the odd
order Voronoi diagrams minus the number of cells in the even order diagram.
Then:
\begin{eqnarray*}
        n \equiv 0(4) & \Rightarrow & A(S) ~{\rm odd}. \\
        n \equiv 2(4)  & \Rightarrow & A(S) ~{\rm even}.
\end{eqnarray*}
\end{lemma}

\begin{proof}
Similar computations as in the proof of Theorem \ref{theuler}.
\end{proof}

Note that as $v_k = c_{k-1} + c_{k-2}$ it follows immediately that:
\begin{eqnarray*}
        \sum_{k=1}^{n-1} (-1)^{k+1}v_k \eis 0,
\end{eqnarray*}
for all $n$, where $v_k$ denotes the number of vertices in the $k$-th order
Voronoi diagram.

\newcommand{\pr}{\mathrm{PosR}}
\newcommand{\dpr}{\mathrm{DPosR}}
\newcommand{\strip}{\mathrm{Strip}}
\hyphenation{pa-ra-me-ter mo-ving cor-res-pond}

\chapter{Limits of Voronoi diagrams.}
\elabel{chlimit}

The Voronoi diagram of a set $S$ of $n$ distinct points in $\R^2$ associates to
a point $p \in S$ that part of the plane that is closer to $p$ than to any other
point in $S$. In this chapter we assume that the position of any of the $n$ points in
$S=S(t)$ is given by a pair of polynomials in one parameter $t$, such that no two
points are represented by the same pair of polynomials. We call the elements of
$S(t)$ polynomial sites. In this setup it is possible that sites coincide at,
say, $t=0$. 
We define a Voronoi diagram $V(S(0)) := \lim_{~t \downarrow 0} V(S(t))$.
That is, $V(S(0))$ is defined as a limit diagram in a particular sense
of  the ordinary Voronoi diagram of the positions, at small positive $t$,
of the polynomial sites.
We show how to extend the notion of type to polynomial sites.
This enables us to determine the combinatorics of $V(S(0))$.
It turns out that in general some sites can be omitted without changing
the boundary of the Voronoi diagram $V(S(0))$. 
We characterize those sites that can not be omitted
and present an efficient 
algorithm to determine these sites together with the boundary of the 
Voronoi diagram.

\section{Introduction.}

Whenever Voronoi diagrams are studied, it is assumed that all points
defining the diagram are distinct.  We call these points sites.
In the case of dynamic 
Voronoi diagrams, where the sites are moving continuously over time,
the assumption above  means that sites are not allowed to coincide 
at any moment. In this chapter we investigate in one particular setting
what happens when we do allow sites to coincide in the plane. 

We consider a set $$S(t)  =  \{p_1(t), \dots, p_n(t)\},$$ of $n$ sites in 
the plane such that
the position of site $p_i$ at time $t$ is given by a pair of
polynomials in $t$, one for every coordinate. 
That is, both the movement of the  $x$- and the $y$-component
are described by a polynomial in $\R[t]$. 

In fact we do not have to restrict ourselves to sites
described by pairs of polynomials: the theory developed in this
chapter works whenever the movement of both the $x$- and $y$-component
of all sites is described by functions that can be expanded
as a convergent Taylor series. For the sake of simplicity we restrict ourselves
to polynomials however. 

In Chapter~1 we have described the Voronoi
diagram of a set of distinct sites both in terms of half-planes and
of empty circles. To check whether some point $q$ is in a given half-plane
$h(p_1,p_2)$ or inside a circle $c(p_1,p_2,p_3)$ for $p_1,p_2,p_3 \in \R^2$
boils down to evaluating the sign of an easy polynomial expression 
$f_h(p_1,p_2,q)$ or $f_c(p_1,p_2,p_3,q)$. 
The main idea of this chapter is to replace $p_i$ for $i =1,2,3$ and $q$ by
$p_i(t),q(t) \in \R[t]^2$ and to evaluate the sign of the 
coefficient of the lowest degree term in 
$f_h(p_1(t),p_2(t),q(t))$ and $f_c(p_1(t),p_2(t),p_3(t),q(t))$,
at $t=0$.
In this way we define a Voronoi diagram $V(S(t))$ at $t=0$,  even if sites do 
coincide at $t=0$. For a set of sites that coincide at $t=0$
we  regard $V(S(0))$ as the limit diagram
$\lim_{\epsilon \downarrow 0}V(S(\epsilon))$. We develop
this arithmetic for polynomial sites in all details in Section
\ref{sprelim}.

Here we add a warning: we can define a limit Voronoi diagram at $t=0$ that is 
consistent with small positive $t$ or small negative $t$. In general these
 two approaches give distinct results. Examples will be presented in the main
 text of this chapter.

There is a similarity between our approach and a well-known 
technique in computational geometry that is used to
avoid computations involving degenerate input. In this technique,
a set of distinct points that is not in general position is perturbed
slightly in such a way that it is in general position after perturbation.
For example, three points that are collinear before perturbation
will not be collinear anymore after perturbation.  As a consequence
it is possible for a {\it generic} algorithm, that is, an
algorithm that can only handle points in general position, 
to handle the input and compute the wished geometric structure, e.g.\ a Voronoi
diagram. Such generic algorithms are in general  much easier than
algorithms for arbitrary input as a lot of degenerated cases can be
ignored. In the so-called {\it Simulation of Simplicity} technique,
cf.\ \cite{EM, Ed2}, polynomials in the variable $\epsilon$ are added
to points of a set of distinct, degenerated points  in such a way  that 
evaluation of the new points for small enough $\epsilon$ produces
a point set that is in general position. An overview
on robust geometric computation is given in \cite{Ya}.

In Section \ref{stype} we extend  the
definition of type, see Definition \ref{dtype}, to
a set of polynomial sites $S(t)$ at $t=0$, provided that $S(0)$ 
fulfills some general position assumptions.  
This type gives us a complete combinatorial structure that matches with
the combinatorial structure of the Voronoi diagram $V(S(t))$ for small
enough positive~$t$.

In Section \ref{svorpoly} we define a Voronoi diagram for $n$ not
necessarily distinct points in the plane and $\binom{n}{2}$ angles
between those points. We apply this definition to introduce 
Voronoi diagrams for a set of polynomial sites $S(t)$ at $t=0$.
The resulting  Voronoi diagram or polynomial sites 
matches with the ordinary Voronoi diagram of the positions of the sites 
at small positive $t$.
In this section we also give the connection between half-planes and Voronoi 
circles defined by $S(t)$ at $t=0$, thereby connecting the notion of type
to polynomial sites diagrams.

The shape of a  Voronoi diagram is defined as the union of the boundaries
of the Voronoi cells.  
It turns out that in general some
sites in a set $S(t)$ of polynomial sites can be omitted without changing
the shape of the Voronoi diagram of $S(t)$. A question that we pose is how
to determine this shape efficiently  and how to characterize
those sites that do determine the shape.
We answer this question by splitting up the problem into two parts. 

First we consider in Section \ref{sposcell}
the shape of the  Voronoi  diagram of one cluster of polynomial sites at $t=0$.
By an $l$-cluster  we mean  a set of polynomial sites
such that the positions of all sites   coincide at one location $l$
at $t=0$.
Lemma \ref{lposarea} fully classifies the sites $p_i(t)$ in the cluster
such that $\text{area}(V(p_i(0)))>0$ at $t=0$. 

The second part of the question is solved by Lemma \ref{lplug}.
This lemma  states that the shape of the Voronoi diagram at $t=0$ of 
an arbitrary set of polynomial sites can be found as follows.
First, compute the ordinary Voronoi diagram of all distinct 
locations at $t=0$.   Second, plug in the cell of location $l$,
the shape of the Voronoi diagram of the $l$-cluster.

Section \ref{sposedge}  shows that a lot of combinatorics can be hidden
in the edges of the shape of the  Voronoi diagram. In Section \ref{ex20poly} 
we demonstrate the theory developed in the former section by a 
somewhat bigger example. We conclude this chapter in Section \ref{spolyk} with some 
remarks on generalizations  to $k$-th order Voronoi diagrams and shortcomings
of this setup.

\section{Preliminaries.}
\elabel{sprelim}

We define half-planes and circles
for polynomial sites. A \bfindex{polynomial site}
$p(t) = (p(t)_x, p(t)_y)$ consists of a pair of
polynomials $p(t)_x, p(t)_y \in \R[t]$.
Two polynomial sites $u(t)$ and $v(t)$ are called
\bfindex{distinct} if they represent distinct elements
in $\R[t] \times \R[t]$. Throughout this section assume that
$S(t)   =  \{ p_1(t), \dots, p_n(t) \}$,
where every $p_i(t)$ a polynomial site
and $p_i(t)$ distinct from  $p_j(t)$ as a polynomial site 
whenever $i \neq j$.

\begin{remark}
The notions in this section are introduced so that they
match the situation for ordinary points in the plane
obtained by substituting a very small 
 {\it positive} value of $t$ in the set $S(t)$ of polynomial sites.
\end{remark}

\subsection{Polynomial lines and their directions.}

Let $u(t)$ and $v(t)$ be two distinct polynomial sites.
\index{polynomial!line}
The {\bf polynomial line} $l_{uv}(t)$ is defined by
\begin{displaymath}
	l_{uv}(t) \quad := \quad   \left|
        \begin{array}{ccc}
                1       &       u(t)_x        &       u(t)_y \\
                1       &       v(t)_x        &       v(t)_y \\
                1       &       x             &       y
        \end{array}
\right| \quad \in \quad \R[t,x,y].
\end{displaymath}
The \bfindex{ruling coefficient} $\text{rc}(f(t))$ of
a polynomial $f(t) \in \R[t]$ is the coefficient
in its lowest degree term. 
The \bfindex{ruling sign} $\text{rs}(f(t))$ is the sign 
of the ruling coefficient $\text{rc}(f(t))$.
We define the \bfindex{direction} $\phi_{uv}$ of the polynomial 
line $l_{uv}(t)$ at $t=0$ as 
the argument of the point
\beq
(\cos \phi_{uv}, \sin \phi_{uv}) & = & \lim_{t\downarrow 0}(\frac{d(t)_x}{\mid d(t)_x \mid},
	\frac{d(t)_y}{\mid d(t)_y \mid})
\eeq
where $d(t) = (d(t)_x, d(t)_y) = v(t) - u(t)$. Note that $\phi_{uv}$ is determined up to
multiples of $2 \pi$. We often use that value of $\phi_{uv}$ that lives in  $(-\pi,\pi]$.

\begin{eexample}
\elabel{exdirection}
Let $u(t)= ( t, t)$ and $v(t)=(- t, t^2)$. Then
$l_{uv}(t) = t^2 + t^3 + t x - t^2 x - 2 t y$.
Putting $l_{uv}(t)$ equal to zero yields 
an ordinary line for every $t \neq 0$. As
$d(t) = v(t) - u(t) = (-2 t, -t + t^2)$, we have
$\lim_{t \to 0} \frac{d(t)_y}{d(t)_x}=\frac{1}{2}$. The
direction of $l_{uv}(t)$ at $t=0$ is given by
$\text{arctan}(\frac{1}{2}) - \pi$, however, as the ruling sign 
$\text{rs}(d(t)_x)$ equals $-1$. 
\end{eexample}

\begin{remark}
Note that this definition of $\phi_{uv}$ indeed  matches with the direction of the directed 
line that passes first through $u(t)$ and then through $v(t)$ for small positive $t$.
For a definition that would match with small negative $t$,
we should take into account the odness or evenness of the 
power of the lowest degree term of $l_{uv}(t)$. 
\end{remark}

\begin{eexample}
Let $u$ and $v$ be as in Example \ref{exdirection}. For small negative $t$,
the direction of the line $l_{uv}$ that passes first through $u(t)$ and 
then through $v(t)$ is close to $\text{arctan}(\frac{1}{2})$. If we multiply
both $u(t)$ and $v(t)$ by $t$, that is change $u(t)$ into $u(t)=(t^2,t^2)$
and $v(t)$ into $v(t)=(-t^2,t^3)$, then the direction of $l_{uv}$ for small
negative $t$ and for small positive $t$  both are close to
$\text{arctan}(\frac{1}{2}) - \pi$. 
\end{eexample}

\subsection{Collinearity.}

Let $u(t), v(t)$ and $w(t)$ be three polynomial sites.
As long as the sites do not coincide we 
can analyze whether they are collinear. Consider 
the determinant $D(t)$ given by
\begin{displaymath}
D(t) \quad  = \quad   D_{u,v,w}(t) \quad  = \quad   \left|
        \begin{array}{ccc}
                1       &       u(t)_x        &       u(t)_y \\
                1       &       v(t)_x        &       v(t)_y \\
                1       &       w(t)_x        &       w(t)_y
        \end{array}
\right|.
\end{displaymath}
We call the three sites \bfindex{collinear} if $D(t)$ is the null polynomial.  
Otherwise let  $\text{rs}(D(t))$ denote the ruling sign of $D(t)$.
If $\text{rs}(D(t))=1$, we say that $w(t)$ is on the 
\bfindex{left} of the polynomial
line $l_{uv}(t)$ at $t=0$, and if $\text{rs}(D(t))=-1$,  we say that $w(t)$ 
is on the \bfindex{right}.

\begin{eexample}
$u(t) = (0,-t)$, $v(t) = (-t,0)$, $w(t) = (t,-2t)$.
Then $D(t) = 0$, so $u,$ $v,$ and $w$ are collinear at $t=0$.
\end{eexample}

\begin{eexample}
\elabel{exontheleft}
Change $w$ to  $w(t) = (t-t^2,-2t)$. Now $D(t) = t^3$,
so $w$ is on the left of the polynomial line  $l_{uv}(t)$ at $t=0$.
\end{eexample}

\begin{eexample}
For small negative $t$, the site $w(t)$ is on the right of $l_{uv}$.
If we change $w(t)$ into $w(t)=(t-t^3,-2 t)$, then $w(t)$ is
on the left of $l_{uv}(t)$ for small negative $t$ while it is still
on the right for small positive $t$.
\end{eexample}

\subsection{The center of a circle.}
\elabel{scirclecenter}

Let $u(t)$, $v(t)$, and $w(t)$ be polynomial sites 
such that $u(t)$, $v(t)$, and $w(t)$ are not collinear at $t=0$.
Let $D=D(t)$ be the determinant of above and let $d=d(t)$, and $e=e(t)$ be the determinants
given by
\begin{displaymath}
d=- \left|
  \begin{array}{lll}\vspace{.7mm}
    {u(t)}_x^2 + {u(t)}_y^2 & u(t)_y & 1 \\[.1cm]
    {v(t)}_x^2 + {v(t)}_y^2 & v(t)_y & 1 \\[.1cm]
    {w(t)}_x^2 + {w(t)}_y^2 & w(t)_y & 1 \\
  \end{array}
\right|, \quad
e= \left|
  \begin{array}{lll}\vspace{.7mm}
   {u(t)}_x^2 + {u(t)}_y^2 & u(t)_x & 1 \\[.1cm]
   {v(t)}_x^2 + {v(t)}_y^2 & v(t)_x & 1 \\[.1cm]
   {w(t)}_x^2 + {w(t)}_y^2 & w(t)_x & 1 \\
  \end{array}
\right|.
\end{displaymath}
We define the {\bf circle center} 
\index{circle!center} of $u(t)$, $v(t)$, and $w(t)$ at $t=0$
as the point $c$ given by the coordinates
\begin{displaymath}
                c_x  =  -\lim_{t\to0} \frac{d}{2 D}, \qquad
                c_y  =  -\lim_{t\to0} \frac{e}{2 D}.
\end{displaymath}
We allow a circle center to be located at infinity.
When leaving out the time dependency, these are of course the ordinary
formulas describing a circle center. 
The circle $C(u(t),v(t),w(t))$ is \bfindex{oriented clockwise} at $t=0$
iff $w(t)$ is on the right of $l_{uv}(t)$ at $t=0$. If $w(t)$ is on the
left of $l_{uv}(t)$ at $t=0$, then $C(u(t),v(t),w(t))$ is 
\bfindex{oriented counterclockwise}.
 
If $u(0)=v(0)=w(0)=(0,0)$, then we say that the circle
defined by $u(t)$,  $v(t)$ and $w(t)$ has \bfindex{positive radius} at $t=0$
if and only if
\begin{eqnarray}
\label{radius}
        {\rm maximum}~(~|c_x|~,~ |c_y|~) \quad > \quad  0.
\end{eqnarray}

\begin{eexample}
\elabel{exradius}
Let again $u(t) = (0,-t)$, $v(t) = (-t,0)$, and  $w(t) = (t-t^2,-2t)$.
Then $a = t^3$ and  $d = e = 4 t^3 - 2 t^4 + t^5$. So $c = (-2,-2)$
and we can conclude that $u,v$ and $w$ define a circle of positive radius 
at $t=0$. As $w(t)$ is on the left of $l_{uv}(t)$ at $t=0$,
the circle $C(u(t),v(t),w(t))$ is oriented counterclockwise
at $t=0$.
\end{eexample}

\begin{eexample}
Let $u(t) = (-t^3,2t)$, $w(t)=(t^3,-2 t)$, and $v(t)=(-t^4, 3 t^2)$.
Then $C(u,v,w)$ is oriented clockwise at $t=0$,
while the circle center $c$ of $C(u,v,w)$
is situated at infinity.
\elabel{exinfradius}
\end{eexample}

\subsection{Cocircularity.} 
\elabel{sscocircular}

Let $u(t)$, $v(t)$, and $w(t)$ be distinct polynomial sites.
Assume that at $t=0$ the sites $u(t)$, $v(t)$, and $w(t)$ are not 
collinear and that the circle $C(u(t),v(t),w(t))$ is oriented clockwise.
As long as the sites do not coincide, we can analyze if  a fourth 
polynomial site $q(t)$ is inside or outside of the circle defined by 
$u(t)$, $v(t)$, and $w(t)$. Consider the polynomial $I(t) \in \R[t]$
defined by 
\begin{displaymath}
I(t) \eis
\left|
\begin{array}{cccc}
        u(t)_x & u(t)_y & u(t)_x^2+ u(t)_y^2 & 1   \\[.1cm]
        v(t)_x & v(t)_y & v(t)_x^2+ v(t)_y^2 & 1   \\[.1cm]
        w(t)_x & w(t)_y & w(t)_x^2+ w(t)_y^2 & 1   \\[.1cm]
        q(t)_x & q(t)_y & q(t)_x^2+ q(t)_y^2 & 1   \\
\end{array}
\right|.
\end{displaymath}
We call the four sites \bfindex{cocircular} at $t=0$ if
$I(t)$ is the null polynomial.
Let otherwise  $\text{rs}(I(t))$ denote the ruling sign of $I(t)$.
If $\text{rs}(I(t))=1$, we say that $q(t)$ is \bfindex{outside} of 
the polynomial
circle $c_{uvw}(t)$ at $t=0$, and if $\text{rs}(D(t))=-1$,  we say that $q(t)$
is \bfindex{inside} the circle.
A set of polynomial sites $S(t)$ is in {\bf general position}
\index{general position!polynomial sites} at $t=0$,
iff no three sites are collinear at $t=0$ and no four sites
are cocircular at $t=0$.

\begin{eexample} 
Let $u(t)$, $v(t)$, and $w(t)$ be as in Example \ref{exradius}. Let
$q(t) = (0,0)$. We check whether $q(t)$ is inside or outside of the 
clockwise oriented circle $C(u(t),w(t),v(t))$ at $t=0$.
As $I(t) = 4 t^4 + O(t^5)$, we conclude that
$q(t)$ is outside of $C(u(t),w(t),v(t))$ at $t=0$.
\end{eexample}

\begin{remark}
For small negative $t$, the circle $C(u(t),w(t),v(t))$ is oriented
counterclockwise. So the orientation swaps at $t=0$. The site
$q(t)$ is inside the clockwise oriented circle $C(u(t),v(t),w(t))$
for small negative $t$.
\end{remark}

\subsection{Ordering the sites.}

The \bfindex{lexicographic ordering} for polynomial sites is as follows.
Let  $u(t)$ and $v(t)$ be two polynomial sites.
First consider the polynomials $u_x$ and $v_x$ that give the 
$x$-component. If $u_x \neq v_x$, order $u(t)$ and $v(t)$ according to 
the coefficients of the term of lowest degree of $u_x$ and $v_x$ that are distinct. 
If $u_x=v_x$, then compare the coefficients of lowest degree that
are distinct of $u_y$ and $v_y$.  We denote this ordering by $O_<$.

\begin{eexample}
Let $u(t) = (2 t^2, t)$ and $v(t) = (t,t)$. 
The first order coefficient of $u_x$ equals $0$ while the first order coefficient of
$v_x$ equals $1$. As this is the lowest order coefficient of $u_x$ and $v_x$ 
that is distinct, it follows
that $u(t) < v(t)$ with respect to $O_<$.
\end{eexample}
\section{The type of a set of polynomial sites.}
\elabel{stype}

Let $S(t) = \{p_1(t), \dots, p_n(t)\}$ be a list of 
distinct polynomial sites that is in general position
at $t=0$. A clockwise polynomial circle $c_{uvw}(t)$  with
$u(t)$, $v(t)$, and $w(t)$ in $S(t)$ is a \bfindex{Voronoi circle}
at $t=0$ iff there are no polynomial sites from $S(t)$ 
inside $c_{uvw}(t)$ at $t=0$. Represent a 
Voronoi circle $c(t)$ by the ordered list of labels of the polynomial
sites $u(t)$, $v(t)$, and $w(t)$ that define $c(t)$. The order
of the labels corresponds to the cyclic, clockwise order of
the defining sites on $c(t)$  at $t=0$, compare with Section 
\ref{scirclecenter}. The set of
all these lists is the \bfindex{type} of $S(t)$ at $t=0$.

\begin{figure}[!ht]
\begin{center}
\setlength{\unitlength}{.6em}
\begin{picture}(10.1,10.2)
\put(5.1,4.95){\makebox(0,0)[cc]{
        \leavevmode\epsfxsize=9\unitlength\epsfbox{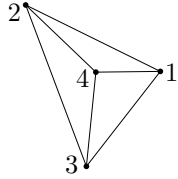}}}
 \put(9.9,5.8){\makebox(0,0)[l]{$1$}}
 \put(0.0,9.6){\makebox(0,0)[l]{$2$}}
 \put(3.6,0.0){\makebox(0,0)[l]{$3$}}
 \put(4.3,5.4){\makebox(0,0)[l]{$4$}}
\end{picture}
\caption{\elabel{fconfig}The Delaunay triangulation of the sites of
Example \ref{expolsites} at $t=0.1$.}
\end{center}
\end{figure}

\begin{eexample}
Let $S(t)$ be the set of polynomial sites given by
\elabel{expolsites}
$$
\begin{array}{lll}
q_1 & = & ( 2 t, \,  2 t^3 + t^4),\\[.05cm]
q_2 & = & ( -2 t - t^3,\,  2 t ),\\[.05cm]
q_3 & = & (-2t^2 - 2 t^4,\,  -3t + 2 t^2 - t^3 + t^4), \\[.05cm]
q_4 & = & (t^2 - 2t^3 + t^5,\,  -t^4).
\end{array}
$$ 
The relative position of the sites at $t=0.1$ is shown
in Figure \ref{fconfig}.
The set $S(t)$ defines at $t=0$ the clockwise oriented circles
$134, 142, 132$ and $243$, where we list only the labels.
The polynomial site $q_4$ is inside the circle $132$ 
at $t=0$, while the other circles are empty at $t=0$.
Therefore the type of $S(t)$ is given by $\{142,134,243\}$.
\end{eexample}

\begin{remark}
Note that this definition of type again represents the situation
for small positive~$t$. This is because the type is defined in terms
of the `inside circle relation', compare Section \ref{sscocircular}.

\end{remark}

\subsection{Abstract Delaunay graph.}

The notion of type enables us to define an 
\bfindex{abstract Delaunay graph} $aD(S(t))$ of $S(t)$ at $t=0$. In fact,
we use Lemma \ref{ltype} as a definition for $aD(S(t))$.
The {\bf vertices} of $aD(S(t))$ at $t=0$ are the labels occurring in the
type. Two labels are connected
by an {\bf edge} if both occur in the same Voronoi circle at $t=0$. 
Note that the  closed paths of length three in $aD(S(0))$
correspond with the Voronoi circles at $t=0$.
The \bfindex{multiplicity} of an edge is the number of distinct Voronoi circles
that contain both the vertices incident to the edge.

\begin{eexample}
Consider the type at $t=0$  
of the set of polynomial sites from Example \ref{expolsites}.
It defines a abstract Delaunay graph $aD(S(t))$ with vertex set
$\{1,2,3,4\}$  and  edge set $\{12,13,14,23,24,34\}$.
The multiplicities of the edges $12$, $23$, and $13$ equal one, those of
the other edges equal two. Compare this to the Delaunay triangulation
shown in Figure \ref{fconfig}.
\end{eexample}

\begin{property}
Let $p_a(t)$ and $p_b(t)$ be two distinct sites of $S(t)$.
\elabel{pmult}
\begin{lijst}
\item If $p_a(t)p_b(t)$ forms an edge of the convex hull $CH(S(t))$ for 
small enough positive $t$, then 
the multiplicity of the edge $ab$ in  $aD(S(0))$ equals one at $t=0$.
\item If $p_a(t)p_b(t)$ forms an edge of the Delaunay triangulation $D(S(t))$
but not of  $CH(S(t))$ for small enough positive $t$, then
the multiplicity of the edge $ab$ in $aD(S(0))$ equals two at $t=0$. 
\end{lijst}
\end{property}

\begin{proof}
Combine the definition of type and Lemma \ref{lunbound}. 
\end{proof}

\subsection{Combinatorial convex hull.}

The convex hull of a set of distinct points in the plane can be defined
as an intersection of half-planes, see Section \ref{sconvexhull}.
But if a set of sites shrinks to one point, the convex hull of
the sites also shrinks to one point.   
However, we can define a  combinatorial convex hull of a set of
polynomial sites at $t=0$ by making  use of half-planes. 

Let $S(t)$ be a set of distinct sites in general position.
It is not necessary that all sites coincide in $(0,0)$ at $t=0$.
Some $i$ is a {\bf vertex} of the \bfindex{combinatorial convex hull}
$cCH(S(t))$ of $S(t)$ at $t=0$  if there exists a polynomial
site $p_j(t) \neq p_i(t)$ such that all sites $p_k(t)$ with
$k \neq i,j$ are on the right of the polynomial line 
$l_{ij}(t)$ at $t=0$. We call such a line a \bfindex{bounding line}.
Two vertices $i$ and $j$ are connected by a directed {\bf edge}, notation
$i \rightarrow j$,  iff
$l_{ij}(t)$ is a bounding line of $cCH(S(t))$ at $t=0$.

\begin{property}
\elabel{pcombi}
Let $S(t)$ be a set of distinct polynomial sites $S(t)$ in general position.
$C_0$ denotes the combinatorial convex hull $cCH(S(t))$ at $t=0$. Let
$a,b \in C_0$ and $p_c(t) \in S(t)$.  
\begin{lijst}
\item $C_0$ is an abstract directed circuit graph.
\item The edges of $C_0$  are exactly the edges of
      multiplicity $1$ of the abstract Delaunay graph
      $aD(S(t))$ at $t=0$.
\item The direction of an edge $ab$ of $C_0$ equals $a\rightarrow b$
iff the Voronoi circle containing both $a$ and $b$ is
oriented like $abc$.
\end{lijst}
\end{property}

\begin{proof} We prove the claims separately.
\begin{lijst}
\item A circuit graph is a connected graph that is regular of
degree two. For $t$ positive, small enough, the boundary of  $CH(S(t))$
is an oriented polygon where the vertices have the same labels as
the vertices of $cCH(S(t))$ at $t=0$. The polygon is  oriented so that
all points of $S(t)$ not lying on an edge $p_i(t)p_j(t)$ of $CH(S(t))$ are on
the right of the directed line defined by 
$p_i(t)$ and $p_j(t)$  that passes first through $p_i(t)$. This orientation
and the connectedness are inherited when $t$ vanishes.
\item This follows from Property \ref{pmult}.
\item If $abc$ is the clockwise orientation of the Voronoi circle defined 
      by $p_a(t)$, $p_b(t)$, and $p_c(t)$, then $p_c(t)(t)$ is on the right 
      of $l_{ab}(t)$ at $t=0$.
\qedhere
\end{lijst}
\end{proof}

\begin{eexample}
\elabel{ex132}
Consider the ordered set $S(t) = \{u(t),v(t),w(t)\}$ of polynomial 
sites introduced in Example \ref{exontheleft}. As $w(t)$ is on the left of $l_{uv}(t)$
at $t=0$, the  combinatorial convex hull $cCH(S(t))$ of $S(t)$
at $t=0$ is given by $132$.
\end{eexample}

\begin{eexample}
 Let $S(t)=\{q_1(t),q_2(t),q_3(t),q_4(t)\}$ as in Example \ref{expolsites}.
Then $cCH(S(t))$ at $t=0$ equals $132$.
\end{eexample}

\section{The Voronoi diagram of a set of 
points and angles between the points.}
\elabel{svorpoly}

In this section we define the Voronoi diagram of a set of points
and angles between the points. This enables us to introduce
a Voronoi diagram for a set of  polynomial sites. 
Moreover, we connect the notions of polynomial lines and polynomial
circles in the context of Voronoi diagrams.

\subsection{Definition.}
\elabel{svorpoldef}

Suppose $S=\{p_1, \dots, p_n\}$ is a set of $n$ not-necessarily distinct 
points in the plane. To every pair of points $(p_i, p_j)$
we add an angle $\alpha_{ij}$ such that the following rule holds:  if two 
points $p_i$ and $p_j$ are distinct
then there is a  unique line though $p_i$ and $p_j$ that makes
some angle $\alpha_{ij} \in \R / 2 \pi \Z$ with the $x$-axis, directed
from $-\infty$ to $\infty$; if $p_i$ and $p_j$ do coincide
then we allow every value in $\R / 2 \pi \Z$ for $\alpha_{ij}$.

Let $\gamma_n$ be a set of $n$ points $p_1(\gamma_n), \dots, 
p_n(\gamma_n)$  and $\binom{n}{2}$ angles 
$\alpha_{12}, \dots, \alpha_{(n-1)n}$ that obey aboves rule.
Fix two points $p_i=p_i(\gamma_n)$ and $p_j=p_j(\gamma_n)$.
The \bfindex{bisection point} $b(p_i, p_j)$ is the point
        $\frac{1}{2}(p_i + p_j)$.
If $p_i\neq p_j$, the bisection point is just the middle of the line segment
$p_ip_j$. If $p_i=p_j$ then  $b(p_i, p_j)$ coincides with the double point
$p_i=p_j$.
The \bfindex{perpendicular bisector} $B(p_i,p_j)$ is the line through
$b(p_i,p_j)$ perpendicular to the angle $\alpha_{ij}=\alpha_{ij}(\gamma_n)$. Let
$\mathbf{n}$ be any non-zero vector, pointing in the direction $\alpha_{ij}$.
The \bfindex{Voronoi half-plane} $vh(p_i,p_j)$ is the half-plane defined by
\begin{eqnarray}
\elabel{eqbiseqn}
        \mathbf{n} \cdot (x-b_x , y-b_y ) & \leq & 0.
\end{eqnarray}
The {\bf Voronoi cell} $V(p_i)$ is defined as
\beq
        V(p_i) & = & \bigcap_{j\neq i} vh(p_i,p_j).
\eeq
A point $x$ is on the {\bf Voronoi edge} $e(p_i,p_j)$ iff
it is on the intersection of the Voronoi cells $V(p_i)$ and $V(p_j)$, that is
\beq
        x \in e(p_i,p_j) & \desda & x \in V(p_i) \cap V(p_j).
\eeq
The {\bf Voronoi diagram} is the  family of subsets of $\R^2$ consisting of
the Voronoi cells $V(p_i)$ and all of their intersections.
The \bfindex{shape} or \bfindex{boundary} of the Voronoi diagram is the union of
the boundaries of the Voronoi cells. 

\begin{remark}
We will always use this definition  in a geometric context which imposes
restrictions on the angles $\alpha_{ij}$. 
\end{remark}

\begin{remark}
We show how to get Equation \ref{eqbiseqn}.
The line given by $\mathbf{n} \cdot (x,y)=0$
defines the line through the origin that is perpendicular to $\mathbf{n}$.
We have to translate this line over the bisection point in order
to get the bisector $B(p_i,p_j)$.
As $\mathbf{n}\cdot\mathbf{n}>0$, the Voronoi half-plane
$vh(p_i,p_j)$ is the half-plane
bounded by the line $B(p_i,p_j)$ in the direction of the tail of $\mathbf{n}$.
\end{remark}

\begin{remark}
Taking for $\mathbf{n}$ the vector $(\cos \alpha_{ij}, \sin \alpha_{ij})$
of length $1$ in the direction $\alpha_{ij}$ gives the inequality
$(y-b_y) \sin \alpha_{ij} \leq (b_x -x) \cos \alpha_{ij}$
for the Voronoi half-plane $vh(p_i,p_j)$.
\end{remark}

\subsection{The Voronoi diagram of a set of polynomial sites.}

Let $S(t)=\{p_1(t), \dots, p_n(t)\}$ be a set of polynomial
sites that is in general position at $t=0$. The
angle $\alpha_{ij}(0)$ for any two sites $p_i(t)$ and $p_j(t)$
is just the direction $\phi_{ij}$ as defined in Section \ref{sprelim}.
We define the Voronoi diagram
of the set of polynomial sites at $t=0$  as the Voronoi diagram of
the points $(p_1(0), \dots, p_n(0))$ and the angles
$\alpha_{12}, \dots, \alpha_{(n-1)n}$.

\begin{remark}
\elabel{rvertex}
At this stage, we do not define Voronoi vertices: if a 
Voronoi vertex is defined as an intersection of three 
Voronoi half-planes, it can become a complete line, 
compare Example \ref{exlinevor}.
\end{remark}

\begin{eexample}
\elabel{exvorpic}
Let $S(t)$ be the set of polynomial sites, introduced in Example \ref{expolsites}.
In Table \ref{thalf-planes} we list the Voronoi half-planes, 
the directions of the bisectors that
bound the half-planes, and the inequalities defining the half-planes.
This results in the Voronoi diagram $V(S(0))$, depicted in Figure \ref{fvorint}.
The Voronoi regions $V(q_1), V(q_2)$ and $V(q_3)$ all have positive
area, while $V(q_4) = \{(0,0)\}$. This shows that the shape of
the Voronoi diagram of $q_1$, $q_2$, $q_3$, and $q_4$ equals the
shape of the diagram of $q_1$, $q_2$, and $q_3$, with $q_4$ omitted.
\end{eexample}

\begin{table}[!ht]
$$
\begin{array}{lll}
\text{half-plane}\qquad\qquad & \text{direction}\qquad\qquad & \text{inequality}\quad\\
& & \\
vh_{1;2} & \pi -\text{arctan} \frac{1}{2}  &  2 x \geq y\\
vh_{1;3} & \text{arctan} \frac{3}{2} - \pi & -\frac{2}{3} x \leq y\\
vh_{1;4} & \pi & x \geq 0\\
vh_{2;3} & -\text{arctan} \frac{5}{2} & \frac{2}{5} x \leq y \\
vh_{2;4} & -\frac{\pi}{4} & x \leq y \\
vh_{3;4} & \frac{\pi}{2} & y \leq 0
\end{array}
$$
\caption{\elabel{thalf-planes}The half-planes, the directions of the 
lines bounding the half-planes,
and the inequalities of the half-planes defined by $q_1(t),q_2(t),q_3(t)$ and
$q_4(t)$ at $t=0$.}
\end{table}

\begin{figure}[!ht]
\begin{center}
\setlength{\unitlength}{0.8em}
\begin{picture}(10,8.4)
\put(5,4.4){\makebox(0,0)[cc]{
        \leavevmode\epsfxsize=10\unitlength\epsfbox{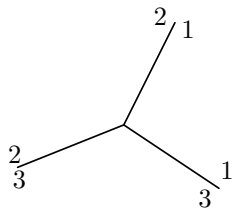}}}
 \put(6.9,8.6){\makebox(0,0)[l]{$2$}}
\put(8.2,8.0){\makebox(0,0)[l]{$1$}}
\put(10.05,1.3){\makebox(0,0)[l]{$1$}}
 \put(9.0,0.0){\makebox(0,0)[l]{$3$}}
 \put(0.23,0.9){\makebox(0,0)[l]{$3$}}
 \put(0.0,2.1){\makebox(0,0)[l]{$2$}}
\end{picture}
\caption{\elabel{fvorint}The Voronoi diagram of the points $q_1(t), q_2(t), q_3(t)$
                and $q_4(t)$ at $t=0$.}
\end{center}
\end{figure} 

\begin{remark}
Given the particular sites in Example \ref{exvorpic}, it is easy to check
that the picture of the Voronoi diagram for small negative values of
$t$ is close to the image under the point reflection in $(0,0)$ of the picture in
Figure \ref{fvorint}.
\end{remark}

\begin{eexample}
\elabel{exlinevor}
Let $S(t)$ be the set of polynomial sites  of Example \ref{exontheleft}
and Example \ref{exradius}, that is, $u(t)=(0,-t)$, $v(t)=(-t,0)$,
and $w(t)=(t-t^2,-2t)$. Then $vh_{uv}$ is given by $y \leq x$,
the half-plane $vh_{u,w}$ by $y \geq x$, and $vh_{v,w}$ by $y \geq x$ as well. Therefore,
$V(u)$ is just the line $y=x$, while $V(v)$ is given by $y \geq x$,
and $V(w)$ by $y \leq x$, all at $t=0$.
\end{eexample}

\subsection{Polynomial bisector.}

Let $u(t)$ and $v(t)$ be two distinct polynomial sites.
We have  introduced the bisector of $u$ and $v$ at $t=0$
as the line passing through the bisection point $b(u,v)$ and
perpendicular to the direction $\phi_{uv}$.  We can derive directly  an 
equation $B_t(u,v)$ for the bisector  depending on $t$ as follows.
We construct two polynomial sites $b(t)$ and $a(t)$ that lie 
on the ordinary bisector for every $t$. The point 
$b(t)$ is just the bisection point $\frac{1}{2}(u(t)+v(t))$.
The other point $a(t)$ is the image of $u(t)$ under the $\frac{\pi}{2}$
rotation around $b(t)$, that is,
\beq	
	a(t) \qquad = \qquad \left(
	\begin{array}{rr}
		0 & -1 \\
		1 & 0 
	\end{array}
	\right) \cdot (u-b)^T + b.
\eeq
Finally, put $B_t(u,v) := l_{ba}(t)$. We call $B_t(u,v)$ the
\index{polynomial!bisector}
{\bf polynomial bisector} of $u(t)$ and $v(t)$.

\begin{lemma}
\elabel{lpolyb}
Let $u(t)$ and $v(t)$ be two polynomial sites.
Let $B(u,v)$ be the perpendicular bisector of $u(t)$ and $v(t)$
at $t=0$. Then
\beq
	B(u,v) & = & \lim_{t\downarrow 0} B_t(u,v).
\eeq
\end{lemma}

\begin{proof}
By construction, the direction of
$B_t(u,v)$ at $t=0$ is perpendicular to $\phi_{uv}$. As 
$B(u,v)$ and $\lim_{t\downarrow 0} B_t(u,v)$ both pass through 
the bisection point $b(0)$, the two lines coincide.
\end{proof}

\begin{eexample}
Let $S(t)=\{u(t),v(t),w(t)\}$ be as in Example \ref{exlinevor}.
The polynomial bisector $B_t(u,v)$ is given by
$$B_t(u,v)=\{(x,y) \in \R^2 \,\mid\,y=x\},$$ while 
$$B_t(u,w)=\{ (x,y) \in \R^2 \,\mid\, y=-2 t + t^2- \frac{1}{2}t^3+x -t x\}.$$
\elabel{exbisec}
\end{eexample}

\subsection{Half-planes and Voronoi circles.}

Let $S(t)$ be a set of polynomial sites that is in general position
at $t=0$.  In Section \ref{svorpoldef} we have introduced the Voronoi diagram 
$V(S(t))$ of $S(t)$ by means of half-planes,
while in Section \ref{stype} we have defined the type of $S(t)$ at $t=0$
in terms of Voronoi circles. The following lemma gives a connection
between polynomial bisectors and circle centers.  

\begin{lemma}
\elabel{lcenter}
Suppose that $p_i(t)$, $p_j(t)$, and $p_k(t)$ are distinct 
non-collinear polynomial sites. 
Denote by  $c$ the circle center of $p_i(t)$, $p_j(t)$, and $p_k(t)$
at $t=0$. Let $B_t(p_i,p_j)$ be the polynomial bisector
of  $p_i(t)$ and $p_j(t)$. Then
\beq
        c  & = & \lim_{t \rightarrow 0}\, 
		(\,B_t(p_i,p_j)\, \cap \, B_t(p_i,p_k)\,).
\eeq
\end{lemma}

\begin{proof}
For $t$ positive and small enough, $p_i(t)$, $p_j(t)$, and $p_k(t)$
are distinct points in general position. For such points 
Lemma \ref{lbicenter} states that the intersection of the bisectors of
the points is the center of the circle defined by the points. 
But then it holds for $t=0$ as well.
\end{proof}

Now we are able to define a Voronoi vertex so that a vertex is
exactly one point in the plane, compare with Remark \ref{rvertex}.
Let $p_i(t)$, $p_j(t)$, and $p_k(t)$ be distinct non-collinear polynomial sites.
Let $c=c_{ijk}(0)$ be the clockwise oriented polynomial circle through
$p_i(t)$, $p_j(t)$, and $p_k(t)$ at $t=0$.
Let $x$ be the circle center of $p_i(t)$, $p_j(t)$, and $p_k(t)$ at $t=0$.
Then $x$ is a \bfindex{Voronoi vertex} $v(p_i,p_j,p_k)$ at $t=0$ if 
$c$ is a Voronoi circle at $t=0$.

Let $u(t)$, $v(t)$, and $w(t)$ be distinct polynomial sites that coincide
in $(0,0)$ at $t=0$. Then the circle center
of $u(t)$, $v(t)$, and $w(t)$ is not equal to $(0,0)$ 
only for special configurations.

\begin{lemma}
\elabel{lisorigin}
Let $u(t)$, $v(t)$, and $w(t)$ be distinct polynomial sites that coincide
in $(0,0)$ at $t=0$. If 
$\phi_{uv} \,\not\equiv\, \phi_{uw} \mod \pi$ or
	$\phi_{uv} \,\not\equiv\, \phi_{vw} \mod \pi$,
then the circle center $c$ of $u(t)$, $v(t)$, and $w(t)$ at $t=0$ equals
$(0,0)$.
\end{lemma}

\begin{proof}
We may assume without loss of generality that 
	$\phi_{uv} ~\not\equiv~ \phi_{uw} \mod \pi$.
Lemma \ref{lcenter} states that $c=\lim_{t \to 0} (B_t(u,v) \cap B_t(u,w))$.
From Lemma \ref{lpolyb} it follows that $c \in B(u,v) \cap B(u,w)$ at $t=0$. 
But then $c$ is the unique intersection point $(0,0)$ of two non-parallel
lines $B(u,v)$ and $B(u,w)$. 
\end{proof}

\begin{eexample} 
Let $u(t)=(0,-t)$, $v(t)=(-t,0)$, and $w(t)=(t-t^2,-2t)$. 
In Example \ref{exradius}
we have shown that $u$, $v$, and $w$ define a circle with center $c = (-2,-2)$ at
$t=0$. We did so by evaluating the determinants $a$, $d$, and $e$ introduced
in Section \ref{scirclecenter}.  Alternatively, we can compute $c(t)$ as
an intersection of two bisectors, say $B_t(u,v)$ and $B_t(u,w)$.
This gives:
\beq
	c(t) & = & B_t(u,v) \cap B_t(u,w) \\
	  & = & (x,x) ~\cap~ (x,-2 t + t^2- \frac{1}{2}t^3+x -t x) \\
	  & = & \frac{1}{2} (-4 + 2 t - t^2, -4 + 2 t - t^2 ).
\eeq 
So $c=c(0)=(-2,-2)$ as expected.  We have already computed the equations 
for the two bisectors in Example \ref{exbisec}. 
\end{eexample}

\section{The positive area cells for one cluster.}
\elabel{sposcell}

In this section we assume we have a set $S$ of polynomial sites that
all coincide in the origin at $t=0$. We show that in general 
some sites can be omitted without changing the shape of the Voronoi
diagram of $S$. We characterize those sites that determine this
shape and present an efficient algorithm to compute the
shape of the Voronoi diagram.

\subsection{The direction hull and positive area cells.}

A \bfindex{zero cluster} at $t=0$ is a set $S(t)$ of 
polynomial sites in general position such that $p_i(0) = (0,0)$
for every site $p_i(t)$ in $S(t)$.
We want to know the shape of the  Voronoi diagram of a zero cluster.
For this purpose it is enough to determine the boundary of the Voronoi
cells of positive area at $t=0$. 

Assume that $S(t)$ is a zero cluster at $t=0$.
Fix some site $v(t)$ in the combinatorial convex hull $cCH(S(t))$ at $t=0$.
There is one incoming edge and one outgoing edge in $cCH(S(t))$ 
at $v(t)$  at $t=0$.
Suppose the incoming edge comes from $u(t)$, while the outgoing edges
goes to $w(t)$.
Site $v(t)$ is called  a \bfindex{corner site} 
of $S(t)$ at $t=0$ if the direction of the incoming edge is distinct from
the direction of the outgoing edge, that is, if  $\phi_{uv} \neq  \phi_{vw}$.
The corner sites of $S(t)$ at $t=0$ are by definition  the vertices
of the \bfindex{direction hull} $DH(S(t))$ at $t=0$. 
Two vertices $u(t)$ and $v(t)$ are connected 
by a directed edge $\stackrel{\longrightarrow}{uv}$ if there 
are no corner sites
on the path from $u(t)$ to $v(t)$ in $cCH(S(t))$.

\begin{eexample}
Let $S(t) = \{q_1(t),\dots, q_4(t)\}$ be as introduced in Example \ref{expolsites}.
From Example \ref{ex132} we know that $cCH(S(t))$ equals $132$.
Compare Table \ref{thalf-planes}. The directions
$\phi_{13} = \text{arctan}\frac{3}{2} $, $\phi_{32}=\pi-\text{arctan}\frac{5}{2}$
and $\phi_{21}=-\text{arctan}\frac{1}{2}$ are all distinct. 
It follows that $DH(S(t))=cCH(S(t))=132$ at $t=0$.
\end{eexample}

\begin{eexample}
\elabel{e123to23}
Let $S(t)=\{q_1(t),q_2(t),q_3(t)\}$ with $q_1(t)=(-t,0)$, $q_2(t)=(0,t^2)$,
and $q_3(t)=(t,-t^2)$.  Then $cCH(S(t)) = 123$, while
$DH(S(t))=13$, as $\phi_{12}=\phi_{23}=0$.
\end{eexample}

\begin{lemma}
\elabel{lbiszero}
Let $S(t) = \{ p_1(t), \dots, p_n(t) \}$ be a zero cluster at $t=0$.
 Then any bisector $B(p_i,p_j)$ passes through $(0,0)$ at $t=0$.
\end{lemma}

\begin{proof}
The bisector $B(p_i,p_j)$ passes through the bisection point 
$b(p_i,p_j) = (0,0)$ at $t=0$.
\end{proof}

\begin{lemma}
\elabel{lposarea}
Let $S(t) = \{ p_1(t), \dots, p_n(t) \}$ be a zero cluster at $t=0$.
Then $V(p_i(0))$ has positive area if and only if $p_i \in DH(S(0))$.
\end{lemma}

\begin{proof}
Assume that $p_i \in DH(S(0))$. As $p_i \in cCH(S(0))$, it follows that
$p_i \in CH(S(t))$ for $t>0$, small enough. That means that there
exists $t_0$ such that for all $0<t<t_0$ the Voronoi cell
$V(p_i(t))$ is unbounded. Moreover, 
$p_i$ is incident with two edges $e(p_{i-1},p_i)$ and $e(p_i,p_{i+1})$,
 say, that are unbounded for all $0<t<t_0$. But this implies that
these edges are unbounded at $t=0$ as well.  
The sites $p_{i-1}$ and $p_{i+1}$ are
the direct predecessor and direct successor of $p_i$ on
$cCH(S(0))$. As $p_i \in DH(S(0))$, the directions 
$\phi_{i-1;i}$ and $\phi_{i;i+1}$ are distinct at $t=0$. But this implies
that $V(p_i(0))$ has two unbounded edges of distinct direction on its 
boundary and by the convexity of $V(p_i(0))$  this implies that
$V(p_i(0))$ has positive area.

For the other direction, assume that $V(p_i(0))$ has positive area at $t=0$. 
It follows from Lemma \ref{lbiszero} that $V(p_i(0))$ is unbounded,
and therefore has two unbounded edges on its boundary of distinct direction.
This implies that $p_i \in DH(S(t))$ at $t=0$.
\end{proof}

The Voronoi diagram $V(S(0))$ of a zero cluster  at $t=0$
looks as follows:
\begin{looplijst}
\item If $p_i(t) \not \in DH(S(t))$ at $t=0$, then $\text{area}(V(S(t)))=0$
      by Lemma \ref{lposarea}.
\item If $p_i(t) \in DH(S(t))$ and if $p_{i-1}(t)$ and $p_{i+1}(t)$
      are its direct predecessor and its direct successor at $t=0$
      in $DH(S(t))$, then
\beq
	V(p_i(0)) & = & vh_0(p_i,p_{i-1}) \cap vh_0(p_i,p_{i+1}). \\
\eeq
\end{looplijst}

\begin{corollary}
\elabel{cposarea}
The shape of the Voronoi diagram $\delta V(S(0))$ of a zero cluster $S(t)$ 
equals the shape of the  Voronoi diagram of $DH(S(0))$.
\end{corollary}

We conclude that for computing the shape of the Voronoi diagram
$\delta V(S(0))$ of a zero cluster, it is enough to consider 
sites on the direction hull $DH(S(0)) \subset S(t)$ only.

\subsection{Determining the direction hull at $t=0$.}

In Section \ref{sprelim} we have defined the combinatorial
convex hull $cCH(S(t))$ of a set of polynomial sites
by imitating the characterization of the ordinary convex hull
as an intersection of half-planes. Exploiting this 
similarity we show how to 
compute $cCH(S(t))$ by an adapted version of the  convex hull algorithm
as presented in \cite{BKOS}, page 6. 

The original algorithm is an incremental algorithm: first it sorts the
$n$ points in the input and next it handles the points one by one.
The original algorithm sorts the input points lexicographically.
This results in a sequence of points that is ordered from left to right
and then from bottom to top. The algorithm determines the convex hull
in two steps. In the first step the \bfindex{upper hull} is determined.
The upper hull is basically that part of the convex hull that is running from the
leftmost point to the rightmost point of the sorted points.   
It constructs the upper hull by adding one point at a time and checking if the
sequence of points that will be the upper hull in the end  keeps going right. 
The lower hull is defined and constructed in a similar way.

Our ordering $O_<$ mimics the lexicographic ordering. Just think of
substituting a very small positive value of $t$ in the set
of polynomial sites and sorting the resulting points lexicographically.
An example is displayed in Figure \ref{figtraj} where the position of
some polynomial sites at $t=0.59$ is indicated by the labels of the sites.
In a similar fashion we define right turns for polynomial sites:
we say that three sites $u(t)$, $v(t)$, and $w(t)$ ordered by $O_<$
make a \bfindex{right turn} at $t=0$, if $w(t)$  is on the right of
$l_{uv}(t)$ at $t=0$.

\begin{alg} \elabel{acCH} Combinatorial convex hull for polynomial sites.
\end{alg}

\begin{algorithmic}[1]
\REQUIRE   set $S(t)$ of $n$ polynomial sites in general position.
\ENSURE   list $L$ containing the vertices of $cCH(S(0))$ in 
clockwise order.
\STATE  Sort the sites by $O_<$, resulting in 
      a sequence $p_1(t), \dots, p_n(t)$.
\STATE   Create list $L_{\text{upper}}=\{p_1,p_2\}$.
\FOR{every $i \in \{3,\dots,n\}$} 
\STATE  Append $p_i$ to $L_{\text{upper}}$.
\WHILE{$L_{\text{upper}}$ contains more than two sites {\bf and} the
          last three sites in $L_{\text{upper}}$ do not make a right turn,}
\STATE Delete the middle of the last three sites from $L_{\text{upper}}$.
\ENDWHILE
\ENDFOR
\STATE  Construct $L_{\text{lower}}$ in a similar way and append the result to
      $L_{\text{upper}}$. Call the resulting list $L$;.
\STATE \return{} $L$.
\end{algorithmic}

The running time of the algorithm is $O( n~\text{log}~n)$,
cf.\ \cite{BKOS}.
Let $S(t)$ be a zero cluster again.  Given  its combinatorial convex 
hull $cCH(S(t))$ at $t=0$, we 
compute in $O( n )$ time the direction hull $DH(S(t))$ at $t=0$.

\newpage
\begin{alg}
Direction hull.
\end{alg}

\begin{algorithmic}[1]
\REQUIRE combinatorial convex hull $cCH(S(0))$ of a zero cluster $S(t)$.
\ENSURE direction hull $DH(S(t))$ at $t=0$.
\STATE Compute the direction  $\phi_{uv}$ for every edge $uv$ in $cCH(S(t))$
      at $t=0$. 
\IF{two consecutive edges $uv$ and $vw$ have the same direction, that
      is, $\phi_{uv}=\phi_{vw}$,} 
\STATE Delete $w$ and its incident edges.
\STATE Create a new edge $uw$.
\ENDIF
\end{algorithmic}

\section{Degenerate Voronoi cells and edges.}
\elabel{sposedge}

Let $S(t)$ be a zero cluster at $t=0$. We have seen before that only sites
in the direction hull $DH(S(0))$  
generate a Voronoi cell of positive area. 
It is possible however that for a site $p_i \in S(t)\setminus DH(S(0))$ 
there exists an edge $e(p_i,p_j)$
that does not collapse on $(0,0)$ at $t=0$. In this section we consider sites
$p_i(t) \in S(t)$ such that there exists, at $t=0$, an edge $e(p_i,p_j)$ with 
$x \in e(p_i,p_j)$, $x \neq (0,0)$. We first give an explicit example
of  an edge $e(p_i,p_j)$ where even $(0,0) \not\in e(p_i,p_j)$, while
$p_i \not\in DH(S(0))$.
Next we present an algorithm that determines all edges $e(p_i,p_j)$ that
contain some $x\in\R^2$ with $x\neq (0,0)$
for a zero cluster  $S(t)$ at $t=0$. 

\begin{figure}[!ht]
\begin{center}
\setlength{\unitlength}{.8em}
\begin{picture}(14,14)
\put(7,7){\makebox(0,0)[cc]{
        \leavevmode\epsfxsize=14\unitlength\epsfbox{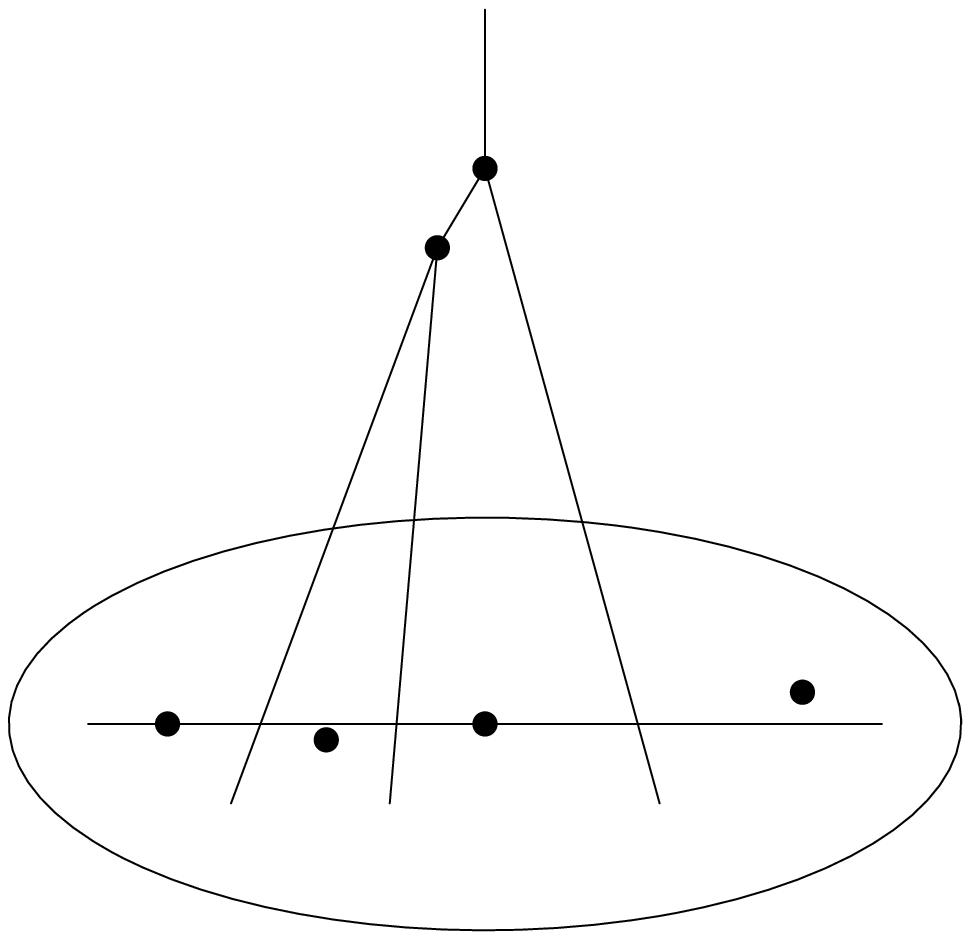}}}
\put(2.2,2.5){\makebox(0,0)[l]{p}}
\put(4.6,2.4){\makebox(0,0)[l]{q}}
\put(7,2.5){\makebox(0,0)[l]{r}}
\put(11.5,2.7){\makebox(0,0)[l]{s}}
\put(4.3,10.5){\makebox(0,0)[l]{pqr}}
\put(8,11.6){\makebox(0,0)[l]{psr}}
\end{picture}
\caption{\elabel{fnoedge}Edge of zero length outside $(0,0)$.}
\end{center}
\end{figure}

\begin{eexample}
\elabel{exshortedge}
Let $S(t)=(p(t),q(t),r(t),s(t))$, where $p(t)=(-2 t, 0)$, 
$q(t)=(-t,-\frac{1}{4} t^2)$, $r(t)=(0,0)$, and
$s(t)=(2 t, 2 t^2)$.  The direction hull at $t=0$ consists of the two
sites  $p(t)$ and $q(t)$, while $\phi_{pq}=0$. Therefore the shape of the
Voronoi diagram $\delta V(S(0))$ is just a vertical line through the origin.
There are two Voronoi cells of positive area at $t=0$: 
$V(p)$, whose cell is the left half-plane, and $V(s)$,
whose cell is the right half-plane.

The type of $V(S(t))$ however is $prq,psr$. A schematic picture
of the situation is given in Figure \ref{fnoedge}. From the type 
it follows that 
combinatorially, $V(S(t))$ consists of two vertices or circle centers,
$prq$ and $psr$ that are connected by an edge $e(p,r)$. 
Besides, there are unbounded edges $e(p,q)$ an $e(q,r)$, incident to 
the circle center $prq$ and edges $e(r,s)$ and $e(p,s)$,
incident to $psr$. Both circle centers are situated at
the point $(0,2)$. This implies that the length of
$e(p,r)$ at $t=0$ equals zero, while the four other edges involved,
have infinite length.
\end{eexample}

Let $S(t)$ be a zero cluster at $t=0$.
An edge $e(p_i,p_j)$ with $p_i,p_j \in S(t)$ is an 
\bfindex{outside edge} at $t=0$ if there exists $x \in e(p_i,p_j)$
such that $x \neq (0,0)$. Here we allow $x$ to be any point
of the form $x=(a,b)$, $x=(\pm \infty, b)$, $x=(a, \pm \infty)$
or $x=(\pm \infty, \pm \infty)$, where $a,b \in \R$. Allowing
these possibilities is 
motivated by Example \ref{exradius}, which demonstrates that circle
centers defined by polynomial sites can be at infinity.

Algorithm \ref{aoutside} determines all outside edges 
of a set of polynomial sites in general position
at $t=0$.

\begin{alg} 
Determining outside edges and their positions.
\elabel{aoutside}
\end{alg}

\begin{algorithmic}[1]
\REQUIRE a zero cluster $S(t)$.
\ENSURE all outside edges with their endpoints at $t=0$.
\STATE Determine the type $T$ of $S(t)$ at $t=0$.
\STATE 
      Any edge $uv$ of multiplicity 1 in the abstract Delaunay
      graph $aD(S(0))$ corresponds to an 
      unbounded edge $e(u,v)$ in the Voronoi diagram $V(S(t))$.
      This edge $e(u,v)$ starts at the circle center of the unique circle
      listed in $T$ that has both $u$ and $v$ on its boundary. The direction
      of $e(u,v)$ is perpendicular to $\phi_{uv}$. Edge $e(u,v)$
      is oriented so that it makes an angle of $\frac{\pi}{2}$
      with $\phi_{uv}$;
\STATE Determine those circles in $T$ that have positive 
      radius; 
\FOR{$a$ and $b$ two labels occurring in some 
      circle of positive radius}
\STATE (We show how to `draw' the edge $ab$).
\STATE Check if we have not drawn $ab$ before in step~2
      as an unbounded edge.
\IF{the combination of the labels $a$ and $b$ occurs twice in the list of
      the circles of positive radius}
\STATE there are two circle centers 
      distinct from $(0,0)$. We connect the two circle centers by an edge $ab$;
\ELSIF{the combination of the labels occurs only once} 
\STATE we connect the corresponding circle center $c_{ab}$ to the origin. 
\ENDIF
\ENDFOR
\end{algorithmic}

Determining the shape of the  Voronoi diagram  
is relatively easy, compare Section \ref{sposcell}. We have seen in
this section, however, that the combinatorial structure on the boundary
itself can get rather complicated. 

\section{Example: 20 polynomial sites.}
\elabel{ex20poly}

In this section we analyze a slightly bigger example, illustrating the
techniques and theory developed in the former sections. We 
consider a set $S(t)$ of 20 polynomial sites.   
We determine
the shape of the Voronoi diagram $V(S(t))$ at $t=0$. Next we focus on
the outside edges of $V(S(t))$ at $t=0$. We visualize the polynomial
sites and compare the outcome of the computation to the visualization.
First we introduce the polynomial sites. $S(t)$ is a set  of 20 distinct sites $(p_1(t), \dots, p_{20}(t))$
such that $p_i(t) = (0,0)$ for $i=1, \dots, 20$:

\begin{displaymath}
\begin{array}{llll}
p_1 & (0, t^5), & p_{11} & (-2 t^2 - 2 t^4, -3 t + 2 t^2 - t^3 + t^4), \\
p_2 & (0, 2 t + t^4 + 2 t^5),  & p_{12} & (-2 t - t^4, 3 t^2 + 3 t^3 - t^5), \\
p_3 & (2 t, 2 t^3 + t^4), &  p_{13} &  (-2 t^2 - t^4, -2 t^2 + 3 t^3), \\
p_4 & (t^2, 2 t + 2 t^2 - t^3),  & p_{14} & (-t + 3 t^2 - 2 t^3 + 3 t^4, t + 3 t^2 - 3 t^5),\\
p_5 & (-2 t^4, -t^5), & p_{15} &  (-t - 2 t^5, -3 t^3 + 2 t^5),\\
p_6 & (3 t^5, t^2 + 2 t^3 + 3 t^5),  & p_{16} & ( t^2 - 2 t^3 + t^5, -t^4),\\
p_7 & (2 t + 3 t^2 - 2 t^3, 2 t - 2 t^2 - t^3),  & p_{17} & (3 t^3 + 2 t^5, -t^4),\\
p_8 & (-2 t - t^3, 2 t), &  p_{18} &  ( 3 t^3 + 3 t^4 + 2 t^5, 2 t - 2 t^3 - t^4),\\
p_9 & (-2 t + t^3, -2 t - t^3),  & p_{19} & ( t^2 + t^3 + 3 t^5, 0),\\
p_{10} & (t + 2 t^2 + 2 t^3, -t + 3 t^3),  & p_{20} &  (3 t^2 + 2 t^3 + 3 t^5,  t^2 + 3 t^3 + 2 t^4).\
\end{array}
\end{displaymath}

We  visualize $S(t)$ by considering any $p_i(t) \in S(t)$ as a plane curve, 
represented by a parametrization $(p_{i,x}(t), p_{i,y}(t))$.
In Figure \ref{figtraj} these twenty curves are plotted
for $t \in [0,1]$. The position of a curve at $t=0.59$ is marked by its label.  

\begin{figure}[!ht]
\begin{center}

\leavevmode\epsfxsize=11cm\epsfbox{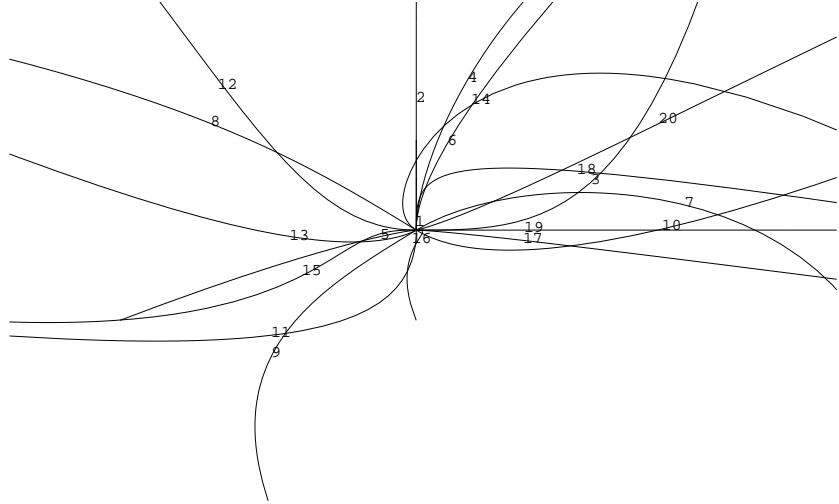}

\end{center}
\caption{The polynomial sites for $t \in [0,1]$. \elabel{figtraj} }
\end{figure}

\subsection{The positive area diagram.}

Recall that the type of $S(t)$ is defined a the set of all Voronoi circles at $t=0$. 
We have computed it using brute force. That is, for every  $\binom{n}{3}$  unordered triples
of polynomial sites $p_i(t)$, $p_j(t)$ and $p_k(t)$ we first determine the orientation
at $t=0$ of the circle $C$ passing though $p_i(t)$, $p_j(t)$ and $p_k(t)$. Next we check
if any of the $n-3$ remaining sites of $S(t)$ is contained in $C$. If not, we conclude
that $C$ is a Voronoi circle and add the labels $i$, $j$ and $k$ ordered with respect
to the orientation of $C$ to the type in the making. It turns out that the type of
$S(t)$ is given by:

\begin{quote}
((1, 5, 6), (1, 17, 5), (1, 6, 17), (2, 8, 4), (2, 4, 18), (2, 14, 8), (2, 18, 14), 
	(3, 20, 7), (3, 11, 10), (3, 10, 20), (4, 7, 18), (5, 13, 6), (5, 17,
     13), (6, 13, 15), (6, 15, 14), (6, 14, 20), (6, 16, 17), (6, 19, 16), (6,
     20, 19), (7, 20, 18), (8, 12, 9), (8, 14, 12), (9, 13, 11), (9, 12, 
    15), (9, 15, 13), (10, 11, 13), (10, 13, 20), (12, 14, 15), (13, 17, 
    16), (13, 16, 19), (13, 19, 20), (14, 18, 20)).
\end{quote}

Counting multiplicities of the edges shows that all edges have multiplicity
$2$ except for the edges
$$ (4,8), (3,7), (3,11), (4,7), (8,9), (9,11). $$
As site $p_7$ is on the right of $l_{4,8}$ at $t=0$ the combinatorial convex
hull $cCH(S(0))$ is given by $(8,4,7,3,11,9)$.
The directions of the edges of $cCH(S(0))$ are 
$\phi_{8;4} = 0$, $\phi_{4;7} = 0$, $\phi_{7;3} = - \frac{\pi}{2}$, 
$\phi_{3;11} = - \text{arctan}(\frac{2}{3}) - \frac{\pi}{2}$,
$\phi_{11;9} =\frac{\pi}{2}+ \text{arctan}~2$, and $\phi_{9;8} = \frac{\pi}{2}$.
The directions $\phi_{8;4}$ and $\phi_{4;7}$ are equal, so site $p_4$ is not a corner
site. We conclude that the direction hull $DH(S(0))$ is given by
$(8,7,3,11,9)$. 
This outcome can almost be `guessed' from Figure \ref{figtraj} by considering
the convex hull of the points $p_1(t)$ to $p_{20}(t)$ for small $t$. 
Knowing the direction hull and the directions of the edges of the direction hull,
we know in fact the shape of the  Voronoi diagram $\delta V(S(0))$. It is depicted on the
left in Figure~\ref{figvlim}.

\subsection{Combinatorics of the edges outside $(0,0)$.}

The only Voronoi circle that has a circle center outside of $(0,0)$ at $t=0$ is the
clockwise oriented circle $(8,12,9)$. Its center is at infinity.
This means that the unbounded edge $e(8,9)$ starts at the Voronoi vertex
at infinity, while two edges $e(8,12)$ and $e(9,12)$ run between 
vertices situated in $(0,0)$ and this vertex at infinity. We have
seen before that $p_4 \in cCH(S(0))$ with $p_8$ its direct predecessor
and $p_7$ its direct successor in $cCH(S(0))$. 
As there are no further circle centers situated outside the origin, 
 there are two unbounded edges $e(4,8)$ and $e(4,7)$ that, 
starting in $(0,0)$, both run upwards, but never meet. A schematic picture
of the combinatorics of the edges outside $(0,0)$ is given on the right
in Figure \ref{figvlim}.

 \begin{figure}[!ht]
\begin{center}

\leavevmode\epsfxsize=8cm\epsfbox{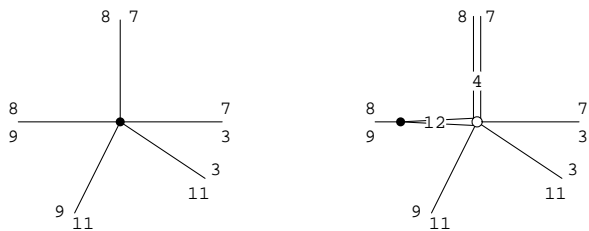}

\end{center}
\caption{The positive area diagram at $t=0$ on the left. The combinatorics outside 
$(0,0)$ included on the right.  \elabel{figvlim} }
\end{figure}

\section{The shape of a  Voronoi diagram at $t=0$.} 
\elabel{swholedia}

In Section \ref{sposcell} we have given a method that determines the
shape of the  Voronoi diagram of a zero cluster at $t=0$. That is, we assumed
that all sites coincide in the origin at $t=0$.
In this section we apply this method
in order to determine the shape of a  Voronoi diagram 
at $t=0$
for an arbitrary set of polynomial sites in general position. 

Let $S(t)$ be a set of polynomial sites of size $n$ that is in general position
at $t=0$.   The {\bf cluster locations}  \index{cluster!location}
$l(0) = \{l_1(0), \dots, l_m(0)\}$ of 
$S(t)$ at $t=0$ are the distinct positions of the sites in $S(t)$ at $t=0$.
An $\mathbf{l_j}${\bf -cluster}, for $l_j \in \R^2$ is a collection of polynomial sites
that coincides at $l_j$ at $t=0$.   Write 
	$S(t)  =  \cup_{j=1}^m S_j(t) $,
where $S_j(t) \subset S(t)$ is the set of sites in $S(t)$ such that
$p_{j_i}(0) = l_j(0)$ for all $p_{j_i} \in S_j(t)$.  That is,
we divide $S(t)$ with respect to the $l_j$-clusters at $t=0$.

\begin{lemma} 
\elabel{lplug}
The shape of the  Voronoi diagram  $\delta V(S(t))$ at $t=0$ of $S(t)$ is 
given by
\beq
        \delta V(S(0)) & = & \bigcup_{j=1,\dots ,m} (\, V(\,l_j(0)\,) 
		\,\cap\, \delta V(\,S_j(0)\,) \,).
\eeq
\end{lemma}

\begin{proof}
First of all, suppose that $|S_j(0)|=1$, for  all $j =1, \dots, m$. 
In this case the lemma just states that the Voronoi diagram of distinct points
coincides with the subdivision of the plane into Voronoi cells, see
Section \ref{svordia}. Next consider an arbitrary point $x$ in the interior
of, say, $V(l_j)(0)$. Then $x$ is, by definition, closer to $l_j$ than
to any other location $l_k$ for $k \neq j$. The lemma follows from
Corollary \ref{cposarea} adapted for an $l_j$-cluster $S_j(t)$.
\end{proof}

\begin{figure}[!ht]
\begin{center}
\setlength{\unitlength}{.9cm}
\begin{picture}(12,6.2)
\put(2,3.3){\makebox(0,0)[cc]{
        \leavevmode\epsfxsize=3\unitlength\epsfbox{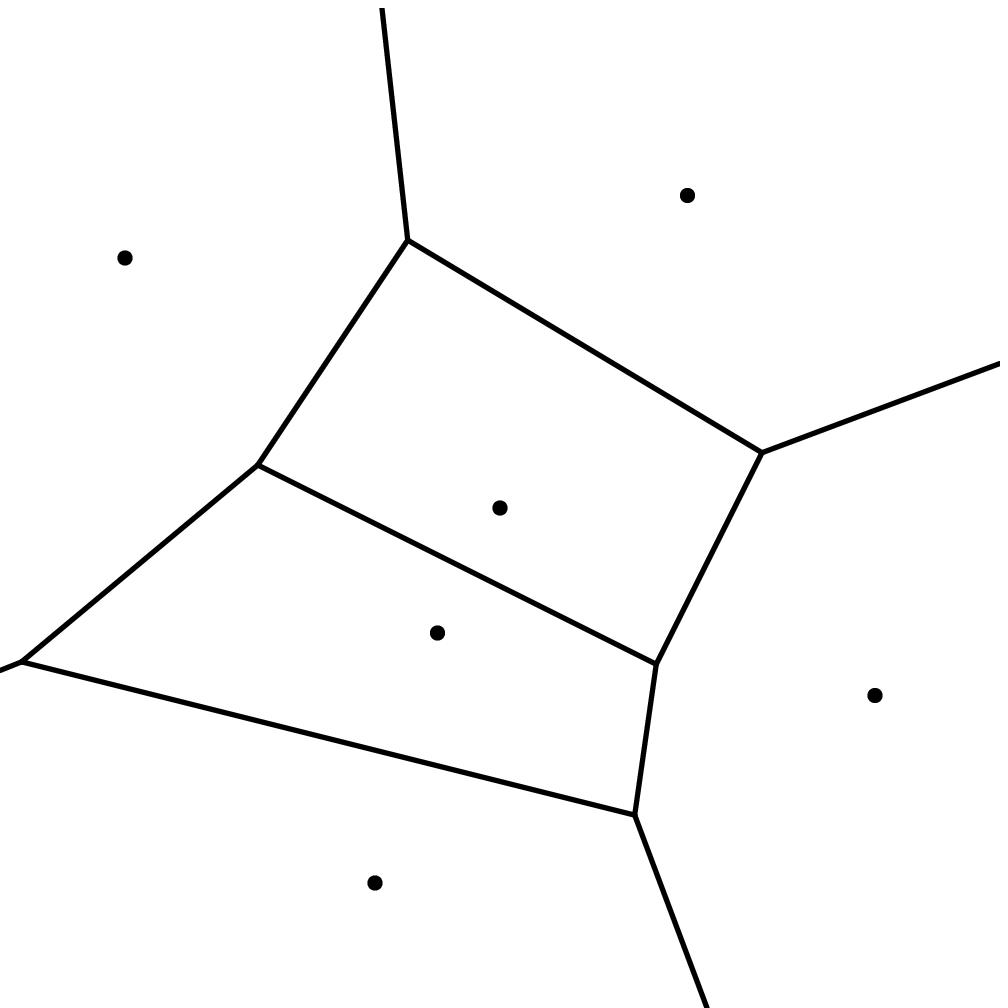}}}
\put(6,5.3){\makebox(0,0)[cc]{
        \leavevmode\epsfxsize=2\unitlength\epsfbox{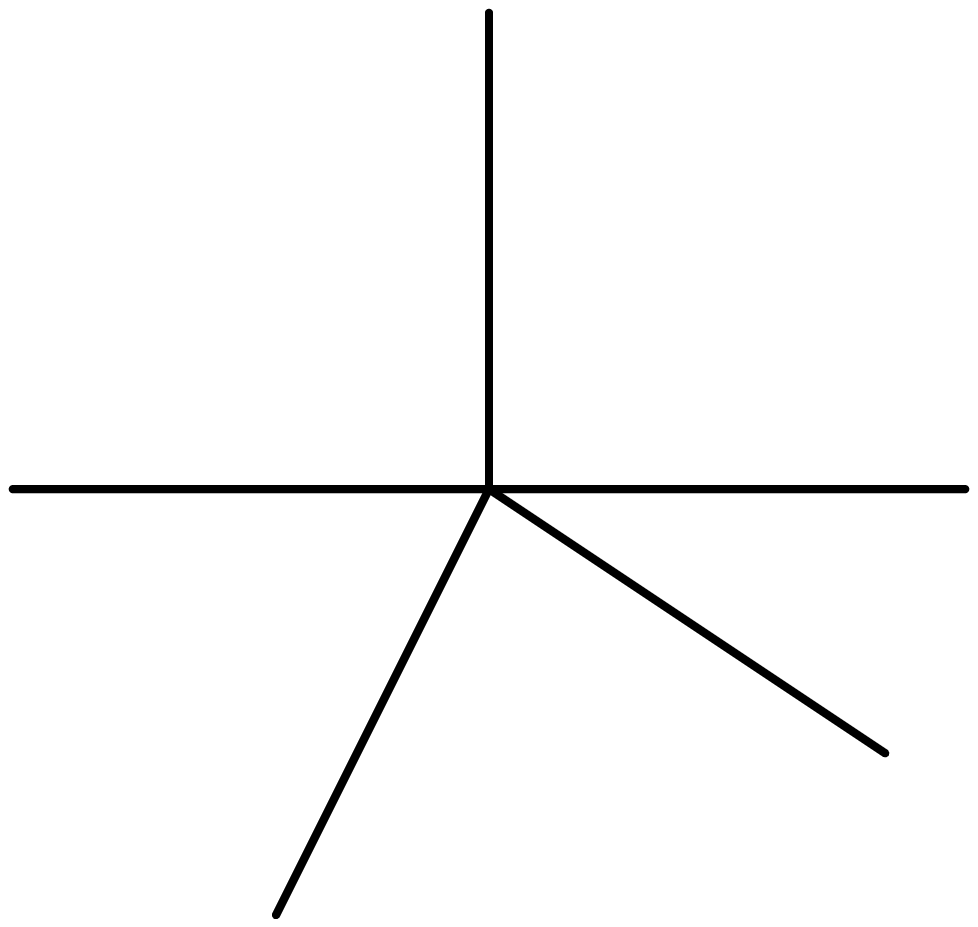}}}
\put(6,1.3){\makebox(0,0)[cc]{
        \leavevmode\epsfxsize=2\unitlength\epsfbox{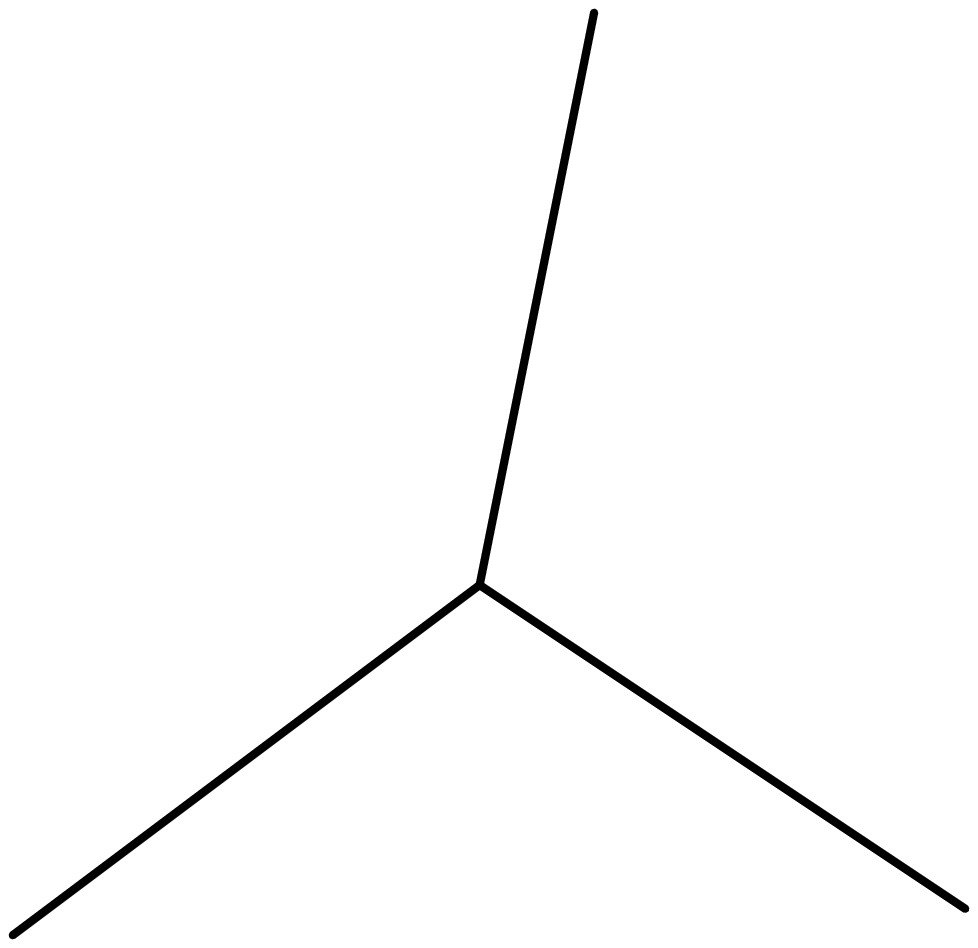}}}
\put(10,3.3){\makebox(0,0)[cc]{
        \leavevmode\epsfxsize=3\unitlength\epsfbox{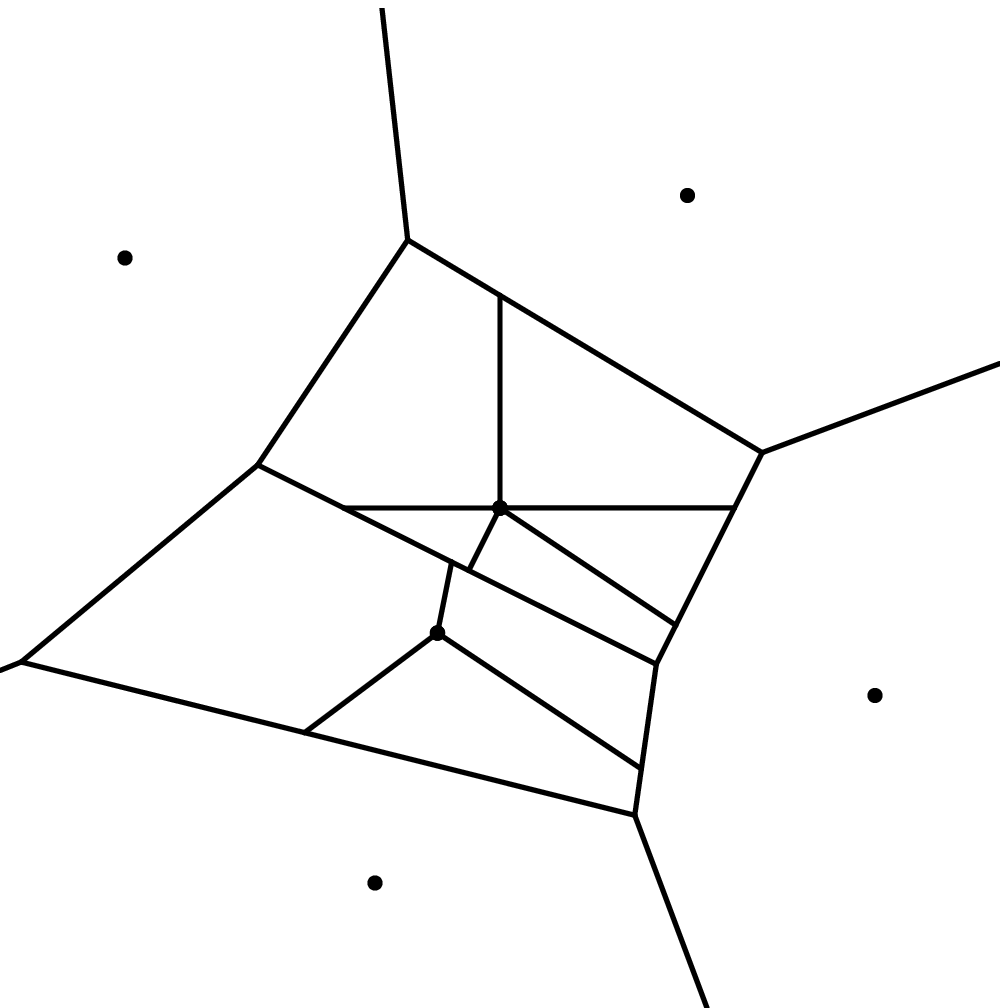}}}
\put(4,3.3){\makebox(0,0)[l]{+}}
\put(8,3.3){\makebox(0,0)[l]{=}}
\end{picture}
 \caption{\elabel{fwithwithout}Plugging diagrams.}
\end{center}
\end{figure}

\begin{eexample}
\elabel{explug}
Suppose we are given the following sites:
\beq
 S(t) & = & ((-6, 4), (-2, -6), (-1, -2 - 3 t), (0, -3 t), (3, 
    5), (6, -3),\\ 
      & & (-1 - 3 t, -2 + t),
      (-2 t, -2 t), (-2 t, 2 t), (2 t, 0), (2 t, 2 t), (-1 + 2 t, -2)).
\eeq
The cluster locations at $t=0$ are given by
\beq
l(0) & = &  ((-6, 4), (-2, -6), (-1,-2), (0,0), (3,5), (6,-3)).
\eeq
The Voronoi diagram $V(l(0))$ of the cluster locations  is presented
on the left in Figure \ref{fwithwithout}. There are two clusters
consisting of more than one point. The shape  of the Voronoi diagram 
of the  zero cluster, given by 
\beq 
S_1(t) & = & ((0, -3 t), (-2 t, -2 t), (-2 t, 2 t), (2 t, 0), (2 t, 2 t)),
\eeq
is presented in the top middle of Figure \ref{fwithwithout}, while
the shape of the  Voronoi diagram of the $(-1,-2)$-cluster 
\beq
	S_2(t) & = & ((-1, -2 - 3 t), (-1 - 3 t, -2 + t), (-1 + 2 t, -2)),
\eeq
is shown in the bottom middle. According to Lemma \ref{lplug}, we can 
plug these two cluster diagrams in the cluster locations diagram in order
to get the shape of the Voronoi diagram of $S(t)$ at $t=0$. 
This final diagram is 
depicted on the right in Figure \ref{fwithwithout}.
\end{eexample}

\section{Generalizations and conclusion.}
\elabel{spolyk}

\subsection{Dropping general position.}

Throughout this chapter we have assumed that sets of polynomial
sites are in general position. General position
is not needed for the methods presented to work.  Dropping general position
just means that a lot of extra cases have to be checked, which distracts from
the main line. Note that the notion of type as introduced in 
Section \ref{stopchange} is defined for point sets that are not 
in general position as well. As we can also compute convex hulls for point
sets that are not in general position, there are no big obstructions
for extending the methods to sets of polynomial sites in arbitrary 
position. 

Recall from Section \ref{sprelim} that
general position for polynomial sites is defined in terms of
cocircularity and collinearity polynomials: if some set $S(t)$
is not in general position at $t=0$, then there is at least
one cocircularity  polynomial $I(t)$ or collinearity polynomial
$D(t)$ that is equal to the null polynomial. This implies that
$S(t)$ is not in general position for any value of $t$.  

\subsection{Generalization to higher order Voronoi diagrams.}

In the $k$-th order Voronoi diagram, the plane is partitioned according
to the $k$ closest sites, see Chapter \ref{chposet}. 
Algorithm \ref{algkvoronoi}  computes the $k$-th order
Voronoi diagram of a set $S$ of distinct points in general position,
using
\begin{itemize}
  \item circles through three points of $S$,
  \item points from $S$ inside these circles,
  \item directions of lines between points of $S$.
\end{itemize}  
We have shown in this chapter that these concepts can be generalized
to sets $S(t)$ of polynomial sites. Therefore, Algorithm \ref{algkvoronoi}
can also be generalized to sets of polynomial sites in general position.
This will give us a combinatorial $k$-th order Voronoi diagram.  

\begin{figure}[!ht]
\begin{center}
\setlength{\unitlength}{1cm}
\begin{picture}(12,4.5)
 \put(2,2.5){\makebox(0,0)[cc]{
        \leavevmode\epsfxsize=4\unitlength\epsfbox{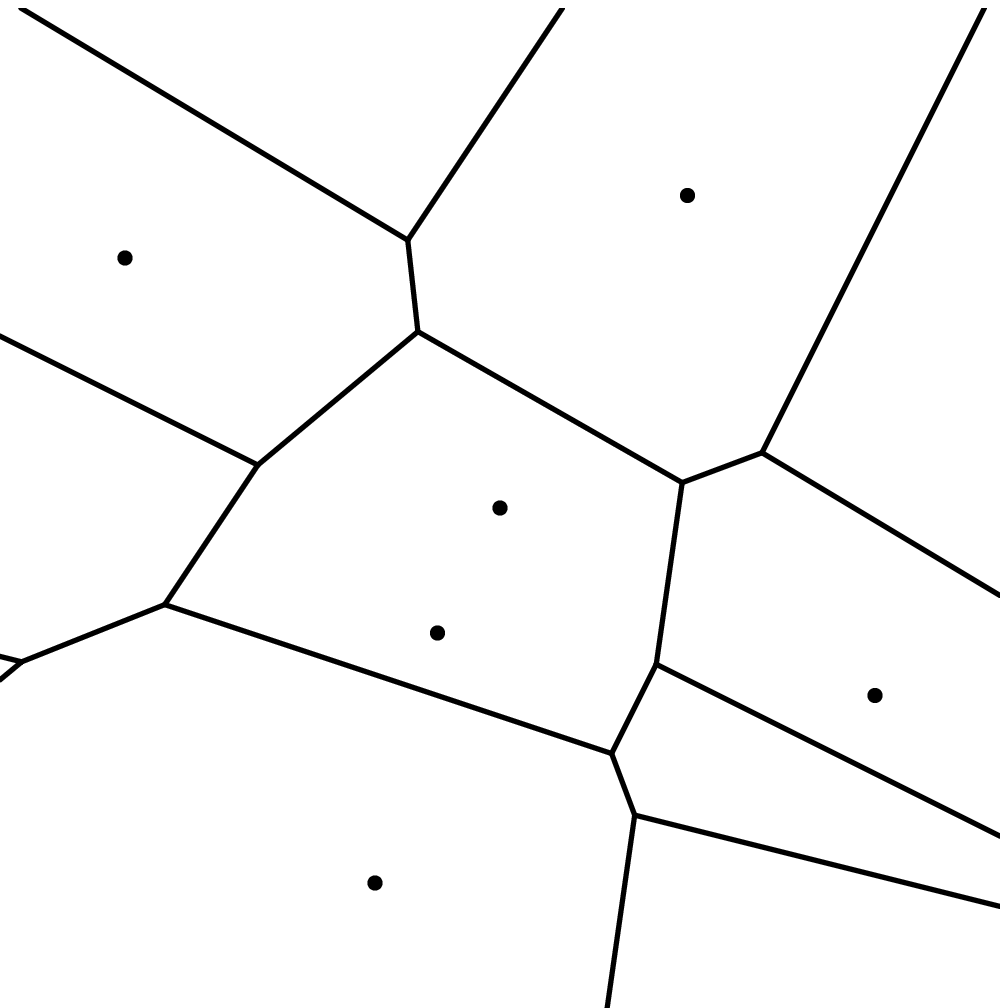}}}
 \put(10,2.5){\makebox(0,0)[cc]{
         \leavevmode\epsfxsize=4\unitlength\epsfbox{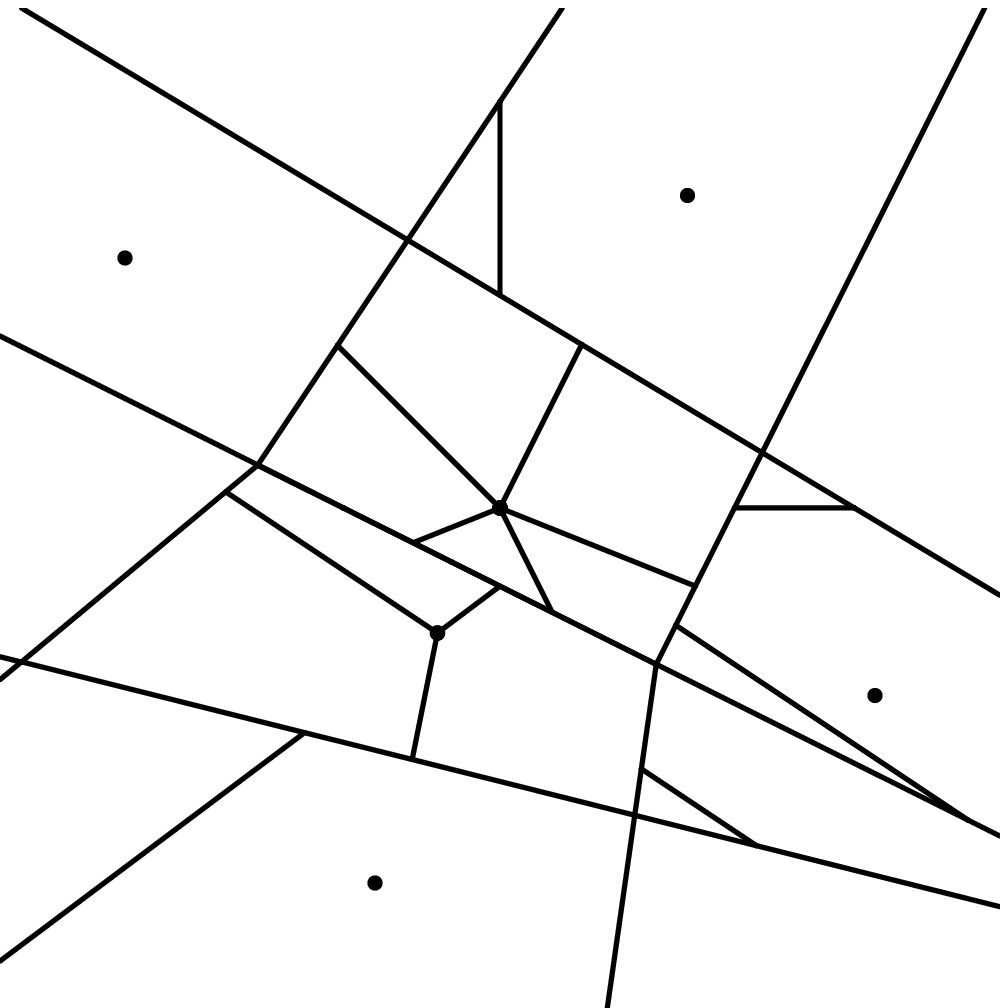}}}
\end{picture}
 \caption{\elabel{fvor2lim}No plugging of second order diagrams.}
\end{center}
\end{figure}

It is also clear that we cannot directly generalize Lemma \ref{lplug}
for obtaining the shape of a $k\text{-th}$ order Voronoi diagrams. 
That is because
the two steps approach, dealing with cluster locations and 
clusters apart, is not allowed
anymore for $k>1$. 

\begin{eexample}
Let $S(t)$ be as in  Example \ref{fvor2lim}.
In Figure \ref{fvor2lim} we show on the
left the second order Voronoi diagrams of the
cluster locations $l(0)$ that we have found in Example \ref{explug}. 
On the right, the second order Voronoi diagram for $S(t)$, with
$t=0.001$ is shown.
\end{eexample}

\subsection{Polynomial arithmetic.}

So far, we did not pay any attention to how to compute efficiently 
with polynomial sites: higher order terms
can often  be skipped without changing the type or the position of the 
Voronoi vertices and edges. But at this point we have no rule 
available that tells us on forehand when it is allowed to skip 
higher order terms.

\begin{eexample} 
Let $S(t)$ be as introduced in Example \ref{expolsites}.
Suppose that we change $q_4$ in Examples \ref{expolsites} and \ref{exvorpic}
from $q_4 = (t^2-2t^3 + t^5, -t^4)$ into $q_4=(t^2,-t^4)$.
Then the type of $S(t)$ at $t=0$ and the shape of the Voronoi diagram
of $S(t)$ at $t=0$ remains the same.
\end{eexample}

This is one reason for considering certain properties of sites or
point sets, rather than the sites themselves. We will do so 
in the following chapters.

\subsection{Conclusion}

The polynomial sites model enables us to introduce limit pictures
of Voronoi diagrams by means of the shape of a Voronoi diagram. 
Moreover, we have seen that this model can be used to
extend the notion of `type' to sets of points  that
sometimes coincide.  The type gives us the complete 
combinatorics of a limit situation. 

The methods presented in this chapter however do only
visualize the shapes of those cells that have positive
area at $t=0$. This demonstrates that not all cells
in the limit Voronoi diagram at $t=0$ are treated  
in an equal way.
The ingredient that is missing is the relation between shape
and scale. Informally, think of dividing the polynomials
defining the polynomial sites by powers of $t$ until
Voronoi cells that had zero area at $t=0$ before dividing
 get positive area.   That is, we want to zoom in 
at a cluster in order to find out the shape of a Voronoi cell of
some arbitrary site $p_i(t)$. In order to solve this problem
we will define 
{\it clickable} Voronoi diagrams in a following chapter. We do so by
exploiting  properties of point sets, like angles between two points 
or hooks between three points.

\chapter{Voronoi diagrams and angle compactifications.}
\elabel{changles}

Given a configuration $c$ of $n$ distinct points, one can determine the set
$a_n$ of angles of lines  through any two  points in $c$.
We distinguish $a_n$, the set of angles up to $2 \pi$,  
and $\overline{a_n}$, the angles up to $\pi$.
We analyze in both cases what sets of angles are possible. We show in which 
cases it is possible to reconstruct the
Voronoi diagram $V(c)$, knowing the angles only. 
We compactify the configuration space of distinct points
by taking the closure of the graph of 
the map that associates the angles to a configuration. 
We present a variety $T_n$ as an  algebraic alternative.
We analyze the connection between
boundary points of the compactification and singularities of $T_n$ for small $n$
and give geometric interpretations.

\section{$\cda$: space of angles on $n$ points.}

In this section we introduce several spaces that will be important 
to us later on. We recall the notion of configuration space 
of $n$ distinct labeled points in the plane. Both an 
introduction on and applications (in robotics!) of
configuration spaces can be found in \cite{AG}. Next we
define two spaces by considering, for $n$ distinct points in the 
plane, all angles between pairs of points. Here we distinguish
angles mod $\pi$ and angles mod $2 \pi$.

\subsection{The angle of two points.}

Given two distinct points $p_i$ and $p_j$, we determine the angle
that the line that passes through $p_i$ and $p_j$ makes with the
positive $x$-axis. We distinguish the directed and undirected line.

\begin{definition}
\elabel{dangle}
See also Figure \ref{falpha12}.
\begin{lijst}
\item For any two distinct points $p_i$ and $p_j$ in the plane 
define the \bfindex{directed angle} $\alpha_{ij} \in \R/ 2 \pi \Z$ as
the argument of the point
\beq
 ( \cos \alpha_{ij}, \sin \alpha_{ij} ) & = & 
		( \frac{d_x}{\|d_x\|}, \frac{d_y}{\|d_y\|} ),
\eeq
where $d = (d_x,d_y) := p_j-p_i$. Note that
$\alpha_{ij}$ is determined up to multiples of $2 \pi$.
\item The \bfindex{undirected angle} $\overline{\alpha_{ij}}$ is defined
as $\overline{\alpha_{ij}} = \overline{\alpha_{ji}}= \alpha_{ij} \mod \pi$.
We often choose $\overline{\alpha_{ij}} \in (-\frac{\pi}{2},\frac{\pi}{2}]$.
\item If $(p_j-p_i)_x \neq 0$, then the \bfindex{slope} $a_{ij}$ is defined as 
$a_{ij}=\frac{(p_j-p_i)_y}{(p_j-p_i)_x}$. If $(p_j-p_i)_x = 0$, then
$a_{ij}=\infty$. It holds that $a_{ij} \in (-\infty, \infty]$.  
\end{lijst}
\end{definition}

\begin{figure}[!ht]
\begin{center}
\setlength{\unitlength}{1cm}
\begin{picture}(6,2.1)
\put(3,1.1){\makebox(0,0)[cc]{
        \leavevmode\epsfxsize=6\unitlength\epsfbox{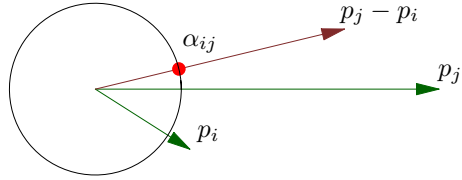}}}
\put(2.5,1.7){\makebox(0,0)[l]{$\alpha_{ij}$}}
\put(2.7,.5){\makebox(0,0)[l]{$p_i$}}
\put(5.9,1.3){\makebox(0,0)[l]{$p_j$}}
\put(4.6,2.1){\makebox(0,0)[l]{$p_j - p_i$}}
\end{picture}
\caption{\elabel{falpha12} The directed angle $\alpha_{ij}$
         between two points $p_i$ and $p_j$.}
\end{center}
\end{figure}

\begin{eexample}
Let $p_i=(1,1)$ and $p_j = (-1,0)$. Then $d=p_j-p_i=(-2,-1)$, so
$\frac{d_y}{d_x}=\frac{1}{2}$, but the angle 
$\alpha_{ij}=\text{arctan}(\frac{1}{2}) + \pi$ as $d_x < 0$.
\end{eexample}

\begin{remark}
\elabel{rtan}
If $(p_j-p_i)_x \neq 0$, then
$a_{ij} = \tan \overline{\alpha_{ij}}$ and 
$\overline{\alpha_{ij}} = \arctan a_{ij}$.
\end{remark}

\subsection{Configuration spaces of $n$ distinct points.}

\begin{definition} 
\elabel{dconf}
The \bfindex{configuration space} $\CONF( \R^2 )$ of $n$ points in $\R^2$
\index{CONF@$\CONF$} is the set
\beq
\CONF( \R^2 ) & := & \{ (p_1, \dots, p_n) \in (\R^2)^n \,\mid\,
	p_i \neq p_j ~\text{if}~ i \neq j\}.
\eeq
Two elements $p,q\in \CONF(\R^2)$
are $\sim_{\{s,t\}}$-equivalent if they only differ
by a scaling, $s$, combined with a translation, $t$. 
The \bfindex{reduced configuration space} is the quotient space
\index{Conf@$\conf$} $\conf =  \CONF / \sim_{\{s,t\}}$.
\end{definition}

$\conf$ is a smooth manifold. We determine its dimension. 
A \bfindex{standard representative} for a class in $\conf$ 
is constructed as follows. Translate
the configuration in such a way that $p_1=(0,0)$. Scale the configuration
in such a way that $p_2$ is at distance one from $p_1$. This shows that we need 
three parameters less than $2 n$ to describe an element in $\conf$.  

\begin{property}
dimension$( \conf )$ $=$  $2 n - 3$.
\end{property}

\begin{remark}
\elabel{rnotcomp}
Note that $\conf$ is not compact: take $n=3$,
then the ratio $\frac{|p_3-p_1|}{|p_2-p_1|}$ is continuous but not bounded.
\end{remark}

\begin{remark}
Let $\hat{c} = \{\hat{p}_1, \dots, \hat{p}_n \}$ be the
      standard representative of a class $[c]$ in $\conf$ and let
      $c=\{p_1, \dots, p_n \}$ be an arbitrary element of $[c]$.
      Then
  	$p_i  = \mid p_2 - p_1 \mid \hat{p_i} + p_1$.
     If we write $T_{\conf}$ for the set of standard
     representatives of classes in $\conf$, then
	$\CONF  =  \{ \hat{c}\R_{>0} + \R^2 \,\mid\, \hat{c} \in T_{\conf} \}$.
\end{remark}

\subsection{Compactification of the graph of the angle map.}

\elabel{scda}
For an element $c \in \CONF$, 
write down for every pair of points $p_i, p_j \in c$, with $i\neq j$,
the angle $\alpha_{ij} \in \R/ 2 \pi Z$ or $\overline{\alpha}_{ij} \in \R/ \pi \Z$.
This gives $\binom{n}{2}$ angles, one for every unordered  pair of labels $i$ and $j$.
\begin{definition}
$\text{DA}_n :=  (\R/  2 \pi \Z)^{\binom{n}{2}}$ \index{DAn@$\DA$}
 is the  space of directed angles; 
$\text{UA}_n :=  (\R/  \pi \Z)^{\binom{n}{2}}$ \index{UAn@$\UA$}
the space of undirected angles. \elabel{defsa}
The \index{directed angle!map}
directed {\bf angle map} $\pda$ is the map
\begin{displaymath}
\begin{array}{rccl}
 \pda: &  \CONF & \rightarrow & \DA, \\
  & (p_1, \dots, p_n)   & \mapsto &  ( \alpha_{ij})_{~1 \leq i < j \leq n}.
\end{array}
\end{displaymath}

The undirected angle map $\pua$ is defined in a similar way.
 \end{definition}

\begin{remark}
As a direct product of circles  $\DA$ and $\UA$ are smooth.
\end{remark}

\begin{remark}
$\pda$ and $\pua$ are well-defined  on $\conf$:
both mapping are constant on classes of $\conf$.
\end{remark}

\begin{definition}
\elabel{dcda}
The \bfindex{graph} of $\pda$ is the set
$ \{(c, \pda(c)) \in (\R^2)^n \times \DA \,\mid\, c \in \CONF \} $.
The \bfindex{compactification} \index{CDA@$\cda$}
$\cda$ of the graph
of the angle map, is the closure of 
$\text{graph}(\pda)$ in $(\R^2)^n \times \DA$. The compactification
$\cua$ is defined in a similar way.
 \index{CUA@$\cua$}
\end{definition}

\begin{remark}
$\cda$ is in fact not a compactification as it is not compact: 
$\CONF$ is not bounded, compare Remark \ref{rnotcomp}. 
The projection map from $\cda$ to $(\R^2)^n$ is proper however.
Recall that
a map is \bfindex{proper} if it is continuous and if the preimage 
of every compact set is compact again. 
\end{remark}

\section{From angles back to point configurations.}
\elabel{secreconstruct}

In this section we consider the mappings 
$\pda: \conf \rightarrow \DA$ and $\pua: \conf \rightarrow \UA$.
We describe the fibers of both mappings by
showing which configurations 
$c$ can be reconstructed from knowing the angles 
$\pda(c)$ or $\pua(c)$ only.

\subsection{Distinct points and angles in $\R/2 \pi \Z$.}

A configuration $c \in \CONF$ is called \bfindex{collinear} iff
all points in $c$ are collinear. A class of configurations
$[c] \in \conf$ is collinear if the class elements are collinear.
Define $\CLCONF  :=  \{ c \in \CONF \,\mid\, c ~\text{collinear} \}$,
and similarly $\clconf$.

\begin{lemma}
\elabel{linj}
The map $\pda: \conf \backslash \clconf  \rightarrow  \DA$
is injective.
\end{lemma}  

\begin{proof}
We have to show that for $[c],[d] \in \conf \backslash \clconf$ with
$[c] \neq [d]$ it holds that $\pda([c]) \neq \pda([d])$.
Let $c$ and $d$ be the standard representative of $[c]$, resp.\ $[d]$.
If $[c] \neq [d]$ then also $c \neq d$. Let $p_i(c)$ be the $i$-th
point of the configuration $c$ and let $\alpha_{ij}(c) \in \R/ 2 \pi \Z$
be the angle between $p_i(c)$ and $p_j(c)$. 
Note that for standard representatives $p_1(c)=p_1(d)=(0,0)$. 
We show that in any case there exists some labels $u$ and $v$
such that $\alpha_{uv}(c) \neq \alpha_{uv}(d)$.

Let $i \in \{ 2, \dots, n \}$ be minimal such
that $p_i(c) \neq p_i(d)$. If $i=2$, then $p_2(c) \neq p_2(d)$,
therefore $\alpha_{12}(c) \neq \alpha_{12}(d)$. Assume that
$i > 2$. In this case 
	$p_2(c) = p_2(d)= (\cos \alpha_{12}, \sin \alpha_{12})$ for some 
$\alpha_{12} \in (-\pi, \pi]$. Let $m$ be the line through 
$p_1$ and $p_i(c)$. If $p_i(d) \not\in m$, then 
$\alpha_{1i}(c)\neq \alpha_{1i}(d)$. Assume that
$p_i(d) \in m$. If $p_2 \not\in m$, then 
$\alpha_{2i}(c) \neq \alpha_{2i}(d)$.

We are left with the case that $p_1, p_2, p_i(c)$ and $p_i(d)$
are all collinear on line $m$. As $c,d \not\in \CLCONF$, there 
exists $j$ such that $p_j(c) \not\in m$. If $p_j(d)=p_j(c)$
then $\alpha_{ji}(c) \neq \alpha_{ji}(d)$. Let
$l$ be the line through $p_1$ and $p_j(c)$. If $p_j(d) \in l$,
then $\alpha_{2j}(c) \neq \alpha_{2j}(d)$. Finally, if
$p_j(d) \not\in l$, then $\alpha_{1j}(c) \neq \alpha_{1j}(d)$.
\end{proof}

\begin{remark}
\elabel{rconstruct}
Suppose we are given an image point $a=\pda([c])$ for some
$[c] \in \conf \backslash \clconf$. We construct the
standard representative from the angles in $a$
as follows. Put $p_1=(0,0)$ and 
$p_2 = (\cos \alpha_{12}, \sin \alpha_{12})$.
A point $p_i$ is on the line $l_{12}$ through $p_1$ and $p_2$ 
if and only if 
\begin{eqnarray}
\elabel{eqcol}
\overline{\alpha_{12}}=\overline{\alpha_{1i}} =
\overline{\alpha_{2i}} \mod \pi.
\end{eqnarray}
As $[c] \in \conf \backslash \clconf$,
we know that there exists $p_i$ that is not on the line $l_{12}$.
We find it by checking Equation \ref{eqcol}. We construct
$p_i$ as the intersection of the lines $l_{1i}$ and
$l_{2i}$. Here $l_{1i}$ is the line that passes
through $p_1$ and has direction $\alpha_{1i}$ and $l_{2i}$ is the line
that passes through $p_2$ with direction $\alpha_{2i}$.
Any other point $p_j$ of the standard representative is
constructed in a similar way: there is always a
pair of vertices from the non-degenerated triangle $p_1p_2p_i$ 
such that the line through these two vertices does not contain $p_j$.
\end{remark}

Let $c \in \CLCONF$. Then all points $p_1, \dots, p_n$ from $c$
are on a common line $l_c$ that makes some undirected angle
$\overline{\alpha_c} \in (-\frac{\pi}{2}, \frac{\pi}{2}]$ 
with the positive $x$-axis.
Order the labels of the points according to the order of the
points on $l_c$ as encountered from left to right, or from bottom
to top, in case $\overline{\alpha_c} = \frac{\pi}{2}$. This defines
an ordered $n$-tuple $\sigma(c)$. The $n$-tuple
$\sigma([c])$ is defined as $\sigma(c)$, for some representative
$c \in [c]$.  Define an 
equivalence class $\sim_{\alpha,\sigma}$ on $\clconf$ as follows:
$[c] \sim_{\alpha,\sigma} [d]$ if and only if both
$\overline{\alpha_{[c]}} = \overline{\alpha_{[d]}}$, and
$\sigma([c])=\sigma([d])$.

\begin{lemma} \elabel{lclconf} Consider $\pda: \CLCONF \rightarrow \DA$. 
\begin{lijst}
\item $\pda$ is constant on classes of 
$\sim_{\alpha, \sigma}$.
\item The map
$\pda: \clconf / \sim_{\alpha,\sigma}  \rightarrow  \DA$
is injective
\end{lijst}
\end{lemma}

\begin{proof} We prove the two claims.
\begin{lijst}
\item Let $[c], [d] \in \clconf$ such that $[c] \sim_{\alpha,\sigma} [d]$.
Let $p_i(c), p_i(d)$ and $p_j(c), p_j(d)$ denote the $i$-th and $j$-th
points of some representatives $c$ and $d$ of the classes $[c]$
and $[d]$. Then
$\overline{\alpha_{ij}(c)} = \overline{\alpha_{ij}(d)}$, as 
all four points are collinear. Moreover, $p_j(c)$ is on the right of
$p_i(c)$ whenever $p_j(d)$ on the right of $p_i(d)$, as
$\sigma(c)=\sigma(d)$. So, $\alpha_{ij}(c) = \alpha_{ij}(d)$.
\item Suppose that $[c] \not\sim_{\alpha, \sigma} [d]$. Let $p_i(c), p_i(d),
p_j(c)$ and $p_j(d)$ be as above. 
If $\overline{\alpha_c} \neq \overline{\alpha_d}$, then clearly 
$\pda([c]) \neq \pda([d])$. Assume $\overline{\alpha_c} =\overline{\alpha_d}$.
Then $\sigma(c) \neq \sigma(d)$. Let $i$ be minimal such that
$k := \sigma(c)_i \neq \sigma( d)_i =: l$
Then $\alpha_{kl}(c) = \alpha_{kl}(d) + \pi$. Therefore
$\pda( [c] ) \neq \pda( [d] )$.
\qedhere
\end{lijst}
\end{proof}

We use Lemma \ref{linj} and Lemma \ref{lclconf} to answer the following
question. Is it possible, given some $a = \pda( c )$, with 
$c \in \CONF$, to determine the Voronoi diagram $V(c)$? 
Let the \bfindex{reduced Voronoi diagram}, notation $\tilde{V}(c)$,
of a point $c \in \CONF$ be the Voronoi diagram $V(c)$ up to scaling and
translation. A reduced configuration of points $[c]$ or a reduced Voronoi
diagram $\tilde{V}(c)$ is \bfindex{reconstructible} from a list of angles $a$
iff $[c]$ resp.\ $\tilde{V}(c)$ is uniquely determined by $a$. 

\begin{property}
\elabel{pconstruct}
$\tilde{V}(c)$ is independent of the choice of $c \in [c]$.
\end{property}

Notation: By $\tilde{V}([c])$, for $[c] \in \conf$, the reduced Voronoi
diagram of some representative $c \in [c]$ is indicated.
Note that if one can reconstruct the class $[c] \in \conf$ from
$\pda(c)$, for $c \in \CONF$, then one can also construct the reduced
Voronoi diagram $\tilde{V}([c])$.

\begin{corollary}
\elabel{crec}
Let $a = \pda( [c] )$, with $[c] \in \conf$.
\begin{lijst}
\item If $[c] \in \conf \setminus \clconf$, then
$[c]$ is reconstructible.
\item If $[c] \in \clconf$, then $[c]$ is not reconstructible,
except for $n=2$. 
\end{lijst}
\end{corollary}

\begin{proof}
The first claim follows from
Lemma \ref{linj} and Remark \ref{rconstruct}. The preimage
 $\pda^{-1}([c])$ consists of distinct classes of $\conf$,
      except for $n=2$, consult Lemma \ref{lclconf}.
This proves the second claim.
\end{proof}

\begin{corollary}
$\tilde{V}([c])$ is reconstructible if and only if 
either $[c] \in \conf \setminus \clconf$, or 
$[c] \in \clconf$ and $n=2$ or $n=3$.
\end{corollary}

\begin{proof} 
From Corollary \ref{crec} it follows that $[c]$ is reconstructible in the 
cases mentioned except
for $[c] \in \clconf$ in case $n=3$. But in this case,
$\tilde{V}([c])$ consists just of two parallel bisectors perpendicular
to $\overline{\alpha_c}$. From $a_3 = \pda( [c] )$
we can determine the order of the three points on the line, which
gives us the labels of the bisector. Changing the relative positions 
of the three points in $[c]$ has no influence  on
$\tilde{V}(c)$ as long as the order of the points is maintained.
Suppose, on the other hand that $[c] \in \clconf$ for $n>3$.
It is easy to show that in this case mutually distinct classes
of configurations in the fiber $\pda^{-1}([c])$ correspond with
mutually distinct reduced Voronoi diagrams. 
\end{proof}

\begin{figure}[!ht]
\begin{center}
\setlength{\unitlength}{.8em}
\begin{picture}(24,5.5)
\put(12,3){\makebox(0,0)[cc]{
        \leavevmode\epsfxsize=24\unitlength\epsfbox{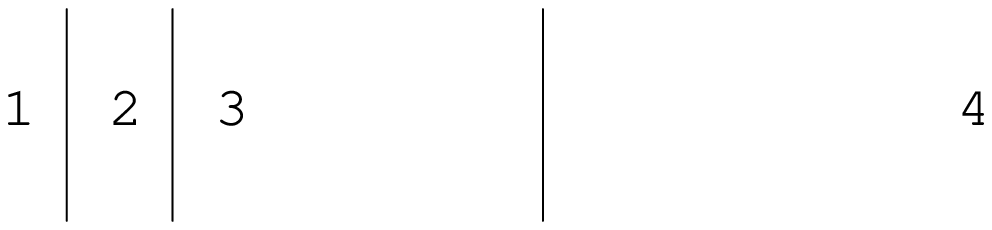}}}
\end{picture}
\caption{\elabel{ffourlinear}The Voronoi diagram of
 four collinear points.}
\end{center}
\end{figure}

\begin{eexample}
\elabel{exfourcol}
Consider the element $d \in {\mbox {\sl DA}_4}$ defined by
$\alpha_{12}=\alpha_{13}=\alpha_{14}=\alpha_{23}=\alpha_{24}=0$.
Suppose that we try to construct a configuration $S$ of labeled points
having those angles. It is clear that any such configuration 
consists of four points on a common horizontal line.
An example is given in Figure \ref{ffourlinear}.
Moreover, it is possible to reconstruct the order of the points
on the line: from left to right we encounter the points
$p_1$, $p_2$, $p_3$ and $p_4$ in that particular order.
This shows however that from the information present in $d$
we cannot determine the ratio $|p_1p_2|/|p_1p_3|$ of the lengths of the
line segments $p_1p_2$ and $p_1p_3$: whatever this ratio is, 
the angles $\alpha_{12}$, $\alpha_{13}$ and $\alpha_{23}$
do not change. The Voronoi diagram  of $p_1$ to $p_4$ consists
of the three vertical bisectors $B(p_1,p_2)$, $B(p_1, p_3)$
and $B(p_2,p_3)$. Again it is impossible to determine the ratio
$|B(p_1,p_2)-B(p_2,p_3)| / |B(p_2,p_3)-B(p_3,p_4)|$ of the distances
between the bisectors.
This shows that it is  impossible to associate a unique
reduced Voronoi diagram to $d \in {\mbox {\sl DA}_4}$.  
\end{eexample}

\subsection{Distinct points and angles in $\R/\pi \Z$.}

Suppose we start with some $c \in \CONF$. Assume we 
have determined $\overline{a} = \pua(c)$, that is, the set
of all angles $\overline{\alpha_{ij}} \mod \pi$ between pairs of
points $(p_i,p_j) \in c$.
We introduce another equivalence class on $\conf$. Let $[c],[d]$ in $\conf$.
The classes $[c]$ and $[d]$ are \bfindex{reflection equivalent},
notation $[c] \sim_R [d]$ if and only if the standard representative $c$ of
$[c]$ equals the standard representative $d$ of $[d]$ up to a reflection
in $p_1(c)=p_1(d)=(0,0)$.

\begin{lemma} The map $\pua:\, (\conf \setminus \clconf) / \sim_R  \,\rightarrow\,  \UA$
is injective
\end{lemma}

\begin{proof}
Proceed as in the proof of Lemma \ref{linj}, but take those 
standard representatives that have $\alpha(c), \alpha(d) 
	\in (-\frac{\pi}{2}, \frac{\pi}{2}]$. It is possible to put mod $\pi$
bars on all $\alpha_{ij}$'s that occur in the proof. 
\end{proof}

If we take angles between points in $\R / \pi \Z$, 
all collinear configurations that have their points on a line with angle 
$\overline{\alpha} \in (-\frac{\pi}{2},\frac{\pi}{2}]$ are mapped
to the same $\overline{a} = (\overline{\alpha}, \dots, \overline{\alpha})$
in $\UA$. So it is impossible to reconstruct the original order of the 
points on the line. 

\begin{corollary}
 Let $\overline{a} = \pua( [c] )$, for $[c] \in \conf$.
\begin{lijst}
\item If $[c] \in \conf \backslash \clconf$, then $[c]$ is reconstructible
up to a point reflection. 
\item If[ $[c] \in \clconf$ and $n=2$, then $[c]$ is reconstructible up to 
      a point reflection.
\item If[ $[c] \in \clconf$ and $n=3$, then $[c]$ is not reconstructible but the
       reduced Voronoi diagram $\tilde{V}([c])$ is reconstructible 
\item If[ $[c] \in \clconf$ and $n>3$, then $\tilde{V}([c])$ 
      is not reconstructible.
\end{lijst}
\end{corollary}

\begin{figure}[!ht]
\begin{center}
\setlength{\unitlength}{1.1em}
\begin{picture}(15,5.2)
\put(7.5,3){\makebox(0,0)[cc]{
        \leavevmode\epsfxsize=15\unitlength\epsfbox{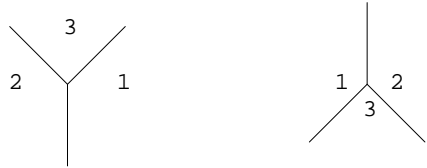}}}
\end{picture}
\caption{\elabel{frotate}$V(c)$ reconstructed from angles mod $\pi$.}
\end{center}
\end{figure}

\begin{eexample}
Let $c=\{p_1,p_2,p_3\}$, with $p_1=(2,0)$, $p_2=(0,0)$ and $p_3=(1,1)$.
The Voronoi diagram $V(c)$ is shown on the left in Figure \ref{frotate}.
Then $\alpha_{12}=\pi$, $\alpha_{13}=\frac{3 \pi}{4}$ and
$\alpha_{23}=\frac{\pi}{4}$. So $\overline{\alpha_{12}}=0$, while
$\overline{\alpha_{13}}=\alpha_{13}$ and $\overline{\alpha_{23}}=\alpha_{23}$.
We get a `reconstructed' configuration $c'$ consisting of $p_1'=(0,0)$, $p_2'=(1,0)$
and $p_3'=(-\frac{1}{2},-\frac{1}{2})$. This configuration $c'$ together with
its Voronoi diagram $V(c')$ is shown on the right in Figure \ref{frotate}.
\end{eexample}

\section{Angle models for small $n$.}

In this section, we analyze the compactifications $\cda$ and $\cua$
for $n=2$. Moreover, we give a complete description of the
image  $\pdat( \CONFt )$. For points on the boundary 
of $\pdat( \CONFt )$, we give a geometric interpretation in terms
of coinciding points.

\begin{remark}
D. G. Kendal has introduced the `The theory of shape' in a statistical
context. As a specific example, the shapes of triangles are analyzed.
To compare with this result, consult \cite{Ke1} and \cite{Ke2}. 
\end{remark}

\subsection{n=2.}

In the directed case,  we have the following diagram:

\begin{displaymath}
\begin{array}{ccccc}
 {\mbox {\sl CONF}}_2 & & \quad\subset\quad &  &
	 \R^2 \times \R^2\\[.1cm] 
	\cap & & &  & \cap \\[.1cm]
 \text{graph}( \psi_{DA_2} ) & \subset & {\mbox {\sl CDA}_2} & \subset  &
	\R^2 \times \R^2 \times \R/2 \pi \Z.
\end{array}
\end{displaymath}

$\psi_{DA_2}$ maps a configuration $c = (p_1,p_2) \in \CONF$ 
to the directed angle $\alpha_{12}$. 
For simplicity we assume  that $p_1 = (0,0)$.
As $c \in {\mbox {\sl CONF}}_2$, this implies that 
$p_2 \in \R^2 \setminus \{(0,0)\}$.
This space $\R^2 \setminus \{(0,0)\}$ is homeomorphic
to a doubly open cylinder: imagine the omitted point $(0,0)$
as one side of the cylinder and infinity as the other side.

A classic construction is the following:
The \bfindex{blow-up} of $\R^2$ at $(0,0)$ is
by definition the closed subset $X \subset \R^2 \times \P^1$ of
all points $(x_1,x_2,y_1,y_2) \in \R^2 \times \P^1$ such that 
$x_1 y_2 = x_2 y_1$. The projection $\pi: X \rightarrow \R^2$
onto the first factor has the following properties, compare \cite{Ha}.

\begin{looplijst}
\item $\pi^{-1}( p_2 )$ consists of a single point, 
        if $p_2 \in R^2 \setminus \{(0,0)\}$.
\item $\pi^{-1}( (0,0) ) = \P^1$.
\item The points of $\pi^{-1}( (0,0) )$ are in 1-1 correspondence
      to the set of (undirected) lines through $(0,0)$.
\end{looplijst}

We consider ${\mbox {\sl CDA}_2}$, the closure of 
$\text{graph}( \psi_{{\text{DA}}_2} )$.

\begin{lemma}
\elabel{cda2}
${\mbox {\sl CDA}_2}$ is homeomorphic to a half-open cylinder times a plane.
\end{lemma}

\begin{proof}[{\bf Proof. (sketch)}]
The plane comes from varying $p_1$, so 
assume $p_1 = (0,0)$. Then $p_2 \in \R^2 \setminus \{(0,0)\}$,
provided that $(p_1, p_2) \in {\mbox {\sl CONF}}_2$. 
Write $p_2  = r ( \cos \alpha, \sin \alpha )$,
with $r \in \R_{>0}$, and $\alpha \in \R / 2 \pi \Z$.
Any configuration $((0,0),\alpha)$, is added exactly once in order
to obtain the closure of $\text{graph}( \psi_{{\text{DA}}_2} )$,
as $\lim_{r \to 0} r ( \cos \alpha, \sin \alpha)$.
The punctured plane $\R^2 \setminus \{(0,0)\}$, is homeomorphic to
a doubly open cylinder. Adding all points of the form
$((0,0),\alpha_{12})$, with $\alpha_{12} \in \R / 2 \pi \Z$
means that we attach one full circle to that end of
the open cylinder that corresponds to $p_2 = (0,0)$. 
\end{proof}

\begin{remark}
Lemma \ref{cda2} demonstrates that $\cda$ can have a `boundary'. 
\end{remark}

\begin{figure}[!ht]
\begin{center}
\setlength{\unitlength}{1.0cm}
\begin{picture}(10,3.6)
\put(5,1.8){\makebox(0,0)[cc]{
        \leavevmode\epsfxsize=10\unitlength\epsfbox{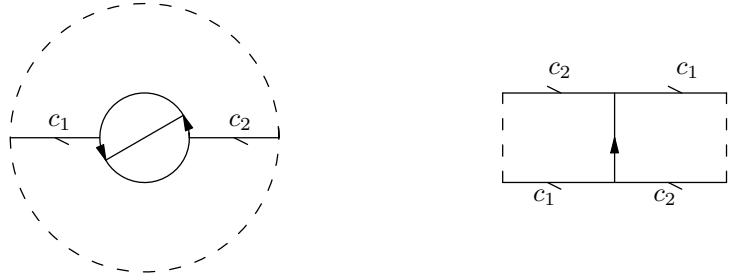}}}
\put(3.2,2.0){\makebox(0,0)[l]{$c_2$}}
\put(.8,2.0){\makebox(0,0)[l]{$c_1$}}
\put(7.25,1.0){\makebox(0,0)[l]{$c_1$}}
\put(8.85,1.0){\makebox(0,0)[l]{$c_2$}}
\put(7.45,2.7){\makebox(0,0)[l]{$c_2$}}
\put(9.15,2.7){\makebox(0,0)[l]{$c_1$}}
\end{picture}
\caption{\elabel{fantipodal} Antipodal identification.}
\end{center}
\end{figure}

In the undirected case, we consider ${\mbox {\sl CUA}_2}$, the closure of
$\text{graph}( \psi_{{\text{UA}}_2} )$.

\begin{lemma}
\elabel{cua2}
${\mbox {\sl CUA}_2}$ is homeomorphic to a M\"obius strip times a plane.
\end{lemma}

\begin{proof}[{\bf Proof. (sketch)}]
Start as in the proof of Lemma \ref{cda2}. 
Write $p_2$ in polar coordinates, that is, $p_2 = r ( \cos \alpha, \sin \alpha )$,
with $r \in \R_{>0}$, and $\alpha \in \R / 2 \pi \Z$.
Any undirected angle $\alpha_{12}$ is added twice, once
for $p_2 = \lim_{r \to 0} r(\cos \alpha, \sin \alpha)$ and once for
$p_2 = \lim_{r \to 0} r(\cos \alpha + \pi, \sin \alpha + \pi)$. 
This gives the antipodal identification as indicated on the left in Figure 
\ref{fantipodal} by the little solid arrows. 
It is allowed to cut the space as long as it is eventually
pasted back together in the same way as it was cut. We make two cuts,
$c_1$ and $c_2$, through the side of the cylinder. This enables us to
perform the antipodal paste, resulting in the rectangle on the
right of Figure \ref{fantipodal}. Pasting back the cuts we have made,
results in the M\"obius strip, compare e.g.\ \cite{Mu}. 
\end{proof} 

\begin{remark}
Note that the construction in the proof of Lemma \ref{cua2} is exactly
the blow-up of $\R^2$ at $(0,0)$ that we have discussed above.
\end{remark}

\subsection{n=3.}

We consider the possible angles mod $2 \pi$  between three distinct points 
$p_1, p_2$ and $p_3$ in the plane. That is, we determine the image
$\psi_{{\text{DA}}_3}({\mbox {\sl CONF}_3})$ in ${\text{DA}}_3$.
  Let $l_{ij}$
denote the directed line that passes first through $p_i$ and last 
through $p_j$. The triangle with vertices 
$p_1$, $p_2$ and $p_3$ is \bfindex{oriented clockwise}
if $p_3$ is on the right of $l_{12}$. 
The triangle is oriented anti-clockwise
if $p_3$ is on the left of $l_{12}$.

\begin{lemma}
\elabel{lposangle}
Let $c=(p_1,p_2,p_3) \in {\mbox {\sl CONF}_3}$.
Let $\triangle$ denote the triangle with vertices $p_1$, $p_2$ and $p_3$.
For $(\alpha_{12}, \alpha_{13}, \alpha_{23}) = \psi_{{\text{DA}}_3}(c)$ 
the following holds.
\begin{lijst}
\item \label{big}
       If $\triangle$ oriented clockwise, then 
       $0 < \alpha_{13}-\alpha_{23}
	 < \alpha_{13} - \alpha_{21} < \pi$  $\mod 2 \pi.$ 
\item If $\triangle$ oriented anti-clockwise, then 
$ 0 < \alpha_{21}-\alpha_{23} < \alpha_{21} - \alpha_{13} < \pi$ $\mod 2 \pi.$ 
\item If $p_1$, $p_2$ and $p_3$ collinear then 
	$$ \alpha_{12} = \alpha_{13} = \alpha_{23} ~\vee~ 
	   \alpha_{12} = \alpha_{13} = \alpha_{32} ~\vee~
	   \alpha_{21} = \alpha_{13} = \alpha_{23}. $$
\end{lijst}
Moreover, for any image point $(\alpha_{12},\alpha_{13}, \alpha_{23}) 
	\in \psi_{{\text{DA}}_3}({\mbox {\sl CONF}_3})$ either
        (i), (ii), or (iii) holds. 
\end{lemma}

\begin{figure}[!ht]
\begin{center}
\setlength{\unitlength}{1.0cm}
\begin{picture}(12,2.4)
\put(6,1.3){\makebox(0,0)[cc]{
        \leavevmode\epsfxsize=12\unitlength\epsfbox{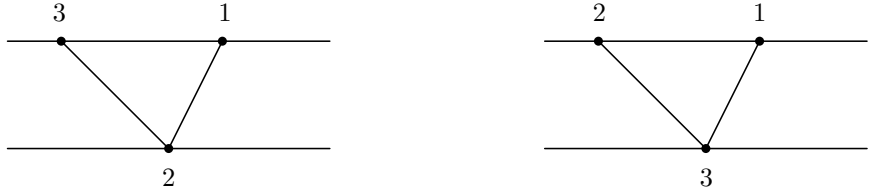}}}
\put(3.15,2.4){\makebox(0,0)[l]{$1$}}
\put(2.4,.2){\makebox(0,0)[l]{$2$}}
\put(0.95,2.4){\makebox(0,0)[l]{$3$}}
\put(10.26,2.4){\makebox(0,0)[l]{$1$}}
\put(8.13,2.4){\makebox(0,0)[l]{$2$}}
\put(9.55,.2){\makebox(0,0)[l]{$3$}}
\end{picture}
\caption{\elabel{ftricompare} The clockwise and anti-clockwise case.}
\end{center}
\end{figure}

\begin{proof}
Note that the differences of angles as they occur in the statements in the lemma
does not change if $\triangle$ is rotated. 
\begin{lijst}
\item Assume $\triangle$ is oriented clockwise. Rotate $\triangle$ such that
$\alpha_{13}=\pi$, compare the triangle on the left in Figure \ref{ftricompare}.
As $p_1$ is on the right of
$l_{23}$ it follows that $\alpha_{21} < \alpha_{23} < \pi$.
\item Assume $\triangle$ is oriented anti-clockwise. Rotate $\triangle$ such that
$\alpha_{12}=\pi$,  compare the triangle on the right in Figure \ref{ftricompare}.
As $p_1$ is on the right of $l_{32}$ it follows that $\alpha_{31}<\alpha_{32}<\pi$.
It follows that $\alpha_{12}-\alpha_{32}<\alpha_{12}-\alpha_{31}<\pi \mod 2 \pi$, or
equivalently, that $\alpha_{21}-\alpha_{23} < \alpha_{21}-\alpha_{13}<\pi \mod 2 \pi$.
\item Proof by inspection.
\end{lijst}
As any triangle on the vertices $p_1, p_2$ and $p_3$ is either oriented 
anti-clockwise, or oriented clockwise or degenerate, the lemma follows.
\end{proof}

Boundary points of
$\psi_{{\text{DA}}_3}({\mbox {\sl CONF}_3})$ are characterized as follows.

\begin{corollary} 
Let $a_3 = (\alpha_{12}, \alpha_{13},\alpha_{23})$ be a boundary 
point of $\psi_{{\text{DA}}_3}({\mbox {\sl CONF}_3})$.
\elabel{cboundconf3}
\begin{lijst}
\item If $a_3 = \psi_{{\text{DA}}_3}(c)$ for $c \in {\mbox {\sl CONF}_3}$,
then $\overline{\alpha_{12}}=\overline{\alpha_{13}}= \overline{\alpha_{23}}$.
\item If $a_3 \not\in \psi_{{\text{DA}}_3}({\mbox {\sl CONF}_3})$, then
$ \alpha_{12} = \alpha_{13}$ $\vee$ $\alpha_{13} = \alpha_{23}$ $\vee$
                $\alpha_{23} = \alpha_{12} + \pi.$
\end{lijst}
\end{corollary}

\begin{figure}[!ht]
\begin{center}
\setlength{\unitlength}{1.06cm}
\begin{picture}(12,4.4)
\put(6,2.2){\makebox(0,0)[cc]{
        \leavevmode\epsfxsize=12\unitlength\epsfbox{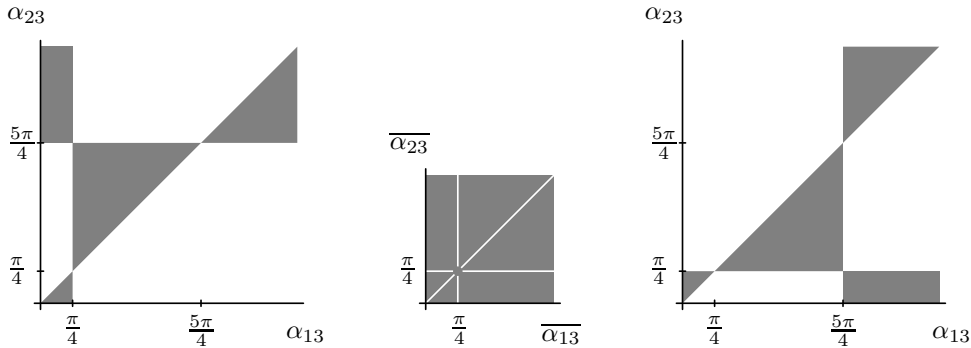}}}
\put(0,1.0){\makebox(0,0)[l]{$\frac{\pi}{4}$}}
\put(4.9,1.0){\makebox(0,0)[l]{$\frac{\pi}{4}$}}
\put(8.05,1.0){\makebox(0,0)[l]{$\frac{\pi}{4}$}}
\put(0,2.6){\makebox(0,0)[l]{$\frac{5\pi}{4}$}}
\put(8.05,2.6){\makebox(0,0)[l]{$\frac{5\pi}{4}$}}
\put(.7,0.25){\makebox(0,0)[l]{$\frac{\pi}{4}$}}
\put(2.25,0.25){\makebox(0,0)[l]{$\frac{5\pi}{4}$}}
\put(5.55,0.25){\makebox(0,0)[l]{$\frac{\pi}{4}$}}
\put(8.75,0.25){\makebox(0,0)[l]{$\frac{\pi}{4}$}}
\put(10.3,0.25){\makebox(0,0)[l]{$\frac{5\pi}{4}$}}
\put(3.5,0.2){\makebox(0,0)[l]{$\alpha_{13}$}}
\put(6.7,0.2){\makebox(0,0)[l]{$\overline{\alpha_{13}}$}}
\put(11.6,0.2){\makebox(0,0)[l]{$\alpha_{13}$}}
\put(0,4.2){\makebox(0,0)[l]{$\alpha_{23}$}}
\put(4.8,2.6){\makebox(0,0)[l]{$\overline{\alpha_{23}}$}}
\put(8,4.2){\makebox(0,0)[l]{$\alpha_{23}$}}
\end{picture}
\caption{\elabel{fmod2pmod1p}Possible values for $\alpha_{13}$,
$\overline{\alpha_{13}}$, $\alpha_{23}$ and $\overline{\alpha_{23}}$, 
in case $\alpha_{12}=\frac{\pi}{4}$ and $\alpha_{12}=\frac{5 \pi}{4}$.}
\end{center}
\end{figure}

\begin{eexample}
\elabel{exposcda3}
Figure $\ref{fmod2pmod1p}$ shows the possible values of $\alpha_{13}$ and
$\alpha_{23}$, for $\alpha_{12} = \frac{\pi}{4}$, on the left,
and, for $\alpha_{12}=\frac{3\pi}{4}$, on the right. Note that
\begin{lijst}
\item Configurations of points that are not collinear
      are mapped by $\psi_{{\text{DA}}_3}$ 
      to the interior of the triangles and rectangles
      in the figure.
\item Configurations that are collinear are mapped to the 
      vertices of the big central triangle on the left and on the right.
\item Other points on the boundary of the triangles and rectangles  do not 
      correspond to configurations in $\text{{\sl {\bf conf}}}_3$.
\end{lijst}
\end{eexample}

\begin{eexample}
\elabel{exua3}
The picture in the middle of Figure \ref{fmod2pmod1p} shows
 the possible values of
$\overline{\alpha_{13}}$ and $\overline{\alpha_{23}}$ for
$\overline{\alpha_{12}}=\frac{\pi}{4}$. The only configurations
in ${{\mbox {\sl UA}_3}}$ that are not image points of
$\psi_{{\text{UA}}_n}( {\mbox {\sl CONF}_3} )$ are the points 
where two angles mod $\pi$  coincide  but not all three.
\end{eexample}

Recall that $\confd \setminus \text{{\sl {\bf clconf}}}_3$
consists of the configurations of three  distinct, non-collinear points
up to scalings and transformations.

\begin{corollary}
Let $a_3 = (x, y, z)$. Then
$a_3 = \pdat([c])$ for some  
$[c] \in \confd \setminus \text{{\sl {\bf clconf}}}_3$
if and only if
 \begin{displaymath}
 \begin{array}{llllll}
 &0<x<\pi; &  y<x; &  y>z; &  z>0, \\
 \vee &  0<x<\pi; &  y>x; &  z<x+\pi; &  y<z, \\
 \vee &    0<x<\pi; &  y>0; &  y<x; &  z>x+\pi; &  z<2 \pi, \\
 \vee &    0<x<\pi; &  z>\pi+x; &  y>z; &  y<2 \pi, \\
 \vee &    \pi < x < 2 \pi; &  y>0; &  y<z; &  z<x-\pi, \\
 \vee &    \pi < x < 2 \pi; &  y>x; &  z>0; &  y<2 \pi; &  y<x-\pi, \\
 \vee &    \pi < x < 2 \pi; &  y>z; &  z>x-\pi; &  y<x, \\
 \vee &    \pi < x < 2 \pi; &  y>x; &  y<z; &  z<2 \pi. &
 \end{array}
 \end{displaymath}
The boundary of $\pdat([c])$ with 
$[c] \in \confd \setminus \text{{\sl {\bf clconf}}}_3$
is depicted in Figure \ref{fthreed}.
\end{corollary}

\begin{figure}[!ht]
\begin{center}
\setlength{\unitlength}{1.05cm}
\begin{picture}(5.6,6.6)
\put(2.8,3){\makebox(0,0)[cc]{
        \leavevmode\epsfxsize=5\unitlength\epsfbox{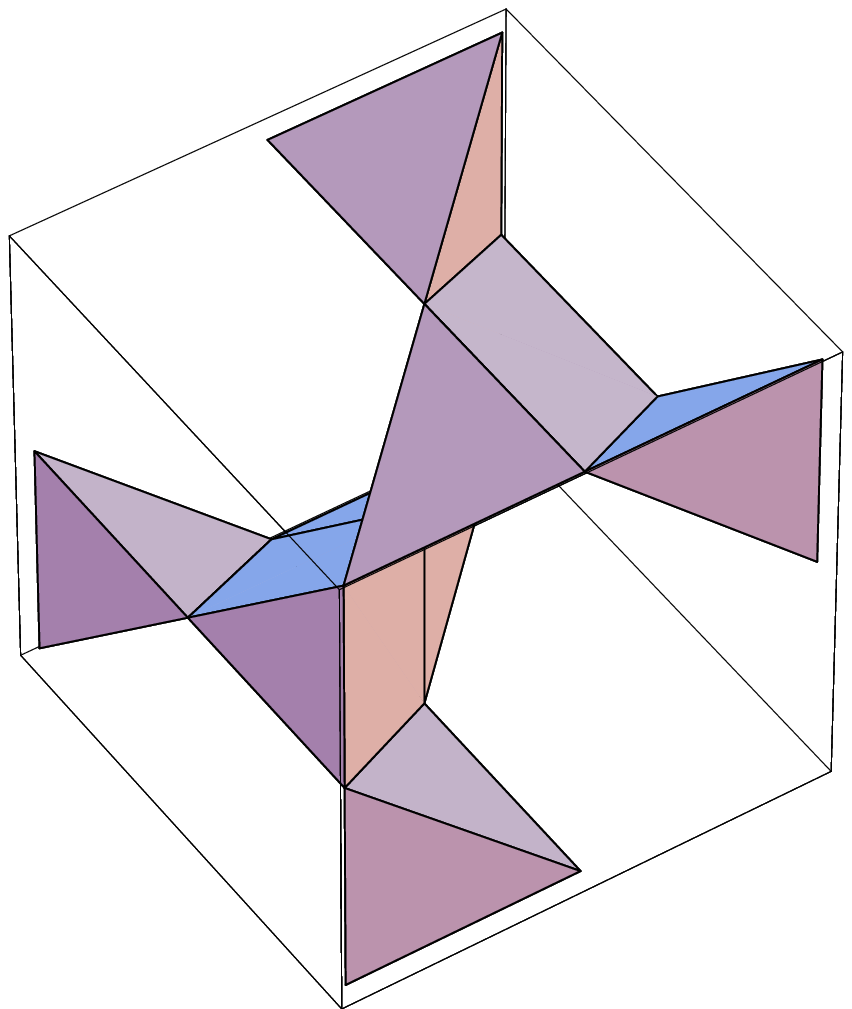}}}
\put(0.05,5.05){\makebox(0,0)[l]{$2 \pi$}}
\put(1.5,5.55){\makebox(0,0)[l]{$\pi$}}
\put(3.2,6.35){\makebox(0,0)[l]{$0$}}
\put(4.4,5.1){\makebox(0,0)[l]{$\pi$}}
\put(5.45,4.0){\makebox(0,0)[l]{$2\pi$}}
\put(5.45,2.7){\makebox(0,0)[l]{$\pi$}}
\put(5.45,1.4){\makebox(0,0)[l]{$0$}}
\put(1.,6.0){\makebox(0,0)[l]{$\alpha_{12}$}}
\put(4.8,5.4){\makebox(0,0)[l]{$\alpha_{13}$}}
\put(5.9,2.8){\makebox(0,0)[l]{$\alpha_{23}$}}
\end{picture}
\caption{
The boundary of $\pdat([c])$ with $[c] \in 
\confd \setminus \text{{\sl {\bf clconf}}}_3$.
\elabel{fthreed}}
\end{center}
\end{figure}

\begin{remark} In Corollary \ref{cboundconf3} we have classified
those points $a_3 \in {{\mbox {\sl DA}_3}}$ that are on the boundary 
of $\pdat(\CONFt)$ but not in $\pdat(\CONFt)$ itself. These 
triples of angles can be interpreted as the triples of angles
that correspond to configurations of three points, such that
exactly two points coincide. Consider the
angles between the points in the three configurations
presented in  Figure \ref{ftogether}. In the leftmost 
configuration, $\alpha_{12} = \alpha_{13}$, in the configuration
in the middle, $\alpha_{21}=\alpha_{23}$, while in the rightmost
configuration, $\alpha_{31}=\alpha_{32}$. This interpretation holds,
as long as the three points are not collinear, compare the first 
statement of Corollary \ref{cboundconf3}.
\end{remark}

\begin{figure}[!ht]
\begin{center}
\setlength{\unitlength}{1.05cm}
\begin{picture}(12,1.8)
\put(6,0.9){\makebox(0,0)[cc]{
        \leavevmode\epsfxsize=12\unitlength\epsfbox{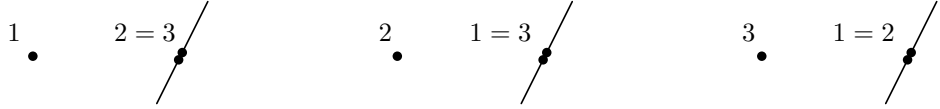}}}
\put(.2,1.2){\makebox(0,0)[l]{$1$}}
\put(4.9,1.2){\makebox(0,0)[l]{$2$}}
\put(9.5,1.2){\makebox(0,0)[l]{$3$}}
\put(1.55,1.2){\makebox(0,0)[l]{$2=3$}}
\put(6.05,1.2){\makebox(0,0)[l]{$1=3$}}
\put(10.65,1.2){\makebox(0,0)[l]{$1=2$}}
\end{picture}
\caption{\elabel{ftogether} Configurations for two coinciding points.}
\end{center}
\end{figure}

\begin{eexample}
 Fix $p_1 = (0,0)$ and $p_2 = (1,1)$. This assures that 
$\alpha_{12} = \frac{\pi}{4}$. Suppose that the position of
$p_3$ is given by $p_3(t)=p_2 + 1.1( \cos t, \sin t)$, for
$t \in [0, 2 \pi)$. So $\alpha_{23} = t$. 
Then $(\alpha_{13}(t), t)$ 
is a curve, that is defined for any  $t$. It is depicted in
Figure \ref{ftravelp3}, on the left. 
\end{eexample}

\begin{eexample}
 Fix $p_1 = (0,0)$ and $p_2 = (1,1)$ again. Let 
$p_3(t) = \sqrt{2} (\cos t, \sin t)$. In this case, the curve
$(\alpha_{13}(t),t)$ is built up out of straight line segments,
see Figure \ref{ftravelp3}, on the right. 
At $t= \frac{5 \pi}{4}$, the point $p_3$ passes through $p_1$.
At this $t$, the curve is not defined.
This $t$ coincides with a swap $\alpha_{13}=\alpha_{13} + \pi$.
\end{eexample}

\begin{figure}[!ht]
\begin{center}
\setlength{\unitlength}{1.03cm}
\begin{picture}(12.3,3.7)
\put(2.25,2){\makebox(0,0)[cc]{
        \leavevmode\epsfxsize=3.5\unitlength\epsfbox{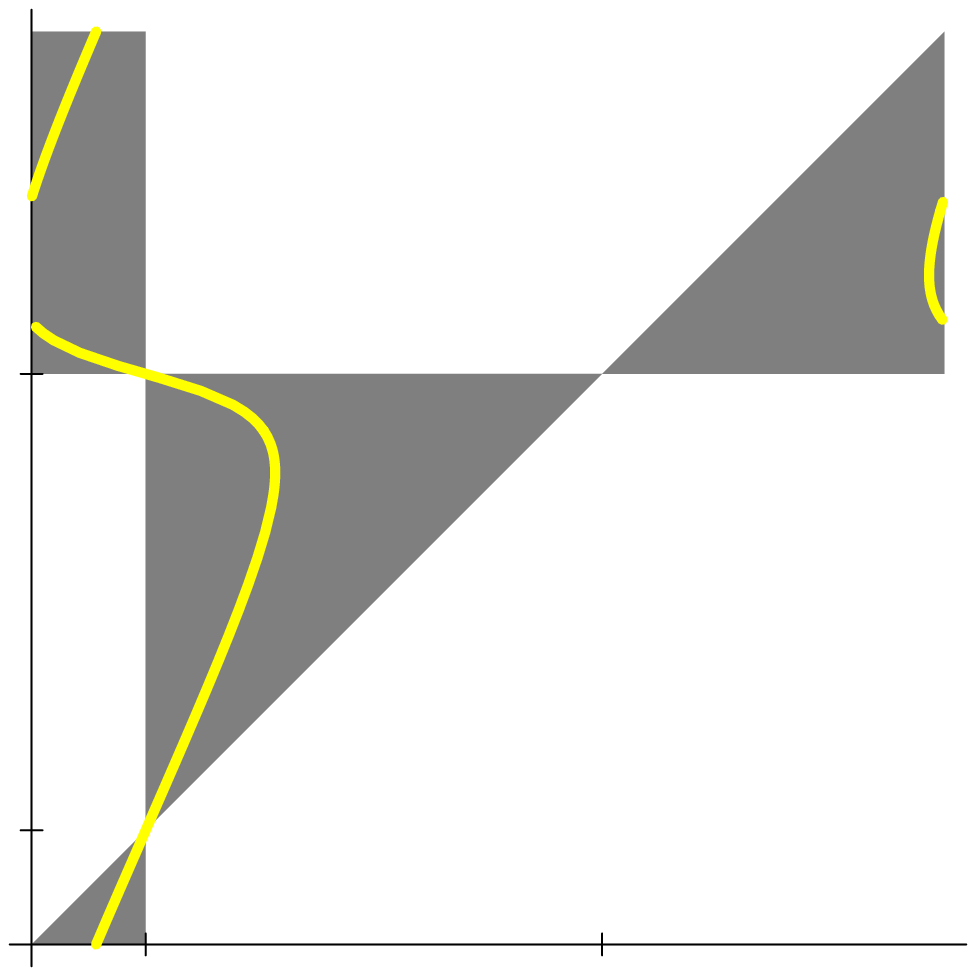}}}
\put(9.25,2){\makebox(0,0)[cc]{
        \leavevmode\epsfxsize=3.5\unitlength\epsfbox{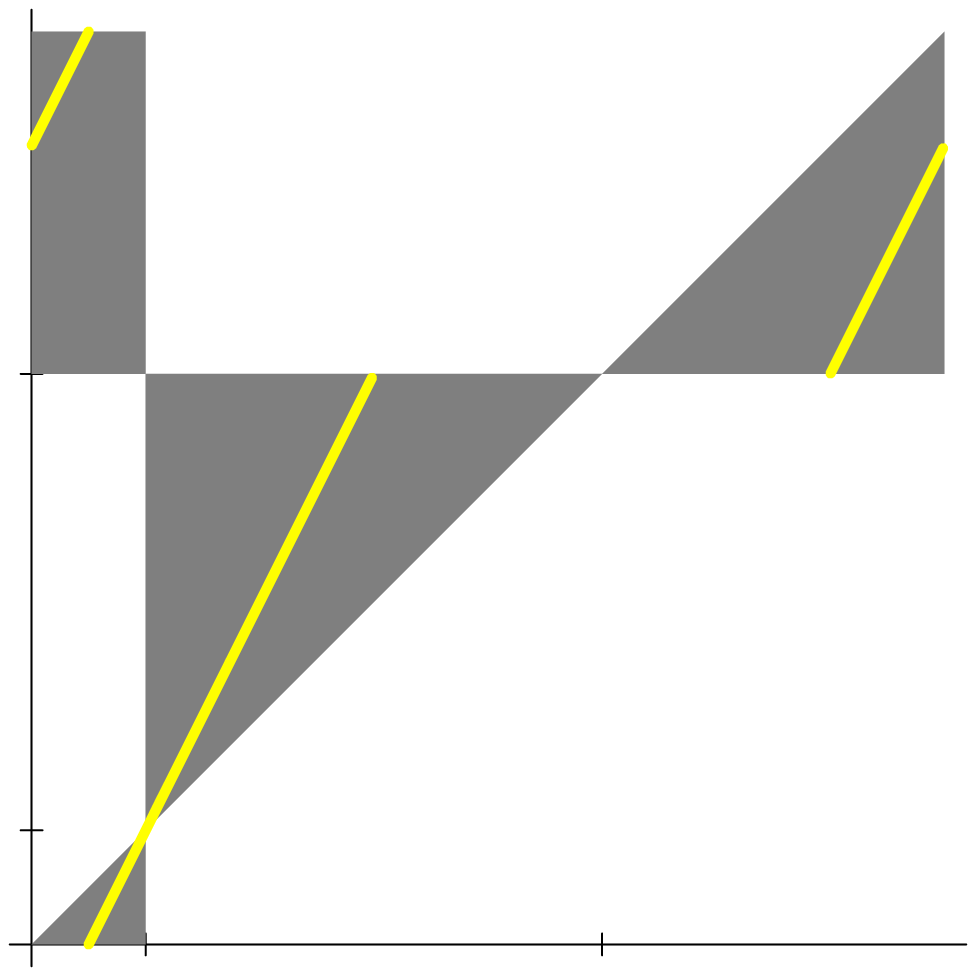}}}
\put(0.2,0.75){\makebox(0,0)[l]{$\frac{\pi}{4}$}}
\put(0.2,2.4){\makebox(0,0)[l]{$\frac{5\pi}{4}$}}
\put(1.0,0){\makebox(0,0)[l]{$\frac{\pi}{4}$}}
\put(2.5,0.){\makebox(0,0)[l]{$\frac{5\pi}{4}$}}
\put(0.1,3.7){\makebox(0,0)[l]{$\alpha_{23}$}}
\put(4.3,0.4){\makebox(0,0)[l]{$\alpha_{13}$}}
\put(7.2,0.75){\makebox(0,0)[l]{$\frac{\pi}{4}$}}
\put(7.2,2.4){\makebox(0,0)[l]{$\frac{5\pi}{4}$}}
\put(8.0,0){\makebox(0,0)[l]{$\frac{\pi}{4}$}}
\put(9.5,0.){\makebox(0,0)[l]{$\frac{5\pi}{4}$}}
\put(7.1,3.7){\makebox(0,0)[l]{$\alpha_{23}$}}
\put(11.3,0.4){\makebox(0,0)[l]{$\alpha_{13}$}}
\end{picture}
\caption{\elabel{ftravelp3} Configurations for moving third point.}
\end{center}
\end{figure}

There are some obvious geometric transformations that relate points
in the  image $\pdat(\CONFt)$. We use these relations for
analyzing the structure of $\pdat(\CONFt)$. Here we consider
$\DAt$ as a quotient of $\R^3$. That is, identify
points $a=(a_x,a_y,a_z)$ and $b=(b_x,b_y,b_z)$  in $\R^3$  iff 
$a_x \equiv b_x \mod 2 \pi$ etcetera. 

\begin{lemma}
\elabel{lperm}
Let $(x,y,z) = \pdat(c)$, with $c=(p_1,p_2,p_3) \in \CONFt$.
Let $\sigma$ be an element of the symmetric group $S_3$.
The consequence of a permutation $\sigma$  of the labels of the
points in $c$ is as follows. 
\begin{displaymath}
\begin{array}{llll}
\sigma \in S_3 \qquad & \sigma(c) \in \pdat(\CONFt) \qquad & 
	{\rm geometric~action~in~} \DAt \\[.2cm]
()      & (x,~ y,~z) \\
(12)    & (x+\pi,~ z,~y) &  T_{(\pi,0,0)} \circ S_{\alpha_{13}=\alpha_{23}}\\
(13)    & (z + \pi,~ y + \pi,~x+\pi) &
T_{(\pi,\pi,\pi)} \circ S_{\alpha_{12}=\alpha_{23}}\\
(23)    & (y,~x,~ z + \pi) &
        T_{(0,0,\pi)} \circ S_{\alpha_{12}=\alpha_{13}}\\
(123)   & (z,~x+\pi,~y+\pi) & R_{-\frac{2 \pi}{3}} \circ T_{(0,\pi,\pi)}\\
(132)   & (y + \pi,~z+\pi,~x) &  R_{\frac{2 \pi}{3}} \circ T_{(\pi,\pi,0)}
\end{array}
\end{displaymath}
 A `$T$' indicates a translation over a given
vector, an `$S$' a reflection in a given plane, while `$R$' indicates a rotation
with respect to the axis through $(0,0,0)$ and $(1,1,1)$ over a given angle.
\end{lemma}

\begin{proof}
A sketch of the six permutations of the labels of the components
of $p$ suffices to determine $\sigma(p)$. The geometric action
that corresponds to applying the three order two elements on
the labels is clear.
The following formula, see \cite{Go}, gives the image $r'$ of the rotation
of a vector $r$ through an angle $\phi$ about an axis $\hat{n}$.
\beq
         r'  & = &  r \cos \phi +   \hat{n} (\hat{n} \cdot r)(1 - \cos \phi) +
                (r \times \hat{n} \sin \phi).
\eeq
In our case, $\hat{n}$ is the normalized vector $\frac{1}{3}\sqrt{3}(1,1,1)$,
while $r=\{x,y,z\}$.
Applying this formula for $\phi=-\frac{2\pi}{3}$ gives $r'=\{z,x,y\}$, 
while $\phi=\frac{2\pi}{3}$ leads to $r'=\{y,z,x\}$. This explains the action 
corresponding to the order three elements of $S_3$.
\end{proof}

\begin{lemma}
\elabel{lrotate}
If $(x,y,z) \in \pdat(\CONFt)$, then also
$(x+\rho, y + \rho, z + \rho) \in \pdat(\CONFt)$ for any $\rho \in [0, 2 \pi)$.
\end{lemma}

\begin{proof}
A rotation of $\rho \in [0, 2 \pi)$ of the baseline, that is the axis with 
respect to which angles are measured, corresponds to adding $\rho$
to all angles $\alpha_{12}$, $\alpha_{13}$ and $\alpha_{23}$.
\end{proof}

Lemmata \ref{lperm} and \ref{lrotate} suggest to consider the set of all 
possible angles between three distinct points up to
the permutations and baseline rotations introduced in the lemmata.
Define $ \Phi_{p,r}^3 := \fpdat$, the \bfindex{fundamental image} of $\pdat( \CONFt )$,
as the image $\pdat( \CONFt )$ up to permutations and rotations. 

\begin{corollary}
$\pdat( \CONFt )$ up to baseline rotations is obtained by projecting
$\pdat( \CONFt )$ on the plane $\alpha_{12} + \alpha_{13} + \alpha_{23}=0$.
\end{corollary}

\begin{proof}
This follows directly from Lemma \ref{lrotate}, as any vector $(\rho, \rho, \rho)$
is orthogonal to the plane $\alpha_{12} + \alpha_{13} + \alpha_{23}=0$.
\end{proof}

We can actually construct the projection $\pi_r$ on the plane 
$P:\,\alpha_{12} + \alpha_{13} + \alpha_{23}=0$ as follows.
$u_1 = (1,-1,0)$ and $u_2=(\frac{1}{2}, \frac{1}{2},-1)$ form an orthogonal
basis of $P$. Let $e_1$ and $e_2$ be the corresponding unit length vectors.
Then $\pi_r$ is given by
\begin{displaymath}
\begin{array}{rlcl}
	\pi_r: &  \R^3 & \rightarrow & \R^2, \\
	       & c     & \mapsto     & ( <c,e_1>, <c,e_2> ).
\end{array}
\end{displaymath}
$\pi_r( \pdat( \CONFt ) )$ is shown on the left in Figure \ref{fprojeps}.

\begin{remark}
\elabel{rsmaller}
$\pi_r$ maps $( \pi, 0, 0)$ to
$p=(\pi A,\pi B)$, maps
$(0,  \pi, 0)$ to $q = (- \pi A, \pi B)$ and maps $(0,0, \pi)$ to 
$(0,2 \pi B)$, where $A = \frac{1}{2} \sqrt{2}$ and $B=\frac{1}{6}\sqrt{6}$. 
As a consequence, any $a_3 \in \pdat(\CONFt)$ has a representative mod 
$2 \pi$ that is mapped by $\pi_r$ in the area shown on the right 
in  Figure \ref{fprojeps}.
\end{remark}

\begin{figure}[!ht]
\begin{center}
\setlength{\unitlength}{1cm}
\begin{picture}(12,5.7)
\put(3,3){\makebox(0,0)[cc]{
        \leavevmode\epsfxsize=4.2\unitlength\epsfbox{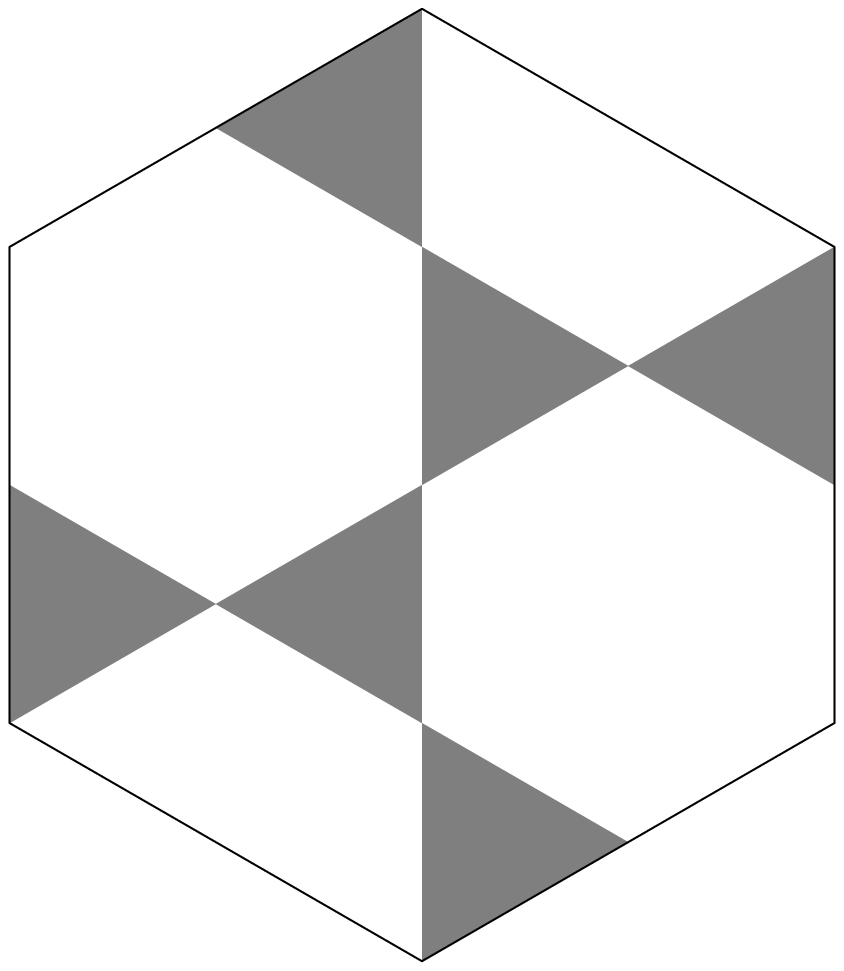}}}
\put(9,3){\makebox(0,0)[cc]{
        \leavevmode\epsfxsize=4.2\unitlength\epsfbox{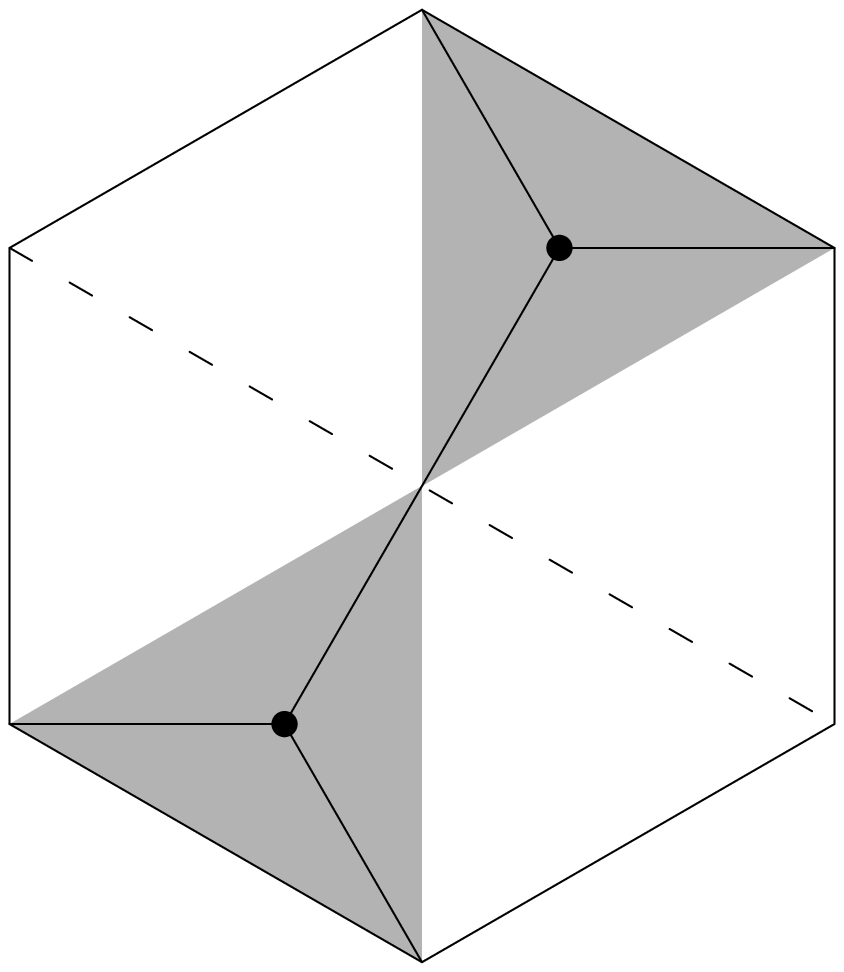}}}
\put(11.25,4.2){\makebox(0,0)[l]{$p$}}
\put(6.8,4.2){\makebox(0,0)[l]{$q$}}
\put(9,5.6){\makebox(0,0)[l]{$r$}}
\put(9.7,4.35){\makebox(0,0)[l]{$-m$}}
\put(8.65,1.8){\makebox(0,0)[l]{$m$}}
\put(8.3,0.4){\makebox(0,0)[l]{$\alpha_{12}=\alpha_{13}$}}
\put(10.5,1.4){\makebox(0,0)[l]{$\alpha_{12}=\alpha_{23}$}}
\put(5.95,1.4){\makebox(0,0)[l]{$\alpha_{13}=\alpha_{23}$}}
\end{picture}
 \caption{\elabel{fprojeps}$\pi_r(\pdat( \CONFt ))$ and six times $\Phi_{p,r}^3$. }
\end{center}
\end{figure}

Next lemma is the `projected' version of Lemma \ref{lperm}.

\begin{lemma}
\elabel{projperm}
Let $(x,y,z) = \pdat((p_1,p_2,p_3))$, with $c=(p_1,p_2,p_3) \in 
\CONFt$. The orbit of $\pi_r(x,y,z))$
under the action of the symmetric group $S_3$ on the labels of 
$p_1$, $p_2$ and $p_3$ is given in the following table.
\begin{displaymath}
\begin{array}{llll}
\sigma \in S_3, & {\rm geometric~action,} & {\rm coordinates.}\\[.2cm]
()  & -  & (A(x-y), & B(x+y-2z))
\\
(12)  & T_{\delta p} + S_{\alpha_{13}=\alpha_{23}} &
(A(\pi+x-z),& B(\pi+ x-2y+z)) \\
(13)  & S_{\alpha_{12}=\alpha_{23}} &
(A(z-y),&B(z+y-2 x)) \\
(23)  & T_{\delta r} +S_{\alpha_{12}=\alpha_{13}} &
(A(y-x),&B(x+y-2z-2\pi)) \\
(123) & R_{\delta m,-\frac{2 \pi}{3}} &
(A(z-x-\pi),&B(x+z-2 y-\pi)) \\
(132) & R_{\delta m,\frac{2 \pi}{3}} &
(A(y-z),&B(2\pi+y+z-2 x))
\end{array}
\end{displaymath}
Here $p,q,r$ and $A, B$ are as in Remark \ref{rsmaller} and $m$ is given
by $m=  (-\frac{\pi A}{3},-\pi B)$. Moreover, $\delta=1$, if the
second coordinate of $\pi_r(x,y,x)<0$ and $\delta=-1$,
if the second coordinate of $\pi_r(x,y,x)>0$.
\end{lemma}

\begin{proof}
We can easily compute the coordinates of the points in the orbit
by projecting the orbit in $\pdat(\CONFt)$ that we found in
Lemma \ref{lperm}.
Concerning the geometric action: the reflection matrix $S_a$
for a reflection of a point in the  line $y=ax$  is given by
\begin{eqnarray*}
        S_a & = & \frac{1}{1+a^2} \left(\begin{array}{cc}
                1-a^2 & 2a \\
               2a    & a^2 -1
       \end{array}\right).
\end{eqnarray*}
It is easy to check that the action of $(12)$ boils down to a
translation over $p$ followed by  a reflection in the line
$\alpha_{13}=\alpha_{23}$. Similarly for the action
of the other two involutions.
The rotation matrix $R_\alpha$ for a rotation of a point around the origin
is given by
\begin{eqnarray*}
        R_\alpha & = &  \left(\begin{array}{cc}
                \cos \alpha & -\sin \alpha \\
                \sin \alpha    &\cos \alpha
       \end{array}\right).
\end{eqnarray*}
A rotation of a vector $v$ around an arbitrary center $m$ is given
by $R_{\alpha}(v-m) + m$. This allows us to check the action
of $(123)$ and $(132)$.
\end{proof}

\begin{corollary}
$\Phi_{p,r}^3$ is given by the triangle in the picture on the
right of Figure \ref{fprojeps} with vertices $(0,0)$, $-m$ and $p$.
\end{corollary}

\begin{proof}
Every element of the orbit of $S_3$ acting on a $\pi_r(x,y,z)$
for some $(x,y,z) = \pdat(c)$, with $c \in \CONFt$
lives `at the same place'  in its own triangle in the picture on the right of
Figure \ref{fprojeps}: the three order $2$ elements map  $\pi_r(x,y,z)$
on one of the three gray triangles on the other site of the dotted 
line $\alpha_{12}=\alpha_{23}$. The two order three elements rotate
$\pi_r(x,y,z)$ clockwise or anti-clockwise into the two 
adjacent gray triangles.  
\end{proof}

\begin{corollary}Let $P$ be the plane orthogonal to the vector
$(1,1,1)$. Then the $2$-dimensional wallpaper group $p6m$
is acting on $P \cap \pdat(\CONFt)$.
\end{corollary}

\begin{proof}
There are order $6$, order $3$ and order $2$ rotations, and $6$
axes of reflection. This characterizes the $17$-th wallpaper
group, see for example \cite{CM}.
\end{proof}

\section{The triangle variety $\mathbf{T_n}$.}

We want to answer the question:
\begin{quote}
are $\cua$ and $\cda$ smooth manifolds?
\end{quote}
Moreover, we are interested in their algebraic counterparts. In this
section we describe an algebraic variety that is very similar to $\cua$.
We know that $\CONF$ is contained in $\cua$:
by definition, $\cua$ equals the closure of the graph of the
undirected angle map:
$$ \CONF ~\subset~ \cua ~\subset~ (\R^2)^n \times (\R/ \pi \Z)^{\binom{n}{2}}.$$
Therefore, $\cua$ contains a `$\CONF$ part' that is smooth. 
The remaining points of $\cua$ lie above the \bfindex{diagonal} $\Delta \subset (\R^2)^n$
consisting of configurations with
at least two coinciding points $p_i$ and $p_j$.
We make an algebraic description for $\R/\pi \Z = \P^1$ by taking coordinates
$(a_{ij}:1)$ and $(1:b_{ij})$, where
$$ a_{ij} = \tan \overline{\alpha}_{ij}; \qquad 
	b_{ij} = \frac{1}{\tan \overline{\alpha}_{ij}}. $$
For simplicity, we consider only  the case where 
$b_{ij} \neq 0$ on each $\P^1$, so we work 
on the $(a_{ij}:1)$-chart. We have transformed $\pua$ in a rational map
$\psi_{\text{slope}}$ given by
\beq
	\psi_{\text{slope}}( (x_0, y_0), \dots, (x_{n-1}, y_{n-1}))
	&=& \{(\frac{y_j-y_i}{x_j-x_i})\}_{0\leq i<j\leq n-1} ,
\eeq
where $((x_0,y_0), \dots, (x_{n-1}, y_{n-1})) \in \CONF$. 
Without loss of generality  we
assume throughout this section that $x_0=y_0 =0$. That is, we consider 
configurations up to translation.  The dimension of $\CONF(\R^2)$ 
up to translations equals $2 n - 2$.
The slope $a_{0i}$, for $i \in \{1, \dots, n-1\}$ is denoted short as $a_i$.
The triangle $T_{ij}$ is the triangle with vertices $(x_0, y_0)$,
$(x_i, y_i)$, and $(x_j, y_j)$. 
The following lemma shows that there exists  a relation between the 
$x$-coordinates of the vertices of $T_{ij}$ and the slopes of
the lines bounding $T_{ij}$. Let
\beq
	t_{ij} & = & a_i x_i - a_j x_j - a_{ij} x_i + a_{ij} x_j.
\eeq

\begin{lemma}
\elabel{ltij}
$t_{ij} = 0$ on  the $(a_{ij}:1)$-chart of $\CONF$. 
\end{lemma}

\begin{proof}
It holds that $y_i  =  a_i x_i$, $y_j  =  a_j x_j$, and
	$y_i - y_j  =  a_{ij} (x_i - x_j)$.
Substitute the former two equations in the last equation.
\end{proof}

\begin{corollary}
\elabel{ccuatij}
On the $(a_{ij}:1)$-chart we have that
$\cua \subset \{ t_{ij} = 0 \}$ for $1 \leq i < j \leq n-1$.
\end{corollary}

A question is if equality holds. That is, if the closed algebraic set 
$\{ t_{ij} = 0 \}_{1 \leq i < j \leq n-1}$ is contained 
in $\cua$. The answer is no. We prove this later on by means of
the six-slopes formula:

\begin{figure}[!ht]
\begin{center}
\setlength{\unitlength}{1cm}
\begin{picture}(6,2.4)
\put(3,1.2){\makebox(0,0)[cc]{
        \leavevmode\epsfxsize=6\unitlength\epsfbox{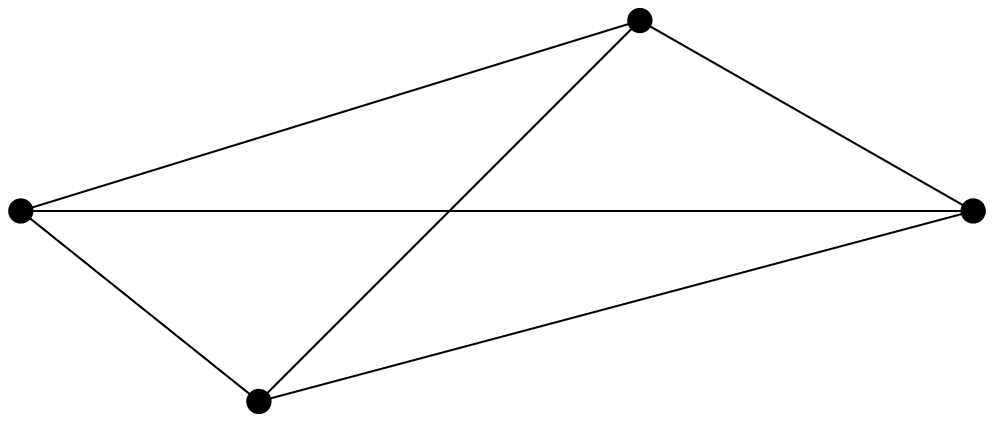}}}
\put(0.5,0.5){\makebox(0,0)[l]{$a_{01}$}}
\put(1.8,1.4){\makebox(0,0)[l]{$a_{02}$}}
\put(1.9,2.1){\makebox(0,0)[l]{$a_{03}$}}
\put(3.9,0.45){\makebox(0,0)[l]{$a_{12}$}}
\put(3.4,1.7){\makebox(0,0)[l]{$a_{13}$}}
\put(4.75,2.0){\makebox(0,0)[l]{$a_{23}$}}
\end{picture}
\caption{\elabel{fsixslopes2}
For four distinct points the six-slopes formula holds.}
\end{center}
\end{figure}

\begin{lemma}[six-slopes formula]
\elabel{a23formula}
\index{six-slopes formula}
Let $p_0$, $p_1$, $p_2$ and $p_3$ be distinct points in the plane.
Then $\Delta=\Delta_{0123} = 0$, where  $\Delta$ is given by
 \begin{eqnarray}
 \elabel{eqa23}
\Delta & = &  (a_1 - a_{12})(a_2-a_{23})(a_3-a_{13}) - 
        (a_1 - a_{13})(a_2-a_{12})(a_3-a_{23}). 
 \end{eqnarray}
\end{lemma}

\begin{proof}
Assume that $p_0 = (0,0)$ and $p_1=(1,a_1)$.
We compute coordinates for the points
$p_2$ and $p_3$. Let $l_{ij}$ be the line through the points $p_i$
and $p_j$.
The lines $l_{02}$ and $l_{03}$ are given by
\beq
        l_{02}:y_2 ~=~ a_2 x,         & \qquad & l_{03}: y_3 ~=~ a_3 x,
\eeq
and $l_{12}$ and $l_{13}$ by
\beq
        l_{12}:y-a_1 ~=~ a_{12}(x-1), & \qquad & l_{13}:y-a_1 ~=~ a_{13}(x-1).
\eeq
We compute the intersections $(x_2, y_2)=l_{02}\cap l_{12}$, 
and $(x_3,y_3)=l_{03}\cap l_{13}$: 
\begin{displaymath}
\begin{array}{lll}
  x_2~=~\frac{a_1-a_{12}}{a_2-a_{12}}  & \qquad &
                x_3~=~\frac{a_1-a_{13}}{a_3-a_{13}} \\
  y_2~=~a_2 x_2 & & y_3~=~a_3 x_3
\end{array}
\end{displaymath}
After some formula manipulation the expression follows as
        $a_{23} = \frac{y_3-y_2}{x_3-x_2}$.
\end{proof}

\begin{remark} Some remarks on Equation \ref{eqa23}.
\begin{lijst}
\item
Interchanging indices $1\leftrightarrow 2$, etcetera, changes the appearance of
the expression for $\Delta_{0123}$, but does not change the expression itself. 
\item By $\Delta_{ijkl}$ we denote $\Delta_{0123}$ with $0$, $1$, $2$ and $3$ replaced
by $i$, $j$, $k$ and $l$.
\item
If no three points of $p_0$, $p_1$, $p_2$ and $p_3$ are collinear, then
we can express $a_{23}$ as follows in terms of $a_1$, $a_2$, $a_3$,
$a_{12}$ and $a_{13}$:
\begin{eqnarray}
\elabel{ra23}
a_{23} & = &
        \frac{ a_3 (a_1 - a_{13})(a_2-a_{12}) - a_2 (a_1 - a_{12})(a_3-a_{13})}
        {(a_1 - a_{13})(a_2-a_{12})-(a_1 - a_{12})(a_3-a_{13})}.
\end{eqnarray}
\end{lijst}
\end{remark}

\begin{corollary}
$\Delta_{ijkl}=0$ on $\cua$. 
\end{corollary}

\begin{proof}
This follows from Lemma \ref{a23formula} and the definition of $\cua$.
\end{proof}

Instead of just looking at the zeros of $t_{ij}=0$, we add the condition
that all $\Delta_{ijkl}$ equal zero as well. This leads to the 
following definition.

\begin{definition}
\elabel{dtriav}
The \bfindex{triangle variety} $T_n$ \index{Tn@$T_n$}
 is  the set of common zeroes of the 
polynomials $t_{ij}$ for $1 \leq i < j  \leq (n-1)$ and 
$\Delta_{ijkl}$ for $0 \leq i < j < k < l \leq (n-1)$.
Any variable, of the form $a_{ij}$ or $x_i$, takes value in $\R$.
\end{definition}

\begin{eexample}
The zero set of a collection of polynomials $f_1. \dots, f_n$  is 
indicated by $V(f_1, \dots, f_n)$. So
$T_4 = V(t_{12}, t_{13}, t_{23}, \Delta_{123})$.
\end{eexample}

Note that the ideal of the triangle variety $T_n$ contains: 
one polynomial $t_{ij}$ for every triangle with vertices
$p_0$, $p_i$ and $p_j$;
one polynomial $\Delta_{ijkl}$ for every quadrilateral
with vertices $p_i$, $p_j$, $p_k$ and $p_l$.
We could consider explicitly triangles with vertices $p_i$, $p_j$
and $p_k$ by including:
\beq
 t_{ijk} & = & a_{ik} (x_i-x_k) - a_{ij} (x_i-x_j) -a_{jk} (x_j-x_k).
\eeq
The following lemma shows that this is not necessary however.

\begin{lemma}
$t_{ijk} = t_{ij}-t_{ik}+t_{jk}$.
\end{lemma}

\begin{proof}
This follows directly from the following equations:
\begin{align*}
t_{ij}&=a_i x_i - a_j x_j - a_{ij} x_i+ a_{ij} x_j, \\
t_{ik}&=a_i x_i - a_k x_k - a_{ik} x_i+ a_{ik} x_k, \\
t_{jk}&=a_j x_j - a_k x_k - a_{jk} x_j+ a_{jk} x_k. 
\qedhere
\end{align*}
\end{proof}

\subsection{Singularities of $T_n$ for small $n$.}

In this section we determine singularities of the
triangle variety $T_n$ for $n=3$ and $n=4$.
For this purpose recall the definition of
\bfindex{singularity}. 

\begin{definition} 
\elabel{dsingular}
See \cite{Ha}. Let $Y \subset A^n$ be an affine variety, and
let $f_1,\dots,f_t$ $\in A=$  $k[x_1,\dots,x_n]$ be a set of generators
for the ideal of $Y$. Then $Y$ is {\bf nonsingular at a point} $P\in Y$
if the rank of the Jacobian matrix $((\partial f_i / \partial x_j)(P))$ is
$n-r$, where $r$ is the dimension of $Y$. If $Y$ is  nonsingular
at every point, then $Y$ is {\bf nonsingular}.
\end{definition}

\subsection{n=3.}

\begin{figure}[ht]
 \begin{center}
\setlength{\unitlength}{1cm}
\begin{picture}(6.2,2.3)
\put(3.5,1.3){\makebox(0,0)[cc]{
        \leavevmode\epsfxsize=6.2\unitlength\epsfbox{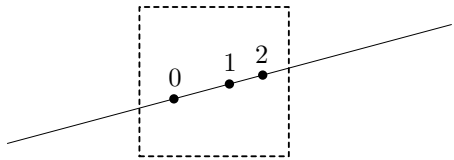}}}
\put(2.75,1.35){\makebox(0,0)[l]{$0$}}
\put(3.48,1.55){\makebox(0,0)[l]{$1$}}
\put(3.9,1.66){\makebox(0,0)[l]{$2$}}
\end{picture}
 \caption{\label{threeonline} Three coinciding points with coinciding
  directions.}
\end{center}
\end{figure}

\begin{lemma}
\elabel{lt3}
An element $c=(x_0,x_1,x_2,a_1,a_2,a_{12}) \in T_3$ is singular iff 
all points and all slopes coincide. That is:
$$
	x_0 = x_1 = x_2, \qquad \text{and} \qquad 
		a_1 = a_2 = a_{12}.
$$
The type of this singularity is $A_{\infty}$.
\end{lemma}

\begin{proof}
$T_3 = V(t_{12})$, with $t_{12} = a_1 x_1 - a_2 x_2 - a_{12} x_1 + a_{12}x_2$. 
$T_3$ is singular if the rank of the Jacobian of $t_3$ is smaller
than $5-4 =1$, that is if all partial derivatives of $t_{12}$ vanish.
\begin{displaymath}
\begin{array}{llllll}
{\rm variable} &x_1 & x_2 &  a_1 & a_2 & a_{12}, \\
{\rm partial~derivative} &a_1-a_{12} & a_{12}-a_2 & x_1 & -x_2 & x_2-x_1.
\end{array}
\end{displaymath}
Solving this system gives the solutions $x_1=x_2=0$ $\wedge$ $a_1=a_2= a_{12}$. 
For the type of the singularity: note that
$t_{12}=(a_1 - a_{12})x_1 + (a_{12} - a_2) x_2$. From this it follows that
the singular set is given by $a_i=a_j=a_{ij}$, and that the singularity is
of type $A_{\infty}$, cf.\ \cite{Si}.
\end{proof}

\begin{remark}
Geometrically, the singularities of $T_3$ correspond to degenerated
 configurations where both the  three points coincide and the
directions between the points coincide, see Figure \ref{threeonline}.
\end{remark}

\subsection{n=4.}

\begin{lemma}
\elabel{lAi}
Define 
\beq
	A_1 & = & (a_2-a_{12})(a_3-a_{23})-(a_2 - a_{23})(a_3-a_{13}), \\
	A_2 & = & (a_1-a_{12})(a_3-a_{13})-
		(a_1 - a_{13})(a_3-a_{23}), \\
	A_3 & = & (a_2-a_{23})(a_1-a_{12})-
		(a_2 - a_{12})(a_1-a_{13}),\\
	A_{12} & = & (a_1-a_{13})(a_3-a_{23})-(a_2-a_{23})(a_3-a_{13}),\\
	A_{13} & = & (a_2-a_{12})(a_3-a_{23})-(a_1-a_{12})(a_2-a_{23}),\\
	A_{23} & = & (a_1-a_{13})(a_2-a_{12})-(a_1-a_{12})(a_3-a_{13}).
\eeq
\begin{lijst}
\item
If $A_1 \neq 0$, then 
$$ x_2  =  x_1 \frac{(a_{12} - a_{13}) (a_{23} - a_3)}{A_1}; \qquad
		x_3  =   x_1 \frac{(a_{12} - a_{13}) (a_{23} - a_2)}{A_1}.$$
If $A_2 \neq 0$, then 
$$ x_1  =  x_2 \frac{(a_{12} - a_{23}) (a_{13} - a_3)}{A_2}; \qquad
		x_3  =   x_2 \frac{(a_{12} - a_{23}) (a_{13} - a_1)}{A_2}.$$
If $A_3 \neq 0$, then
$$ x_2  =  x_3 \frac{(a_{23} - a_{13}) (a_{12} - a_1)}{A_3}; \qquad
                x_1  =   x_3 \frac{(a_{23} - a_{13}) (a_{12} - a_2)}{A_3}.$$
\item $\frac{\partial \Delta}{\partial a_{ij}} = A_{ij}$. 
\end{lijst}
\end{lemma}

\begin{proof} 
As $t_{23}=0$ on $T_n$, it follows that $(a_{23}-a_3)x_3 = (a_{23}-a_2)x_2$.
As $t_{12}-t_{13}=0$, we arrive at
\beq
(a_{12}-a_2) x_2 - (a_{13}-a_3) x_3 & = & a_{12} x_1 - a _{13} x_1.
\eeq
Add $(a_{23} - a_3)$ times this last equation to $f_{23}$ and substitute
$(a_{23}-a_2)x_2$ for  $(a_{23}-a_3)x_3$.  This gives that
$A_1 x_2 = (a_{12}-a_{13})(a_{23}-a_3) x_1$, resulting in the formula for $x_2$.
The formula for $x_3$ is obtained by adding $(a_{23}-a_2)(f_{12}-f_{13})$ to $f_{23}$
and substituting $(a_{23}-a_3)x_3$ for $(a_{23}-a_2)x_2$.
Substitute the expression for $x_2$ in $f_{12}$, multiply by $A_1$ and reorder:
\beq
	x_1 a_1 A_1 + x_1(a_{12}-a_2)(a_{12}-a_{13})(a_{23}-a_3) & = & x_1 a_{12} A_1,\\
	x_1 \Delta & = & 0.
\eeq
The last equation holds in two cases. In the first case, $x_1=0$. But from the
expressions for $x_2$ and $x_3$ it follows that then  $x_2=0$ and 
$x_3=0$ as well, while no conditions are set on the slopes. In the other case 
$\Delta = 0$, while for $x_2$ and $x_3$ the two formulas in the Lemma hold.
The equations in case  of $A_2 \neq  0$ and $A_3 \neq 0$ follow
by relabeling. This proves the first claim.
The second claim  follows by inspection.
\end{proof}

\begin{corollary}
\elabel{ctwoc}
$\cuaf$ $\neq$ $V(t_{12},t_{13},t_{23})$.
\end{corollary}

\begin{proof}
This follows from the proof of Lemma \ref{lAi}: The variety $V(t_{12},t_{13},t_{23})$
contains $0=x_1=x_2=x_3$ as a component, {\it without} the condition that
$\Delta_{0123}=0$. 
\end{proof}

\begin{remark}
Corollary \ref{ctwoc} and its proof explain why we have added the $\Delta_{ijkl}$'s
in the definition of triangle variety $T_n$, cf.\ Definition \ref{dtriav}.
\end{remark}

\begin{lemma}
\elabel{lsingt4}
A configuration $c \in T_4$ is singular iff, 
up to relabeling, both $x_0=x_1=x_2$ and $a_1=a_2=a_{12}$. 
\end{lemma}

\begin{proof}  
Consider the Jacobian $J_4$ of $T_4$ and apply Definition \ref{dsingular}.
The Jacobian of $T_4$ where the nine variables are in the order
$x_1$, $x_2$, $x_3$, $a_1$, $a_2$, $a_3$, $a_{12}$, $a_{13}$, $a_{23}$ is given
by:
\begin{displaymath}
\left(
\begin{array}{ccccccccc}
\!\!\!\!a_1 - a_{12} & a_{12}-a_2 & 0 & x_1 & -x_2 & 0 & x_2 -x_1 & 0 & 0 \\
\!\!\!\!a_1 - a_{13} & 0 & a_{13}-a_3 & x_1 & 0 & -x_3 & 0 & x_3 -x_1 &  0 \\
\!\!\!\!0 & a_2 -a_{23} & a_{23}-a_3 & 0 & x_2 & -x_3 & 0 & 0 & x_3 - x_2 \\
\!\!\!\!0 & 0 & 0 & -A_1 & -A_2 & -A_3 & \pm A_{12} & \pm A_{13} &  \pm A_{23} 
\end{array}
\!\!\!\!\right).
\end{displaymath}

The dimension of $T_4$ equals $2 \cdot 4 -2=6$. Therefore, a configuration is singular
iff the rank of its Jacobian is equal to or smaller than $2$. 
Denote by $m(i,j,k)$ the submatrix of $J_4$ consisting of columns
$i$, $j$ and $k$ where the fourth row is deleted.
We distinguish several cases, by considering the number of distinct clusters
of coinciding points in $c$.
\begin{lijst}
\item {\bf [$\geq 3$ clusters]}
We may assume that $x_1 \neq 0$, $x_2 \neq 0$, $x_3 \neq 0$ and $x_1 \neq x_2$. 
Consider $m(5,6,7)$: 
\begin{displaymath}
m(5,6,7) \quad = \quad 
\left(\begin{array}{ccc}
	-x_2 & 0    & x_2 -x_1 \\
	0    & -x_3 & 0 \\
	x_2  & -x_3 & 0
\end{array}
\right).
\end{displaymath}
The determinant of $m(5,6,7)$ equals $x_2x_3(x_2 -x_1) \neq 0$. This implies
non-singularity.
\item {\bf [$2$ clusters: (2,2)]}  We may assume that $x_1 = 0$ and that $x_2 = x_3 \neq 0$.
       In this case the determinant of $m(5,6,7)$ equals $(x_2)^3$, which implies 
	non-singularity	.
\item {\bf [$2$ clusters: (3,1)]} See also Figure \ref{fsingT4}.
Assume that for $c \in T_4$ it holds that $x_0=x_1=x_2$;
that $a_1=a_2=a_{12}$ and that $x_3 \neq x_0$.
As a consequence $a_3=a_{13}=a_{23}$. It follows that the first and last row
$J_4$ vanishes. Therefore, $c$ is singular.
If we only assume that $x_0=x_1=x_2$ and that $x_3 \neq x_0$, then
the determinants $\text{det}(m(1,6,8))=(x_3)^2(a_1 - a_{12})$
and $\text{det}(m(2,6,8))=(x_3)^2(a_{12} - a_2)$ show that
the condition $a_1 = a_2 = a_{12}$ is necessary for a singularity. 
\item {\bf [$1$ cluster]} It holds that $x_0=x_1=x_2=x_3=0$.  
Consider $J_4$. First assume that the last row of $J_4$
      equals zero. As $\Delta=0$ on $T_4$ this means that we
      are looking for solutions of the system of equations
$$ A_1 ~=~ A_2 ~=~ A_3 ~=~ A_{12} ~=~ A_{13} ~=~ A_{23} ~=~ 
      \Delta ~=~ 0.$$
      Applying Mathematica's {\tt Reduce} function results in 10
      reductions. Each of these 10 reductions contains 
      equations of the form $a_1 = a_2 = a_{12}$.
We are left with the case that some $A_{ij} \neq 0$. If
$\text{rank}(J_4) \leq 2$, then the determinants 
of all $2 \times 2$ submatrices of $m(1,2,3)$ vanish. 
This gives a system of nine equations that can be reduced
to the six equations of the form 
$ a_1 = a_2 = a_3 = a_{12} = a_{13}$,
that is, five slopes are equal. Again, each of these
equations contains a condition of the form $a_1=a_2=a_{12}$.
\qedhere
\end{lijst}
\end{proof}

\begin{figure}[ht]
 \begin{center}
\setlength{\unitlength}{1cm}
\begin{picture}(7,1.9)
\put(3.5,1){\makebox(0,0)[cc]{
        \leavevmode\epsfxsize=7\unitlength\epsfbox{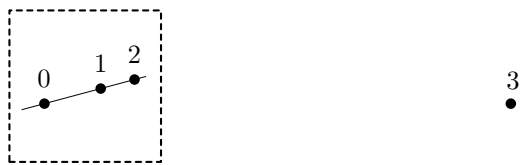}}}
\put(0.6,1.1){\makebox(0,0)[l]{$0$}}
\put(1.35,1.3){\makebox(0,0)[l]{$1$}}
\put(1.8,1.42){\makebox(0,0)[l]{$2$}}
\put(6.83,1.07){\makebox(0,0)[l]{$3$}}
\end{picture}
\caption{\elabel{fsingT4} Typical singular configuration: three
coinciding, collinear points with a fourth point. }
\end{center}
\end{figure}

\begin{remark} 
The singularities of $T_4$ are closely related to the singularities
of $T_3$. The singularity in the (3,1) case is of type $A_\infty$
just as for $T_3$. But if all points coincide a more complicated 
singularity occurs as several $A_\infty$ singularities `meet': 
one can move any of the four points away in such a way that the
three remaining points are as in the configuration of Lemma \ref{lt3}.
\end{remark}

It is still an open question whether $\cuaf = T_4$.  Maybe we 
need to add some relations or inequalities to $T_4$ to obtain equality?

\chapter{Continuity of the Voronoi map. \elabel{chcont}}

In Chapter 4 we have introduced the compactification $\cda$
of the graph of the angle map, applied to $n$ distinct point in the plane.
By means of the extended definition of Voronoi diagram, as introduced in
Chapter 3, we are able to consider the Voronoi diagram $V(\gamma_n)$ of a data set
$\gamma_n \in \cda$. The main result of this chapter is Theorem \ref{hausdorffcont}. 
It states, up to a compactness condition, that two data sets $\gamma_n$ and $\eta_n$
in $\cda$ that are Euclidean close, have Voronoi diagrams that are Hausdorff close.
That is, the Voronoi diagram generated by a set of not necessarily distinct 
points and the angles between those points does not change 
dramatically if we perturb both the points and the angles only slightly.

\section{Introduction.}

In Chapters 3 and 4 we have extended the notion of
Voronoi diagram to configurations of points together
with the angles between the points. This allows
us to consider Voronoi diagrams of point sets that include
coinciding points. More precisely, in Section \ref{scda},
we have introduced a compactification $\cda$, as the closure of
the graph of the angle map
\beq
	\pda: \CONF & \rightarrow & \DA.
\eeq
The angle map $\pda$ maps a set $c$ of $n$ distinct points in the plane
to the angles in $\R/2 \pi \Z$ between the points in $c$.

In this chapter we call the elements $\gamma_n \in \cda$, consisting of
$n$ points in $\R^2$ and $\binom{n}{2}$ angles in $\R/2 \pi \Z$, data sets.
With any data set $\gamma_n \in \cda$, we associate a Voronoi diagram $V(\gamma_n)$ 
in Section \ref{scprelim}.

One can wonder whether the data sets are a robust way of storing
Voronoi diagrams. That is, we consider the following 
question: 
if we perturb a data set $\gamma_n \in \cda$ slightly, how does the 
corresponding Voronoi diagram $V(\gamma_n)$ change? By perturbing
a data set $\gamma_n \in \cda$ we mean that both the points components
and the angle components of $\gamma_n$ are allowed to be perturbed slightly, 
as long as  the perturbed data set $\tilde{\gamma}_n$ is again in $\cda$.
A first result is given in Theorem \ref{onepcont}, that is a kind of one point
continuity theorem. It states that a point $x$ that was on the one-skeleton
of a Voronoi diagram $V(\gamma_n)$ before perturbing, cannot be too far from the
one skeleton of the Voronoi diagram of the perturbed data set.

Before we are able to give a more general continuity result, we first need
a suitable metric on the set of Voronoi diagrams. A metric that is often used
to compare pictures, is the so-called Hausdorff metric. Two subsets $A$ and $B$ 
of $\R^2$ are close to each other in the Hausdorff metric, iff the maximal Euclidean
distance from any point $b \in B$ to the set $A$ is small and vice versa.
For a precise definition, see Section \ref{shausdorff}.

In order to prove a continuity theorem with respect to the
Hausdorff metric, we have to add some  restrictions  to 
our underlying point configurations.
We add four so-called camera points, points  that are very far away around our 
configuration.
Now we restrict the configurations we are interested in to a relatively
small bounded subset of the plane. If we only change the configurations
within this subset, the camera points guarantee us that the Voronoi
diagrams of the configurations are not changing outside a second, larger
bounded subset of the plane.
But then we are ready for proving the main 
result of this chapter: Theorem \ref{hausdorffcont} states that,
up to the compactness condition, two data sets that are Euclidean close,
have Voronoi diagrams that are Hausdorff close.
Note that it is essential that the continuity theorem holds on the whole
of the chosen compact subset of $\cda$: that is, the points 
added to the graph of $\pda$ by taking the closure are essential
as there availability is
used in the proof of Theorem \ref{hausdorffcont}. 
\section{Preliminaries.}
\elabel{scprelim}

Throughout this chapter we compare data sets $\gamma_n,
\eta_n \in \cda$. Write 
\beq 
\gamma_n & = & (p_1(\gamma_n), \dots,
p_n(\gamma_n), \alpha_{12}(\gamma_n), \dots, \alpha_{(n-1)n}(\gamma_n)),
\eeq
and similarly $\eta_n$. 
We say that $\eta_n$ is within \bfindex{distance} $\delta$ of $\gamma_n$,
notation $d(\gamma_n, \eta_n) < \delta$, iff
\beq
	\underset{1 \leq k \leq n}{\text{max}}
		 \| p_k(\gamma_n) - p_k(\eta_n) \| < \delta
		& \wedge &
	\underset{1 \leq i<j\leq n}{\text{max}} 
		\| \alpha_{ij}(\gamma_n) -
		\alpha_{ij}(\eta_n) \| < \delta.
\eeq
Here, $\|.\|$ denotes the ordinary 
Euclidean distance on $\R^2$ and $\R/2 \pi \Z$.

From Section \ref{svorpoly} we recall the definition of Voronoi
diagram in terms of not necessarily non-coinciding points in the plane
and directed angles between those points. This enables
us to define for any $\gamma_n \in \cda$ a Voronoi diagram $V(\gamma_n)$. 
In Section \ref{svordia} we have seen that the Voronoi cell of a generator
$p_i$  can be written as the intersection of all Voronoi
half-planes $vh(p_i,p_j)$. We use this characterization to
introduce the Voronoi cell of a point $p_i(\gamma_n)$.

Fix two points $p_i=p_i(\gamma_n)$ and $p_j=p_j(\gamma_n)$. 
The \bfindex{bisection point} $b(p_i, p_j)$ is the point 
	$\frac{1}{2}(p_i + p_j)$.
If $p_i\neq p_j$, the bisection point is just the middle of the line segment
$p_ip_j$. If $p_i=p_j$ then  $b(p_i, p_j)$ coincides with the double point
$p_i=p_j$.
The \bfindex{perpendicular bisector} $B(p_i,p_j)$ is the line through
$b(p_i,p_j)$ perpendicular to the angle $\alpha_{ij}=\alpha_{ij}(\gamma_n)$. Let
$\mathbf{n}$ be any non-zero vector, pointing in the direction $\alpha_{ij}$.
The \bfindex{Voronoi half-plane} $vh(p_i,p_j)$ is the half-plane defined by
\begin{eqnarray*}
\elabel{eqbiseqn2}
        \mathbf{n} \cdot (x-b_x , y-b_y ) & \leq & 0.
\end{eqnarray*}
The {\bf Voronoi cell} $V(p_i)$ is defined as
\beq
        V(p_i) & = & \bigcap_{j\neq i} vh(p_i,p_j).
\eeq
A point $x$ is on the {\bf Voronoi edge} $e(p_i,p_j)$ iff
it is on the intersection of the Voronoi cells $V(p_i)$ and $V(p_j)$, that is
\beq
        x \in e(p_i,p_j) & \desda & x \in V(p_i) \cap V(p_j).
\eeq
The {\bf Voronoi diagram} is the  family of subsets of $\R^2$ consisting of
the Voronoi cells $V(p_i)$ and all of their intersections.
The \bfindex{one-skeleton} of the Voronoi diagram is the union of
the boundaries of the Voronoi cells. Note that the notion of one-skeleton equals
the notion of shape of a Voronoi diagram, introduced in Section \ref{svorpoldef}.

\section{One point continuity.}

In this section we prove the following theorem. 

\begin{theorem}
\elabel{onepcont}
Let $\gamma_n \in \cda$ be a set of data representing $n$ generators in $\R^2$
and $\binom{n}{2}$ angles in $\R/2 \pi \Z$ between these generators.
Suppose that $x$ is on the one-skeleton of the corresponding Voronoi
diagram $V(\gamma_n)$. 
Then for every $\epsilon > 0$ there exists $\delta > 0$,
such that when we perturb $\gamma_n$ by not more than $\delta$, 
the one-skeleton of any perturbed diagram is within distance  $\epsilon$
of $x$.
\end{theorem}

In order to prove this theorem, we perform the following steps:
\begin{looplijst}
\item We decrease $\epsilon$ a finite number of times: for every
      pair of generators $(p,q)$  whose bisector $B(p,q)$ does not
      pass through $x$ we decrease $\epsilon$, if necessary, in such a way that 
\beq
	d(x, B(p,q)) & > & 2 \epsilon.
\eeq
\item For every pair of generators $(p,q)$ we show how
to choose $\delta_{pq}$. 
 If $x \in B(p,q)$  before perturbing then after perturbing by $\delta_{pq}$,
\beq
	B(p,q) \cap B_{\epsilon}(x) & \neq & \emptyset.
\eeq
If  $x \notin B(p,q)$  before perturbing then after perturbing by $\delta_{pq}$,
\beq
        B(p,q) \cap B_{\epsilon}(x) & = & \emptyset.
\eeq
Of course this last statement can only hold
if we adjust $\epsilon$ as indicated above.
\item Set 
$\delta   :=  \min_{p,q} \delta_{pq}$.
\end{looplijst}

These steps  gives us enough control on the bisectors to complete
 the proof.
We first show that the bisection point of two generators $p$ and $q$ can not
get more perturbed than the generators themselves.

\begin{lemma}
\elabel{bisection} Let $p$ and $q$ be two points in the plane.
Let $b(p,q) = (p+q) /2$.
If $\max(~ \| \bar{p} - p \|,~ \|\bar{q} - q \|~) < \delta$, then
\beq
	\| b(\bar{p},\bar{q}) - b(p,q) \| & < & \delta.
\eeq
\end{lemma}

\begin{proof}
$\| b(\bar{p},\bar{q}) - b(p,q)\|$ = $\|\frac{1}{2}(\bar{p}+\bar{q}) - \frac{1}{2}(p+q)\|$
= $\| \frac{1}{2}(\bar{p}-p) + \frac{1}{2}(\bar{q}-q)\| $
	$\leq$ $\frac{1}{2}\|\bar{p}-p\| + \frac{1}{2}\|\bar{q}-q\|$ $<$ $\delta$.
\end{proof}

We apply Lemma \ref{bisection} in choosing $\delta_{pq}$ 
such that a bisector $B(p,q)$
that passes through $x$ before perturbation still passes through
$B_{\epsilon}(x)$ after perturbation.

\begin{lemma}
\elabel{inin}
Let $p$, $q$ and $\alpha_{pq}$  be data from $\gamma_n \in \cda$ such that
        $x  \in  B(p,q)$.
Write $r = \| x - b(p,q)\|$. 
Let 
\beq
	\delta_{pq} & = & \min( \frac{\epsilon}{1+r},1).
\eeq
If $\max(~\|\bar{p}-p\|,~\|\bar{q}-q\|,~\|\bar{\alpha}_{pq}-\alpha_{pq}\|)
	< \delta_{pq}$ then
the bisector $B(\bar{p},\bar{q})$ 
passes through $B_{\epsilon}(x)$.
\end{lemma}

\begin{proof}
We treat the perturbation of the  bisector angle $\alpha_{pq}$
and the bisection point $b(p,q)$ separately and add up the maximal effect.
Let $l_{\bar{\alpha}_{pq}}$ be the image of a rotation
over $\pm \delta_{pq}$ around $b(p,q)$ of the line
$B(p,q)$. Then $x$ is at distance $r \sin \delta_{pq}$ of $l_{\bar{\alpha}_{pq}}$.
Moreover, from Lemma \ref{bisection} it follows that 
$\|b(\bar{p},\bar{q}) - b(p,q)\| < \delta_{pq}$.
Let $d = \| x - B(\bar{p},\bar{q})\|$.  Adding up the two effects gives that
$d < \delta_{pq} + r \sin \delta_{pq}$. We want $d$ to be smaller than 
$\epsilon$. As $\delta_{pq}$ is a majorant
for $\sin\delta_{pq}$, we are safe if we ensure, as claimed,  that
        $\delta_{pq}  <  \epsilon/(1 + r)$.
\end{proof}

The next lemma gives a value of $\delta$ that ensures 
that a bisector that misses $x$ before perturbation, stays away far enough
from $x$ after perturbation.

\begin{lemma}
\elabel{notnot}
Let $p,q$ and $\alpha_{pq}$ be data from $\gamma_n \in \cda$.  Let 
$x \in \R^2$ be such that
\beq
	\| x -  B(p,q) \| & > & 2 \epsilon.
\eeq
Write $r = \| x - b(p,q)\|$.
The angle between the bisector $B(p,q)$ and the line through $b(p,q)$ and
$x$ is denoted by $\gamma \in (0,\frac{\pi}{ 2}]$. 
Let 
\beq
        \delta_{pq} & = & \frac{ \gamma r - 2 \epsilon}{r +2}.
\eeq
If $\max(\,\|\bar{p}-p\|,\,\|\bar{q}-q\|,\,\|\bar{\alpha}_{pq}-\alpha_{pq}\|\,)
        < \delta_{pq}$ then
$\| x - B(\bar{p},\bar{q}) \|  >  \epsilon$. 
\end{lemma}

\begin{figure}[!ht]
 \begin{center}
\setlength{\unitlength}{1cm}
\begin{picture}(5,3.6)
\put(2.5,2){\makebox(0,0)[cc]{
        \leavevmode\epsfxsize=5\unitlength\epsfbox{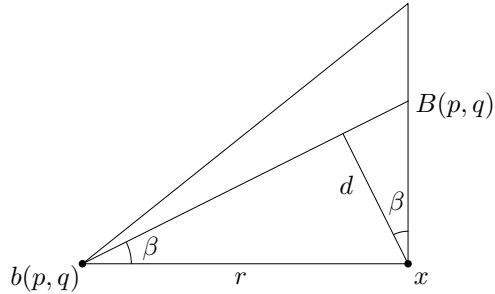}}}
\put(-0.1,0.2){\makebox(0,0)[cc]{$b(p,q)$}}
\put(2.5,0.2){\makebox(0,0)[cc]{$r$}}
\put(4.9,0.2){\makebox(0,0)[cc]{$x$}}
\put(1.3,.6){\makebox(0,0)[cc]{$\beta$}}
\put(3.9,1.44){\makebox(0,0)[cc]{$d$}}
\put(4.57,1.2){\makebox(0,0)[cc]{$\beta$}}
\put(5.35,2.5){\makebox(0,0)[cc]{$B(p,q)$}}
\end{picture}
 \end{center}
 \caption{\elabel{angledist} $d$ is the distance from the bisector $B$ to $x$.}
\end{figure}

\begin{proof}
First we show that indeed $\delta_{pq} > 0$. As
	$2 \epsilon < r \sin \gamma < r \gamma$ it follows that 
	$\gamma r - 2 \epsilon > 0$.
We treat the perturbation of the bisector angle $\alpha_{pq}$  and 
the bisection point $b(p,q)$ separately and add up the maximal
effect.
Let $l_{\bar{\alpha}_{pq}}$ be the image of a rotation
over $\pm \delta_{pq}$ around $b(p,q)$ of the line $B(p,q)$.
The 
distance $d$ of $x$ to $l_{\bar{\alpha}_{pq}}$ s given by $r \sin \beta$. 
Here $\beta := \gamma - \delta_{pq}$, see also Figure \ref{angledist}.
Moreover, from Lemma \ref{bisection} it follows that
$\| b(\bar{p}, \bar{q}) - b(p,q) \| < \delta_{pq}$
This implies that the composite minimal distance, that is the distance of 
$x$ to $B(\bar{p},\bar{q})$ 
could become as small as $d_1 :=  r \sin (\gamma - \delta_{pq}) - \delta{pq}$.
We want this distance $d_1$ to stay bigger than $\epsilon$. As
$\frac{1}{2}(\gamma-\delta_{pq})$ is a minorant of $\sin (\gamma - \delta_{pq})$
we are safe if we ensure, as claimed,  that
\begin{equation*}
        \delta_{pq} ~ < ~ \frac{ \gamma r - 2 \epsilon}{r +2}. \qedhere
\end{equation*}
\end{proof}

The following corollary is obvious but useful.

\begin{corollary}
\elabel{hpq}
Let $x, p, q, \alpha_{pq}$ and $\delta_{pq}$ be as in Lemma \ref{notnot}.
Suppose that $(p,q,\alpha_{pq})$ is perturbed by a vector of 
length at most $\delta_{pq}$. 
For every point $y \in B_{\epsilon}(x)$ it holds that
\beq
y \in vh(p,q)~ {\mbox {\rm before perturbation}} & \desda &
y \in vh(p,q)~ {\mbox {\rm after perturbation,}}
\eeq
where $vh(p,q)$ denotes the Voronoi half-plane of generators $p$ and $q$
containing $p$.
\end{corollary}

We are ready to prove Theorem \ref{onepcont}.

\begin{proof}[{\bf Proof of Theorem \ref{onepcont}.}]
\elabel{finproof} 
Let $x$ be on the $1$-skeleton of $V(\gamma_n)$ before perturbation. Let
$\epsilon > 0$.
Adjust $\epsilon$ a finite number of times: for every
pair of generators $(p,q)$  whose bisector $B(p,q)$ does not
pass through $x$ we set
\beq
	\epsilon & = & \min (\epsilon, \frac{1}{3}  d(x, B(p,q)) ).
\eeq
This ensures that for every such pair $(p,q)$ it holds that 
$d(x, B(p,q))  >  2 \epsilon$.

Next determine a $\delta_{pq}$ for every  pair of generators in such a way
that i) 
if $x \in B(p,q)$, then the condition in Lemma \ref{inin} is fulfilled; 
ii) if $x \notin B(p,q)$, then  the condition of Lemma \ref{notnot} is fulfilled. 
Let $\delta$ be the minimum of all these $\delta_{pq}$, that is,
$\delta~:=~\min_{(p,q)} \delta_{pq}$.
Define
\beq
	P & := & \{ p \in \gamma_n \mid V(p) ~\cap~ B_{\epsilon}(x)  \neq  
	\emptyset, ~\mbox{{\rm before perturbation}}\}.
\eeq
Consequently, $P$ has at least two elements. Say 
$P = \{p_1, \dots, p_m \}$, for some $m \geq 2$.
Before perturbation it holds that
\beq
         B_{\epsilon}(x) & = & \bigcup_{i=1}^{m}
		 ( V(p_i) \cap B_{\epsilon}(x) ). 
\eeq
Suppose that, after perturbation, $B_{\epsilon}(x)  \subset  V(s)$,
for some generator $s$. We show that this leads to a contradiction.
First note that $s \in P$ leads to a contradiction immediately: 
we can apply Lemma \ref{inin}.
This implies that after perturbation some bisector
$b(s, p_i)$ passes through $B_{\epsilon}(x)$. 
We are left with the case that $s \notin P$. Fix a point
$y \in B_{\epsilon}(x)$ such that $y$ is not on any bisector
before or after perturbation. Then
\begin{center}
\begin{tabular}{ll}
        $y\in V(p_1)$ &  before perturbation, say,\\
        $y \in V(s)$  &  after perturbation.
\end{tabular}
\end{center}
We concentrate on the bisector $B(p_1,s)$ now. If $x \not\in B(p_1,s)$
before perturbation, we apply Corollary \ref{hpq}:
any point in $B_{\epsilon}(x)$
stays in the same half-plane $vh(p_1,s)$. 
As a consequence, $y  \in  vh(p_1,s)$
before and after perturbation. This is in contradiction with
our assumption that $y\in V(s)$ after perturbation. So $x \in B(p_1,s)$
before perturbation. This implies that
\[
        B(p_1,s) \cap B_{\epsilon}(x) \neq \emptyset,
\]
after perturbation: a contradiction with our assumption as well.
\end{proof}
\section{The distance between two Voronoi diagrams.}

\elabel{shausdorff}
\subsection{Hausdorff distance.}
If we want to compare two Voronoi diagrams, a suitable notion of distance 
is the Hausdorff distance: 
two sets $A$ and $B$ are within \bfindex{Hausdorff distance} r 
iff $r$ is the smallest number such that  any point of $A$ is within
distance $r$ from some point of $B$ and vice versa.
Let us give this definition more formally. 
Suppose we have a metric space $(X,d)$. 
For $A \subset X$ and $r >0$ we define the open neighborhood
$N_r(A)$ as the set
\beq 
N_r(A) & := & \{~ y ~|~ d(x,y) < r, ~{\rm for~some~} x \in A \}.
\eeq 
Think of $N_r(A)$ as the territorial border of $A$.

\begin{definition}
Let $(X,d)$ be a metric space.  The {\bf Hausdorff distance}
$h(A<B)$ between
two subsets $A,B \subset X$ is defined as
$h(A,B)  :=  \inf\{r: A \subset N_r(B)~{\rm and}~ B \subset N_r(A) \}$.
\end{definition}

The Hausdorff distance $h$ defines a metric on the set of nonempty 
compact subsets of $(X,d)$, see \cite{Ca}.  Note that
\beq	
h(A,B) < r & \desda & A \subset N_r(B) ~\wedge~ B \subset N_r(A).
\eeq

\subsection{A compactness condition.}
\elabel{compactness}

In this section we put a compactness condition on the configuration
space of data sets. It consists of two ingredients. Suppose that
$\{p_1, \dots, p_n\}$ denotes the underlying point set of some
data set $\gamma_n \in \cda$. First we restrict the domain of $\{p_1, \dots, p_n\}$
to a closed disk $\U$. Moreover, we add four camera points
$c_1, c_2, c_3$ and $c_4$ far away outside  $\U$.
Lemma \ref{complement} shows that we can choose $\U$ and
$c_1, c_2, c_3$ and $c_4$ in such a way that it is guaranteed that the 
Voronoi diagram of $\gamma_n$, extended with $c_1, c_2, c_3$ and $c_4$,
changes inside some compact subset of the plane only, provided that
we perturb $\{p_1, \dots, p_n\}$ within $\U$.

For $\U$ we choose  the unit disk, that is
$\U  :=  \{ x \in \R^2 ~\mid~ \|x\| \leq 1 \}$.

\begin{definition} 
Let $\gamma_n \in \cda$ be a set of data representing $n$ generators
$\{p_1,\dots,p_n\}$ in $\U$
and $\binom{n}{2}$ angles in $\R/2 \pi \Z$ between those generators. Let furthermore
$$
\begin{array}{ll}
        c_1 = (-N,0), & c_2 = (0,N), \\
        c_3 = (N,0), & c_4 = (0,-N).
\end{array}
$$
for some $N \gg 1$ be the so-called \bfindex{camera points}. The 
{\bf Voronoi diagram with camera points} $V_N(\gamma_n)$ is 
the Voronoi diagram of the data set $\overline{\gamma}_n$, which consists
of the generators $ \{p_1,\dots,p_n\} ~\cup ~ \{ c_1, \dots, c_4 \}$
and the $\binom{n+4}{2}$ angles between those generators. 
\end{definition}  

\begin{lemma}
\elabel{complement}
Let $N > 2 + 2 \sqrt{2}$ and let $B_N^C := \{x \in \R^2 \mid d(x,0) > N\}$.
Then $V_N(\gamma_n) \cap B_N^C$ does not change, provided that $\gamma_n$
is perturbed inside $\U$.
\end{lemma}

\begin{proof}
Let $z=R(\cos \phi, \sin \phi)$ be an arbitrary point in $B_N^C$.
It is enough to show that  $z$ is always closer to at least one camera point
than to $\U$. 
Using the symmetry we can assume, without loss of generality, that
$0  ~< ~\phi ~ \leq ~ \frac{\pi}{4}$.
So it is enough to show that $ d(z, c_3) ~ < ~ d(z, \U)$
or, equivalently, that $d(z, \U)^2 - d(z,c_3)^2 ~ > ~ 0$.
Now,
\beq
	d(z, c_3)^2 & = & R^2 \sin^2 \phi + (N-R \cos \phi)^2, \\
		  & = & R^2 + N^2 - 2 N R \cos \phi,
\eeq
is maximal when $\phi = \frac{\pi}{4}$. From now on assume that
\beq
	d(z, c_3)^2 & = & R^2 + N^2 - \sqrt{2} N R .
\eeq
As $d(z,\U)=R-1$, we have to show that
\beq
	d(z,\U)^2 - d(z,c_3)^2 \quad 
	= \quad (R^2 - 2 R + 1) - (R^2 + N^2 - \sqrt{2} N R) \quad > \quad 0.
\eeq
Substituting $R = N + X$, where $X>0$, gives that
\beq 
	d(z,c_3)^2 - d(z,\U)^2 & = & N((\sqrt{2}-1)N-2) + (\sqrt{2}N-2) X + 1.
\eeq
Now, $ (\sqrt{2}-1)N-2 > 0$, whenever $N > 2 + 2 \sqrt{2}$, while
$(\sqrt{2} N - 2)  > 0$ whenever $N > \sqrt{2}$.
This proves the lemma.
\end{proof}

\section{Continuity of the Voronoi map.}

\subsection{The Voronoi map.}

In our original configuration space, generators live in $\R^2$.
But under the compactness condition above, generators all
live on the unit disk $\U$ and we have four additional camera points,
$c_1$, $c_2$, $c_3$ and $c_4$.
That is, we have a restricted configuration space
\begin{eqnarray*}
\cdau & := & \{c_1,\dots,c_4; p_1, \dots, p_n;
	\alpha_{12},\dots,\alpha_{(n+3)(n+4)} \},
	\quad p_i \in \U; ~\alpha_{ij} \in \R/ 2 \pi \Z.
\end{eqnarray*}
Let $V(n) \subset P(\R^2)$ be the set of one-skeletons of
Voronoi diagrams, defined by data sets $\gamma_n \in \cdau$. Denote by 
$d$ be the Euclidean metric on $\cdau$ and by $h$ the Hausdorff
metric on  $P(\R^2)$. We show that the map $f_V$,
\bmap
f_V:& (\cdau,d)   & \rightarrow & (V(n),h), \\
    &    \gamma_n & \mapsto & V(\gamma_n),
\emap
that maps a data set $\gamma_n$ to its Voronoi diagram, is
continuous. By definition, this means  that we have to show that
\begin{multline}
\elabel{ucont}
\forall \gamma_n \in \cdau,~ \forall \epsilon > 0, ~\exists
 \delta > 0, ~\forall  \eta_n \in \cdau   \\
d(\gamma_n,\eta_n) < \delta \quad \Rightarrow \quad
h(V(\gamma_n),V(\eta_n)) < \epsilon. \qquad
\end{multline}
In Theorem \ref{onepcont} we have proved that
\[
\forall (\gamma_n,x \in V(\gamma_n)),~ \forall 
	\epsilon > 0,~ \exists \delta >0:~ 
     d(\gamma_n,\eta_n) \leq \delta ~\Rightarrow~ d(x,V(\eta_n)) \leq \epsilon.
\]
In this formula, $\delta$ really depends on both the particular position
of $x$ on $V(\gamma_n)$ and the particular diagram $V(\gamma_n)$. 
The uniform version of this claim is given by 
\begin{multline}
\elabel{uni}
\forall \epsilon >0, \quad \exists \delta > 0,\quad \forall (\gamma_n,x \in V(\gamma_n)):~\\
	d(\gamma_n,\eta_n) \leq \delta \quad\Rightarrow\quad d(x,V(\eta_n)) \leq \epsilon.
\qquad
\end{multline}

\begin{remark}
It is clear that (\ref{uni}) does not hold if we regard
the one-skeleton of Voronoi diagrams that correspond with data sets
that do not contain the camera points.
Suppose for example that our data consist of two points 
$p$ and $q$ that are close to
each other. If we fix $p$, move $q$ slightly, thereby changing 
$\alpha_{pq}$ slightly as well, some point that was far away on the bisector 
$B(p,q)$ before moving $q$, is on a big distance from $B(p,q)$
after disturbance. 
\end{remark}

\subsection{$\{\mathbf{(\gamma_n, x)~|~ x \in V(\gamma_n))\}}$ is closed
	in $\mathbf{\cda \times \R^2}$.}

In this section we show that the set
$ \{(\gamma_n, x) ~|~ x \in V(\gamma_n)\}$ is closed in
$\cda \times \R^2$.
Or, equivalently, that the complement 
$ \{(\gamma_n, x) ~|~ x \not\in V(\gamma_n)\} $ 
is open. But this means exactly that we have to prove the following:
$$ \forall \gamma_n,~x \not\in V(\gamma_n),~\exists \epsilon >0,~\exists \delta>0:~
	d(\gamma_n,\eta_n) ~<~ \delta ~\Rightarrow~ d(x, V(\eta_n))~>~\epsilon.
$$
We do so in the following lemma.

\begin{lemma}
\elabel{notonlemma}
Let $\gamma_n$ be a data set representing $n$ generators in the plane
and $\binom{n}{2}$ angles between those generators.
Suppose that $x$ is 
in the interior of Voronoi cell $V(p)$, for some generator $p\in \gamma_n$.
Then there exists $\epsilon > 0$ and  $\delta > 0$
such that the following holds:  if we perturb $\gamma_n$ by not 
more then $\delta$,
the one-skeleton of any perturbed diagram is at distance at least 
$\epsilon$ of $x$.
\end{lemma}

\begin{proof}
{\bf [Determine $\mathbf{\epsilon}$.]} 
Start with some $\epsilon > 0$. 
For every pair of generators $(p,q)$  whose bisector $B(p,q)$ does not
pass through $x$ we set
\beq
        \epsilon & = & \min (\epsilon, \fr{3}  d(x, B(p,q)) ).
\eeq
This ensures that $d(x, B(p,q)) ~>~ 2 \epsilon$
for every such pair $(p,q)$.

{\bf [Determine $\mathbf{\delta(\epsilon)}$.]} 
As a consequence, we can apply Lemma \ref{notnot} for every pair of generators
$(p,q)$ whose bisector $B(p,q)$ does not pass through $x$. The lemma gives us 
$\delta_{pq}$ such that any perturbation of the data by $\delta_{pq}$
implies that $d(x,B(p,q)) ~ > ~ \epsilon$
after perturbation. Let $\delta~:=~\min_{(p,q)} \delta_{pq}$ be the minimum 
of all these $\delta_{pq}$.

{\bf [Show that this $\delta$ works.]}
Suppose that the claim of the  lemma is not true. 
This means that after perturbation
there are generators $q_1$ and $q_2$ such that both
$$
	B_{\epsilon} (x) ~\cap~ V(q_1) \neq \emptyset,
\qquad \text{and} \qquad
	B_{\epsilon} (x) ~\cap~ V(q_2) \neq \emptyset,
$$
hold. It follows that, after perturbation, 
	$B(q_1,q_2) ~\cap~ B_{\epsilon}(x) ~\neq~ \emptyset$.
Because of our choice of $\delta$ and $\epsilon$ this implies though that
$x  \in  B(q_1,q_2)$
before perturbation. As $x \in V(p)$, before perturbation, we conclude that 
	$q_1 \neq p$ and $q_2 \neq p$.
But then there exists $y \in B_{\epsilon}(x)$ such that
$$	
	y \in h(p,q_1),
$$
{\it before} perturbation, and
$$ 
	y \in h(q_1,p),	
$$
{\it after} perturbation. This in contradiction with Corollary \ref{hpq}
however.
\end{proof}

\subsection{Proof of the continuity.} 

In this section we prove the following theorem:

\begin{theorem}
\elabel{hausdorffcont}
Let $\gamma_n \in \cdau$. 
Then the Voronoi diagram $V(\eta_n)$ of any data set
$\eta_n \in  \cdau$ that is Euclidean-close to $\gamma_n$, 
is Hausdorff-close to the Voronoi diagram $V(\gamma_n)$:
\beq
\forall \epsilon >0,~ \exists \delta > 0,~ :~
        d(\gamma_n,\eta_n) \leq \delta & \Rightarrow & h(V(\gamma_n),V(\eta_n)) 
	\leq \epsilon
\eeq
\end{theorem}

\begin{proof}
The proof is divided into a number of steps.\setcounter{proof}{0}

\step
\begin{bf}
It is enough to prove the following assertion.
\end{bf}
\begin{multline}
\elabel{cuni}
 \forall \epsilon > 0, \quad \exists \delta > 0,\quad \forall \gamma_n \in \U,\quad
\forall x \in V(\gamma_n),\quad \forall \eta_n \in \U: \\
        d(\gamma_n,\eta_n) \leq \delta \quad \Rightarrow \quad d(x,V(\eta_n)) 
	\leq \epsilon.
\qquad
\end{multline}
In (\ref{cuni})
the implication holds for all choices of $x$ on $V(\gamma_n)$. This implies
that in fact
\beq
        d(\gamma_n,\eta_n) \leq \delta & \Rightarrow & d(x, V(\eta_n))
	 \leq \epsilon \qquad \forall x.
\eeq
This together with the fact that we can interchange $\gamma_n$ and $\eta_n$ defines
the Hausdorff distance and therefore implies the claim of the theorem.

\step
\begin{bf}
Construct a convergent sequence.
\end{bf}\\
Suppose that Assertion \ref{cuni} does not hold. That means that
the negation must be true, where we replace $\forall \delta$ by 
$\forall m$:
\begin{equation}
\begin{split}
\elabel{neguni} \exists \epsilon>0: ~\forall m \in \N,~ \exists (\gamma_n, x \in V(\gamma_n)),~
\exists \eta_n: \\
 ( d(\gamma_n, \eta_n) \leq \frac{1}{m} \quad\wedge\quad 
		d(x, V(\eta_n)) > \epsilon ).
\end{split}
\end{equation}
Fix such $\epsilon$ and call it $\epsilon_0$.
We can find for every $m \in \N$ a triple $t_m$
\begin{eqnarray}
\elabel{triple}
	t_m & := & ( \gamma_n^{m},~ x^{m} \in V(\gamma_n^{m}),~ \eta_n^{m} ),
\end{eqnarray}
such that Assertion \ref{neguni} is true.

\step $\mathbf{x^m \in B_N,~\forall m}$.\\
Recall that $B_N$ denotes the disk $B_N=\{x \in \R^2~|~d(x,0) \leq N\}$.
Suppose that $x^m \not\in B_N$ for some  $x^m \in V(\gamma_n^m)$. Then Lemma
\ref{complement} tells us that
$ x^m \in V(\eta_n)$, for all $\eta_n \in \cdau $.
So Assertion \ref{neguni} can never be true.

\step {\bf Remarks on compactness.}\\
$\gamma_n$ and $\eta_n$ both live in the compact set
$$
    	\U^n \times \{c_1,c_2,c_3,c_4\}\times (\R/2 \pi \Z)^{\binom{n+4}{2}}.
$$
We have proved in  Lemma \ref{notonlemma} that the set 
$ \{(\gamma_n, x) ~|~ x  \in V(\gamma_n) \}$
is closed. This means that  we get a sequence $(t_m)$  of 
triples of the form (\ref{triple}) that
live on a compact set. So, $(t_m)$ has some convergent subsequence
$t_m(k)$ that converges to, say,
$(\tilde{\gamma}_n, \tilde{x}, \tilde{\eta}_n)$.
Note that $\tilde{x}  \in  V(\tilde{\gamma}_n)$.

\step
\begin{bf}
A contradiction by combining Assertion \ref{neguni} and Theorem \ref{onepcont}.
\end{bf}\\
Because of convergence, $V(\tilde{\gamma}_n)=V(\tilde{\eta}_n)$, so 
$\tilde{x}  \in  V(\tilde{\eta}_n)$.
	 Apply Theorem \ref{onepcont} with
	\begin{itemize}
		\item $\tilde{x} \in V(\tilde{\eta}_n)$,
		\item $\epsilon ~=~\fr{3} \epsilon_0$.
	\end{itemize}
      The theorem gives us a $\delta_0$ such that  
\begin{eqnarray}
        \elabel{fromth}
	d( \tilde{\eta}_n, \theta_n) ~<~ \delta_0 & \Rightarrow &
	d( V(\theta_n), \tilde{x} ) ~<~ \fr{3}{\epsilon_0}.
\end{eqnarray}
Choose $m$ so big that the following condition hold:
  \begin{eqnarray} 
	\elabel{climit}
	  d(~(\gamma_n^m,x^m,\eta_n^m),(\tilde{\gamma}_n, \tilde{x},\tilde{\eta}_n)~)
		& < & \min( \delta_0, \fr{3}\epsilon_0).
  \end{eqnarray}
We show that this leads to a contradiction on $d(x^m, V(\eta_n^m) )$.
\begin{lijst}
	\item The sequence $t_m$ is constructed in such a way from
	 Assertion \ref{neguni} that 
	\begin{eqnarray}
	\elabel{bigger}
	d(x^m, V(\eta_n^m) ) & > & \epsilon_0.
        \end{eqnarray}

	\item
We have chosen $m$ so big, (\ref{climit}), that
\begin{eqnarray}
	\elabel{sm1}
	d(x^m, \tilde{x}) & < & \fr{3} \epsilon_0.
\end{eqnarray}
\item But $m$ is also big enough such that $d(\tilde{\eta}_n,\eta_n^m) ~<~ \delta_0$.
Therefore we conclude from (\ref{fromth}), with $\theta_n := \eta_n^m$,  that
\begin{eqnarray}
	\elabel{sm2}
	d(V(\eta_n^m), \tilde{x}) & < & \fr{3} \epsilon_0.
\end{eqnarray}
Now (\ref{sm1}) combined with (\ref{sm2}) gives a contradiction  with (\ref{bigger}).
\end{lijst}
This proves the theorem.
\end{proof}

\newpage
\thispagestyle{empty}
\chapter{Clickable Voronoi diagrams and hook compactifications.}

\label{chmodel}
We define a compactification space $\XAH$ of the configuration 
space of $n$ distinct points in the plane, by
considering data elements of pairs of points and triples of points.
For every pair of points we write down the angle mod $\pi$
between the two points, and for every ordered triple $(p_i,p_j,p_k)$
of points, we specify a hook $\ho{ik}{ij}$. This hook expresses how 
to construct the point $p_k$ given $p_i$ and $p_j$.
Now $\XAH$ is defined as the closure of the image space of
all angles and hooks on $n$ distinct points. 
We show that configurations that are added by taking the closure 
have a natural nested structure, easily revealed by analyzing
the hooks. The main result of this chapter is an explicit 
construction establishing $\XAH$ as the graph of a function.
This construction shows that $\XAH$ is a smooth manifold. 
If we replace the angles mod $\pi$ by angles mod $2 \pi$
we get a compactification space $\xedah$. This space is isomorphic
to the manifold with corners $\FMt(n)$, introduced by
Kontsevich-Soibelman. On $\xedah$ we define clickable
Voronoi diagrams.

\section{Fulton-MacPherson related models.}

In this chapter we describe a compactification $\XAH$ of the
configuration space of $n$ distinct points in the plane.
The basic idea of the compactification consists 
of considering the geometry of all pairs and all
triples of distinct points. In section \ref{sfmintro}
we give an overview of the contents of this chapter: it
describes in an informal way the construction and some properties
of this compactification $\XAH$. And it explains how clickable
Voronoi diagrams are related to $\XAH$. Sections \ref{sah}
to \ref{sclick} contain the mathematical constructions
referred to in Section \ref{sfmintro}, while Section 
\ref{sfmconclusion} is devoted to concluding remarks.
  We start however in this section by recalling the
famous compactification due to Fulton an MacPherson,
and by presenting a related compactification as described by
Kontsevich and Soibelman. These two compactifications have some
important features that we mimic and extend in our compactification
$\XAH$.

\subsection{The Fulton-MacPherson compactification.}

The following compactification is defined in terms of
a nonsingular algebraic variety. The authors note however,
cf. \cite{FM}, page 188/189, that the same constructions
work for complex manifolds, as well as for real manifolds.
So, let $X$ be a nonsingular algebraic variety, and let
\beq
	 \CONF(X) & = & \{ (p_1, \dots, p_n) \in X^n ~\mid~
		p_i \neq p_j ~\text{if}~ i \neq j \},
\eeq
be the configuration space of $n$ distinct labeled points in $X$,
cf.\ Definition \ref{dconf}.
For a subset $S$ of the set of labels $\{1, \dots, n\}$,
the small diagonal $\Delta_S$ is given by
\beq
        \Delta_S & = & \{ (p_1, \dots, p_n) \in \CONF(X) ~|~
                p_{i_1} = \dots = p_{i_m} \text{ for any }
                i_j \in S \}.
\eeq
For any such $S$ with at least two points, 
let $\text{Bl}_{\Delta}(X^S)$ be the blow-up of $X^S$ along $\Delta_S$.
For the definition of blow-up, consult \cite{Ha}.
There is a natural embedding 
\beq
	\CONF(X) & \subset & X^n \times \prod_{|S| \geq 2} \text{Bl}_{\Delta}(X^S).
\eeq
The \index{Fulton-MacPherson!compactification}
{\bf Fulton-MacPherson compactification} $X[n]$ is defined as the closure
of $\CONF(X)$ in this product.
We list some properties that are stated respectively as 
Theorems 1, 2 and 3 in \cite{FM}.
The definition of divisor can be found in \cite{Ha}.
\begin{looplijst}
\item $X[n]$ is nonsingular.
\item For $n\geq 2$, $X[n]$ is the closure of $\CONF(X)$ in the product
      of the $\text{Bl}_{\Delta}(X^S)$ for $S\subset \{1, \dots,n\}$
      of cardinality 2 and 3.
\item For each $S \subset \{1, \dots,n\}$ with at least two elements,
      there is a nonsingular divisor $D(S) \subset X[n]$ such that 
\begin{lijst}
	\item The union of these divisors is $X[n] \backslash \CONF(X)$.
	\item An intersection of divisors $D(S_1) \cap \dots \cap D(S_r)$ 
        is nonempty iff the sets $S_k$  are \bfindex{nested} in the
        sense that each pair $S_i$ and $S_j$ is either disjoint or one is
        contained in the other.
     \end{lijst} 
\end{looplijst}
Important for us is the `first geometric description' of $X[n]$, as presented
in Section~1 of \cite{FM}. A configuration in $X[n]$ is called
\bfindex{degenerate}, if not all points in the configuration are distinct.
Any such configuration $c$ is described by a set of data given the 
locations of the points, and if two or more points coincide at say $x$,
then a \bfindex{screen} is specified
for the set $S \subset \{1, \dots, n\}$ of labels of the coinciding points.
The data elements describing a screen for $S$ at $x$ consist of a labeled
set of points $x_a$ in the tangent space $T_x$ of $x$, such that
not all points in $x_a$ are equal. This process is repeated for $x_a$
until all points in the configuration $c$ are separated in some screen.
Think of repeatedly zooming in at clusters of points until all points are 
separated, cf.\ Figure \ref{screens.eps}.
The divisor $D(S)$ of $X[n]$, mentioned in the property above
consists of all configurations $c$ in $X[n]$ that have a screen that 
contains exactly the points $x_a$ for $a$ in $S$.

\subsection{The Fulton-MacPherson operad.}

\elabel{sKS}
In \cite{Ko} and \cite{KS} the 
\index{Fulton-MacPherson!operad}
{\bf Fulton-MacPherson operad} $\FMd$ is discussed.
A definition of operad can be found in \cite{Ko}.
As we consider only point configurations in $\R^2$,
we set $d=2$ in the sequel. 
The operad $\FMt=\{\FMt(n)\}_{\geq0}$ is defined as follows, cf.\ \cite{KS}:
\begin{looplijst}
\item $\FMt(0) := \emptyset$.
\item $\FMt(1) :=$ point.
\item $\FMt(2) := \conft = S^1$.
\item For $n \geq 3$, the space $\FMt(n)$ 
   \index{FMt@$\FMt(n)$} is a manifold with corners,
      its interior is $\conf$, and all boundary strata are certain products
      of copies of $\confp$ for $n'<n$.
\end{looplijst} 
A \bfindex{manifold with corners} is  a topological space that is
locally homeomorphic with $\R^n_{\geq 0}$. The space $\FMt(n)$, $n \geq2$ 
is defined explicitly as follows.
\begin{definition}
For $n\geq2$, the manifold with corners $\FMt(n)$ is the closure of the image of
$\conf$ in the compact manifold $(S^1)^{\binom{n}{2}} \times 
	[0, + \infty]^{6\binom{n}{3}}$ under the map
\beq
	[(p_1, \dots, p_n)] & \mapsto & ( (\alpha_{ij})_{1\leq i<j\leq n}, 
	\be{ik}{ij}),
\eeq
where $i$, $j$, and $k$ are pairwise distinct indices,
$\alpha_{ij} \in \R / 2 \pi \Z$,  and $\be{ik}{ij}=
	\frac{|p_i-p_k|}{|p_i-p_j|}$.
\end{definition}

That is, Kontsevich and Soibelman define the space $\FMt$ in terms of data elements
for pairs and triples of points in any reduced configuration
$c \in \conf$: for any pair of points they write down the directed angle
between the two points, while
for every ordered triple $(p_i,p_j,p_k)$ of points they specify the ratio
$\be{ik}{ij}$ of the line segments $p_ip_k$ and $p_ip_j$.

\subsection{Combining the models.}

From the Fulton-MacPherson compactification we borrow the screen model, 
while from Kontsevich-Soibelman we borrow the angles and ratios. We extend, 
however, the ratios to hooks, by marking both the ratio $\be{ik}{ij}$
of the line segment $p_ip_k$ and $p_ip_j$ and the angle $\al{ik}{ij}$
between the line segments. In this way every hook has two
representatives. One with a positive ratio $\be{ik}{ij}$ and one 
with a negative `ratio' $\be{ik}{ij}$. In the latter case we just
add $\pi$ to the hook angle $\al{ik}{ij}$. These two distinct 
representatives are identified by an equivalence relation $\sim_k$. 
This will be important for getting a smooth model. 
Another adjustment that we make for obtaining smoothness is that we 
write down  the angle between two distinct points up to multiples of $\pi$
instead of up to $2 \pi$. To summarize, we get a map
\beq
	\psi_{\ah}: \conf(\R^2) & \rightarrow &
		\ah ~:=~ (\R / \pi \Z)^{\binom{n}{2}} \times
			( ([-\infty, \infty] \times \R/ 2 \pi \Z) / \sim_k)
				^{6\binom{n}{3}},
\eeq
and define $\XAH$ as the closure of the image of $\conf$ in $\ah$.
The data elements and spaces introduced so far are discussed in detail
in Section \ref{sah}. In section \ref{nestsscreens} we show that we can associate
with any $x \in \XAH$ a nested set of subsets of $\{1, \dots, n\}$.
This is done by analyzing the ratios. It turns out that for degenerate
configurations certain ratios equal zero. If 
$\be{ik}{ij}$ is close to zero, then the length of the line segment $p_ip_k$ is
very small compared to the length of the line segment $p_ip_j$.
We then say that $p_i$ and $p_k$ coincide with respect to
$p_j$.

Combinations of compactifications of configuration spaces and combinatorics  
receive a lot of attention in recent years. Elegant examples can be found in 
\cite{De1} and \cite{De2}.

\section{Informal introduction.}
\elabel{sfmintro}

In this section we introduce most of the terminology that we use
in the rest of the chapter, in a more intuitive setting.

\subsection{Screens, clusters and nests.}

Suppose we are looking at a plane 
containing $n\geq 2$ {\bf sites}. By sites we mean points, but in order
to be able to use points in other contexts, we prefer
the word sites.
 When we see all $n$ sites
at first glance we are happy, but when we seem to see
fewer sites we tend to look more closely.
That is, we zoom in at a site in order to see if it is really
one site, or in fact consists of a number of sites. At the
same time we will lose sight of other sites. We repeat
this procedure until we have found all $n$ sites.

Imagine that the plane
is in fact a computer {\bf screen}. If some sites coincide, that
is, if they form a {\bf cluster}, we click on the cluster.
Another screen pops up that contains exactly the sites in the 
cluster. If  we click on a site and nothing happens, we
know that at the level we are looking at, this site does not coincide with another
site. So every non-trivial cluster corresponds to a screen.
By continuing to click we always find all $n$ sites.
Label the sites with labels from $1$ to $n$ and
write down the labels that are visible in every screen.
This produces a {\bf nest},  a nested set of subsets of $\{1, \dots, n\}$,
i.e.\ the clusters.  If no sites coincide the nest consists 
of the set $\{1, \dots, n\}$ only.

The {\bf top screen} is the first screen we see. It contains all sites,
but maybe not all sites are separated.
We organize things in such a way that every screen
contains at least $2$ non-coinciding sites. That is, in every screen
at least two sites are separated.  We are
not interested where in the plane a screen focuses. But we do
want to know the exact relative position of the sites.  We conclude that
we consider site sets up to {\bf scaling} and {\bf translation}.

\subsection{Pinpointing sites.}

We want to describe such a family of  screens filled by sites.
A naive approach is by  listing the  screens together with the 
coordinates of the sites that occur in each screen. This may describe  
the situation well when the sites do not move, but assume that two sites
$p_i$ and $p_j$ that did coincide first in some screen start moving apart.
 This could change the relative positions of the sites in all screens
were both $p_i$ and $p_j$ occur, but it is not at all clear how
a change in coordinates in one screen should influence coordinates
in other screens. 

\begin{figure}[!ht]
\begin{center}
\setlength{\unitlength}{1.cm}
\begin{picture}(5,3.2)
\put(2.5,1.6){\makebox(0,0)[cc]{
        \leavevmode\epsfxsize=4\unitlength\epsfbox{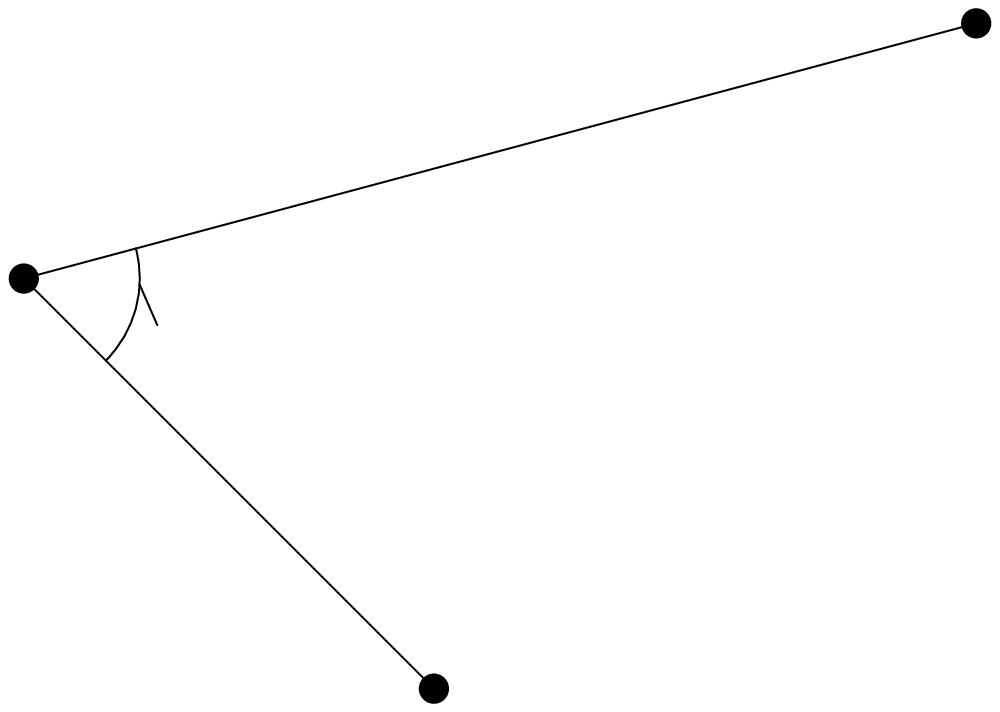}}}
\put(1.35,1.8){\makebox(0,0)[l]{$\al{ik}{ij}$}}
\put(2.4,2.8){\makebox(0,0)[l]{$\be{ik}{ij}$}}
\put(0.2,2){\makebox(0,0)[l]{$p_i$}}
\put(2.5,0.5){\makebox(0,0)[l]{$p_j$}}
\put(4.7,2.9){\makebox(0,0)[l]{$p_k$}}
\end{picture}
 \caption{\elabel{fhook}
The hook $\ho{ik}{ij}=(\be{ik}{ij},\al{ik}{ij})$ hinged at $p_i$ from $p_j$ to $p_k$.
The little arrow indicates the positive direction.}
\end{center}
\end{figure}

Instead of marking coordinates of the sites involved, we mark the {\bf angle}
$\alpha_{ij}$ between any two sites $p_i$ and $p_j$ and the {\bf hook} $\ho{ik}{ij}$
between any ordered triple $(p_i,p_j,p_k)$ of sites. A hook $\ho{ik}{ij}$
consists of an angle $\al{ik}{ij}$ and a ratio $\be{ik}{ij}$. Consider 
Figure \ref{fhook}. Suppose that the sites $p_i$ and $p_j$, thus the leg
$p_ip_j$ have been constructed. If we rotate $p_j$ by angle $\al{ik}{ij}$
with respect to center of rotation $p_i$ and multiply the image of the
rotation by $\be{ik}{ij}$, we get $p_k$.

Using these hooks and angles we can fill a screen with sites. Suppose that for a 
set of distinct sites $p_1, \dots, p_n$ we know all hooks and angles between
the sites. We show how to fill a screen with the sites. We put $p_1$ in
the origin, as we do not care about translation. As we do not care about scaling
either, we put $p_2$ at distance $1$ of $p_1$. We use the angle $\alpha_{12}$
to assure that the line through $p_1$ and $p_2$ has the right direction, that is,
we set $p_2 = (\cos \alpha_{12}, \sin \alpha_{12})$. All other sites $p_i$,
for $i \in 3, \dots,n$, can be constructed as in Figure \ref{fhook}
by a hook $\ho{1i}{12}$.

\subsection{Degenerate configurations.}

Let $\CONF$ be the configuration space of $n$ distinct sites in the plane
and let $\psi_{\ah}$ be the map that maps an element $c \in \CONF$
to the $\binom{n}{2}$ angles and the $6\binom{n}{3}$ hooks between 
the sites. In fact we embed the image space in a slightly bigger
space $\ah$, allowing also negative ratios. The space that we are really
interested in is the closure of the image of $\CONF$ in $\ah$.
We call this space $\XAH$. Most points $x \in \XAH$ can be realized
as image points of a configuration $c \in \CONF$, that is
$x = \psi_{\ah}(c)$. We call a point $x \in \XAH$ {\bf degenerate} if
$x \in \overline{\psi_{\ah}(\CONF)} \backslash \psi_{\ah}(\CONF)$.
It turns out that the degenerate points are exactly those points $x \in  \XAH$
 such that at least one ratio coordinate $\be{ik}{ij}(x)=0$. 

We go back to the nests and screens. By analyzing the set of all 
$\be{ik}{ij}(x)$ in some point $x \in \XAH$, we can associate a nest $C(x)$
with $x$. If no $\be{ik}{ij}(x)=0$, we think of $x$ as the image of a 
configuration where all sites are distinct, and in this case we associate 
the nest $\{\{1, \dots, n\}\}$ with $x$. And with this nest we associate
exactly one screen where all $n$ sites occur and are distinct. In the 
other case at least one $\be{ik}{ij}(x)=0$. As we can introduce 
$\be{ik}{ij}=\frac{|p_i-p_k|}{|p_i-p_j|}$ for distinct sites, we think of
$\be{ik}{ij}(x)=0$ as if $p_i$ is very close to $p_k$ as seen from $p_j$.
But this means exactly that $p_i$ and $p_k$ form a cluster with respect 
to $p_j$. In fact, one can prove that the set $\mathcal{C}(x)$ 
consisting basically  of all sets 
$C_{ij} = \{k ~|~ \be{ik}{ij} \neq 0 \} \cup \{i,j\}$ defines a nest
on $\{1, \dots, n\}$. So $\mathcal{C}(x)$ is the nest that corresponds 
to $x \in \XAH$.  For every cluster in $\mathcal{C}(x)$ with at least 
two elements we define a screen as a copy of the plane. 

\subsection{Factorizing $\XAH$.}

Fix some $x \in \XAH$. By analyzing the $C_{ij}$ as explained above,
we associate a family of screens with $x$ denoted by $\mathbf{x}$-{\bf screens}.
Our next goal is to fill the $x$-screens. The main idea of working 
with screens is as follows: the site set 
$p_1, \dots ,p_n$ that we think of as corresponding to some
degenerate configuration $q \in \XAH$, should be separated completely,
somewhere in the hierarchy of $x$-screens. This puts a requirement
on those points $q\in \XAH$ that we can use to fill the $x$-screens 
with: the point $q$ should not be more nested than $x$ itself, in the following
sense: if $\be{ik}{ij}(x) \neq 0,\infty$, then also 
$\be{ik}{ij}(q) \neq 0,\infty$. Another way to put this: we are not allowed to
move down in the stratification with respect to the coincidence of sites.
We do not need all $\binom{n}{2} + 6 \binom{n}{3}$ data elements 
$$ (\alpha_{ij}(q))_{1 \leq i< j \leq n}, \ho{ik}{ij}(q) ),
	\qquad i,j, k,~\text{pairwise distinct indices},$$
for filling the $x$-screens with sites $p_1(q), \dots, p_n(q)$,
given some suitable $q\in \xah$. We just need enough data elements
to pinpoint each site once. By exploiting the nested structure
of $x$, we can avoid bad choices in picking the data elements.
We do not want to be confronted with legs of length $0 \cdot \infty$ 
for example.
These prerequisites result in a factor $\dom(x)$ of data elements
of $\xah$ that are both suited and enough to fill the $x$-screens,
given some suitable $q \in \xah$. Note that the choice of the data 
elements depends on $x$, although we use of course the
data elements of $q$ itself in order to fill the $x$-screens
with sites $p_1(q), \dots, p_n(q)$. 

\subsection{The hooked tree.}

Again by using the nested structure of our fixed $x \in \xah$,
we can define a partial order on the labels $1, \dots, n$.
This partial order tells in which order exactly 
we have to add the $q$-sites $p_1(q), \dots, p_n(q)$ in the
$x$-screens. It turns out that with every label $i \in 1, \dots, n$
we can associate an $\mathbf{x}$-{\bf tag} $l_x(p_i)$. Such an $x$-tag is
a pointer to one coordinate of $\dom(x) \subset \XAH$. 
As any such coordinate is either an angle or a hook, we conclude that 
every label $i$ corresponds to a geometric construction
(construct $p_i$ as the end of a line segment of a given angle, 
or construct $p_i$ as the end point of a given hook) that tells
how to construct $p_i$. 

\begin{figure}[!ht]
\begin{center}
\setlength{\unitlength}{1em}
\begin{picture}(24,14)
\put(12,8){\makebox(0,0)[cc]{
        \leavevmode\epsfxsize=20\unitlength
\epsfbox{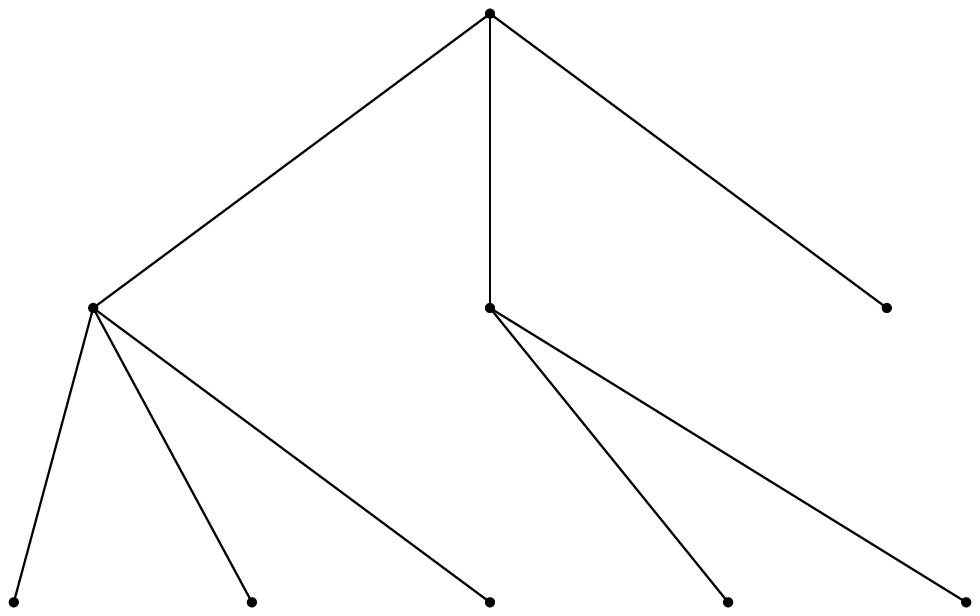}}}
\put(20,7.2){\makebox(0,0)[l]{4}}
\put(2.3,0.9){\makebox(0,0)[l]{1}}
\put(7.1,0.9){\makebox(0,0)[l]{2}}
\put(12.0,0.9){\makebox(0,0)[l]{6}}
\put(16.8,0.9){\makebox(0,0)[l]{3}}
\put(21.7,0.9){\makebox(0,0)[l]{5}}
\put(7.7,10){\makebox(0,0)[l]{top}}
\put(12.3, 10){\makebox(0,0)[l]{$\alpha_{1,3}$}}
\put(18.1, 10){\makebox(0,0)[l]{$h^{14}_{13}$}}
\put(7.2, 3.5){\makebox(0,0)[l]{$h^{12}_{13}$}}
\put(11.0, 3.5){\makebox(0,0)[l]{$h^{16}_{12}$}}
\put(20.5, 3.5){\makebox(0,0)[l]{$h^{35}_{31}$}}
\end{picture}
\caption{\elabel{ht1.eps} The hooked tree for the nest $<\{1,2,6\},\{3,5\}>$.}
\end{center}
\end{figure}

We store the set of $x$-clusters in a tree. The vertices of the tree 
are exactly the $x$-clusters. Two $x$-clusters are connected by
an edge if one cluster is maximal in the other. Some edges are
labeled by $x$-tags. This is done in such a way that if going up in 
the tree from leaf $i$, the first $x$-tag encountered points exactly
to that hook or angle that is used to construct $p_i$ in the
set of $x$-screens. This tree is called the {\bf hooked tree}.

\subsection{Filling screens.}

We discuss the use of the partial order a bit further. Our setup is such
that in each $x$-screen separation occurs. Fix some $x$-screen $S$.
As mentioned above, any $x$-screen $S$ corresponds to some $x$-cluster
$C_S$. Assume that all $x$-screens $T$, such that the cluster $C_T$
contains $C_S$, have already been filled. We say that such screens $T$
are {\bf above} $S$. The $x$-screen $S$ should exactly contain the
sites with labels in the $x$-cluster $C_S$ after filling. Moreover,
the sites corresponding to the maximal subclusters of $S$ are separated
in the $x$-screen $S$. Let $C_1$ be the maximal subcluster of $S$
with smallest minimal label $i$ and let $C_2$ be the maximal subcluster 
with second smallest minimal label. The first site we put in $S$
is the site $p_i$, that we place at the origin. The second site to 
construct is the site $p_j$, that we place at distance $1$ of $p_i$.
The angle that the line segment $p_ip_j$ makes with the positive axis
is the so-called {\bf screen orientation} $\mathcal{O}_S$ of $S$.
It is defined recursively in terms of the screens above $S$. 
Note that by now we have constructed one leg $p_ip_j$ in $S$. All other
sites in $S$ can be constructed using hooks.   

Suppose, for instance, that we want to construct all sites in some maximal 
subcluster $C_k \neq C_1,C_2$, where $k$ is the minimal label of $C_k$ and therefore
$k>i,j$. We construct $p_k$ by means of a hook $\ho{ik}{ij}$. One can think
of this hook as a {\bf hook on scale}, as $p_i,p_j$ and $p_k$ are separated
at the same level. Note that $k$ is automatically minimal in the maximal
subcluster of $C_k$ that has smallest minimal label.   Let $m \in C_k$, with $m>k$, 
be minimal in the second maximal subcluster of $C_k$.  We construct $p_m$ 
using $\ho{km}{ki}$. This is a different type of hook
as $p_k$ and $p_m$ are separated one step deeper in the nest structure. Think of
this type of hook as a {\bf explosion hook}: for $x$ itself $|\be{km}{ki}(x)|=0$,
as $p_k$ and $p_m$ `coincide with respect to $p_i$'. If we take an arbitrary
point $q \in \dom$, then $|\be{km}{ki}(q)| \geq 0$, and if $|\be{km}{ki}(q)|>0$
we could say that an explosion occurs, as suddenly $p_k$ and $p_m$ do
not coincide anymore with respect to~$p_i$. 

\subsection{$\XAH$ is locally the graph of a function.}

One purpose of this setup is to minimize the dimension of the 
factor $\dom(x)$, that is, we want to use as few
data elements of $\XAH$ for filling the $x$-screens as possible.
That we cannot do better
is shown by Lemma \ref{dimdom} that proves that the dimension of the factor 
$\dom(x)$ equals the dimension of the reduced configuration space $\conf$.
Given a set of filled $x$-screens, we can just read off angles and
hooks from the screens. Theorem \ref{thconsist}
proves that filling screens using data elements is consistent with reading off 
data elements, in the sense that we get back the data elements that we started with.
Note that besides the data elements in the $\dom(x)$ factor,
we can also read off the data elements in the complementary factor $\rng(x)$ of
$\dom(x)$ in the space of all angles and hooks on $n$ points, $\ah$.
This shows that $\XAH$ can be written as the graph of some
function $\readoff: \dom(x) \rightarrow \rng(x)$. As a consequence,
we can prove in Theorem \ref{tsubmani} that $\XAH$ is a  smooth manifold.
In fact there are some extra requirements on the structure presented
to obtain smoothness that we did not discuss here. Details can 
be found in Sections \ref{sah} to \ref{ssmooth}.

\subsection{Connection with Kontsevich-Soibelman.}

In \cite{KS}, see also Section \ref{sKS}, Kontsevich and Soibelman
describe a manifold with corners $\FMt(n)$ that is closely related
to our smooth manifold $\xah$. The exact relation is
given in Theorem \ref{tfiber} that describes a map
$f: \FMt(n) \rightarrow \xah$ together with the fibers of $f$.
Due to this close relation we can use the construction of the
filled $x$-screens for points $x \in \FMt(n)$ as well. This is
important in the following application. For an overview of the 
compactifications introduced so far, consult Table \ref{tspace}.

\begin{table}[!ht]
\begin{center}
\begin{tabular}{lclclcl}
& & angles  & & angles  \\
&& mod $2 \pi$ && mod $\pi$\\\\
clickable & $\longleftarrow$ &
$\FMt(n)$  & $\longrightarrow$ & $\xah$ & $\phantom{\longrightarrow}$
        & angles \\    
Voronoi &&&&&& and \\
diagrams &&&&&& hooks\\\\
&&$\Big\downarrow$ & & $\Big\downarrow$ \\ \\
Voronoi & $\longleftarrow$  &  $\cda$ & $\longrightarrow$ & $\cua$ &
         $\phantom{\longrightarrow}$ & angles  \\
diagrams &&&&&& and\\
&&&&&& points
\end{tabular}
\end{center}
\caption{\elabel{tspace}Overview of compactifications.}
\end{table}

\subsection{Adding Voronoi diagrams in the $x$-screens.}

We apply the construction of filled $x$-screens to the analysis
of Voronoi diagrams of points sets that may include
coinciding points. We proceed by associating a Voronoi diagram
to any point $x \in \FMt(n)$. First, consider some 
non-degenerate point $x \in \FMt(n)$. Then by definition, as in 
the $\xah$ case, $x$ can be expressed as image of some 
configuration $c \in \CONF$. As no $\be{ik}{ij}(x)=0$ for a
non-degenerate $x \in \FMt(n)$, there is exactly one $x$-screen
associated with $x$, the top screen, that contains all
sites $p_1(x), \dots, p_n(x)$ in such a way that all sites
are distinct. The Voronoi diagram $V_{FM}(x)$ associated to
$x$ in this case is the Voronoi diagram of the $n$ sites in the
top screen. 

Next assume that $x$ is a degenerate element of $\FMt(n)$.
Suppose that we have filled the $x$-screens with 
the $x$-sites $p_1(x), \dots, p_n(x)$ as described above.
In this case, there is more then one $x$-screen, but it is 
guaranteed that in any $x$-screen $S$ exactly those clusters
are separated that are maximal subclusters of the cluster 
$C_S$. We define a Voronoi diagram $V_{FM}(x)$ in terms of the 
$x$-filled $x$-screens in two steps. In the first step we
add the classic Voronoi diagram of the maximal subclusters
of $C_S$ to any $x$-screen $S$. This is called the
initialization step. In the so-called completion step
we recursively add the completed diagrams of all screens
below $S$ to $S$ for all $x$-screens $S$. This is done
as follows.  A non-trivial cluster is a cluster 
that consists of at least two sites. Any non-trivial
maximal cluster $M$ of $S$ corresponds to 
an $x$-screen $S_M$ that is directly below $S$. We paste the
completed Voronoi diagram $V_{FM}(S_M)$ of $S_M$ into
the Voronoi cell of $M$ in the initialized Voronoi
diagram $V(S)$.

This gives us a \bfindex{clickable} Voronoi diagram $V_{FM}(x)$ for
any $x \in \FMt$. By clicking on a non-trivial 
maximal cluster, if present, a screen appears. This
screen  contains the sites
in the cluster and the Voronoi diagram of the sites.
Recall the polynomial sites from Chapter \ref{chlimit}. 
By computing ratios and angles we can associate an element
$x_S \in \FMt$ to a set $S(t)$ of polynomial sites. 
Then the completed Voronoi diagram 
$V_{FM}(T)$ of the top screen $T$ for $x_S$ matches the 
limit Voronoi diagram of $S(t)$ at $t=0$.

\section{Angles and hooks.}

\elabel{sah}
We start this section by recalling and introducing data elements defined on pairs
and triples of distinct points in $\R^2$. These
data elements form the basis of the compactification of
a suitable quotient
of the configuration space of $n$ distinct points in the plane
that we define later on in this section.

\subsection{Data elements of pairs and triples of points.}

Given two distinct points $p_i$ and $p_j$, we can determine the angle
between  the line through $p_i$ and $p_j$ and the positive
$x$-axis. We have done this in Definition \ref{dangle}. There we have introduced
the directed angle $\alpha_{ij} \in \R/ 2 \pi \Z$ and the undirected angle
$\overline{\alpha_{ij}} \in \R/\pi \Z$. 
We store information on  the geometry of three distinct points in hooks, see 
Figure \ref{hook}.

\begin{figure}[!ht]
\begin{center}
\setlength{\unitlength}{1cm}
\begin{picture}(5,3.2)
\put(2.5,1.6){\makebox(0,0)[cc]{
        \leavevmode\epsfxsize=4\unitlength\epsfbox{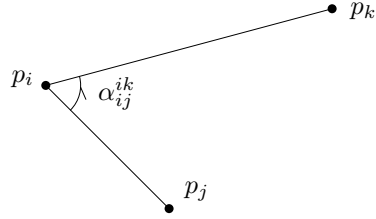}}}
\put(1.35,1.8){\makebox(0,0)[l]{$\al{ik}{ij}$}}
\put(0.2,2){\makebox(0,0)[l]{$p_i$}}
\put(2.5,0.5){\makebox(0,0)[l]{$p_j$}}
\put(4.7,2.9){\makebox(0,0)[l]{$p_k$}}
\end{picture}
 \caption{\elabel{hook}
The angle $\al{ik}{ij}$ from $p_j$ to $p_k$ hinged at $p_i$.
The small arrow indicates the positive direction.}
\end{center}
\end{figure}

\begin{definition}
Let $p_i$, $p_j$, and $p_k$ be three distinct points in the plane. 
The {\bf angle} $\al{ik}{ij}$ from $p_j$ to $p_k$ hinged at 
$p_i$ is given by
$\al{ik}{ij}  :=  \alpha_{ik} - \alpha_{ij} \mod 2 \pi$.
The {\bf ratio} \index{ratio} $\be{ik}{ij}$ from $p_j$ to $p_k$  hinged at $p_i$
is given by
\beq
         \be{ik}{ij} & := & \frac{| p_i - p_k | }{ | p_i - p_j |}, \qquad
        \be{ik}{ij} \in (0, \infty).
\eeq
The {\bf hook} \index{hook} $\ho{ik}{ij}$ from $p_j$ to $p_k$ hinged at $p_i$ 
equals the ratio together with the angle:
$h^{ik}_{ij}  :=  ( \be{ik}{ij}, \al{ik}{ij} )$.
The point $p_i$ is the {\bf hinge point} \index{hinge point}
 of the hook $\ho{ik}{ij}$
while the line segments $p_ip_j$ and $p_ip_k$ form its
{\bf legs} \index{leg}.
\end{definition}

Note that one triple of distinct points gives six hooks.

\begin{remark}
An alternative is to consider hooks in the complex plane. Let $p_i$, $p_j$,
and $p_k$ be distinct points in $\C$.  Note that 
$\alpha_{ij} = (p_i-p_j)/|p_i-p_j|$. It is easy to check that, 
with $I=\sqrt{-1}$, it holds that
\beq
	\frac{p_i-p_k}{p_i-p_j} & = & \be{ik}{ij} e^{I \al{ik}{ij}}.
\eeq
\end{remark}

\subsection{Geometric interpretation of the hook $\ho{ik}{ij}$. \elabel{secgeo}}

\begin{notation}
Let $p$ and $q$ be points in the plane and $\alpha \in [-\pi, \pi)$
an angle. 
By $\rot_{\alpha, p}(q)$ \index{rot@$\rot_{\alpha, p}(q)$}  denote
the image of $q$ under the  anti-clockwise rotation around $p$ by 
angle $\alpha$. 
We abbreviate  $\rot_{\alpha, (0,0)}$  by $\rot_{\alpha}$.
\end{notation}

Suppose we are given  two distinct points $p_i$ and $p_j$,
together with a  hook $\ho{ik}{ij}$.
The hook $\ho{ik}{ij}$ can be seen as 
a prescription that tells how to construct  the point $p_k$.
Using 
$\al{ik}{ij}$, we first construct 
$p_j'$, the image of the rotation $\rot_{\al{ik}{ij}, p_i}$, applied to
$p_j$:
\beq
	p_j' & = & \rot_{\al{ik}{ij}, p_i}(p_j).
\eeq
The half-line that starts in $p_i$  and passes through $p_j'$
contains $p_k$.
The ratio $\be{ik}{ij}$ fixes the distance from $p_k$ to $p_i$,  
thereby fixing $p_k$ itself:
\beq
	p_k & = & \be{ik}{ij} (p_j' - p_i) + p_i.
\eeq
{\bf Summary}: we interpret a hook $\ho{ik}{ij}$ as a point rotation 
followed by a vector multiplication. As a result the leg $p_ip_j$ is
transformed into the leg $p_ip_k$.

\subsection{Negative ratios.}

\begin{figure}[!ht]
\begin{center}
\setlength{\unitlength}{1em}
\begin{picture}(12.8,8)
\put(6.4,4){\makebox(0,0)[cc]{
        \leavevmode\epsfxsize=11.2\unitlength\epsfbox{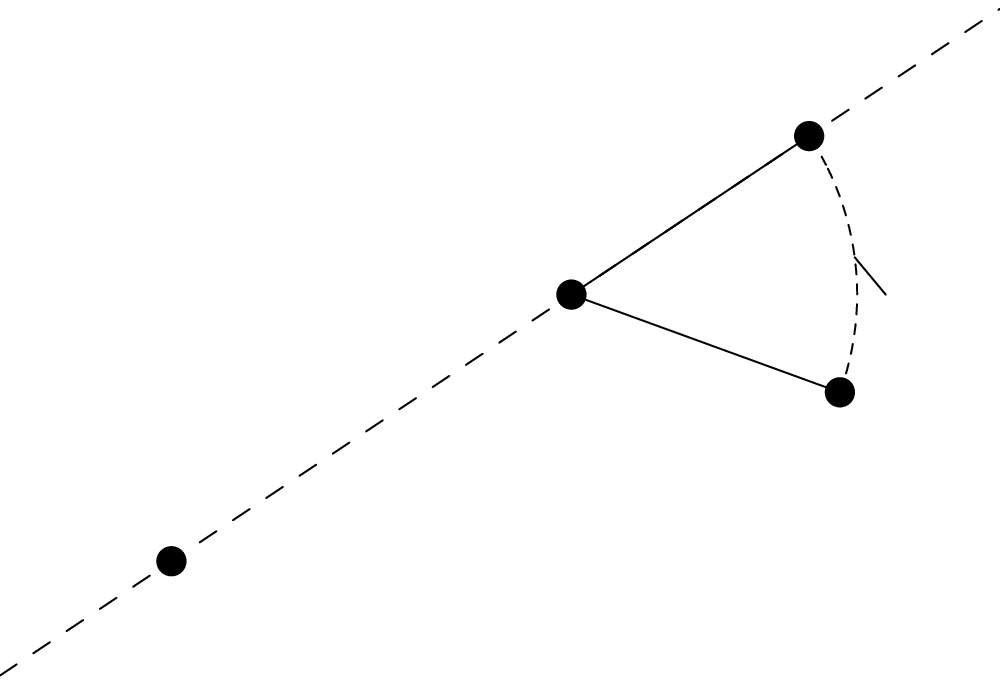}}}
\put(8.3,5.0){\makebox(0,0)[l]{$\al{ik}{ij}$}}
\put(9.6,7.7){\makebox(0,0)[l]{$p_j'$}}
\put(10.4,2.9){\makebox(0,0)[l]{$p_j$}}
\put(6.0,5.3){\makebox(0,0)[l]{$p_i$}}
\put(12.4,7.76){\makebox(0,0)[l]{$A_{ij'}$}}
\put(1.9,2.7){\makebox(0,0)[l]{$p_k$}}
\end{picture}
\caption{\elabel{rotbeta1} The local coordinate axis $A_{ij'}$.}
\end{center}
\end{figure}

In this section we give sense to `negative ratios', as illustrated by
Figure \ref{rotbeta1}.
Think of the line through the points 
$p_i$ and $p_j'$ as a local coordinate axis $A_{ij'}$.
The origin on this axis is $p_i$ and the point with coordinate~$1$
is $p_j'$.  The position of $p_k$ on this axis is given by the coordinate
$\be{ik}{ij}$. 
In this setting, any value of $\be{ik}{ij}$ in the interval
$(-\infty, \infty)$ is meaningful: it just indicates the point
with coordinate $\be{ik}{ij}$ on the $A_{ij'}$-axis.
Suppose we allow the full interval $(-\infty, \infty)$
for values of $\be{ik}{ij}$. Then there are two different ways to 
construct a point $p_k$ by a hook $\ho{ik}{ij}$, given the points
$p_i$ and $p_j$. The hooks   $(\be{ik}{ij}, \al{ik}{ij})$ and
 $(-\be{ik}{ij}, \al{ik}{ij}+\pi)$ both result in the same point $p_k$.
If identifying $-\infty$ with $\infty$, we have by now
described a map
\bmap
	\phi_{ijk}: &  \CONFt(\R^2)
		& \rightarrow & \P^1 \times \R/2 \pi \Z \\
	& (p_i, p_j, p_k) & \mapsto & ( \be{ik}{ij}, \al{ik}{ij} ),
\emap
under the identification
\begin{eqnarray}
\elabel{kleinid}
(\be{ik}{ij}, \al{ik}{ij}) & \sim_k & (-\be{ik}{ij}, \al{ik}{ij}+\pi).
\end{eqnarray}
Recall that $\CONFt(\R^2)$ denotes the configuration space of three
distinct points in the plane. $\P^1$ denotes the \bfindex{projective line}. 
It can be defined
as $S^1$ with antipodal points identified.  A point on the
projective line can be seen as the coordinate $q$
of a third point on an axis with respect to two fixed points $0 \neq1$.
A homeomorphism from the circle to the line $\R \cup \{\infty\}$ given
by \bfindex{stereographic projection} is shown in Figure 
\ref{figstereo}. A line through $P$
intersects the circle in another point, $Q$,  and hits the $x$-axis in a 
point $R$. When $P$ and $Q$ coincide, this line is horizontal,
so $P$ is mapped to $R=\infty$, cf. also \cite{R}.

\begin{figure}[!ht]
\begin{center}
\setlength{\unitlength}{1cm}
\begin{picture}(4,2.5)
\put(2,1){\makebox(0,0)[cc]{
        \leavevmode\epsfxsize=4\unitlength\epsfbox{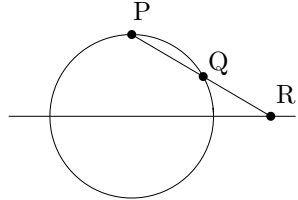}}}
\put(1.8,2.4){\makebox(0,0)[l]{P}}
\put(2.8,1.7){\makebox(0,0)[l]{Q}}
\put(3.7,1.3){\makebox(0,0)[l]{R}}
\end{picture}
\caption{\elabel{figstereo} Stereographic projection yields
a homeomorphism between the circle $S^1$ and the projective line
$\R\P^1=\R \cup \infty$.}
\end{center}
\end{figure}
\subsection{Klein bottle $\K{ik}{ij}$.}

\begin{figure}[!ht]
\begin{center}
\setlength{\unitlength}{1cm}
\begin{picture}(4,2)
\put(2,1){\makebox(0,0)[cc]{
        \leavevmode\epsfxsize=4\unitlength\epsfbox{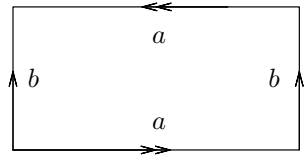}}}
\put(2,1.55){\makebox(0,0)[l]{$a$}}
\put(2,.4){\makebox(0,0)[l]{$a$}}
\put(.35,1){\makebox(0,0)[l]{$b$}}
\put(3.55,1){\makebox(0,0)[l]{$b$}}
\end{picture}
\caption{\elabel{klein2} The Klein bottle is obtained  by identifying sides
$a$ and $b$ according to the arrows.}
\end{center}
\end{figure}

The \bfindex{Klein bottle} is the quotient space obtained
from a rectangle by identifying opposite sides, see Figure \ref{klein2}; 
see also \cite{Mu}. In Figure \ref{kstrip}
we have depicted another rectangle, but this time we add
a description of the axes. For the horizontal axis we take
the projective line $\P^1$.  The vertical axis is $\R/2 \pi \Z$.
Recall the identification 
 $(\be{ik}{ij}, \al{ik}{ij}) ~\sim_k~
        (-\be{ik}{ij}, \al{ik}{ij}+\pi)$  introduced in  \ref{kleinid}. 

\begin{lemma} 
\elabel{lklein}
The quotient space $ \P^1 \times \R/ 2 \pi \Z  ~\slash~ \sim_k $
defines a Klein bottle.
\end{lemma}

\begin{proof}
Consider Figure \ref{kstrip}. Any class in $\P^1 \times \R / 2 \pi \Z / \sim_k$
has a representative in $[0,\infty] \times [0, 2 \pi]$. On
$[0,\infty] \times [0, 2 \pi]$, there are three identifications.

\begin{lijst}
\item The periodicity of directed angles
in $\R / 2 \pi \Z$ is indicated by $\gg$. 
\item The identification 
$(0,\al{ik}{ij}) \sim_k (0, \al{ik}{ij} + \pi)$ is indicated by $\triangle$.  
\item If $(0,\al{ik}{ij}) \sim_k (0, \al{ik}{ij} + \pi)$ holds on the
compact set $\P^1 \times \R / 2 \pi \Z / \sim_k$, then also
$(\infty,\al{ik}{ij}) \sim_k (\infty, \al{ik}{ij} + \pi)$, as 
any $\be{ik}{ij}=0$ corresponds to $\be{ij}{ik}=\infty$.
This gives the third identification, indicated in Figure \ref{kstrip}
by $\wedge$. 
\end{lijst}
For identifying the quotient space we are allowed
to cut along $\al{ik}{ij}=\pi$, if we eventually paste back again
along the same cut.
This cut is indicated by $\propto$. We do cut and paste the bottom
half on the top half, applying 
$\triangle$ and $\wedge$. This results in the cylinder
shown on the right of Figure \ref{kstrip}. On the top of the cylinder,
$\gg$ and $\propto$ are directed clockwise, while on the 
bottom they are directed counterclockwise. This shows
that we are in the situation of Figure \ref{klein2}. 
\end{proof}

\begin{figure}[!ht]
\begin{center}
\setlength{\unitlength}{1cm}
\begin{picture}(12,4)
\put(6,2){\makebox(0,0)[cc]{
        \leavevmode\epsfxsize=11\unitlength\epsfbox{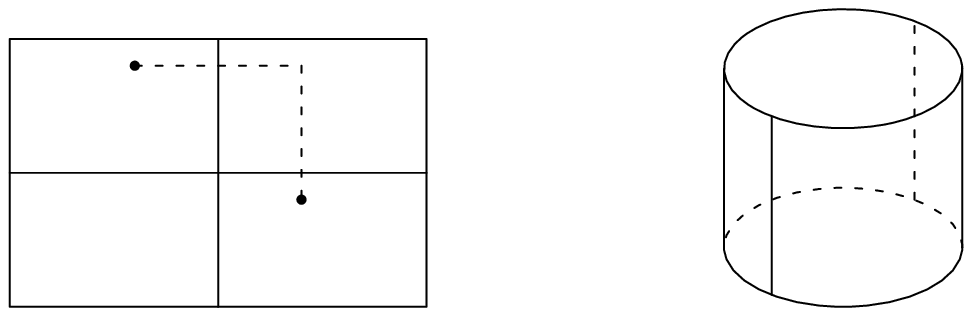}}}
\put(.75,3.8){\makebox(0,0)[l]{$\uparrow\al{ik}{ij}$}}
\put(5.6,.4){\makebox(0,0)[l]{$\rightarrow\be{ik}{ij}$}}
\put(0.1,0.6){\makebox(0,0)[l]{0}}
\put(0.1,1.9){\makebox(0,0)[l]{$\pi$}}
\put(0.1,3.35 ){\makebox(0,0)[l]{$2 \pi$}}
\put(.4,0){\makebox(0,0)[l]{$-\infty$}}
\put(3,0){\makebox(0,0)[l]{0}}
\put(5.2,0){\makebox(0,0)[l]{$\infty$}}
\put(5.29,1.1){\makebox(0,0)[l]{$\wedge$}}
\put(5.29,2.6){\makebox(0,0)[l]{$\wedge$}}
\put(2.96,1.1){\makebox(0,0)[l]{$\triangle$}}
\put(2.96,2.6){\makebox(0,0)[l]{$\triangle$}}
\put(4.2,.38){\makebox(0,0)[l]{$\gg$}}
\put(4.2,3.32){\makebox(0,0)[l]{$\gg$}}
\put(4.2,1.84){\makebox(0,0)[l]{$\propto$}}
\put(10.66,2.9){\makebox(0,0)[l]{$\wedge$}}
\put(9.05,2.0){\makebox(0,0)[l]{$\triangle$}}
\put(9.75,3.65){\makebox(0,0)[l]{$\gg$}}
\put(9.99,0.39){\makebox(0,0)[l]{$\gg$}}
\put(10.0,2.34){\makebox(0,0)[l]{$\propto$}}
\put(10.0,1.67){\makebox(0,0)[l]{$\propto$}}
\end{picture}
 \caption{\elabel{kstrip} The identification
        $(\be{ik}{ij}, \al{ik}{ij}) ~\sim_k~
        (-\be{ik}{ij}, \al{ik}{ij}+\pi)$. }
\end{center}
\end{figure}

What is essential for us is that the Klein bottle is a smooth manifold.
Note that with every ordered triple $(i, j, k)$ we have associated a Klein 
bottle $\K{ik}{ij}$ as the target space of $\phi_{ijk}$.
We can write  $\phi_{ijk}$ as a composition: 
$\phi_{ijk}= c_{ijk} \circ h_{ijk}$.
The map $h_{ijk}$ maps an ordered triple of distinct
points to its hook $(\be{ik}{ij}, \al{ik}{ij})$, where $\be{ik}{ij}>0$. 
The second map, $c_{ijk}$, maps a hook $(\be{ik}{ij}, \al{ik}{ij})$ to its
class $[(\be{ik}{ij}, \al{ik}{ij})]$ under the equivalence relation $\sim_k$.
We have:
\begin{displaymath}
\begin{array}{clcccc}
\phi_{ijk}: & (\R^2)^3 \backslash \Delta_{ijk} & \stackrel{h_{ijk}}{\rightarrow} &
	[0,\infty] \times \R / 2 \pi \Z  & \stackrel{c_{ijk}}{\rightarrow} &
	( [-\infty,\infty] \times \R / 2 \pi \Z ) / \sim_k, \\
	& (p_i,p_j,p_k) & \mapsto & (\be{ik}{ij}, \al{ik}{ij}) & \mapsto &
	[(\be{ik}{ij}, \al{ik}{ij})].
\end{array}
\end{displaymath}

\begin{proposition}
The map 
\beq
c_{ijk}^{-1}:~  ( [-\infty,\infty] \times \R / 2 \pi \Z ) / \sim_k
	&\rightarrow&  [0,\infty] \times \R / 2 \pi \Z,
\eeq
is $2$ to $1$ if $\be{ik}{ij}$ equals $0$ or $\pm \infty$ and $1$ to $1$ else.
\end{proposition}

\begin{proof}
 The inverse images of $c_{ijk}$ are as follows:
\begin{lijst}
\item $c_{ijk}^{-1}(~[(0, \al{ik}{ij})]~) = \{(0, \al{ik}{ij}),
	(0, \al{ik}{ij} + \pi)\}$;
\item $c_{ijk}^{-1}(~[(\infty, \al{ik}{ij})]~) = \{(\infty, \al{ik}{ij}),
        (\infty, \al{ik}{ij} + \pi)\}$;
\item $c_{ijk}^{-1}(~[(-\infty, \al{ik}{ij})]~) = \{(\infty, \al{ik}{ij}),
        (\infty, \al{ik}{ij} + \pi)\}$;
\item $c_{ijk}^{-1}(~[(\be{ik}{ij}, \al{ik}{ij})]~) = 
	(\be{ik}{ij}, \al{ik}{ij})$, if $0 < \be{ik}{ij} < \infty$;
\item $c_{ijk}^{-1}(~[(\be{ik}{ij}, \al{ik}{ij})]~) = 
	(-\be{ik}{ij}, \al{ik}{ij} + \pi)$, if $-\infty < \be{ik}{ij} < 0$.
\qedhere\end{lijst}
\end{proof}

\begin{remark}
Consider the image on the right in Figure \ref{kstrip}.
If we cut along
$\be{ik}{ij}=0$, that is, along $\triangle$,
we obtain a M\"obius strip with $\be{ik}{ij}=0$ on the boundary.
Cutting along $\be{ik}{ij} = \infty$, that is,
along $\wedge$,  leaves us
with two pieces joined along $\propto$. This results in
a cylinder with $\be{ik}{ij}=0$ on one boundary component and
$\be{ik}{ij}=\infty$ on the
other component.
\end{remark}

\subsection{The data map from $\conf$ to $\ah$.}

In this section we describe a map from a set of distinct points
in the plane to all angles and hooks between pairs and triples
of those  points.
Recall that in Definition \ref{dconf} we have introduced the reduced
configuration space $\conf$ as the space of $n$ distinct points in the 
plane up to scaling and translations.
For a representative $c$ of a
class $[c] \in \conf$, write the following data: 
for every pair of points $p_i, p_j \in c, i\neq j$,
 the undirected angle $\overline{\alpha}_{ij}=\overline{\alpha}_{ji}$;
for every triple $p_i, p_j, p_k$ of points in $c$, the 
six hooks $\ho{ik}{ij}$, $\ho{ij}{ik}$, $\ho{jk}{ji}$, $\ho{ji}{jk}$, 
$\ho{ki}{kj}$, and $\ho{kj}{ki}$.
This gives $\binom{n}{2}$ unordered angles and $6 \binom{n}{3}$
hooks. Note that the angles and hooks that we obtain are independent of
the choice of a representative of a class $c \in \conf$.

\begin{definition}
\elabel{defah}
$\ah$ \index{AH@$\ah$}
 is the space of hooks and angles on $n$  points:
\beq
\ah & := &  (\R/ \pi \Z)^{\binom{n}{2}} \quad  \times \quad
	 (~\P^1 ~\times~ (\R/  2 \pi \Z)~/\sim_k)^{ 6 \binom{n}{3}},
\eeq 
where $\sim_k$ denotes the identification defined in \ref{kleinid}.
\end{definition}

\begin{remark}
Being a direct product of circles and Klein bottles, $\ah$ is smooth.
\end{remark}

In the next definition, we introduce a 
\bfindex{compactification} of the reduced configuration space
$\conf$.

\begin{definition}
\elabel{defdatamap}
The \bfindex{data map} $\psi_{\ah}$ is the map
\begin{displaymath}
\begin{array}{rccl}
 \psi_{\ah}:& \conf & \rightarrow & \ah, \\
 \psi_{\ah}:& \conf & \rightarrow & 
	(\R/ \pi \Z)^{\binom{n}{2}} \quad  \times \quad  
		(~\P^1 ~\times~ (\R/  2 \pi \Z)/\sim_k~)^{ 6 \binom{n}{3}}, \\
 & {}[ (p_1, \dots, p_n ) ] & \mapsto &  ( (\overline{\alpha_{ij}}) 
_{1 \leq i < j \leq n}
	,\,(\be{ik}{ij}, \al{ik}{ij}) \,\mid
 i,j,k~ \text{pairwise distinct} ),
\end{array}
\end{displaymath}
and $\XAH$ \index{XAH@$\XAH$} is the closure of the image of $\conf$ in $\ah$.
\end{definition}

Alternatively we could allow non-negative ratios only,
 construct a data map from
$\conf$ to $
 (\R/ \pi \Z)^{\binom{n}{2}} \times 
         (~[0,\infty] ~\times~ (\R/  2 \pi \Z)~)^{ 6 \binom{n}{3}}$,
and define $\XAH$ as the closure of $\conf$ in this product.
In this approach, however, `jumps' occur in every hook $\ho{ik}{ij}$.
Indeed, suppose that $p_j \neq p_k$ are fixed while
$p_i$ moves `through' $p_k$. Then the ratio $\be{ik}{ij}$ stays close to zero,
but $\al{ik}{ij}$ jumps to $\al{ik}{ij}+\pi$ at collision, as $\alpha_{ik}$
changes direction, while $\alpha_{ij}$ stays the same. 

\subsection{Example: a smooth compactification of $\confd$.}

Consider three points $p_1$, $p_2$, and $p_3$ in the plane. Their configuration
space is given by
$\confd  =  \{ (p_1, p_2, p_3) ~|~  
p_1 \neq p_2, p_1  \neq p_3, p_2 \neq p_3 \}.$ Then $\confd$ is mapped
by $\psi_{\ah}$  to a 
rather big space:
\beq
	\psi_{\ah}:~ \confd & \rightarrow & 
		S_{12} \times S_{13} \times S_{23} \times
			\K{13}{12} \times  \K{12}{13} \times \K{23}{21}
				   \times \K{21}{23} \times
		\K{32}{31} \times \K{31}{32}.
\eeq

Indeed, we can make a smooth compactification of $\confd$ by mapping it
to a much smaller space: 

\begin{lemma}
Let $\psi$ be the map
	$\psi:~ \confd  \rightarrow S_{12} \times \K{13}{12}$
given by $(p_1, p_2, p_3) \mapsto (\alpha_{12}, \be{13}{12})$.
Then $\overline{\psi[{\mbox {\sl {\bf conf}}_3}]} =  S_{12} \times \K{13}{12}$.
Moreover, $S_{12} \times \K{13}{12}$ is smooth.
\end{lemma}

\begin{proof}
We have to show that for every point $y \in S_{12} \times \K{13}{12}$
and every neighbourhood $U$ of $y$ it holds that 
$U \cap \psi[{\mbox {\sl {\bf conf}}_3}] \neq \emptyset$. 
Any point $y \in  S_{12} \times \K{13}{12}$ has coordinates:
\beq
y & = &  (\alpha_{12}(y),\be{13}{12}(y),\al{13}{12}(y)).
\eeq
If 
$0 < |\be{13}{12}(y)| < \infty$, then 
	$(\alpha_{12}(y), \be{13}{12}(y), \al{13}{12}(y))  =  \psi(p_1,p_2,p_3)$,
where $p_1 = (0,0)$, $p_2 = ( \cos \alpha_{12}, \sin \alpha_{12})$, and
$p_3 = \be{13}{12} ( \cos (\alpha_{12} + \al{13}{12}), 
		\sin (\alpha_{12} + \al{13}{12})).$
If ratio $\be{13}{12}(y)=0, \pm \infty$, then in any neighbourhood of 
$y$ there exists a point $z$ such that 
$0 < |\be{13}{12}(z)| < \infty$, so $z \in \psi[\text{{\sl {\bf conf}}}_3]$.
\end{proof}

Later on we prove that ${\mbox {\sl XAH}[3]}$ is another smooth 
compactification  of ${\mbox {\sl {\bf conf}}_3}$. 
Mapping of $\confd$ into
$\ah$ has the advantage that there are no preferred labels.

\section{Nests and screens.\elabel{nestsscreens}}
The main result of this section is Theorem \ref{mathcalC}
that associates a nested set of subsets of $\{1,\dots,n\}$
to an element $x \in \xah$. This is done by
analyzing those ratios $\be{ik}{ij}(x)$ that are equal to zero.

\subsection{Clusters, nests, and screens.}

\begin{definition}
Let $Z$ be a collection of sets.
\begin{lijst}
\item $Z$ is \bfindex{nested} iff either $A \subset B$ or $B \subset A$ or
$A \cap B = \emptyset$ for any two elements $A,B \in Z$. 
\item A set $U$ is a {\bf maximal subset} of a set $W$, w.r.t.\ some nested
set $Z$ iff $U,W \in Z$, $U \subsetneq W$ and there exists no $V \in Z$ with
$U \subsetneq V \subsetneq W$.
\item A {\bf nest} \index{nest} $Z$ on $\{1\dots n\}$ is a 
nested set of subsets of 
$\{1, \dots, n\}$ that includes
the set $\{1, \dots, n\}$ itself and all singleton sets $\{i\}$, where 
$i=1, \dots, n$.
\item  The elements of a nest are called {\bf clusters}. \index{cluster}
 The singleton sets
       are called the {\bf trivial clusters}. \index{trivial cluster}
\end{lijst}
\end{definition}

\begin{eexample}
\elabel{exZ}
The set $Z  =  \{ \{1,2,3,4,5,6\}, \{1,2,6\}, \{3,5\}, \{1\}, \{2\},
                \{3\}, \{4\}, \{5\}, \{6\} \}$
is a nest on $\{1\dots 6\}$. For simplicity, we often write only the 
non-obvious elements of $Z$. The obvious elements are
the singleton sets and the total set $\{1, \dots, n\}$. 
We indicate this by using the $\langle,\rangle$
brackets. In our example we would write
$Z  =  \langle \{1,2,6\}, \{3,5\} \rangle$.
\end{eexample}

Let $Z(n)$ be a nest.
Associate a copy $S_C$ of the plane to every non-trivial cluster 
$C \in Z(n)$. Every copy is labeled by its corresponding 
cluster. 
In this way we get a family of screens
$\mathcal{S}_Z  = \{ S_C ~|~ C \in Z, ~C ~\text{non-trivial} \}$,
associated to the nest $Z$.
In case there is only one non-trivial cluster, we
have exactly one screen, $S_{[n]}$, where $[n]$ stands for
$\{1, \dots, n \}$. \index{n@$[n]$}
In Example \ref{exZ} we have screens $S_{[6]}$, $S_{1,2,6}$,, and $S_{3,5}$.

In this section we show how to associate a nest $Z$ to a point
$x \in XAH[n]$. In the following sections we describe how to fill
the family of screens $\mathcal{S}_Z$  corresponding to a nest 
$Z$ in a meaningful way.

\subsection{Properties of ratios.}

\begin{lemma}
\elabel{betaconf}
Let $p_i$, $p_j$, $p_k$, and $p_l$ be distinct points in the plane.
Then the following relations hold:
\begin{lijst}
\item $\be{ik}{ij} \in (0,\infty)$;
\item $\be{ij}{ik} = 1 / \be{ik}{ij}$;
\item $\be{ik}{ij} + \be{jk}{ji} \geq 1$;
\item $\be{ij}{il}\cdot \be{ik}{ij} = \be{ik}{il}$;
\item $\be{ik}{ij}\cdot\be{kl}{ki} = \be{il}{ij}\cdot\be{lk}{li}$.
\end{lijst}
\end{lemma}

\begin{proof}{\tiny .}
\begin{lijst}
\item Follows directly from the definition.
\item $1/\be{ik}{ij} ~=~ \dfrac{|p_i-p_k|}{|p_i-p_j|} ~=~ \be{ij}{ik}$.
\item This is in fact the triangle inequality:
$$ \be{ik}{ij} + \be{jk}{ji} ~ = ~ \dfrac{|p_i - p_k|}{|p_i - p_j|}
		+ \dfrac{|p_j - p_k|}{|p_i - p_j|} 
	~ = ~ \dfrac{| p_i - p_k| + | p_j - p_k|}{|p_i - p_j |} 
        ~ \geq ~ 1.
$$
\item $ \be{ij}{il}\cdot \be{ik}{ij} ~ = ~ \dfrac{|p_i-p_j|}{|p_i-p_l|} 
					\cdot \dfrac{|p_i-p_k|}{|p_i-p_j|} 
                            ~ = ~ \dfrac{|p_i-p_k|}{|p_i-p_l|} 
                            ~ = ~ \be{ik}{il} $. 
\item Apply the definition.
\qedhere
\end{lijst}
\end{proof}

\begin{corollary}
\elabel{propbeta}
For  $x \in XAH[n]$  and $i$, $j$, $k$, $l$ distinct labels we have
\begin{lijst}
\item $\abe{ik}{ij} \in [0,\infty]$;
\item $\abe{ij}{ik} = 1 / \abe{ik}{ij}$;
\item $\abe{ik}{ij} + \abe{jk}{ji} \geq 1$;
\item If $\max(~\abe{ij}{il},~\abe{ik}{ij} ) < \infty$, then  
	$\abe{ij}{il}\cdot \abe{ik}{ij} = \abe{ik}{il}$.
\end{lijst}
\end{corollary}

\begin{proof}
The claims follow from Lemma \ref{betaconf} by taking limits, but notice
that for $x \in \XAH$ some $\be{ik}{ij}$ can be negative, see Definition
\ref{defah}. 
\end{proof}

\subsection{Separating clusters.}

\begin{remark}
For the rest of this section, we identify a site $p_i$ with its label $i$.
\end{remark}

Set
\begin{equation}
\elabel{eqcij}
	C_{ij} ~:=~  \{ p_k ~|~ \abe{ij}{ik} > 0 \} \cup \{i,j\}. 
\end{equation}

Intuitively, think of the following. Recall the geometric
interpretation of a hook $\ho{ij}{ik}$ given in Section \ref{secgeo}.
Suppose we use $\ho{ik}{ij}$ to draw a point $p_j$ in the plane, given 
the points $p_i$ and $p_k$. The ratio  $\abe{ij}{ik}$ gives the distance 
between the points $p_i$ and $p_j$, compared to the unit distance $p_ip_k$.
Imagine an observer situated at $p_k$.  When $\abe{ij}{ik}$ gets very small,
$p_j$ moves very close to $p_i$. But as long as $\abe{ij}{ik} > 0$, the
observer at $p_k$ still can distinguish $p_i$ and $p_j$.
Therefore one could say in this setting that $C_{ij}$ consists of all 
points $p_k$ such that $p_i \neq p_j$ w.r.t.\ $p_k$.  We call 
$C_{ij}$ the \bfindex{separating cluster} of $p_i$ and $p_j$.

\begin{figure}[ht]
\begin{center}
\setlength{\unitlength}{1cm}
\begin{picture}(12,3)
\put(6,1.5){\makebox(0,0)[cc]{
        \leavevmode\epsfxsize=3\unitlength
\epsfbox{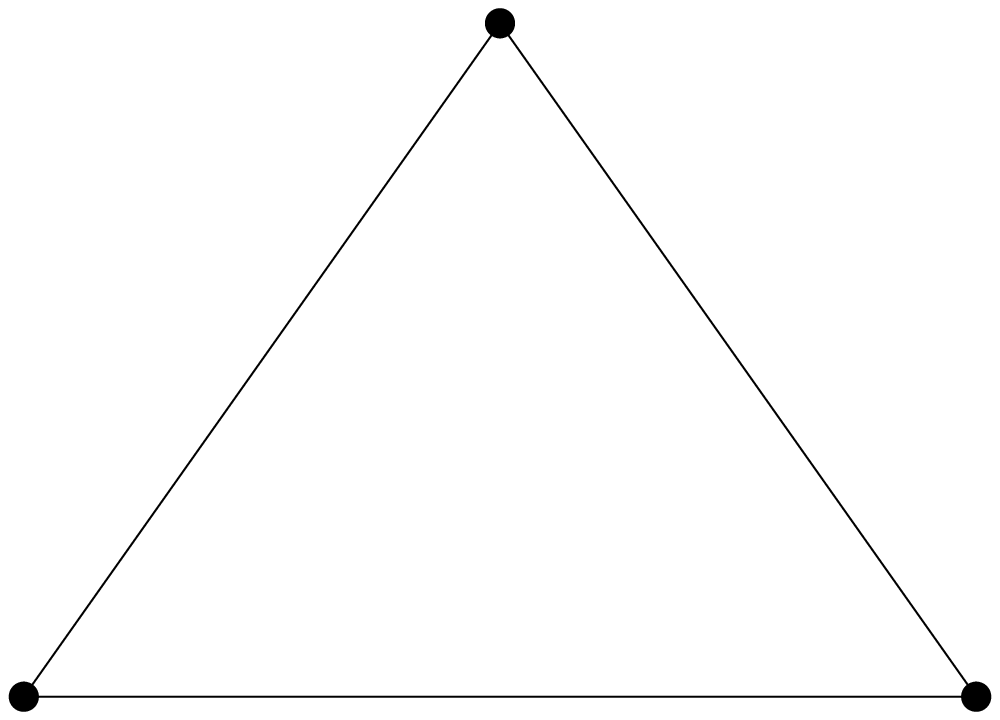}}}
\put(4.3,0.2){\makebox(0,0)[l]{$p_i$}}
\put(7.5,0.2){\makebox(0,0)[l]{$p_k$}}
\put(6,2.8){\makebox(0,0)[l]{$p_j$}}
\put(6.9,1.6){\makebox(0,0)[l]{$I$}}
\put(6,.2){\makebox(0,0)[l]{$J$}}
\put(4.9,1.6){\makebox(0,0)[l]{$K$}}
\end{picture}
\caption{\elabel{triangle} The triangle $p_ip_jp_k$.}
\end{center}
\end{figure}

\begin{lemma}
We have $C_{ij} = C_{ji}$.
\end{lemma}

\begin{proof}
Let us  start with a heuristic proof: consider Figure \ref{triangle}.
Note that $I = |p_k-p_j|$, $J=|p_k-p_i|$ and
$K=|p_j-p_i|$. Therefore $\abe{ij}{ik} = K/J$ and 
	$\abe{ji}{jk}=K/I$. We have to show
that $K/J > 0$ implies that $K/I > 0$ as well. $K/J > 0$ only excludes that $p_i$ and
$p_j$ coincide w.r.t.\ $p_k$. But this is also the only way to make $K/I$ small. \\
More formally: 
we have to prove that $\abe{ij}{ik}>0$ implies that $\abe{ji}{jk}>0$.
First of all: as $\abe{ij}{ik} > 0$ it follows that $\abe{ik}{ij} < \infty$. Now suppose
that $\abe{ji}{jk}=0$. Then also $\abe{ji}{jk}\cdot\abe{ik}{ij} = \abe{ki}{kj} = 0$.
This contradicts the triangle inequality applied to
$\abe{ki}{kj}$ and $\abe{ji}{jk}$. 
\end{proof}

\begin{lemma}
\elabel{pkpl}
If $p_k \in C_{ij}$ and $p_l \not\in C_{ij}$, then $p_l \not\in C_{ik}, C_{jk}$.
\end{lemma}

\begin{proof}
First, $p_k \in C_{ij}$ implies that $\abe{ik}{ij} < \infty$. Secondly, 
as $p_l \not\in C_{ij}$ it follows that $\abe{ij}{il} = 0$.
So  both $\abe{ij}{il}$ and $\abe{ik}{ij}$  are finite, which means that we
can apply Corollary \ref{propbeta} and get
$\abe{ik}{il} ~=~ \abe{ij}{il} \abe{ik}{ij} ~=~ 0$.  
\end{proof}

\begin{lemma}
\elabel{pkCij}
Let $i,j$ and $k$ be distinct labels. Then
\beq
	p_k \in C_{ij} & \Rightarrow & C_{ik} \subset C_{ij},\\
        p_k \not\in C_{ij} & \Rightarrow & C_{ij} \subset C_{ik}.
\eeq
\end{lemma}

\begin{proof} 
If $p_k \in C_{ij}$, then 
$C_{ik} \subset C_{ij}$ by Lemma \ref{pkpl}.  
We are left with the case $p_k \not\in C_{ij}$. We 
prove that in this case $C_{ij} \subset C_{ik}$. For this purpose,
we fix $p_l \in C_{ij}$ and show that $p_l \in C_{ik}$.
Since $p_k \not\in C_{ij}$, it follows that $\abe{ij}{ik} =0$. 
Due to $p_l \in C_{ij}$, one has $|\be{ij}{il}|>0$, and in 
particular $|\be{il}{ij}|$ is finite.
By Corollary \ref{propbeta}, 
$\abe{il}{ik} ~=~ \abe{ij}{ik}\cdot \abe{il}{ij} ~=~ 0.$
As a consequence $\abe{ik}{il} > 0$,
so indeed $p_l \in C_{ik}$.
\end{proof}

\begin{lemma}
\elabel{papb}
If $p_a, p_b \in C_{ij}$, then $C_{ab} \subset C_{ij}$.
\end{lemma}

\begin{proof}
If $p_b \in C_{ia}$, it follows from Lemma \ref{pkCij}
      that $C_{ab} \subset C_{ib} \subset C_{ij}$.
Suppose that $p_b \not\in C_{ia}$. Then Lemma \ref{pkCij} gives
$C_{ia} \subset C_{ib}$. This implies that $p_a \in C_{ib}$. Therefore,
$C_{ab} \subset C_{ib} \subset C_{ij}$.
\end{proof}

\subsection{The clusters form a nest.}

\begin{theorem}
\elabel{mathcalC} Fix $x \in \xah$.
Let $ \mathcal{C}(x)  :=  \{C_{ij}\}_{1 \leq i < j \leq n} ~\cup~
\{\{1\},\dots,\{n\}\}$.  Then $\mathcal{C}(x)$ is a nest on $n$.
\end{theorem}

\begin{proof}
We have to show that for any clusters $C_{st}$ and $C_{ij}$ either 
$C_{st} \cap C_{ij} = \emptyset$ or one of them is contained in the other.
Lemma \ref{pkCij} deals with the case $s=i$: it shows that the nest
condition is fulfilled for $C_{it}$ and $C_{ij}$.
When $C_{st} \cap C_{ij} = \emptyset$ the nest condition is 
fulfilled as well.
So, assume that $s$, $t$, $i$, and $j$ are all distinct and that
$p_m  \in  C_{st} \cap C_{ij}$.
We have proved in Lemma \ref{pkCij} that the following inclusions hold:
\begin{displaymath}
\begin{array}{ccc}
	C_{im} & \star & C_{sm} \\
	\cap   &                & \cap \\
        C_{ij} &                & C_{st}\\
        \cup   &                & \cup \\
        C_{jm} & \ast & C_{tm}
\end{array}
\end{displaymath}
The inclusions $\star$ and $\ast$ are still open. 
There are four possibilities:
\begin{center}
\begin{tabular}{ccccc}
 & 1 & 2 & 3 & 4 \\
$\star$ & $\subset$ & $\supset$ & $\subset$ & $\supset$ \\
$\ast$ & $\subset$ & $\supset$ & $\supset$ & $\subset$
\end{tabular}
\end{center}
In the first case, $p_i$ and $p_j$ are in $C_{st}$. So, 
$C_{ij} \subset C_{st}$ by Lemma \ref{papb}. In the second case,
$C_{st} \subset C_{ij}$. Concentrate on the case 3 now.
\beq
	C_{im} \subset C_{sm} & \Rightarrow & p_i \in C_{st},\\
	C_{jm} \supset C_{tm} & \Rightarrow & p_t \in C_{ij}.
\eeq
If $p_j \in C_{sm}$, then $p_j \in C_{st}$, so both $p_i$ and  $p_j$  are in $C_{st}$.
This  implies that $C_{ij} \subset C_{st}$.
When  $p_j \not\in C_{sm}$, we can apply Lemma \ref{pkpl} as $p_i \in C_{sm}$
and conclude that  $C_{is} \subset C_{ij}$. But this implies that 
$p_s \in C_{ij}$. 
As $p_t \in C_{ij}$ as well, it follows that $C_{st} \subset C_{ij}$. 
The fourth case follows from the third by relabeling.
\end{proof}

From Theorem \ref{mathcalC} we obtain an important corollary.

\begin{corollary}
If $x \in XAH[n]$, then $\mathcal{C}=\mathcal{C}(x)$ forms a nest.
\end{corollary}

For $x \in \XAH$, we call a cluster $C \in \mathcal{C}(x)$ an {\bf $x$-cluster}.
\index{x-cluster@$x$-cluster}
By now, the following statements are obvious.

\begin{corollary}
Let $\mathcal{C}$ be a nest. Then
\elabel{difij}
\begin{lijst}
\item $C_{ij}$ is the smallest cluster of $\mathcal{C}$ that contains both $i$ and $j$;
\item $i$ and $j$ belong to distinct clusters in the screen $S_{C_{ij}}$.
\end{lijst}
\end{corollary}

\begin{definition} 
\elabel{defsepa}
Fix $x \in \xah$.
\begin{lijst}
\item Let $M=\{i_1, \dots, i_m\}$, for $m$ in $2, \dots, n$, be a 
subset of $\{1, \dots, n\}$.
The {\bf separating cluster} $C_M$ of $M$ in $x$ is
the smallest cluster of $\mathcal{C}(x)$ that contains ${i_1}, \dots, {i_m}$.
The {\bf separating screen} \index{separating screen}
$S_{M}$ is the unique screen that contains exactly
the elements of~$M$. 
\item Let $\mathcal{C}(x)$ be as in Theorem \ref{mathcalC}.
The $\mathbf{x}${\bf-screens} are the screens $S_C$ for $C \in \mathcal{C}(x)$. 
\index{x-screens@$x$-screens}
\end{lijst}
\end{definition}

Suppose we are given a point $x \in \XAH$.  By determining the
clusters $C_{ij}$ for $x$, we find the nest $\mathcal{C}(x)$
for $x$. From the nest, we get a list of screens $\mathcal{S}_\mathcal{C}$,
the $x$-screens.
In this section we introduce a data structure
that enables us to fill the $x$-screens, the  hooked tree.
The hooked tree $ht(x)$ encodes both the nest $\mathcal{C}(x)$ of $x$ and
the relative positions of the underlying point set $p_1, \dots, p_n$ 
of $x$.  We define $ht(x)$ recursively.

\subsection{Definition of hooked tree.}

\begin{definition}
\elabel{htree}
The \bfindex{hooked tree} $ht(x)$ of a point 
$x \in \XAH$ is a rooted tree. Some of its edges have tags.
The {\bf vertices} of $ht(x)$ are the clusters of the
nest $\mathcal{C}(x)$. The {\bf root} is the cluster $\{1, \dots, n\}$.
The {\bf leaves} are the clusters $\{1\}, \dots, \{n\}$.
Two clusters are connected by an {\bf edge} if one of the clusters is
maximal in the other.

The $\mathbf{x}${\bf -tags} 
\index{x-tags@$x$-tags}  refer to hooks and angles in $\ah$.
Fix a cluster $S$. 
Order the maximal subclusters in $S$ according to  the smallest label 
in each maximal subcluster. Suppose $S$ has maximal subclusters 
$c_1, \dots ,c_k$ with minimal labels
$$
        j_1 ~<~ j_2 ~<~ \dots ~<~ j_k.
$$
If $S$ is not the top cluster $[n]$, there exists a cluster $T$
that contains $S$ as one of its maximal subclusters. Suppose that $T$
its maximal subclusters have minimal labels 
$$
        i_1 ~<~ i_2 ~<~ \dots ~<~ i_m.
$$
Tag the edges from cluster $S$ to its maximal subclusters
as follows.
\begin{description}
 \addtolength{\itemsep}{-0.5\baselineskip}
\item[Edge to $\mathbf{c_1}$.]
{\em If} $S$ is the top cluster $[n]$ {\em then} the edge from $S$ to $c_1$ is 
tagged \bfindex{top}.  {\em Else}: no tag.
\item[Edge to]$\mathbf{c_2}$.
\begin{lijst}
\item {\em If} $S$ is the top cluster $[n]$ {\em then} the edge 
to $c_2$ is tagged $\alpha_{1j_2}$.
         Such tag is called an {\bf $x$-type 2.a} tag.
	\index{x-type@$x$-type!2.a}
\item {\em Else if}  $j_1 >  i_1$ {\em then} tag  the edge from $S$ 
      to $c_2$ by 
        $h^{\,j_1j_2}_{\,j_1i_1}$. This is an
      {\bf $x$-type 2.b} tag. \index{x-type@$x$-type!2.b}
\item {\em Else if} $j_1 = i_1$ {\em then} tag the edge to $c_2$
      by $h^{\,j_1j_2}_{\,j_1i_2}$.
      This is an {\bf $x$-type 2.c} tag.
\end{lijst}
\item[Edge to $\mathbf{c_t}$, for  $\mathbf{t\geq 3}$.]
      The edge from $S$ to
      $c_t$, for $t>3$, is  tagged $ \ho{j_1j_t}{j_1j_2}$.
      This is an {\bf $x$-type 3}\index{x-type@$x$-type!3} tag.
\end{description}
\end{definition}

The tag $\alpha_{ij}$ in $x$-tags does not stand for the value $\alpha_{ij}(x)$,
but rather for a map $\ah \rightarrow \R/\pi \Z$. The indices $i$ and 
$j$ indicate which particular factor of the form  $\R/\pi \Z$ in $\ah$ is meant.

\begin{definition}
Fix $x \in \XAH$. Let $q$ be an arbitrary point in $\ah$.
A $q$-coordinate $\ho{ik}{ij}(q)$, respectively, $\alpha_{ij}(q)$, $\al{ik}{ij}(q)$,
$\be{ik}{ij}(q)$ is a {\bf $\mathbf{x}$-type 2}  hook, respectively,  angle,  angle, 
{\bf ratio} iff   $\ho{ik}{ij}$,  $\alpha_{ij}$, $\al{ik}{ij}$ or $\be{ik}{ij}$ 
occurs in $x$-tags;  {\bf $\mathbf{x}$-type 3}  hooks,  angles, and ratios
are defined in a similar way.
\end{definition}

\subsection{The hooked path.}

In this section we associate a unique $x$-tag 
to every site $p_i$. Moreover, we show that the hooked tree $ht(x)$ induces
a partial ordering on the labels $1, \dots, n$. 

\begin{definition}
\elabel{defhookedpath}
Fix $x \in \XAH$. 
\begin{lijst}
\item To any site $p_i$, for $i =1,\dots,n$ we associate an
$x$-tag $l_x(p_i)$: it is the first $x$-tag that one encounters
      when going up in $ht(x)$, starting from the leaf $i$.
\item The \bfindex{predecessor} $p_x(p_i)$  of $p_i$ with respect to $x$,
      for $i=2, \dots,n$ is defined as follows:
\end{lijst}
\begin{looplijst}
\item $p_1$ has no predecessor.
\item If $l_x(p_i)$ is an $x$-type 2 tag then  $i=j_2$ in
       Definition \ref{htree}.
\begin{itemize}
\item In case of $x$-type 2.a, $p_x(p_i)=p_1$.
\item In case of $x$-type 2.b, $p_x(p_i)=p_{j_1}$.
\item In case of $x$-type 2.c, $p_x(p_i)=p_{i_2}$.
\end{itemize}
\item If $l_x(p_i)$ is an $x$-type 3 tag then  $i=j_t$ in
       Definition \ref{htree} and $p_x(p_i)=p_{j_2}$.
\end{looplijst}
\begin{lijst}
\setcounter{mijnlijst}{2}
\item To a site $p_i$, for $i=1, \dots,n$ we associate a \bfindex{hooked path}
 $L_x(p_i)$: 
\end{lijst}
\begin{looplijst}
\item $L_x(p_1) = ()$.
\item $L_x(p_i)$ equals $L_x(p_x(p_i))$ with $p_i$ prepended.
\end{looplijst}
If a site $p_j$ occurs in $L_x(p_i)$, then we say that $p_j$ is $x$-{\bf above}
$p_i$. \index{x-above@$x$-above} 
If $p_j = p_x(p_i)$ or $p_i = p_x(p_j)$, then we call $p_ip_j$ an $x$-{\bf leg}.
\index{x-leg@$x$-leg}
\end{definition}

\begin{lemma}
The predecessor relation gives a partial order on 
the set $p_1, \dots , p_n$.
\end{lemma}

\begin{proof}
We have to show that there are no loops. For this purpose we define the 
\bfindex{depth} of a site $p_i$  as the number of edges from the root in $ht(x)$ to the
node   where the edge that has $l_x(p_i)$ as its $x$-tag starts.
Consider site $p_i$. Site $p_1$ has no predecessor.
Assume that $i\neq 1$ and that $p_j$ is its  predecessor: $p_j =  p_x(p_i)$. 
We claim that either $ \text{depth}(p_j)$ $<$ $\text{depth}(p_i)$,
or
$\text{depth}(p_j) = \text{depth}(p_i)$ and $j < i$.
Indeed, if $l_x(p_i)$ is of $x$-type 3, then the last case holds.
And  if  $l_x(p_i)$ of $x$-type 2, then 
$\text{depth}(p_i)$ $ =$  $\text{depth}(p_j) + 1$.
\end{proof}

\begin{lemma}
\elabel{lhtree11}
There is a $1$ to $1$ correspondence between the sites $p_1, \dots, p_n$
and the $x$-tags.
\end{lemma}

\begin{proof}
In Definition \ref{defhookedpath} we have associated an $x$-tag with every
site $p_i$: the $x$-tag $l_x(p_i)$ associated with $p_i$ is the first tag
one encounters going up from the leaf $p_i$ in the hooked tree $ht(x)$.
To finish the proof, we turn around the association: 
we show that for every $x$-tag $l$ there is a unique path to a leaf of the
hooked tree such that no other $x$-tag is encountered. If the edge
tagged $l$ is itself incident with a leaf we are done. 
So assume that the edge tagged $l$ does not end in a leaf. Then it ends in a 
vertex of the hooked tree corresponding to a cluster $S$ of at least two elements. 
By Definition \ref{htree} only the edge from $S$ to its maximal subcluster $M$
with minimal label $j_1$  has no tag. $j_1$ will also be the minimal label of every
subcluster of $M$ that contains $j_1$. Therefore there is a unique  path without
labels to $p_{j_1}$.  
\end{proof}

\subsection{$\dom(x)$ and $\rng(x)$.}

The $x$-tags define a factor of $\xah$.

\begin{definition} Fix $x \in \xah$.
\begin{lijst}
\item The factor  $\dom(x)$ is the factor of $\ah$ that 
corresponds to the angles and hooks in $x$-tags.\index{domx@$\dom(x)$}
The factor of $\ah$ complementary to $\dom(x)$ is denoted by
$\rng(x)$. \index{rngx@$\rng(x)$}
\item By $tv_x$,  denote the projection of $\ah$ onto $\dom(x)$;
\index{tv@$tv_x$}
`tv' stands for \bfindex{tree values}.
The projection of $\ah$ onto $\rng(x)$ is denoted by $ctv_x$, 
the complementary tree values. \index{ctv@$ctv_x$}
\end{lijst}
\end{definition}

The following lemma states that $ \dim(\,\dom (x)\,)$  $=$  $\dim(\,\conf )$.
This shows that the hooked tree is optimal in the sense that it refers
to as many coordinates of $\ah$ as the dimension requires.

\begin{lemma}
\elabel{dimdom}
Let $x \in \XAH$. Then $\dim(\, \dom (x) \,)  =  2 n -3$.
\end{lemma}

\begin{proof}
The dimension of $\dom(x)$ is given by the sum of the degrees of freedom in the
hooks and angles referred to  by $x$-tags. 
Every hook is defined by one angle $\al{ik}{ij}$ and one ratio $\be{ik}{ij}$,
so the dimension of one hook equals $2$. The dimension of a single angle
$\alpha_{ij}$ equals~$1$. The dimension of `top' is zero.
From Lemma \ref{lhtree11} it follows that the $x$-tags refer to
$n-2$ hooks and $1$ angle as there are $n$ sites and one `top' tag.
\end{proof}

\subsection{Example of a hooked tree.}

\begin{figure}[ht]
\begin{center}
\setlength{\unitlength}{1.05em}
\begin{picture}(24,15)
\put(12,8){\makebox(0,0)[cc]{
        \leavevmode\epsfxsize=16\unitlength
\epsfbox{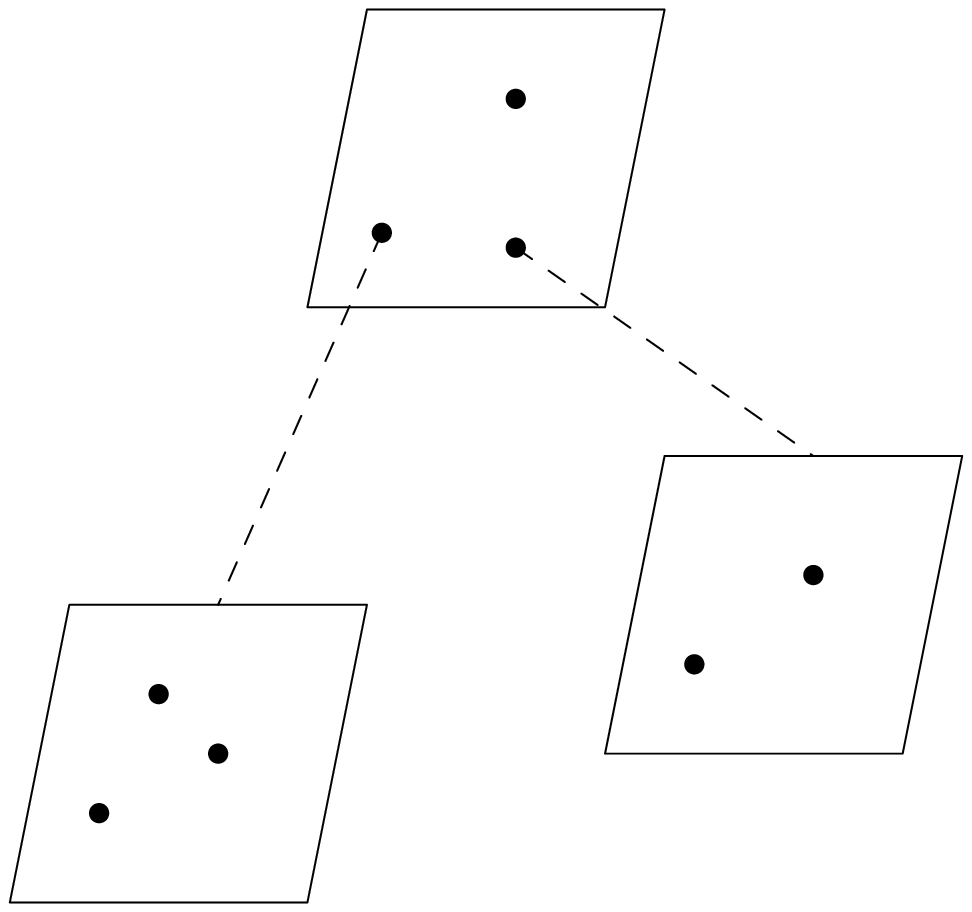}}}
\put(13.2,14){\makebox(0,0)[l]{4}}
\put(5.5,1.5){\makebox(0,0)[l]{1}}
\put(6.2,5){\makebox(0,0)[l]{2}}
\put(8.5,3.5){\makebox(0,0)[l]{6}}
\put(16,5){\makebox(0,0)[l]{3}}
\put(17.7,6.8){\makebox(0,0)[l]{5}}
\end{picture}
\caption{\elabel{screens.eps} The screens for the nest 
	$\langle\{1,2,6\},\{3,5\}\rangle$.}
\setlength{\unitlength}{1.05em}
\begin{picture}(24,15)
\put(12,8){\makebox(0,0)[cc]{
        \leavevmode\epsfxsize=20\unitlength
\epsfbox{ht.eps}}}
\put(20,7.2){\makebox(0,0)[l]{4}}
\put(2.3,0.9){\makebox(0,0)[l]{1}}
\put(7.1,0.9){\makebox(0,0)[l]{2}}
\put(12.0,0.9){\makebox(0,0)[l]{6}}
\put(16.8,0.9){\makebox(0,0)[l]{3}}
\put(21.7,0.9){\makebox(0,0)[l]{5}}
\put(7.7,10){\makebox(0,0)[l]{top}}
\put(12.3, 10){\makebox(0,0)[l]{$\alpha_{1,3}$}}
\put(18.1, 10){\makebox(0,0)[l]{$h^{14}_{13}$}}
\put(7.2, 3.5){\makebox(0,0)[l]{$h^{12}_{13}$}}
\put(11.0, 3.5){\makebox(0,0)[l]{$h^{16}_{12}$}}
\put(20.5, 3.5){\makebox(0,0)[l]{$h^{35}_{31}$}}
\end{picture}
\caption{\elabel{ht.eps} The hooked tree for the nest 
	$\langle\{1,2,6\},\{3,5\}\rangle$.}
\end{center}
\end{figure}

\begin{eexample}
\elabel{exht}
Let us consider again, as in Example \ref{exZ}, the nest 
	$Z = \langle\{1,2,6\},\{3,5\}\rangle$.  
A realization of $Z$ consists 
of the three filled screens $S_{[6]}, S_{1,2,6}$ and $S_{3,5}$.
Recall that $[6] = \{1,\dots,6\}$.
A possible realization is given in 
Figure~\ref{screens.eps}. The hooked tree corresponding to
$Z$ is given in Figure~\ref{ht.eps}. Note that the dimension of $\dom(x)$
is indeed $9$.

The first column of  Table \ref{xqscheme} gives all sites in $Z$.
The $x$-tags corresponding to each site  are shown in the second column. 
The third column gives for every $x$-tag its $x$-type. The fourth column 
gives $tv_x(x)$ coordinates 
corresponding to the way the $x$-screens are filled in Figure \ref{screens.eps}.
How this filling is done exactly,  and the meaning of the
fifth and sixth column will be explained later in this chapter.

From the data for the $tv_x(x)$ coordinates it follows, for example, that the length
of the leg $12$ is zero with respect to the length of the leg $13$.
And as the length of the leg $16$ is of the same order as the length
of the leg $12$, it must vanish with respect to the length
of the leg $13$ too. Continuing in this way reveals the complete
nest structure.
\end{eexample}

\begin{table}[ht]
\begin{center}
\begin{tabular}{ccccccc}
site & $x$-tag & $x$-type & $tv_x(x)$ & $q \in \dom(x)$ & representatives \\
&\\
1 & top & top \\[.1cm]
3 & $\alpha_{1,3}$ & 2.a & $-7^{\circ}$ & 
	 $69^{\circ}$ & same \\[.1cm]
2 & ($\be{12}{13}$, $\al{12}{13}$)  & 2.c  & ($0$, $71^{\circ}$) &
	($.35$, $-91^{\circ}$)	&  ($-.35$, $89^{\circ}$) \\[.1cm]
4 & ($\be{14}{13}$, $\al{14}{13}$)   & 3 & ($1.5$, $50^{\circ}$) &
	  ($-1.92$, $4^{\circ}$) & same \\[.1cm]
5 & ($\be{35}{31}$, $\al{35}{31}$)  & 2.b & ($0$, $44^{\circ}$) &
	($-.34$, $170^{\circ}$) &  ($.34$, $-10^{\circ}$) \\[.1cm]
6 & ($\be{16}{12}$, $\al{16}{12}$)  & 3 & ($1$, $-37^{\circ})$ &
	  ($.91$, $-11^{\circ}$) & same 
\end{tabular}
\end{center}
\caption{\elabel{xqscheme}$x$-tags and values.}
\end{table}

\section{Filling $x$-screens.}

In this section we present a method for filling $x$-screens,
given a suitable point $q \in \dom(x)$.

\subsection{Standard form of $x \in \XAH$.}

First we put a fixed $x \in \XAH$ in standard form:  
if necessary we
change some representatives in Klein bottles $\K{ik}{ij}$
and/or change the representative for the  top angle 
$\alpha_{1j_2}(x)$ modulo $\pi$.
We only use the $\dom(x)$ factor of $\ah$ in order to fill
the $x$-screens. Therefore we only have to choose standard representatives
for those angles and hooks that are referred to in $x$-tags. 

\begin{definition}
\elabel{dstandard}
$x \in \XAH$ is represented in \bfindex{standard form} iff 
\begin{itemize}
\item $\alpha_{1j}$,  $\al{ik}{ij}$ of $x$-type 2 $\Rightarrow$
	$\alpha_{1j}(x)$, $\al{ik}{ij}(x)$ $\in [-\pi/2, \pi/2)$ and
\item $\be{ik}{ij}$ of $x$-type 3 $\Rightarrow$
        $\be{ik}{ij}(x) > 0$.
\end{itemize}
\end{definition}

By the following, we do not have to 
make a sign choice for any $x$-type~2 ratio $\be{ik}{ij}(x)$.

\begin{lemma}
\elabel{clbex}
Fix $x \in \XAH$. 
\begin{lijst}
\item If $\ho{j_1j_2}{j_1i_{1,2}}$ is an $x$-type 2 hook,  
then $\be{j_1j_2}{j_1i_{1,2}}(x)  = 0$.
\item If $\ho{j_1j_t}{j_1j_2}$ is an $x$-type 3 hook, then
$0 < |\be{j_1j_t}{j_1j_2}(x)| < \infty$. 
\end{lijst}
\end{lemma}

\begin{proof}
Follow the notation of Definition \ref{htree}. 
We give a proof for the $x$-type 2.b case. The $x$-type 2.c case
is similar and is omitted.
From Definition \ref{htree} we know that there exists
an $x$-screen $S$ such that  $p_{j_1}, p_{j_2} \in S$ but
$p_{i_1} \not\in S$. The claim follows from the definition
of separating cluster, see in particular Equation \ref{eqcij}.
If $\be{j_1j_t}{j_1j_2}$ is of $x$-type 3, then $p_{j_2} \in C_{j_1j_t}$,
so $|\be{j_1j_t}{j_1j_2}(x)| > 0$. Suppose that 
$|\be{j_1j_t}{j_1j_2}(x)| = \infty$. 
Then $|\be{j_1j_2}{j_1j_t}(x)| = 0$, but this would imply that
$p_{j_2} \not\in C_{j_1j_t}$. That is a contradiction.
\end{proof}

Changing an $x$-type 3 ratio $\be{ik}{ij}(x)$ from negative to positive 
involves of course
a swap $\al{ik}{ij}(x) = \al{ik}{ij}(x) + \pi$ for the corresponding $x$-type 3 angle:
we have to oblige the Klein bottle identification \ref{kleinid}.

\begin{remark} 
We also refer to $x$-type 3 hooks, angles and ratios as 
\bfindex{hooks on scale}, etc. The unique $x$-type 2.a angle is referred
to as \bfindex{top angle}. Hooks, angles and ratios of $x$-type  2.b or 2.c 
are also called \bfindex{explosion hooks}, etc.
\end{remark}

An explanation for these names can be found in the geometric interpretation
of a hook $\ho{ik}{ij}$, given in Section \ref{secgeo}. Think of what changes
in the hook $\ho{ik}{ij}$ as we start from the fixed value $x \in \XAH$
and move away from this point in some general direction inside
$\ah$. If $\ho{ik}{ij}$ is of $x$-type 2, then $\be{ik}{ij}(x)=0$. Perturbing
$x$ to $\tilde{x}$ implies in general  that $|\be{ik}{ij}(\tilde{x})|> 0$.
In the geometric interpretation, this has the effect of an explosion. 
If $\ho{ik}{ij}$ is of $x$-type 3, then $0 < |\be{ik}{ij}(x)| < \infty$ and
shifting $x$ in general direction  does not change this. In our geometric 
interpretation this means that the transformed leg, before and after perturbation,
 has finite but positive length. This explains the name `hook on scale'.
The name `top angle' refers to the place of the $x$-tag `top' in the hooked tree 
$ht(x)$.

\begin{definition}
\elabel{cond}
Let $x \in \XAH$ be in standard form.
A point $q \in \dom(x)$ is \bfindex{valid} with respect to 
$x$ iff the following two conditions on the coordinates of $q$ hold:
\begin{description}
\item[non-orthogonality:] 
no explosion  angle  is orthogonal to the corresponding 
angle $\al{ij}{ik}(x)$. 
The top angle $\alpha_{1j_2}(q)$  is not orthogonal to the top angle
$\alpha_{1j_2}(x)$.
\item[finiteness:] $|\be{ij}{ik}(q)| < \infty$ for all
ratio coordinates $\be{ij}{ik}(q)$. 
\end{description}
\end{definition}

\begin{remark}
$tv_x(x)$ is valid with respect to $x$..
\end{remark}

\begin{proof}
The non-orthogonality is obvious, while
the finiteness follows from Lemma \ref{clbex}.
\end{proof}

\subsection{Example: ordering the sites for the drawing process.}

\begin{figure}[ht]
\begin{center}
\setlength{\unitlength}{1.15em}
\begin{picture}(16.8,12)
\put(8.4,6){\makebox(0,0)[cc]{
        \leavevmode\epsfxsize=8.3\unitlength
\epsfbox{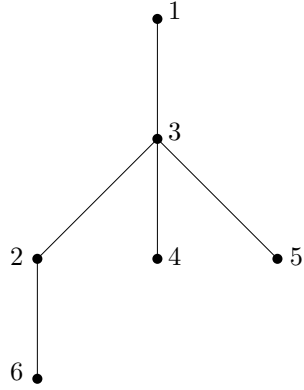}}}
\put(3.7,4.1){\makebox(0,0)[l]{2}}
\put(8.9,4.1){\makebox(0,0)[l]{4}}
\put(12.9,4.1){\makebox(0,0)[l]{5}}
\put(3.7,0.3){\makebox(0,0)[l]{6}}
\put(8.9, 8.2){\makebox(0,0)[l]{3}}
\put(8.9, 12.2){\makebox(0,0)[l]{1}}
\end{picture}
\caption{\elabel{posetht.eps} The poset induced by the predecessor relation.}
\end{center}
\end{figure}

\begin{eexample}
\elabel{exposetht}
Figure \ref{posetht.eps} depicts the poset induced by the predecessor
relation in our running example, see also Examples \ref{exZ} and \ref{exht}.
This poset tells us, for example, that
in order to draw $p_6$ in the top screen, we should draw
$p_1$, $p_3$, and $p_2$ first and in that order and, moreover, that we can ignore 
$p_4$ and $p_5$. So, the hooked path $L_x(p_6)$ is given by
	$L_x(p_6)  =  (p_6, p_2, p_3, p_1)$.

This hooked path $L_x(p_6)$ tells us how to draw $p_6$ in the top 
$x$-screen using $ht(x)$.
First $p_1$ is placed in the origin. Then $p_3$ is drawn using
$l_x(p_3) = \alpha_{13}$. As by now the leg $p_1p_3$ exists,
we can apply $l_x(p_2) = \ho{12}{13}$ in order to add the leg
$p_1p_2$. Finally, the leg $p_1p_6$ is formed by applying
$l_x(p_6) = \ho{16}{12}$.
\end{eexample}

\subsection{The construction $\draw_x(q)$.}

Construction \ref{draw} indicates how to fill $x$-screens with
sites $p_1, \dots, p_n$, given a valid point $q \in \dom(x)$.
First we define an ordering on the $x$-screens. 

\begin{definition}
\elabel{above} An $x$-screen $T$ is \bfindex{above} an $x$-screen $S$ iff the set
of labels of $S$ is a subset of the set of labels of $T$.
\end{definition}

The following lemmas relate the $x$-type of a hook to clusters $C$ in $\mathcal{C}(x)$.

\begin{lemma} 
\elabel{whichscreen}  Let $\ho{ik}{ij}$ be an $x$-tag.
\begin{lijst}
\item If $x$-type$(\,\ho{ik}{ij}\,) = 3$ then $C_{ij} = C_{ijk} =  C_{ik}$. 
\item If $x$-type$(\, \ho{ik}{ij}\,) = 2$ then $C_{ij} \neq C_{ik}$ and 
      moreover, $C_{ij}$ is directly above $C_{ik}$.
\end{lijst}
\end{lemma}

\begin{proof} We prove the two claims in the Lemma.
\begin{lijst}
\item By construction, an $x$-type 3 hook involves the minimal labels of the
first, the second, and another maximal cluster in a fixed cluster. Therefore
$$C_{ij}=C_{ijk}=C_{ik}.$$
\item Analyzing the cases in Definition \ref{htree} shows that $i$ points
to the minimal label in some cluster $S$ and $k$ points to the second
lowest label in $S$. Moreover, $i$ is a minimal label for one of the maximal
clusters of $T$, the cluster directly above $S$, while $j$ is the minimal
label of one of the other maximal clusters in $T$. Therefore
\begin{equation*}
C_{ij} ~=~ T ~\neq~ S ~=~ C_{ik}.
\qedhere
\end{equation*}
\end{lijst}
\end{proof}

\begin{lemma}
\elabel{leminC}
Let $C$ be a $x$-cluster. Let $p_{m_1}$ and $p_{m_2}$ be the sites with minimal 
labels in the first, resp., the second maximal cluster of $C$.
Let $p_k \in C$, but $p_k \neq p_{m_1}, p_{m_2}$, and let $\ho{ik}{ij}=l_x(p_k)$.
Then $p_i,p_j \in C$.
\end{lemma}

\begin{proof}
Let $M_1, M_2, \dots, M_s$ be the maximal subclusters of $C$. Suppose that
$p_k \in M_l$ for l in $\{1,\dots,s\}$.  
Suppose that $x$-type$(\ho{ik}{ij}) = 3$. By Definition \ref{htree}, there
are two possibilities. Either $p_k$ has minimal label in $M_l$, and in this 
case $l > 2$, $p_i =p_{m_1}$, and $p_j=p_{m_2}$ or, $p_k, p_i$ and $p_j$
belong all three to $M_l$.
We are left with the case $x$-type$(\ho{ik}{ij})=2$. According to Definition
\ref{htree},  $p_k$ has minimal label in the second maximal
cluster of some cluster $D$. As $p_k \neq p_{m_2}$, this means that 
either $D$ equals  $M_l$ or $D$ is a subcluster of $M_l$. It follows from Lemma
\ref{whichscreen} that in both cases $p_i, p_j \in C$, 
since these two sites are contained in the cluster directly above $D$.
\end{proof}

\construction {\bf ~(} $\mathbf{\draw_x(q)}${\bf ~)}. 
\elabel{draw}
\index{draw@$\draw_x(q)$}
\\Filling $x$-screens.

Input:  $x$-screens, valid $q \in \dom(x)$.\\
Output:  $x$-screens filled with $q$, prescription $\draw_x(q)$.

\step {\bf Choosing representatives.}\\
Choose representatives for the top angle $\alpha_{1j_2}(q)$ and the 
      explosion angles $\al{ij}{ik}(q)$ in $ht(q)$. 
      Do so by changing the present value, if necessary, to the value 
      modulo $2 \pi$ most close to the corresponding value $\alpha_{1j_2}(x)$
      or $\al{ij}{ik}(x)$ in $tv(x)$.
      This choice is unique because of the non-orthogonality condition in Definition 
      \ref{cond}. See also Figure \ref{pi2}, which demonstrates that there
      always exists a value of $\alpha(q)$ within distance $\pi/2$ of 
      $\alpha(x)$.

\begin{figure}[!ht]
\begin{center}
\setlength{\unitlength}{1.15cm}
\begin{picture}(10,3)
\put(5,1.5){\makebox(0,0)[cc]{
        \leavevmode\epsfxsize=2.5\unitlength
\epsfbox{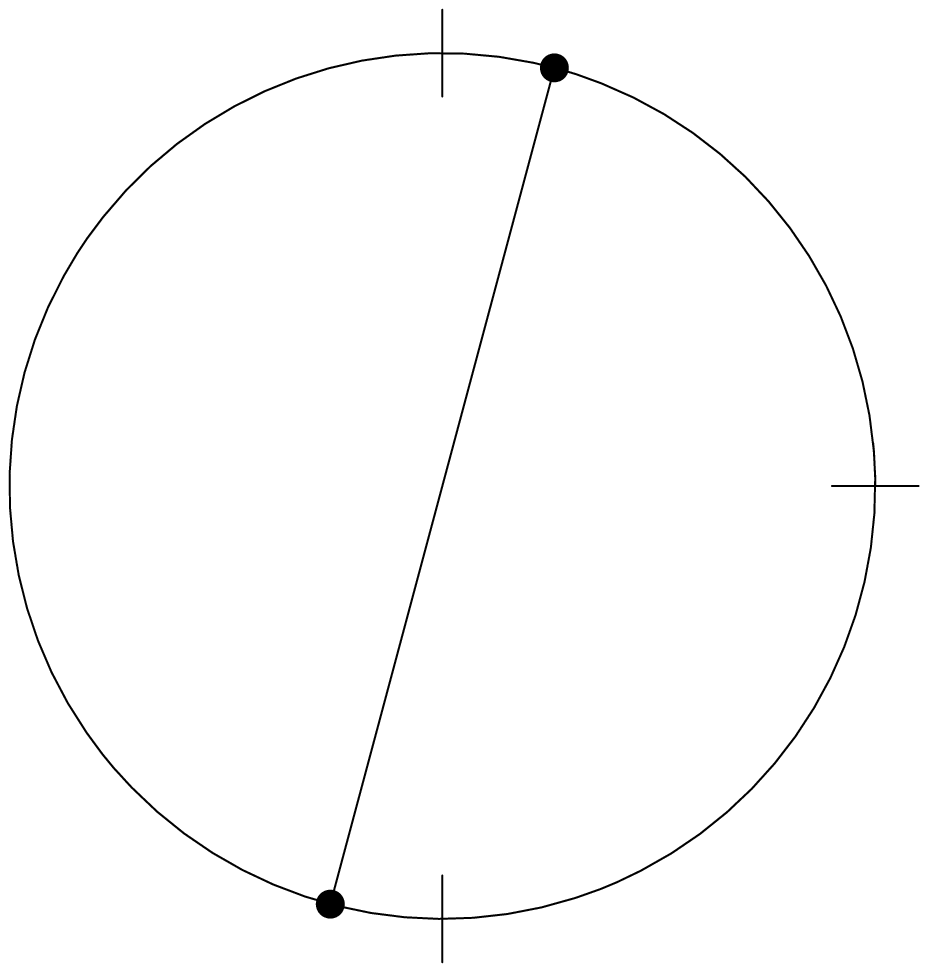}}}
\put(6.3,1.6){\makebox(0,0)[l]{$\alpha(x)$}}
\put(5.3,2.8){\makebox(0,0)[l]{$\alpha(q)$}}
\put(4.3,0.2){\makebox(0,0)[l]{$\alpha(q)$}}
\end{picture}
\caption{\elabel{pi2} Non-orthogonality condition.}
\end{center}
\end{figure}
      Simultaneously, change the sign of an explosion ratio
      $\be{ik}{ij}(q)$ whenever you change the explosion angle
      by $\al{ik}{ij}(q)=\al{ik}{ij}(q)+\pi$. Otherwise the $\K{ik}{ij}$
      Klein bottle identification \ref{kleinid} is not obliged.

\step {\bf Screen orientation.}\\
We define the \bfindex{screen orientation} of an $x$-screen $S=S_C$  by 
induction on the `above' relation introduced in Definition \ref{above}. 
Assume that the screen orientation is defined for all $x$-screens 
above $S$. Let $T$ be the screen directly
above $S$. Following the labeling in Definition \ref{htree}, the screen
orientation $\mathcal{O}_S$ is defined by
        \begin{equation}
	\label{eqscreenos}
   \mathcal{O}_S ~:=~
         \begin{cases}
           \alpha_{1j_2}, & \text{$S$ is a top screen}, \\
           \mathcal{O}_T + \alpha^{j_1j_2}_{j_1i_2},
                         & j_1~=~i_1,\\
           \mathcal{O}_T + \alpha^{j_1j_2}_{j_1i_1} + \pi,
                         & j_1~=~i_2,   \\
           \mathcal{O}_T + \al{i_1j_1}{i_1i_2}+
                         \alpha^{j_1j_2}_{j_1i_1} + \pi,
                         & j_1~=~i_t, ~t=3, \dots, m, \quad
                         \be{i_1j_1}{i_1i_2} ~>~ 0,\\
           \mathcal{O}_T + \al{i_1j_1}{i_1i_2}+
                         \alpha^{j_1j_2}_{j_1i_1},
                         & j_1~=~i_t, ~t=3, \dots, m, \quad
                         \be{i_1j_1}{i_1i_2} ~<~ 0.
         \end{cases}
 \end{equation}  
\step {\bf Drawing the sites.} \\
The recipe for drawing a site $p_i$ in a screen $S=S_C$ for $i \in C$,  
      is given inductively, by induction on the  predecessor relation.
      Below, $T$ denotes the screen directly above $S$.
\begin{algorithmic}[1]
\STATE Site $p_1$ is the first site in every path $L_x(p_i)$
        and is drawn in `the origin' in any screen where $p_1$ occurs.
\STATE Next suppose that $i \neq 1$ and that all sites $p_m$ above $p_i$ with
        $m \in C$ have already been drawn.
\IF{$p_i$ is the first site in the first cluster of $C$}
\STATE $p_i$ is drawn in the origin of $S$.
\ELSIF{$p_i$ is the first site in the second cluster of $C$}
\STATE put $p_i$ at $( \cos \mathcal{O}_S, \sin \mathcal{O}_S )$
\ELSE 
\STATE $l_x(p_i) = \ho{ji}{jk}$, where $p_j$ and
              $p_k$ are above $p_i$. By the induction hypothesis, the leg
              $p_jp_k$ is already present in $S$, see also
              Lemma~\ref{leminC}. The leg $p_jp_i$ is constructed as follows:
\STATE Rotate $p_jp_k$ over $\al{ji}{jk}$ around $p_j$.
\STATE Scale the result by a factor $\be{ji}{jk}$ with respect
                      to $p_j$.
\STATE  To summarize,
                      \begin{eqnarray}
                        \elabel{rotbet}
                        p_i & = & \be{ji}{jk} \rot_{\al{ji}{jk}}(p_k - p_j)
                                + p_j.
                      \end{eqnarray}
\STATE  By induction, the length of the leg $p_jp_k$ is finite.
        The finiteness condition in Definition \ref{cond} assures
         that the resulting leg $p_jp_i$ is finite as well.
\ENDIF
\end{algorithmic}

\begin{remark}
Note that the legs $p_ip_j$ constructed in $\draw_x(q)$ while
filling the screens are exactly the $x$-legs. 
\end{remark}

Some drawings are better than others.

\begin{definition}
\elabel{defaccept}
A valid $q\in \dom(x)$ is \bfindex{accepted} with respect to $x$ 
iff in any $x$-screen $S$ filled by $\draw_x(q)$, $q$-sites 
that belong to distinct maximal $x$-clusters of $S$ do not coincide.
\end{definition}

\begin{remark} 
Suppose that $q \in \dom(x)$ is valid and accepted. 
If type$_x(\be{ik}{ij})=3$, then $|\be{ik}{ij}(q)| > 0$.
\end{remark}

\begin{definition}
Every valid and accepted  $q \in \dom(x)$ gives a filling of the $x$-screens
consisting of $n$ sites.  The set of all filled $x$-screens that can
be obtained in this way is denoted by $\scr(x)$.
\index{scr-x@$\scr(x)$}
\end{definition}

\subsection{Example: filling $x$-screens.}

\begin{figure}[!ht]
\begin{center}
\setlength{\unitlength}{1cm}
\begin{picture}(6,6)
\put(3,5){\makebox(0,0)[cc]{
        \leavevmode\epsfxsize=2\unitlength\epsfbox{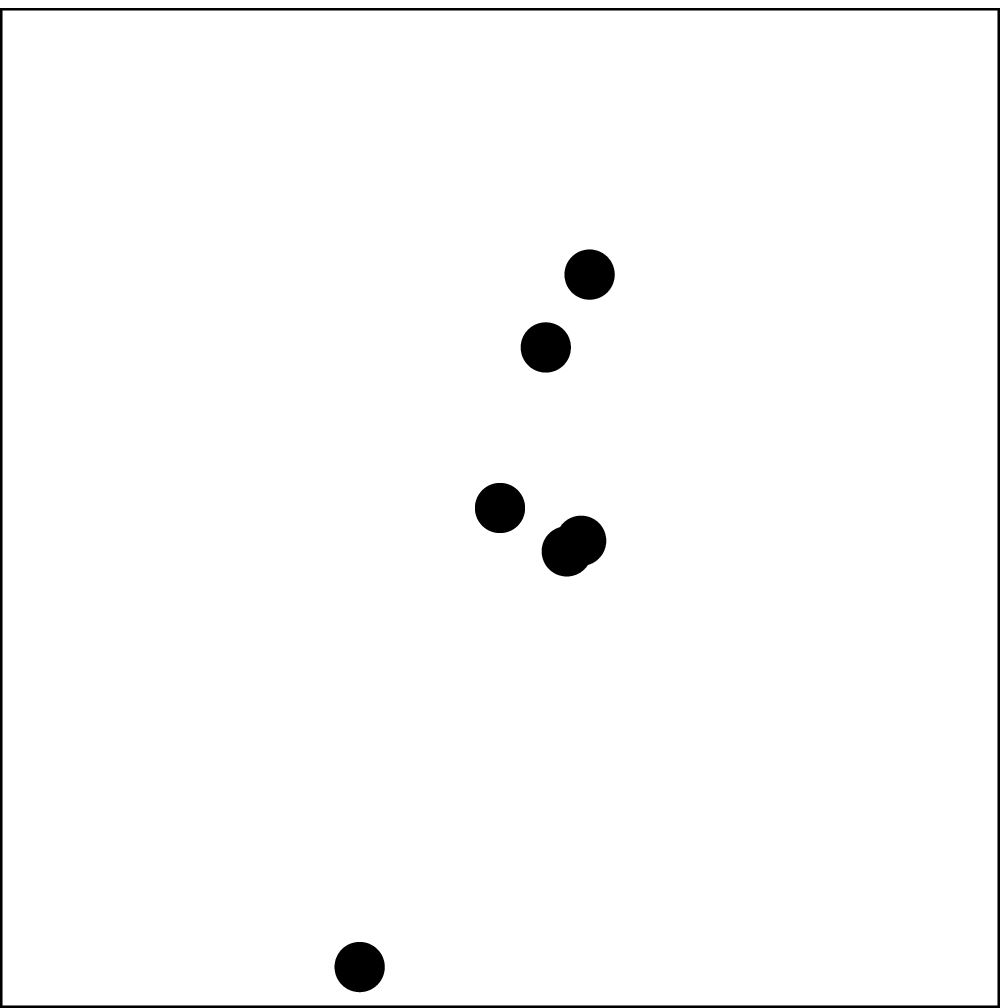}}}
\put(1.5,1){\makebox(0,0)[cc]{
        \leavevmode\epsfxsize=2\unitlength\epsfbox{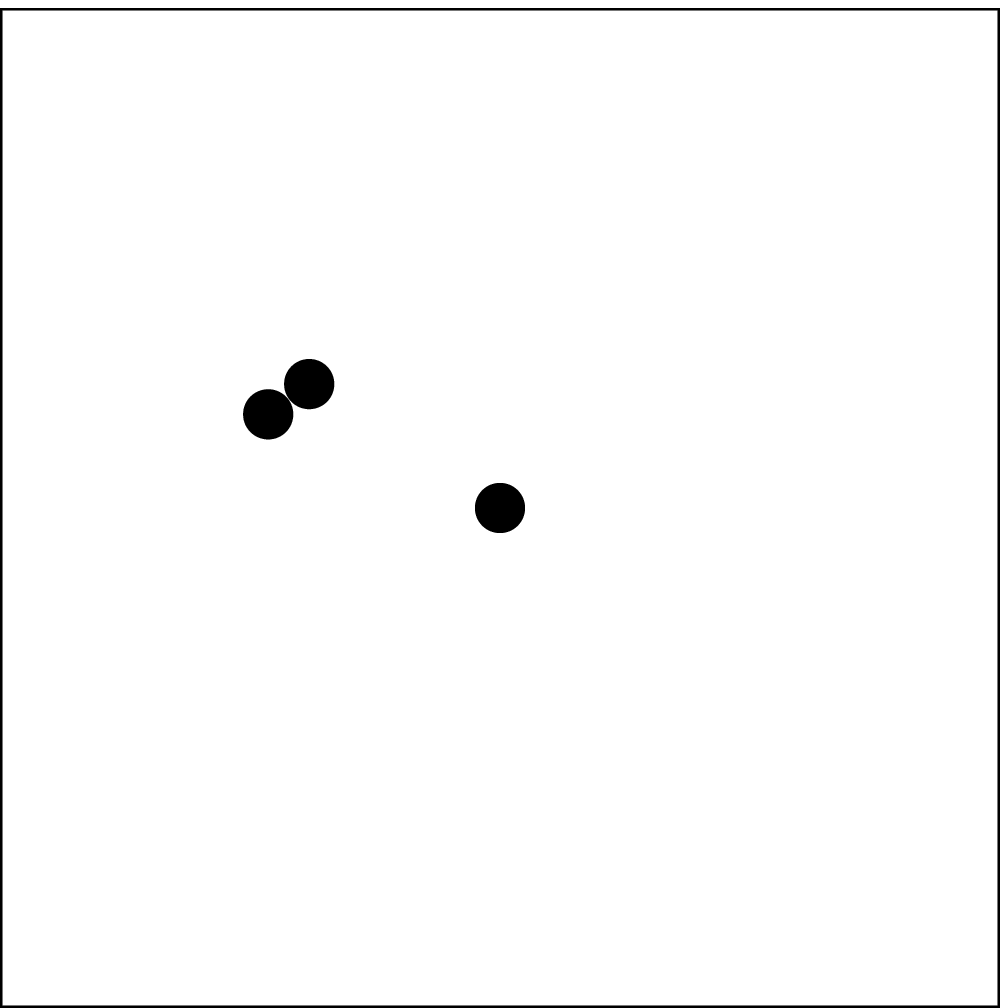}}}
\put(4.5,1){\makebox(0,0)[cc]{
        \leavevmode\epsfxsize=2\unitlength\epsfbox{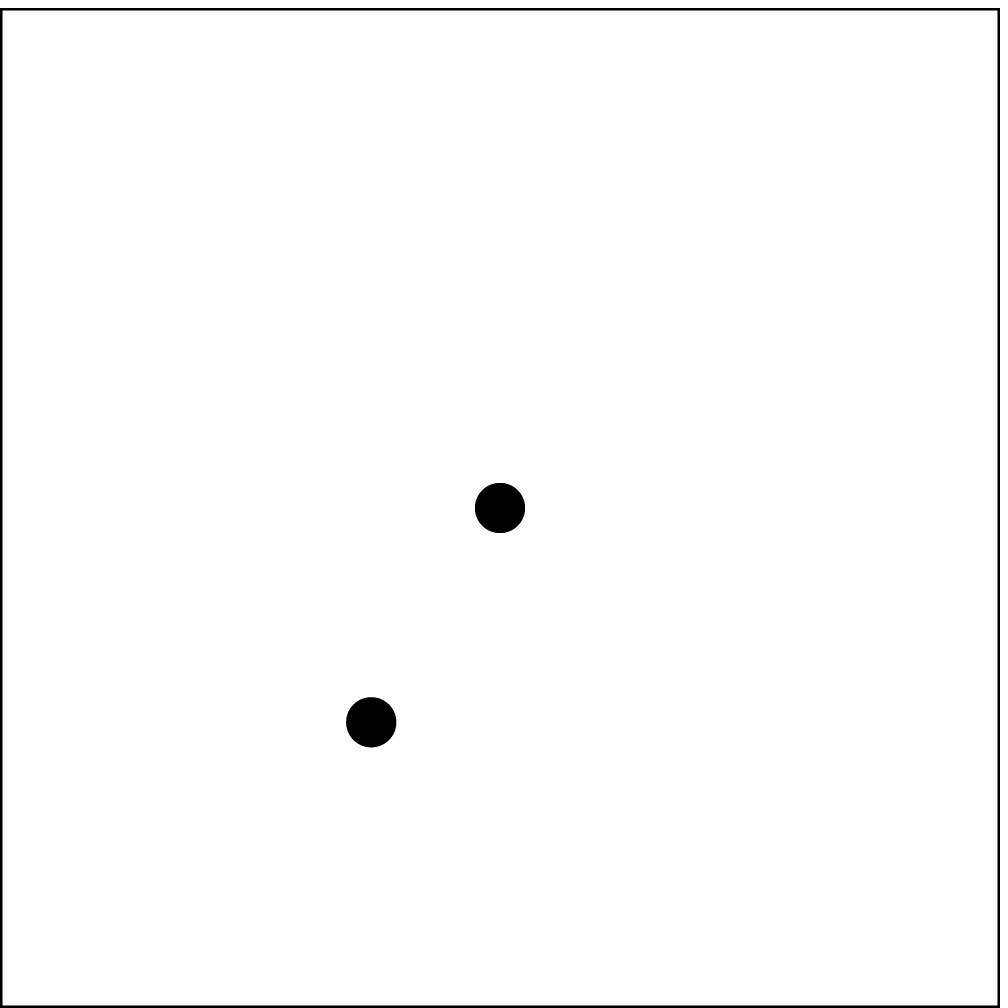}}}
\put(2.8,5){\makebox(0,0)[c,c]{1}}
\put(3.3,4.75){\makebox(0,0)[c,c]{2}}
\put(3.4,5.7){\makebox(0,0)[c,c]{3}}
\put(2.95,4.2){\makebox(0,0)[c,c]{4}}
\put(2.95,5.4){\makebox(0,0)[c,c]{5}}
\put(3.1,4.65){\makebox(0,0)[c,c]{6}}
\put(1.75,1.2){\makebox(0,0)[c,c]{1}}
\put(.9,1.0){\makebox(0,0)[c,c]{2}}
\put(1.3,1.5){\makebox(0,0)[c,c]{6}}
\put(4.75,1.2){\makebox(0,0)[c,c]{3}}
\put(4.1,.4){\makebox(0,0)[c,c]{5}}
\unitlength=1mm
\linethickness{.3mm}
\dottedline{1}(17,20)(31,50)
\dottedline{1}(43,20)(32,55)
\end{picture}
\caption{\elabel{qfilled} $x$-screens filled with some $q\in\dom(x)$.}
\end{center}
\end{figure}

\begin{eexample}
\elabel{exfill}
In Figure \ref{qfilled} we have  filled the $x$-screens of our running example,
see Examples \ref{exZ}, \ref{exht} and \ref{exposetht},
with the $q$-values present in the the fifth column in Table \ref{xqscheme}.
The chosen representatives for the $q$-values are given in the last
column. According to Step 2 in Construction \ref{draw}, we have changed the
representatives for $\ho{12}{13}$ and $\ho{35}{31}$. 
Let us look more carefully at the relative positions of the sites 
$p_1, p_2,$ and $p_6$ in the top screen $S_{[6]}$ and the screen 
$S_{1,2,6}$. It seems like the complete cluster $1,2,6$ has been rotated
over $\pi$ w.r.t.\ site $p_1$ in the bottom screen. 
In order to find an explanation, we compute the coordinates of $p_2$
in both $S_{[6]}$ and $S_{1,2,6}$. The hooked path $L_x(p_2)$
equals $(p_2, p_3, p_1)$.
In both screens $p_1 = (0,0)$. In the top screen, 
$p_3 = ( \cos \mathcal{O}_{S_{[6]}}, \sin \mathcal{O}_{S_{[6]}} )$,
where $\mathcal{O}_{S_{[6]}} = \alpha_{13}$.
We determine the position of  $p_2$ in the top screen, using Equation \ref{rotbet}:
$$
	p_2(S_{[6]}) ~=~ \be{12}{13} \rot_{\al{12}{13}}( p_3 )
	  ~=~ \be{12}{13} ( \cos \alpha_{13} + \al{12}{13},
				\sin \alpha_{13} + \al{12}{13} ).
$$
In $S_{1,2,6}$, however, $p_2$ is the first site in the second cluster,
so 
\beq
p_2(S_{1,2,6}) & = & ( \cos \mathcal{O}_{S_{1,2,6}}, \sin \mathcal{O}_{S_{1,2,6}} ).
\eeq
According to Step~2 in Construction \ref{draw},  
$$
	\mathcal{O}_{S_{1,2,6}} ~=~ \mathcal{O}_{S_{[6]}} + \al{12}{13}
		~=~ \alpha_{13} + \al{12}{13}\text{,}
$$
as $i_1=1$, $i_2=3$ and $j_1=1, j_2 =2$. We conclude that  the only 
difference between the argument of $p_3$ in  $S_{[6]}$ and $S_{1,2,6}$
can be caused by $\be{12}{13}$: this ratio is used in drawing $p_2$ in
$S_{[6]}$ but not used in drawing $p_2$ in $S_{1,2,6}$. Indeed,
sign$(\,\be{12}{13}\,)=-1$, which causes a $\pi$ shift in arg$(p_2)$  from the
top screen to the screen $S_{1,2,6}$.
\end{eexample}
\section{Reading off data elements from filled screens.}

Let $q \in \dom(x)$. Suppose we have filled 
the $x$-screens with $q$-sites by applying $\draw_x(q)$. 
In this section we tell
exactly how to read off angles and hooks coordinates
from the $q$-filled $x$-screens. 
More precisely, we explain how to determine 
every angle $\overline{\alpha_{ij}}$ between two
sites $p_i(q)$ and $p_j(q)$ and
every hook $\ho{ij}{ik}$ between
three sites $p_i(q)$, $p_j(q)$, and $p_k(q)$. 
Basically we read off lacking coordinates in a suitable $x$-screen. 
Recall the notions of separating clusters  and separating screens that
were introduced in Definition \ref{defsepa}. These notions are used
for fixing the $x$-screens used for 
determining the angles and hooks coordinates.

\begin{remark}
The separating screen $S_{ij}(x)$ is the only $x$-screen in which it is guaranteed
that $p_i(q) \neq p_j(q)$ for $q$ being valid and accepted.
\end{remark}

\subsection{Determining hooks and angles from $\draw_x(q)$.}

\begin{definition}
\elabel{readoff}
Let $q \in \dom(x)$ be valid and accepted. From $\draw_x(q)$ we determine
a point $\tilde{q}=\readoff_x(q) \in \ah$ as follows. 
\index{readx@$\readoff_x(q)$}
\begin{lijst}
\item $\alpha_{ij}(\tilde{q}) \in \R/2\pi \Z$ is the angle of the
directed line through $p_i(q)$ and $p_j(q)$ in the 
$q$-filled screen $S_{ij}(x)$.
\item  $\al{ik}{ij}(\tilde{q}) = \alpha_{ik}(\tilde{q})
	 - \alpha_{ij}(\tilde{q})$.
\item Read off all sites  coordinates in screen $S_{ijk}$. 
If $p_i(q) \neq p_j(q)$ in $S_{ijk}$, then
	\begin{eqnarray}
		\elabel{readoffbeta}
		\be{ik}{ij}(\tilde{q}) & = & \frac{(p_k(q)-p_i(q))\cdot(p_j'(q)-p_i(q))}
			{|p_j'(q)-p_i(q)|^2},
	\end{eqnarray}
        where 
	 \begin{eqnarray}
                \elabel{pjprime}
	p_j'(q) & = & \rot_{\al{ik}{ij}(\tilde{q})}(p_j(q) - p_i(q)) + p_i(q).
         \end{eqnarray}
Else, $ \be{ik}{ij}(\tilde{q})  =  1 /  \be{ij}{ik}(\tilde{q})$.
\end{lijst}
\end{definition}

Note that by definition of separating screen, not all three sites
$p_i$, $p_j$, and $p_k$ coincide in $S_{ijk}$. Therefore 
$p_i(q)=p_j(q)$ in $S_{ijk}$ implies that $p_i(q) \neq p_k(q)$.

\begin{figure}[!ht]
\begin{center}
\setlength{\unitlength}{1em}
\begin{picture}(12.8,8)
\put(6.4,4){\makebox(0,0)[cc]{
        \leavevmode\epsfxsize=11.2\unitlength\epsfbox{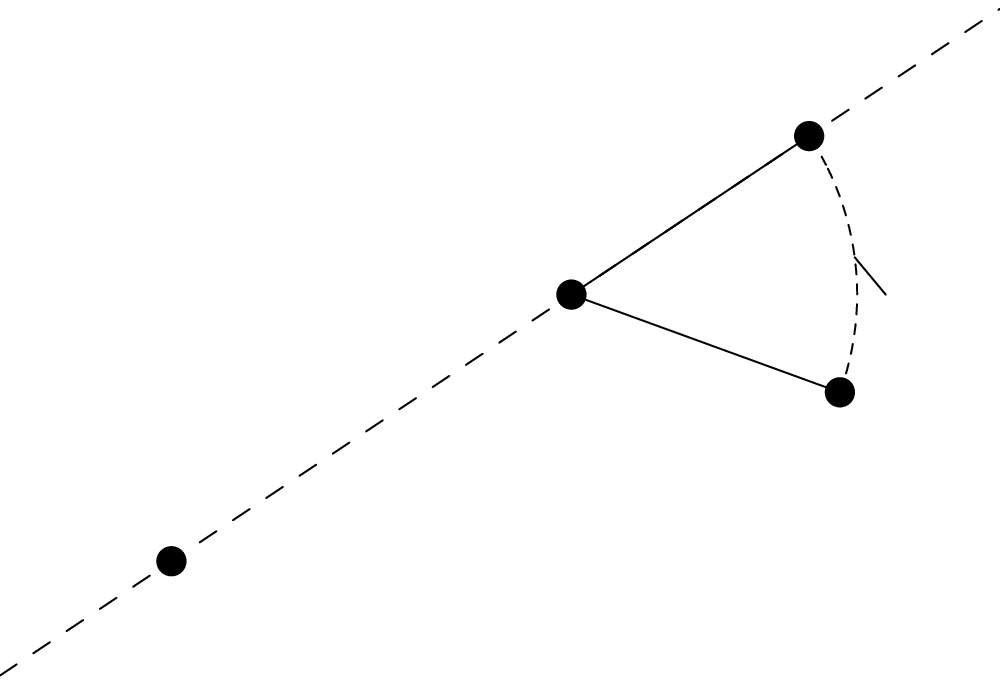}}}
\put(8.3,5.0){\makebox(0,0)[l]{$\al{ik}{ij}$}}
\put(9.6,7.7){\makebox(0,0)[l]{$p_j'$}}
\put(10.4,2.9){\makebox(0,0)[l]{$p_j$}}
\put(6.0,5.3){\makebox(0,0)[l]{$p_i$}}
\put(1.9,2.7){\makebox(0,0)[l]{$p_k$}}
\end{picture}
\caption{\elabel{rotbeta}Determining $\be{ik}{ij}$.}
\end{center}
\end{figure}

\begin{remark}
We give a geometric explanation for Definition \ref{readoff}.
It is rather straightforward to read off $\alpha_{ij}(\tilde{q})$ and
$\al{ik}{ij}(\tilde{q})$, so let us focus on $\be{ik}{ij}(\tilde{q})$.
Go to the separating screen $S_{ijk}(x)$. 
Write down the coordinates of $p_i(q)$, $p_j(q)$, and $p_k(q)$
      in this screen. First consider the case where $p_i(q) \neq p_k(q)$.
      If proceeding through the definition in the right
      order, we know by now $\al{ik}{ij}(\tilde{q})$. See also Figure \ref{rotbeta}.
      Use $\al{ik}{ij}(\tilde{q})$ in order
      to rotate the leg $p_j - p_i$ onto the line $l_{ik}$ through $p_i$
      and $p_k$, thereby fixing a point $p_j'$.
      Finally determine the coordinate of $p_k$ on the $ij'$-axis
      with respect to the `unit vector' $p_j - p_i$ as follows.
      Let $x = p_k-p_i$, $y = p_j'-p_i$, and $\theta$ the angle between
      $x$ and $y$.
      From the definition of dot product
        $x \cdot y  =  |x| |y| \cos \theta$,
      it follows that the quantity $\frac{|x|}{|y|} \cos \theta$
      we look for is given by $(x \cdot y)/|y|^2$.  
      Note that in our case $\theta$ always equals either $0$ or $\pi$.
\end{remark}

\begin{proposition}
The components of $\readoff_x(q)$ are infinitely often differentiable
($C^\infty$) on the set of valid and accepted $q \in \dom(x)$.
\end{proposition}

\begin{proof}
The main problem is the continuity of $\readoff_x(q)$, so this is
stressed in the proof. To make the proof complete, replace every 
occurrence of `continuous' by $C^\infty$.
\begin{lijst}
\item  $\ah$ is smooth as a direct product of circles and Klein 
      bottles. $\dom(x)$ is a factor of $\ah$ and is a 
      product of circles and Klein bottles as well. $\readoff_x(q)$
      determines a point in $\ah$.
\item We restrict ourselves to valid $q \in \dom(x)$. This means that we
      have imposed two conditions on $q$: the non-orthogonality condition,
      which is is an open condition as it assures that one point 
      of a circle $S^1$ is avoided.
      The finiteness condition is open as well, as one point, infinity,
      of $\P^1$ is avoided.
\item For valid $q \in \dom(x)$, the construction $\draw_x(q)$ is defined. 
      The construction consists of three steps. First of all representatives
      are chosen for the top angle and explosion angles. The chosen  representatives
      depend continuously on $q$. No jumps can occur due to the
      non-orthogonality condition. In step 2 screen orientations
      $\mathcal{O}_S$  are defined. The formulas presented
      in Equation \ref{eqscreenos}, depend continuously on the chosen representatives
      for $q$.  In the third step, points are drawn in the screen. The positions
      of these points depend continuously on the screen orientation $\mathcal{O}_S$
      and the $q$ coordinates, compare Equation \ref{rotbet}.
\item  Filled screens that are not accepted are thrown away. A filled
      screen is not accepted when certain clustering between sites occurs. 
      Avoiding clustering is yet another open condition.
\item For valid and accepted $q \in \dom(x)$, the point $\readoff_x(q) \in\ah$
      is defined. The formulas for the coordinates of 
      $\readoff_x(q)$ are given in Definition \ref{readoff}.
      All equations depend continuously on the position of the
      sites in accepted $\draw_x(q)$.
	\qedhere
\end{lijst}
\end{proof}

\subsection{Example: reading off hooked tree elements.}

\begin{eexample}
\elabel{exreadoff}
Applying Definition \ref{readoff}, we read off  the values of the
angle and hooks in $\dom(x)$ from Figure \ref{qfilled}. 
These values are presented
in Table \ref{fromscreens}. Next to the data elements to read off,
the read off screen is given. Compare the values in Table
\ref{fromscreens} with the values of $q \in \dom(x)$ in 
Table \ref{xqscheme} used above 
to fill the $x$-screens.  Note that these two sets of  values are indeed 
the same but only up to identification in the appropriate Klein bottle.
\end{eexample}

\begin{table}[ht]
\begin{center}
\begin{tabular}{cccccccccc}
angle & screen & value & \qquad\qquad & angle & value & \qquad\quad &
         ratio & screen & value \\
&\\
$\alpha_{1,3}$ & $S_{[6]}$ & $69^{\circ}$ &  & & \\[.1cm]
$\alpha_{1,2}$ & $S_{1,2,6}$ & $157^{\circ}$ &
        & $\al{12}{13}$ & $88^{\circ}$  & & $\be{12}{13}$ & $S_{[6]}$ & -.35 \\[.1cm]
$\alpha_{1,4}$ & $S_{[6]}$ & $-107^{\circ}$ &
        &        $\al{14}{13}$ & $-176^{\circ}$  & & $\be{14}{13}$ & $S_{[6]}$ &
 1.92 \\[.1cm]
$\alpha_{3,5}$ & $S_{3,5}$ & $-123^{\circ}$ &
        & $\al{35}{31}$ & $-12^{\circ}$& &  $\be{35}{31}$ & $S_{[6]}$ & .34 \\[.1cm]
$\alpha_{1,6}$ & $S_{1,2,6}$ & $146^{\circ}$ &
        & $\al{16}{12}$ & $-11^{\circ}$ &  & $\be{16}{12}$ & $S_{1,2,6}$ & .91 
\end{tabular}
\end{center}
\caption{\elabel{fromscreens}Read off values.}
\end{table}

\section{Consistency Theorem.}
In this section we show that our construction of 
filling screens by means of hooked tree data is consistent
with the procedure of reading off data from the filled screens.

\subsection{Managing $x$-type 3 hooks.}

\begin{lemma} 
\elabel{readtype3}
Let $q \in \dom(x)$ be valid and accepted and put 
$\tilde{q} = \readoff_x(q) \in \ah$. 
Suppose that type$_x(\ho{ik}{ij}) = 3$. Then
\begin{equation*}
  \al{ik}{ij}( \tilde{q} ) \quad = \quad 
	\begin{cases}
	\al{ik}{ij}( q ), & \be{ik}{ij}( q ) > 0,\\
	\al{ik}{ij}( q ) + \pi, & \be{ik}{ij}( q ) < 0.
	\end{cases}
\end{equation*}
\end{lemma}

\begin{figure}[!ht]
\begin{center}
\setlength{\unitlength}{1cm}
\begin{picture}(12,3)
\put(6,1.5){\makebox(0,0)[cc]{
        \leavevmode\epsfxsize=12\unitlength\epsfbox{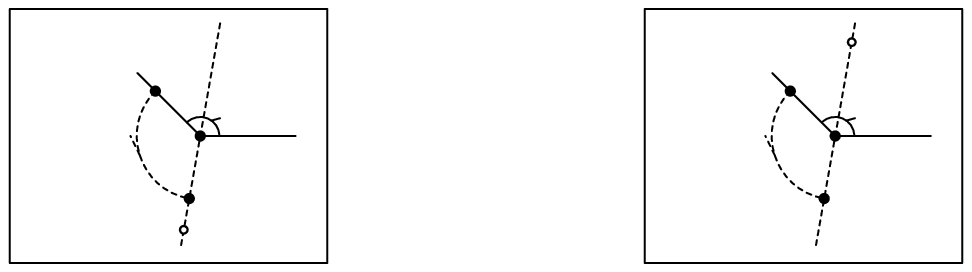}}}
\put(2.7,1.25){\makebox(0,0)[l]{1}}
\put(10.35,1.25){\makebox(0,0)[l]{1}}
\put(2.2,2.3){\makebox(0,0)[l]{2}}
\put(9.75,2.3){\makebox(0,0)[l]{2}}
\put(2.65,0.8){\makebox(0,0)[l]{2'}}
\put(10.3,0.8){\makebox(0,0)[l]{2'}}
\put(2.6,0.3){\makebox(0,0)[l]{t}}
\put(10.15,2.7){\makebox(0,0)[l]{t}}
\put(3.0,1.9){\makebox(0,0)[l]{$\mathcal{O}_S$}}
\put(10.6,1.9){\makebox(0,0)[l]{$\mathcal{O}_S$}}
\put(.7,1.3){\makebox(0,0)[l]{$\al{1t}{12}(q)$}}
\put(8.3,1.3){\makebox(0,0)[l]{$\al{1t}{12}(q)$}}
\end{picture}
\caption{\elabel{readtype$_x$3}$\be{j_1j_t}{j_1j_2}(q) > 0$ 
	and $\be{j_1j_t}{j_1j_2}(q) < 0$.}
\end{center}
\end{figure}

\begin{proof}
Follow the labeling in Definition \ref{htree}:
$i=j_1$, $j = j_2$, and $k=j_t$. For convenience we set
$j_1=1$, $j_2=2$ and $j_t = t$. By Lemma \ref{whichscreen},
all evaluations take place in $x$-screen  $S=S_{12}=S_{12t}=S_{1t}$.
For an illustration of the proof, see Figure \ref{readtype$_x$3}.
As $p_1=(0,0)$ and $p_2 = (\cos \mathcal{O}_S, \sin \mathcal{O}_S)$ it
follows that
$$ \al{1t}{12}(\tilde{q}) ~ = ~ 
		\alpha_{1t}(\tilde{q}) - \alpha_{12}(\tilde{q})
			~ = ~
			\text{arg}(p_t) - \text{arg}(p_2)
			~ = ~ 
			 \text{arg}(p_t) - \mathcal{O}_S.$$
If $\be{1t}{12}(q) > 0$, then
	$\text{arg}(p_t) = \mathcal{O}_S + \al{1t}{12}(q)$,
so $\al{1t}{12}(\tilde{q})  = \al{1t}{12}(q)$.
If, on the other hand, $\be{1t}{12}(q) < 0$, then
        $\text{arg}(p_t)  =  \mathcal{O}_S + \al{1t}{12}(q) + \pi$.
So, $\al{1t}{12}(\tilde{q})  =  \al{1t}{12}(q) + \pi$ as required.
\end{proof}

\begin{lemma}
\elabel{lembe3pos}
Let $q \in \dom(x)$ be valid and accepted. Set 
$\tilde{q} := \readoff_x(q) \in \ah$. 
Suppose that type$_x(\be{ik}{ij}) = 3$. Then
\begin{lijst}
\item $|\be{ik}{ij}(\tilde{q})|=|\be{ik}{ij}(q)|$.
\item $\be{ik}{ij}(\tilde{q}) > 0$.
\end{lijst}
\end{lemma}

\begin{proof}
As in the proof of Lemma \ref{readtype3} set $i=1, j=2, k=t$ and conclude
from Lemma \ref{whichscreen} that the read off screen is given by
screen $S=S_{12}=S_{1t}$.
We use Equation \ref{readoffbeta} for determining $\be{j_1j_t}{j_1j_2}(\tilde{q})$,
while remembering that $p_1=(0,0)$.
\begin{eqnarray}
\elabel{be1t12}
	\be{1t}{12}(\tilde{q}) & = & \frac{p_t \cdot p_2'}{p_2' \cdot p_2'}.
\end{eqnarray}
$p_2'$ is given by
$$
	p_2' ~=~ \rot_{\al{1t}{12}} p_2 ~=~ 
	( \cos [\mathcal{O}_S + \al{1t}{12}(\tilde{q})],
	\sin [\mathcal{O}_S + \al{1t}{12}(\tilde{q})] ).
$$
Therefore $|p_2'|=1$, so $\be{1t}{12}(\tilde{q})$ reduces to
$\be{1t}{12}(\tilde{q}) = p_t \cdot p_2'$. Observe that
$p_2 = (\cos \mathcal{O}_S, \sin \mathcal{O}_S)$  and apply 
Lemma \ref{readtype3}. If $\be{1t}{12}(q)>0$, then
$\al{1t}{12}(\tilde{q})=\al{1t}{12}(q)$ and by evaluating
Equation \ref{be1t12}  it follows that 
	$\be{1t}{12}(\tilde{q}) = \be{1t}{12}(q)$. 
If, on the other hand, $\be{1t}{12}(q)<0$, then
$\al{1t}{12}(\tilde{q})=\al{1t}{12}(q) + \pi$, so
$p_t = - \be{1t}{12}(q) p_2'$, which implies that
$\be{1t}{12}(\tilde{q}) = - \be{1t}{12}(q)$, as claimed.
\end{proof}

\begin{eexample}
\elabel{exconsist3}
In the hooked tree in our running example, see  Examples
\ref{exZ}, \ref{exht}, \ref{exposetht},  \ref{exfill} and
\ref{exreadoff}, two $x$-type 3 hooks
occur, $\ho{14}{13}$ and $\ho{16}{12}$. 
The $q$-values, that are those values that where used to
{\it fill} the $x$-screens, and the $\readoff_x(q)$-values,
the values that where {\it read off} from the filled $x$-screens, are
collected in Table \ref{tabswap}. 
\end{eexample}

\begin{table}[ht]
\begin{center}
\begin{tabular}{ccc}
& $q$ coordinate & $\readoff_x(q)$ coordinate \\
&\\
$\ho{14}{13}$ & $(~-1.92, 4^{\circ}~)$ &  $(~1.92, -176^{\circ}~)$\\
$\ho{16}{12}$ & $(~.91, -11^{\circ}~)$ &  $(~.91, -11^{\circ}~)$
\end{tabular}
\caption{\elabel{tabswap}Swap of representatives for type$_x$ 3 hooks.}
\end{center}
\end{table}

Note that $\be{14}{13}(q) < 0$. According to Lemmas \ref{readtype3} and  
\ref{lembe3pos} this implies that the representative in $\K{14}{13}$
swaps. Indeed, $\be{14}{13}$ changes sign, while
$\al{14}{13}(\tilde{q}) = \al{14}{13}(q) + \pi$. Ratio $\be{16}{12}(q)$
is positive
so here nothing should change and this is reflected in the values
in the table.

\subsection{Managing $x$-type 2 hooks.}

Next lemma shows that in case of a $x$-type 2 hook, 
the representative in the Klein bottle stays the same.

\begin{lemma}
\elabel{readtype2}
Let $q \in \dom(x)$ be valid and accepted. Put
$\tilde{q} = \readoff_x(q) \in \ah$. Suppose that 
$x$-type$(\ho{ik}{ij})=2$.  Then
\beq
	(\,\be{ik}{ij}(\tilde{q}),\,\al{ik}{ij}(\tilde{q})\,) & = & 
		(\,\be{ik}{ij}(q),\,\al{ik}{ij}(q)\,)
\eeq
\end{lemma}

\begin{proof}
Throughout the proof, we label as in Definition \ref{htree}. 
The screen orientations 
$\mathcal{O}_S$ for the different cases are defined in 
Equation \ref{eqscreenos}. Set $k=j_2$ and $i=j_1$.
We distinguish three possibilities for the hinge point $p_{j_1}$.
\begin{lijst}
\item $j_1=i_1$. In this case, $j=i_2$. First we show that 
$\al{i_1j_2}{i_1i_2}(\tilde{q}~) = \al{i_1j_2}{i_1i_2}(q)$.
In order to determine $\alpha_{i_1i_2}(\tilde{q})$, we need
the coordinates
$p_{i_1}(T) = (0,0)$ and 
$p_{i_2}(T)=( \cos \mathcal{O}_T, \sin \mathcal{O}_T)$ in $T=S_{i_1i_2}$. 
It follows that $\alpha_{i_1i_2}(\tilde{q})=\text{arg}(p_{i_2})=
	\mathcal{O}_T$.
For determining $\alpha_{i_1j_2}(\tilde{q})=\alpha_{j_1j_2}(\tilde{q})$ we compute 
in screen~$S$ the coordinates of $p_{j_1}(S)=(0,0)$ and:
        $$p_{j_2}(S) ~ = ~ ( \cos \mathcal{O}_S, \sin \mathcal{O}_S)
	~ = ~  ( \cos \mathcal{O}_T + \al{i_1j_2}{i_1i_2},
		 \sin \mathcal{O}_T + \al{i_1j_2}{i_1i_2}).$$
We get that 
	$\al{i_1j_2}{i_1i_2}(\tilde{q})$
		 $=$  $\mathcal{O}_T + \al{i_1j_2}{i_1i_2}(q) 
			- \mathcal{O}_T$
		 $=$  $\al{i_1j_2}{i_1i_2}(q)$.

Next we determine $\be{i_1j_2}{i_1i_2}(\tilde{q})$. For this purpose
we compute the coordinates of $p_{i_1}$, $p_{i_2}$, $p_{i_2}'$, and
$p_{j_2}$ in $S_{i_1i_2j_2} = T$. We gave the first two above, the
latter two are given by
\beq
	p_{i_2}'(T) & = &  ( \cos \mathcal{O}_T + \al{i_1j_2}{i_1i_2},
                 \sin \mathcal{O}_T + \al{i_1j_2}{i_1i_2}),\\
        p_{j_2}(T)  & = & \be{i_1j_2}{i_1i_2} p_{i_2}'.
\eeq
By substituting those coordinates in Equation \ref{readoffbeta} 
we obtain that $ \be{i_1j_2}{i_1i_2}(\tilde{q}~)  = 
 \be{i_1j_2}{i_1i_2}(q)$, as claimed.
\begin{figure}[!ht]
\begin{center}
\setlength{\unitlength}{1em}
\begin{picture}(20,10)
\put(10,5){\makebox(0,0)[cc]{
        \leavevmode\epsfxsize=20\unitlength\epsfbox{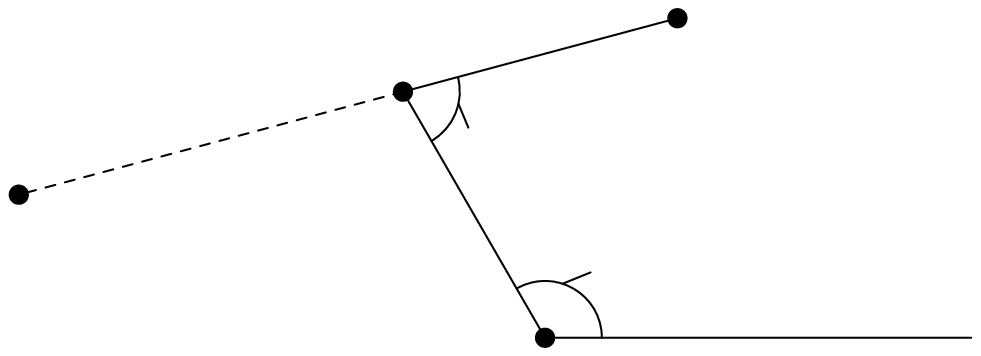}}}
\put(12.7,3){\makebox(0,0)[l]{$\mathcal{O}_T$}}
\put(9.9,6.2){\makebox(0,0)[l]{$\al{i_2j_2}{i_2i_1}(q)$}}
\put(9.5,2){\makebox(0,0)[l]{$i_1$}}
\put(0.3,6.0){\makebox(0,0)[l]{$j_2$}}
\put(8,8){\makebox(0,0)[l]{$i_2$}}
\put(14.2,9){\makebox(0,0)[l]{$i_1'$}}
\end{picture}
\caption{\elabel{jisi2}$\be{i_2i_1}{i_2j_2}(\tilde{q})
                                = \be{i_2i_1}{i_2j_2}(q)$. In the picture,
$\be{i_2i_1}{i_2j_2}(q)<0$.}
\end{center}
\end{figure}
\item $j_1=i_2$. In this case, $j=i_1$. Again we start by showing 
	that $\al{i_2j_2}{i_2i_1}(\tilde{q})
	=\al{i_2j_2}{i_2i_1}(q)$. For an illustration of
      what happens in screen $T=S_{i_1i_2}$, consult Figure \ref{jisi2}.
      As in the case $j_1=i_1$ it holds that
      $\alpha_{i_1i_2}(\tilde{q}) = \mathcal{O}_T$. Therefore,
	$\alpha_{i_2i_1}(\tilde{q})  =  \mathcal{O}_T + \pi$. Angle
      $\alpha_{i_2j_2}(\tilde{q})=\alpha_{j_1j_2}(\tilde{q})$ is determined 
      from $p_{j_1}(S)=(0,0)$ and $p_{j_2}(S)=
		( \cos \mathcal{O}_S, \sin \mathcal{O}_S)$ in $S=S_{j_1j_2} = S$:
screen orientation $\mathcal{O}_S$ is given by
	$\mathcal{O}_S  =  \mathcal{O}_T + \al{j_1j_2}{j_1i_1} + \pi$.
Therefore $\alpha_{j_1j_2}(\tilde{q})  =  \text{arg}(p_{j_2})
			 =  \mathcal{O}_T + \al{j_1j_2}{j_1i_1} + \pi$.
So
	$$ \al{i_2j_2}{i_2i_1}(\tilde{q}) ~ = ~ \alpha_{i_2j_2}(\tilde{q})
		-  \alpha_{i_2i_1}(\tilde{q})
		~ = ~ \mathcal{O}_T + \al{j_1j_2}{j_1i_1} + \pi -\mathcal{O}_T
			+ \pi
		~ = ~ \al{j_1j_2}{j_1i_1}(q).
$$
For determining $\be{i_2j_2}{i_2i_1}(\tilde{q})$, which occurs in the 
$x$-screen $T$, we apply Equations 
\ref{readoffbeta} and \ref{pjprime} directly:
\beq
\be{i_2j_2}{i_2i_1}(\tilde{q}) & = & 
	\frac{ (p_{j_2} - p_{i_2}) \cdot (p_{i_1}' - p_{i_2})}
		{|p_{i_1}' - p_{i_2}|^2}, \\ 
	& = & 
	(p_{j_2} - p_{i_2}) \cdot (p_{i_1}' - p_{i_2}),\\
	& = & ( \be{i_2j_2}{i_2i_1}(q) ( p_{i_1}' - p_{i_2}) + p_{i_2}
		-p_{i_2} ) \cdot (p_{i_1}' - p_{i_2}),\\
	& = & \be{i_2j_2}{i_2i_1}(q) (p_{i_1}' - p_{i_2}) 
		\cdot (p_{i_1}' - p_{i_2}),\\ 
	& = & \be{i_2j_2}{i_2i_1}( q ). 
\eeq
\item $j_1=i_t, t = 3, \dots m$. In this case $j=i_1$.
We start by analyzing $\al{i_tj_2}{i_ti_1}$. For determining 
$\alpha_{i_ti_1}(\tilde{q})$ we compute in $T=S_{i_1i_t}$
coordinates $p_{i_1}(T)=(0,0)$, $p_{i_2}(T)=
	( \cos \mathcal{O}_T, \sin \mathcal{O}_T )$, and 
\beq        
p_{i_t}(T) & = & \be{i_1i_t}{i_1i_2}(q)\rot_{\al{i_1i_t}{i_1i_2}(q)}
			(p_{i_2}),\\
                & = & \be{i_1i_t}{i_1i_2}(q)( 
	\cos( \mathcal{O}_T + \al{i_1i_t}{i_1i_2}(q)),
	\sin( \mathcal{O}_T + \al{i_1i_t}{i_1i_2}(q)) ).
\eeq
Therefore,
\begin{equation*}
  \alpha_{i_ti_1}(\tilde{q})\quad=\quad
         \begin{cases}
  \mathcal{O}_T +  \al{i_1i_t}{i_1i_2}(q), 
		& \be{i_1i_t}{i_1i_2}(q) ~>~ 0, \\
   \mathcal{O}_T +  \al{i_1i_t}{i_1i_2}(q) + \pi,
	 & \be{i_1i_t}{i_1i_2}(q) ~<~ 0.
         \end{cases}
  \end{equation*}
We use the coordinates,  $p_{i_t}(S)=p_{j_1}(S)=(0,0)$ and 
$p_{j_2}(S) = ( \cos \mathcal{O}_S, \sin \mathcal{O}_S)$ in 
	$S=S_{i_tj_2}=S_{j_1j_2}$ for determining 
$\alpha_{i_tj_2}(\tilde{q})$.
From Equation \ref{eqscreenos} it follows that
\begin{equation*}
  \mathcal{O}_S\quad=\quad
         \begin{cases}
         \mathcal{O}_T + \al{i_1i_t}{i_1i_2}(q)
                + \al{j_1j_2}{j_1i_1}(q) + \pi,
                & \be{i_1i_t}{i_1i_2}(q) ~>~ 0 ,\\
	\mathcal{O}_T + \al{i_1i_t}{i_1i_2}(q)
                + \al{j_1j_2}{j_1i_1}(q), 
                & \be{i_1i_t}{i_1i_2}(q) ~<~ 0.
         \end{cases}
  \end{equation*}
After subtracting it follows that 
$\al{i_tj_2}{i_ti_1}(\tilde{q})=\al{i_tj_2}{i_ti_1}(q)$. 
We only have to show that $\be{i_tj_2}{i_ti_1}(\tilde{q})
	=  \be{i_tj_2}{i_ti_1}(q)$. We do so by filling in
coordinates in  Equation \ref{readoffbeta}.  The relevant
coordinates in screen $S_{i_1i_tj_2}=T$ are as follows. 
\beq
	p_{i_1} & = & (0,0) ,\\
	p_{i_2} & = & ( \cos \mathcal{O}_T, \sin \mathcal{O}_T),\\
	p_{i_t} & = & \be{i_1i_t}{i_1i_2}(q) 
		\rot_{\al{i_1i_t}{i_1i_2}}( p_{i_2} - p_{i_1} ) + p_{i_1},\\
                & = & \be{i_1i_t}{i_1i_2}(q)
		\rot_{\al{i_1i_t}{i_1i_2}}( p_{i_2} ) ,\\
	p_{i_1}' & = & \rot_{\al{i_tj_2}{i_ti_1}}( p_{i_1} - p_{i_t} )  
			+ p_{i_t},\\	
		& = & \rot_{\al{i_tj_2}{i_ti_1}}( - p_{i_t} ) + p_{i_t},\\
	p_{j_2} & = & \be{i_tj_2}{i_ti_1} \rot_{\al{i_tj_2}{i_ti_1}}
			( p_{i_1}- p_{i_t} ) + p_{i_t},\\
		& = & \be{i_tj_2}{i_ti_1} \rot_{\al{i_tj_2}{i_ti_1}}
                        ( - p_{i_t} ) + p_{i_t}.
\eeq
Therefore
\begin{equation*}
	\be{i_tj_2}{i_ti_1}(\tilde{q})~  = ~
		\frac{ \be{i_tj_2}{i_ti_1}(q)
			 \rot_{\al{i_tj_2}{i_ti_1}}(-p_{i_t})
			\cdot \rot_{\al{i_tj_2}{i_ti_1}}(-p_{i_t}) }
		{ \rot_{\al{i_tj_2}{i_ti_1}}(-p_{i_t}) \cdot 
			\rot_{\al{i_tj_2}{i_ti_1}}(-p_{i_t}) }
	 ~=~  \be{i_tj_2}{i_ti_1}(q).
\qedhere
\end{equation*}
\end{lijst}
\end{proof}

\begin{eexample}
\elabel{exconsist2}
In our running example we have two type$_x 2$ hooks,
$\ho{12}{13}$ and $\ho{35}{31}$.
The reader is invited to check that the representatives computed
for $\readoff_x(q)$ in $\K{12}{13}$
and $\K{35}{31}$ that are the same as those
for $q$ itself.
See in particular Example \ref{exfill} and Example \ref{exreadoff}.
\end{eexample}

\subsection{Consistency theorem for hooked tree elements.}

\begin{theorem}
\elabel{thconsist}
Fix $x \in \XAH$. Let $q \in \dom(x)$ be valid and accepted. Then
\begin{displaymath}
\begin{array}{ccccccc}
\dom(x) & \rightarrow & \scr(x) & \rightarrow & \ah & \rightarrow & \dom(x) ,\\
&\\
q & \mapsto & \draw_x(q) & \mapsto & \readoff_x(q) & 
	\stackrel{tv_x}{\mapsto} & \tilde{q}.
\end{array}
\end{displaymath}
is the identity map.
\end{theorem}

\begin{proof}
We analyze the three types of data elements present in $ht(x)$,
see Definition \ref{htree}.
The type$_x$ 3 hooks are dealt with in Lemma \ref{readtype3} and the
type$_x$ 2 hooks in Lemma \ref{readtype2}.  We are left with the
type$_x$ 2.a angle between the top screen and its second maximal cluster. 
In Definition \ref{htree} this angle is labeled $\alpha_{1j_2}$.
The separating screen of $p_1$ and $p_{j_2}$ is the top screen $S_{[n]}$,
so we have to read off in this screen. The coordinates of $p_1$ and
$p_{j_2}$ in $S_{[n]}$ are  
	$p_1 = (0,0)$ and $p_{j_2} = ( \cos \alpha_{1j_2}(q),
	\sin \alpha_{1j_2}(q))$. 
The read off value $\alpha_{1j_2}(\tilde{q})$ equals 
$\text{arg}(p_{j_2}) = \alpha_{1j_2}(q)$. 
This completes the proof of the theorem.
\end{proof}

\section{$\XAH$ is smooth. \elabel{ssmooth}}

\begin{lemma}
\elabel{lemconsconf}
Fix $x \in \XAH$. Let $q \in \XAH$  be such that:
\begin{lijst}
\item $tv_x(q)$ is valid and accepted;
\item $q = \psi_{\ah}([c_q])$ for some $[c_q] \in \conf$.
\end{lijst}
Then 
	$q  =  \readoff_x( tv_x(q))$.
\end{lemma}

\begin{proof}
The consistency theorem, Theorem \ref{thconsist},
proves the claim for those $q$ coordinates that belong to
$\dom(x)$.  In the rest of the proof, we analyze 
$\alpha_{ij} ~\text{mod}~ \pi$ and $\ho{ij}{ik}$.
Write $\tilde{q} := \readoff_x( tv_x(q))$. 

First we show that $\alpha_{ab}(q) = \alpha_{ab}(\tilde{q})$ for
$a,b \in 1, \dots,n$, with  $a \neq b$.
Definition \ref{readoff} tells that we have to determine 
$\alpha_{ab}(q)$ in the screen $S_{ab}$. Let $j_1$ and $j_2$ be the first label
in the first cluster and the first label in the second cluster of screen $S_{ab}$.
Then $p_{j_1}$ and $p_{j_2}$ are the sites drawn first and second in 
Construction \ref{draw}. 
If $a \neq j_1,j_2$, then $p_a$ was 
constructed out of the $x$-leg $p_{j_1j_2}$ by repeatedly applying 
hooks from $ht(x)$ in the order indicated by the
hooked path $L_x(p_a)$. The same holds for $p_b$. 
For hooks $\ho{ik}{ij}$ that are referred to in $ht(x)$ one has
$\ho{ik}{ij}(q)  =  \ho{ik}{ij}(\tilde{q})$.  
As a consequence, the positions of $p_a$ and $p_b$ relative to
$p_{j_1}$ and $p_{j_2}$ are the same in $S_{ab}$ as in 
any configuration $c_q \in [c_q]$. 

The angle of the directed line passing first through $p_{j_1}$ and then
through $p_{j_2}$ in $S_{ab}$ is exactly the screen orientation 
$\mathcal{O}_{S_{ab}} = \mathcal{O}_{S_{j_1j_2}}$.
Writing out any screen orientation according to the recursion
in Equation \ref{eqscreenos} gives a
list of data elements corresponding  exactly to the hooked path
$L_x(p_b)$.
As all angles occurring in the screen orientation definition are 
coordinates of $\dom(x)$, Theorem \ref{thconsist} applies again. 
At any type~2 angle $\al{ik}{ij}$ however a $\pi$ switch can
occur: the corresponding type~2 ratio $\be{ik}{ij}$ is not used in the 
definition. Moreover, representatives might swap in Step~1  of 
Construction \ref{draw}.  It follows that
\beq
        \mathcal{O}_{S_{j_1j_2}}  ~=~
                \begin{cases} \alpha_{j_1j_2}(q), \\
                              \alpha_{j_1j_2}(q) + \pi
                \end{cases}
\eeq

To show that $\ho{ik}{ij}(q) =  \ho{ik}{ij}(\tilde{q})$ is easy now.
The proof consists of three observations:
\begin{lijst}
\item Any $q$-filled $x$-screen $S$ is just a copy of 
      the original configuration, due  to the Consistency Theorem. But,
      it also satisfies the following:
      points not belonging to cluster $C(S)$ 
       corresponding to screen $S$ are omitted;
      points in $S$ are scaled such that the distance 
	from $p_{j_1}$ to $p_{j_2}$ equals one;
      points in $S$ are possibly rotated over 
	$\pi$ according to Construction \ref{draw}, Step~1. 
\item By Definition \ref{readoff}, $\al{ik}{ij}(\tilde{q}) = \alpha_{ik}(\tilde{q})
	- \alpha_{ij}(\tilde{q})$. From the proof of part 1.\ of this lemma 
      it follows that the angle $\al{ik}{ij}(\tilde{q})$ equals
       $\al{ik}{ij}(q)$ up to $\pi$. 
\item In case $\al{ik}{ij}(\tilde{q}) = \al{ik}{ij}(q) + \pi$ 
      this $\pi$-switch is automatically corrected by the sign of $\be{ik}{ij}$,
      as  $\al{ik}{ij}(\tilde{q})$ is used in Equation \ref{pjprime}.
\qedhere
\end{lijst} 
\end{proof}

\begin{lemma} 
\elabel{lemgraph}
Fix $x \in \XAH$. Let $U \subset \ah$ be the set
of points in $\ah$ such that $tv_x(u)$ is valid and accepted for all
$u \in U$. Then
\begin{eqnarray}
\elabel{eqgraph}
   \XAH \cap U & = & \text{graph}(ctv_x(\readoff_x)).
\end{eqnarray}
\end{lemma}

\begin{proof} We prove the two inclusions.
First let $y \in \XAH \cap U$. Then its projection 
$q_y := tv_x(y) \in \dom(x)$ is valid and accepted
with respect to $x$. 
There always exist points $z \in \XAH \cap U$ of the form
$z= \psi_{\ah}([c_z])$, for $[c_z] \in \conf$, that are arbitrary close to $y$.
For such points $z$ expression (\ref{eqgraph}) holds by Lemma \ref{lemconsconf}. 
So these points can be expressed as the image  of a configuration
consisting of non-coinciding points.  Therefore these points $z$ belong to 
the graph of $\text{graph}(ctv_x(\readoff_x))$.
This implies that  $y$ itself belongs to the graph 
as well, as $y$ can be  constructed as limit point of points
of the form $z = \psi_{\ah}([c_z])$, for $[c_z] \in \conf$ in $\XAH$.

Secondly, suppose we start with a point $g_q$ on the graph. This implies
that $q=tv_x(g_q) \in \dom(x)$ is valid and accepted. 
As both conditions (valid, accepted) are open conditions, there exists
$\epsilon$ with $ 0  <  \epsilon  <<  1$,
such that if we replace all $q$-ratios that vanish by $\epsilon$ that
then $q(\epsilon)$ is still valid and accepted.  
But this means that $\draw_x(q(\epsilon))$ produces a filling of
the $x$-screens in which all sites in the top screen are distinct. 
Writing down the coordinates gives a  class
$[c_{q(\epsilon)}] \in \conf$. Then
\beq
	\psi_{\ah}([c_{q(\epsilon)}]) & = &  
		(q(\epsilon); ctv_x(\readoff_x(q(\epsilon))) ).
\eeq
The claim follows by observing that $\XAH$ is the closure of 
$\psi_{\ah}(\conf)$ in $\ah$.
\end{proof}

\begin{theorem}
\elabel{tsubmani}
$\XAH$ is a submanifold of $\ah$.
\end{theorem}

\begin{proof}
Being a submanifold of a manifold  can be determined locally, 
see \cite{BG}, section 6.2.  Lemma \ref{lemgraph} shows that 
$\XAH$ can locally be written as
the graph of a function, and this is one of the characterizations of submanifold.
See \cite{BG}, Theorem 2.1.2.\, notably part (iv).
\end{proof}

\section{Voronoi diagrams in the $x$-screens.}
\elabel{sclick}

So far in this chapter we have introduced the space $\XAH$ in
order to model configurations of coinciding sites in the plane.
We have seen that we can associate with every $x \in \XAH$
a family of $x$-screens. Moreover, for almost all $q \in \XAH$,
we can fill the $x$-screens with $q$-sites  so
that  for any two $q$ sites $p_i(q)$ and $p_j(q)$ there
exists an $x$-screen where $p_i(q)$ and $p_j(q)$ are distinct.

In this section we add Voronoi diagrams to the $x$-screens. 
As in any $x$-screen at least two sites are distinct, we could
apply the ordinary definition of Voronoi diagram. But there is
an ambiguity in the screen model, as we 
take angles $\alpha_{ij}$ between pairs of points in $\R/\pi\Z$
and allow negative ratios $\be{ij}{ik}$.

\begin{eexample} Recall Example \ref{exfill}. This example demonstrates
that $q$-clusters can get rotated over $\pi$ in distinct $x$-screens.
\end{eexample}

\subsection{Connection with Kontsevich-Soibelman space.} 

We nail down the clusters in such a way that they cannot get swapped
over $\pi$ in distinct $x$-screens. We do this by
taking angles between pairs of points in $\R/2 \pi \Z$ and
allowing ratios $\be{ik}{ij} \in [0,\infty]$ only.
That is, we consider the map
\bmap
        \psi_{\edah}:& \conf  & \rightarrow &
                (\R / 2 \pi \Z)^{\binom{n}{2}} \times
                         ([0, \infty] \times \R/ 2 \pi \Z)^{6 \binom{6}{3}}, \\
	& {}[(p_1, \dots, p_n)] & \mapsto 
	& ( (\alpha_{ij})_{1 \leq i < j \leq n}, 
			(\be{ik}{ij}, \al{ik}{ij})).
\emap
where $i$, $j$, and $k$ are pairwise distinct indices. Denote
by $\edah$ the product $(\R / 2 \pi \Z)^{\binom{n}{2}} \times
     ([0, \infty] \times \R/ 2 \pi \Z)^{6 \binom{6}{3}}$.
Recall the manifold with corners $\FMt(n)$ introduced by
Kontsevich and Soibelman, see \cite{KS} and 
described in Section \ref{sKS}.

\begin{proposition}
Let $\text{XEDAH}_n$ denote the closure of $\psi_{\edah}( \conf )$
in $\edah$. Then
\beq
\text{XEDAH}_n & \cong & \FMt(n).
\eeq
\end{proposition}

\begin{proof}
$\FMt$ is defined as the closure of the image of $\conf$ 
under the map
\beq
	[\{p_1, \dots, p_m \}] & \mapsto &
		((\alpha_{ij})_{1 \leq i < j \leq n}, \be{ik}{ij}),
\eeq
where $i$, $j$, and $k$ are pairwise distinct indices and $\be{ik}{ij} =
	\frac{|p_i-p_k|}{|p_i-p_j|}$. As $\al{ik}{ij} = \alpha_{ik} - \alpha_{ij}$,
the claim follows directly.
\end{proof}

\begin{figure}[!ht]
 \begin{center}
 \setlength{\unitlength}{1cm}
 \begin{picture}(10,2.5)
 \put(5,1.5){\makebox(0,0)[cc]{
         \leavevmode\epsfxsize=6.67\unitlength\epsfbox{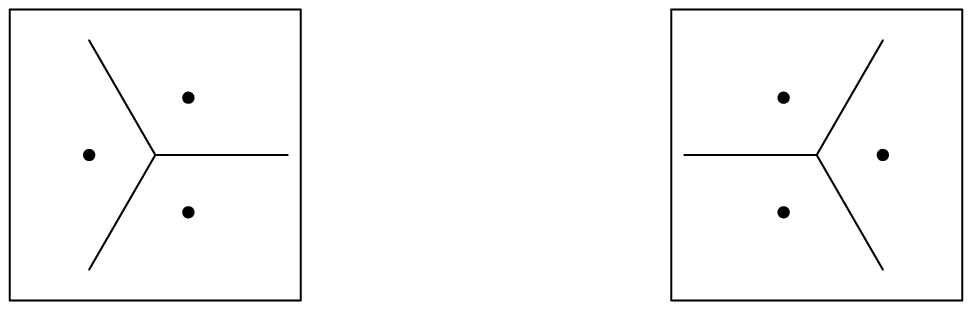}}}
 \put(2.1,1.5){\makebox(0,0)[l]{$1$}}
 \put(3.1,.8){\makebox(0,0)[l]{$2$}}
 \put(3.1,2.2){\makebox(0,0)[l]{$3$}}
 \put(7.9,1.5){\makebox(0,0)[l]{$1$}}
 \put(6.8,2.2){\makebox(0,0)[l]{$2$}}
 \put(6.8,.8){\makebox(0,0)[l]{$3$}}
 \end{picture}
 \caption{\elabel{fcorner}The Voronoi diagram of $p_1(t)$, $p_2(t)$
 and $p_3(t)$. On the left for negative $t$, on the right for
 positive $t$.}
 \end{center}
 \end{figure}

\begin{eexample}
Consider $p_1(t)=t(1,0)$, $p_2(t)= t ( -\frac{1}{2}, \frac{1}{2} \sqrt{3})$,
and $p_3(t)=t(-\frac{1}{2}, -\frac{1}{2} \sqrt{3})$. We compare the
angles associated to $p_1$, $p_2$, and $p_3$ in
${\mbox {\sl XAH}[3]}$ and $\FMt(3)$ for positive and
negative $t$. 
First we consider ${\mbox {\sl XAH}[3]}$.
The undirected angles $\overline{\alpha_{12}}= -\frac{\pi}{3}$,
$\overline{\alpha_{13}}= \frac{\pi}{3}$ and 
$\overline{\alpha_{23}}= \frac{\pi}{2}$ do not change at $t=0$.
In $\FMt(3)$, however, angles $\alpha_{12}$, $\alpha_{13}$ and $\alpha_{23}$
change by $\pi$ at $t=0$.
Think of this as moving into another corner in the manifold
with corners. In Figure \ref{fcorner} the Voronoi diagram of
$p_1(t)$, $p_2(t)$ and $p_3(t)$ is shown. On the left for
any negative value of $t$, on the right for any positive value of $t$.
 Note the jump in the diagram at $t=0$. It demonstrates the need of the
directed  angles in the screen model if we want to use it for
displaying Voronoi diagrams.
\elabel{exswapt}
\end{eexample}

\begin{figure}[!ht]
 \begin{center}
 \setlength{\unitlength}{1.12cm}
 \begin{picture}(10,5.4)
 \put(5,2.5){\makebox(0,0)[cc]{
         \leavevmode\epsfxsize=8\unitlength\epsfbox{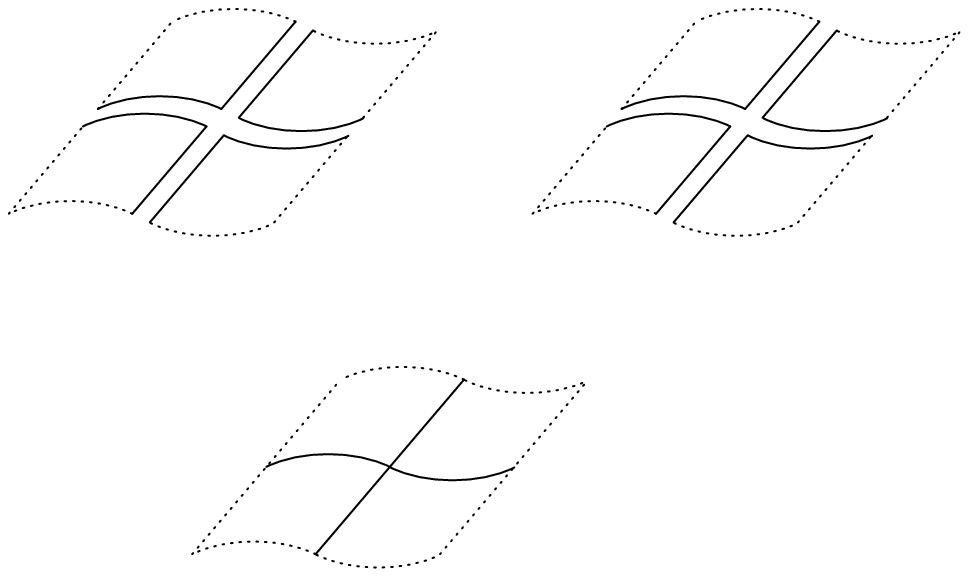}}}
 \put(3.3,0){\makebox(0,0)[l]{$\beta=0$}}
 \put(5.6,1.1){\makebox(0,0)[l]{$\beta'=0$}}
 \put(1.3,0.5){\makebox(0,0)[l]{$\overline{\alpha}=\overline{\alpha}_0$}}
 \put(3.1,5.2){\makebox(0,0)[l]{$\alpha=\alpha_0$}}
 \put(7.0,5.2){\makebox(0,0)[l]{$\alpha=\alpha_0 + \pi$}}
 \put(2.5,2.5){\makebox(0,0)[l]{$\searrow$}}
 \put(5.7,2.5){\makebox(0,0)[l]{$\swarrow$}}
 \put(9.5,2.5){\makebox(0,0)[l]{$\downarrow f$}}
 \put(9.5,1.1){\makebox(0,0)[l]{$\xah$}}
 \put(9.5,3.8){\makebox(0,0)[l]{$\FMt(n)$}}
\end{picture}
 \caption{\elabel{fourcorners}The manifold with corners $\FMt(n)$ above 
	the smooth manifold $\xah$, with $\overline{\alpha}$,
$\beta$, and $\beta '$ in $\dom(x)$. That slice is presented where
all coordinates are fixed except for $\beta$ and $\beta'$.}
 \end{center}
 \end{figure}

\begin{theorem}
\elabel{tfiber}
See also Figure \ref{fourcorners}.
Consider the map
\bmap
	f: & \FMt(n) & \rightarrow & \XAH, \\
& ((\alpha_{ij})_{1 \leq i<j\leq n}, (\be{ik}{ij}, \al{ik}{ij})) & \mapsto &
	( (\overline{\alpha_{ij}}_{1 \leq i < j \leq n}),
		[(\be{ik}{ij}, \al{ik}{ij})]),
\emap
where $[.]$ denotes the class of  $[(\be{ik}{ij}, \al{ik}{ij})]$
in the Klein bottle $\K{ik}{ij}$ with respect to the equivalence
relation $\sim_k$. Fix $x \in \XAH$.
Let $m$ be the
number of $\be{ik}{ij}(x)$ in $\dom(x)$ such that $\be{ik}{ij}(x)=0$.
Then
\beq
|f^{-1}(x)| & = & 2^{m+1}.
\eeq
\end{theorem}

\begin{proof}
Fix a class $[(\be{ik}{ij}(x), \al{ik}{ij}(x))]$ such that
$\be{ik}{ij}(x)=0$. Any such class has two elements,
$(0, \al{ik}{ij})$ and $(0, \al{ik}{ij} + \pi)$.
By  Lemma \ref{lemgraph}, $\XAH$ can be written as the graph of a function
on $\dom(x)$. 
Therefore those $\be{ik}{ij}$ components that are not in $\dom(x)$
are a function of those components that are in $\dom(x)$. Moreover,
$\readoff_x(x)$ even determines the representative for those
$\be{ik}{ij}$ components that are not in $\dom(x)$.
For any $\be{ik}{ij} \in \dom(x)$, the two representatives
correspond to two distinct pieces of the boundary in $\FMt(n)$. 
We conclude that we can use the signvector of 
those $\be{ik}{ij}(x) \in \dom(x)$ with $\be{ik}{ij}(x)=0$ 
as a coordinate system pointing
to the distinct corners in the fiber of the map $f$.
The additional $1$ in $2^{m+1}$ comes from 
the two representatives of the unique $\alpha_{ij}$-component
in $\dom(x)$.
\end{proof}

\begin{eexample}
\elabel{exfiber} 
Consider $x\in \xaht$ given by 
\begin{displaymath}
\begin{array}{cccccc}
\overline{\alpha}_{12} & \overline{\alpha}_{13} &  \overline{\alpha}_{23} &
	(\be{13}{12}, \al{13}{12}) & 
	(\be{23}{21}, \al{23}{21}) & 
	(\be{32}{31}, \al{32}{31})  \\[.1cm]
\frac{\pi}{6} & \frac{\pi}{6} & \frac{\pi}{4} &
(1,0) & (0, -\frac{5 \pi}{12}) & (0,\frac{\pi}{12})
\end{array}
\end{displaymath}
Theorem \ref{mathcalC} associates to $x$ the nest $\langle\{1,2,3\},\{2,3\}\rangle$.
The factor $\domt(x)$ consists of the two components 
$\overline{\alpha}_{12}(x)$ and 
$\ho{23}{21}(x)$. The standard form of $x$ restricted to $\domt(x)$ is given by 
the representatives $\alpha_{12}(x) = \frac{\pi}{6}$ and 
$\al{23}{21}(x) = -\frac{5 \pi}{12}$,  compare Definition \ref{dstandard}.
The fiber $f^{-1}(x)$  consists exactly of those
points in $\FMt(3)$ that are mapped on the same standard form
in $\XAH$. We obtain these points by analyzing the four distinct
representatives for $\alpha_{12}$ and $\al{23}{21}$.
Choose for example $\alpha_{12} = -\frac{5 \pi}{6}$ and
$\al{23}{21}=-\frac{5 \pi}{12}$. It follows from the nest structure
that  $\alpha_{13}=\alpha_{12}$, while 
$\alpha_{23} = \al{23}{21}+\alpha_{21} = 
-\frac{5 \pi}{6} + \frac{\pi}{6} = -\frac{\pi}{4}$.
The four points $f_1, \dots, f_4$ in $\FMt(3)$ obtained in this way with 
coordinates   $(\alpha_{12}, \alpha_{13}, \alpha_{23}, \be{13}{12},
		\be{21}{23}, \be{31}{32})$,
are given by 
\begin{displaymath}
\begin{array}{lll}
f_1 &=& (\frac{\pi}{6}, \frac{\pi}{6}, \frac{3 \pi}{4}, 1, 0, 0),\\[.1cm]
f_2 &=& (\frac{\pi}{6}, \frac{\pi}{6}, -\frac{ \pi}{4}, 1, 0, 0),\\[.1cm]
f_3 &=& (-\frac{5 \pi}{6}, -\frac{5 \pi}{6}, -\frac{\pi}{4}, 1, 0, 0), \\[.1cm]
f_4 &=& (-\frac{5 \pi}{6}, -\frac{5\pi}{6}, \frac{3 \pi}{4}, 1, 0, 0). \\
\end{array}
\end{displaymath}
\end{eexample}

\begin{remark}
In Chapter \ref{changles} we have encountered
a similar situation in comparing the spaces $\cuat$ and $\cdat$ of 
undirected and directed angles on three points. 
Recall particularly Examples \ref{exposcda3} and \ref{exua3}. 
\end{remark}

We can define and fill $x$-sreens
in $\FMt(n)$ as well. Fix a point $x \in \FMt(n)$. The set
of vanishing ratios $\be{ik}{ij}(x)$ for $x$ defines a family
of $x$-clusters and corresponding $x$-screens by applying the methods
of Section \ref{nestsscreens}. A stripped version
of Construction \ref{draw} 
fills the $x$-screens, given some $q \in \dom(x)$. In short,
we discuss the adaptations in Construction \ref{draw}.\\
{\bf Ad Step 1. } We do not have representatives in $\edah$, so we 
skip Step~1.\\
{\bf Ad Step 2. } The definition of $\mathcal{O}_S$ still holds,
but negative ratios $\be{ik}{ij}$ do not occur in $\edah$.\\
{\bf Ad Step 3.} In Step~3 nothing changes.

\subsection{Adding Voronoi diagrams in the $x$-screens.}

In this section we define a Voronoi diagram $V_{FM}(x)$
for $x \in \FMt(n)$.

\begin{definition}
\elabel{dvorfm}
Let $x \in \FMt(n)$. The Voronoi diagram $V_{FM}(x)$ 
is defined in terms of the $x$-screens filled by
$\text{draw}_x(x)$. Fix an $x$-screen $T$ filled
by $\text{draw}_x(x)$. Let $C$ be the cluster corresponding
to $T$. Its maximal subclusters are denoted by
$C_1, \dots, C_m$. Suppose that the coordinates of $C_1, \dots, C_m$ 
in $T$ are given by $t_1, \dots, t_m$.
\begin{description}
\item[Initialization.] $V(T)$ is defined as the Voronoi diagram
$V(\{t_1, \dots, t_m\})$.
\item[Completion.] This is a recursive step from the screens displaying
single sites up to the top screen. Assume that the maximal 
subscreens $S_1, \dots, S_m$ of $T$ have been completed. 
Denote the completed diagram in $S_i$ by $V_{FM}(S_i)$. After 
the initialization step, any cluster $C_i$ lives in a Voronoi cell
$V(t_i)$ in $T$. Define
\beq
	V_{FM}(t_i) & = & V(t_i) \cap V_{FM}(S_i),
\eeq
and
\beq
	V_{FM}(T) & = & \bigcup_{i=1\dots m}V_{FM}(t_i).
\eeq
\end{description}
The Voronoi diagram $V_{FM}(x)$ is the rooted tree of all $V_{FM}(T)$,
where $T$ runs over all $x$-screens. The ordering on $V_{FM}(T)$
is inherited from the ordering on the $x$-screens. 
\end{definition}

\begin{figure}[!ht]
\begin{center}
\setlength{\unitlength}{1cm}
\begin{picture}(12,6)
\put(3,5){\makebox(0,0)[cc]{
        \leavevmode\epsfxsize=2\unitlength\epsfbox{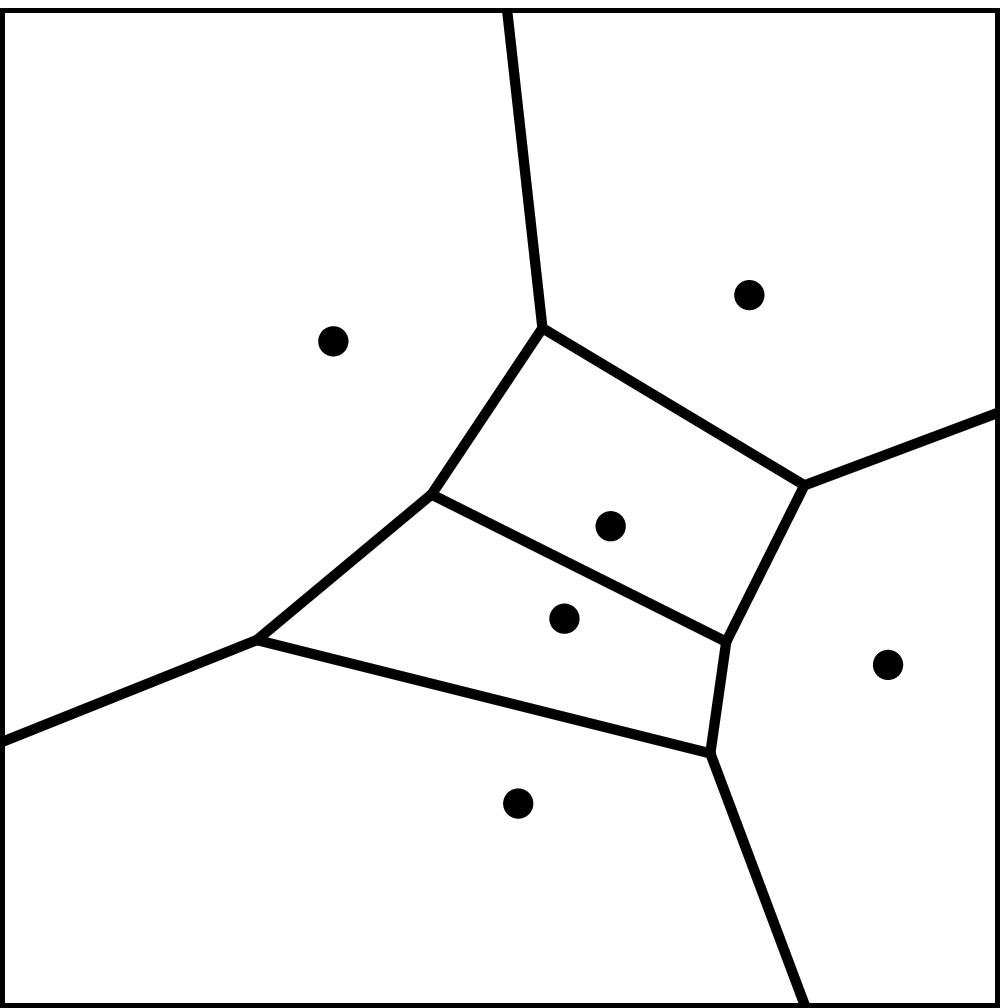}}}
\put(1.5,1){\makebox(0,0)[cc]{
        \leavevmode\epsfxsize=2\unitlength\epsfbox{fmplug2.eps}}}
\put(4.5,1){\makebox(0,0)[cc]{
        \leavevmode\epsfxsize=2\unitlength\epsfbox{fmplug3.eps}}}
\put(9,5){\makebox(0,0)[cc]{
        \leavevmode\epsfxsize=2\unitlength\epsfbox{fmplug4.eps}}}
\put(7.5,1){\makebox(0,0)[cc]{
        \leavevmode\epsfxsize=2\unitlength\epsfbox{fmplug2.eps}}}
\put(10.5,1){\makebox(0,0)[cc]{
        \leavevmode\epsfxsize=2\unitlength\epsfbox{fmplug3.eps}}}
%
%
\put(2.6,5.5){\makebox(0,0)[c,c]{1}}
\put(3.1,4.2){\makebox(0,0)[c,c]{2}}
\put(3.6,5.7){\makebox(0,0)[c,c]{5}}
\put(3.8,4.4){\makebox(0,0)[c,c]{6}}
\put(1.45,.4){\makebox(0,0)[c,c]{3}}
\put(.9,1.3){\makebox(0,0)[c,c]{7}}
\put(2.1,1.3){\makebox(0,0)[c,c]{12}}
\put(4.8,.3){\makebox(0,0)[c,c]{4}}
\put(3.8,.6){\makebox(0,0)[c,c]{8}}
\put(3.9,1.7){\makebox(0,0)[c,c]{9}}
\put(5.1,.8){\makebox(0,0)[c,c]{10}}
\put(5,1.8){\makebox(0,0)[c,c]{11}}
\put(5,1.8){\makebox(0,0)[c,c]{11}}
%
%
\put(8.6,5.5){\makebox(0,0)[c,c]{1}}
\put(9.1,4.2){\makebox(0,0)[c,c]{2}}
\put(9.6,5.7){\makebox(0,0)[c,c]{5}}
\put(9.8,4.4){\makebox(0,0)[c,c]{6}}
\put(7.45,.4){\makebox(0,0)[c,c]{3}}
\put(6.9,1.3){\makebox(0,0)[c,c]{7}}
\put(8.1,1.3){\makebox(0,0)[c,c]{12}}
\put(10.8,.3){\makebox(0,0)[c,c]{4}}
\put(9.8,.6){\makebox(0,0)[c,c]{8}}
\put(9.9,1.7){\makebox(0,0)[c,c]{9}}
\put(11.1,.8){\makebox(0,0)[c,c]{10}}
\put(11,1.8){\makebox(0,0)[c,c]{11}}
\put(5,1.8){\makebox(0,0)[c,c]{11}}
\unitlength=1mm
\linethickness{.3mm}
\dottedline{1}(17,20)(32,48)
\dottedline{1}(43,20)(33,49)
\dottedline{1}(77,20)(92,48)
\dottedline{1}(103,20)(93,49)
\end{picture}
 \caption{\elabel{fmoreplug}Constructing $V_{FM}(x)$: on the left, initialization;
on the right, completion.}
\end{center}
\end{figure}

\begin{eexample}
\elabel{exmoreplug}
Recall Example \ref{explug}. In this example we consider an ordered set of sites
in $\R[t]\times \R[t]$ given by
\beq
 S(t) & = & ((-6, 4), (-2, -6), (-1, -2 - 3 t), (0, -3 t), (3,
    5), (6, -3),\\
      & & (-1 - 3 t, -2 + t),
      (-2 t, -2 t), (-2 t, 2 t), (2 t, 0), (2 t, 2 t), (-1 + 2 t, -2)).
\eeq
For any pair of sites $p_i(t)$ and $p_j(t)$ we have defined an
angle or direction $\alpha_{ij}(t)$ at $t=0$. In a similar
way we can define a ratio 
\beq
	\be{ik}{ij}(0) & = &  \lim_{t \downarrow 0} 
		\frac{|p_i(t)-p_k(t)|}{|p_i(t)-p_j(t)|}.
\eeq 
This  establishes a map $\chi$  from a
set $S(t)$ of $n$ sites in $\R[t]\times\R[t]$ to
a point $x := \chi(S(0)) \in \FMt$. We construct
the Voronoi diagram $V_{FM}(x)$ by applying Definition \ref{dvorfm}. 
The nest of $x$ is given
by $\langle\{3,7,12\},\{4,8,9,10,11\}\rangle$. Figure
\ref{fmoreplug} shows on the left the three filled screens 
together with the Voronoi diagrams added in the initialization step.
On the right the diagram $V_{FM}(x)$ is shown. It is constructed 
by adding the two diagrams 
corresponding to clusters $3,7,12$ and $4,8,9,10,11$ in 
the cells of those clusters in the top screen.
\end{eexample}

\begin{eexample}
\elabel{exdepth2}
We analyze an example where the maximal depth of a site equals $2$. 
We start with the ordered set of sites in $\R[t] \times \R[t]$ given
by $S(t)=\{p_1, p_2, p_3, p_4 \}$, where
\begin{displaymath}
\begin{array}{ll}
p_1 ~= & (0,0),\\[.1cm]
p_2 ~= & t ( \cos \frac{\pi}{3}, \sin \frac{\pi}{3} ),\\[.1cm]
p_3 ~= & p_2 + t^2 ( \cos \frac{\pi}{6}, \sin \frac{\pi}{6} ),\\[.1cm]
p_4 ~= & p_3 + t^3 ( \cos \frac{\pi}{2}, \sin \frac{\pi}{2} ).
\end{array}
\end{displaymath}
As in Example \ref{exmoreplug}, one can construct a point 
$y = \chi(S(0)) \in \FMt$. The nest of $y$ is given by
$\langle 234, 34 \rangle$. In Figure \ref{fmoredepth} the tree
of filled screens for $y$ is shown: on the left after the initialization step
but before the completion step
of Definition \ref{dvorfm}; on the right after the completion
step. In this example one can see clearly how the Voronoi
diagram $V_{FM}(T)$ of the top screen $T$ is built up recursively
out of the screens corresponding with the subclusters of
$\{1,2,3,4\}$.
\end{eexample}

\begin{figure}[!ht]
\begin{center}
\setlength{\unitlength}{1cm}
\begin{picture}(5.4,8.6)
\put(2.7,4.3){\makebox(0,0)[cc]{
        \leavevmode\epsfxsize=5.4\unitlength\epsfbox{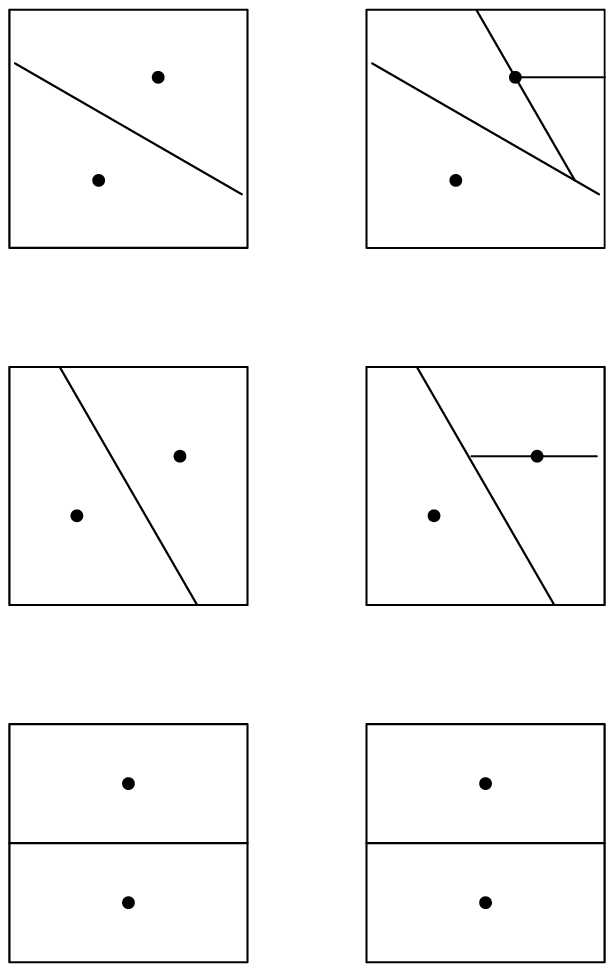}}}
\put(0.83,6.9){\makebox(0,0)[c,c]{1}}
\put(0.9,3.8){\makebox(0,0)[c,c]{2}}
\put(1.4,.8){\makebox(0,0)[c,c]{3}}
\put(1.4,1.8){\makebox(0,0)[c,c]{4}}
\put(3.9,6.9){\makebox(0,0)[c,c]{1}}
\put(3.95,3.8){\makebox(0,0)[c,c]{2}}
\put(4.48,.8){\makebox(0,0)[c,c]{3}}
\put(4.48,1.8){\makebox(0,0)[c,c]{4}}
\unitlength=1mm
\linethickness{.3mm}
\dottedline{1}(12,54)(14.5,77.2)
\dottedline{1}(43,54)(45.5,78)
\dottedline{1}(12,23)(16.5,44.5)
\dottedline{1}(43,23)(47.3,44.5)
\end{picture}
 \caption{\elabel{fmoredepth}Constructing $V_{FM}(y)$: on the left,
initialization; on the right, completion.}
\end{center}
\end{figure}

\begin{remark}
Note the compatibility with the construction of limit 
Voronoi diagram for sites in $\R[t]\times \R[t]$ as constructed
in Chapter~3. Compare also Examples \ref{exmoreplug} and \ref{exdepth2}.
\end{remark}

\section{Conclusion.}
\elabel{sfmconclusion}

\subsection{An easy model to remember.}

In this chapter we have presented a real version of
the Fulton-MacPherson compactification $X[n]$. 
An advantage of our approach is that we do not need
any machinery from algebraic geometry neither in constructing
our compactification space $\XAH$ nor in proving that it is a smooth
manifold. Instead of algebraic blowups we have used 
the angles between two points and the hooks between three points. 
This idea was already proposed by Kontsevich and Soibelman. 
We have shown  how to adapt their approach so that the
resulting modified compactification space is a smooth manifold instead of
just a manifold with corners. Moreover, we have shown that an
explicit analysis of the ratios on triples of points gives 
a combinatorial counterpart of Fulton-MacPherson's 
description of degenerated configurations in terms of screens.

Because of its nature, the construction and description
of the space $\XAH$ can serve as a illustrative example
to the theory of configuration spaces that has attracted a lot
of interest in recent years. The setup in terms of screens can serve
as a general framework for studying properties of point sets in the plane.
Especially for point sets that contain almost coinciding points,
a description by angles and hooks can be a more robust alternative
for storing the relative positions of the points. Another advantage
of working with hooks and angles is that translations
and scalings are already eliminated. That is, one concentrates
on point configuration up to affine transformations.

\cite{FM}, Theorem~2 states that the Fulton-MacPherson compactification $X[n]$
equals the closure of $\CONF$ in the product of those 
$\text{Bl}_{\Delta}(X^S)$ for $S\subset\{1, \dots, n\}$ of cardinality
2 and 3. Note that our model also supports this `three is enough' motto:
we only consider properties of pairs of points (the angles) and
triples of points (the hooks).

\subsection{Relation to earlier chapters.}

A first application is given by studying limits of Voronoi diagrams: 
we have introduced a method to associate a collection of filled $x$-screens
to a point $x \in \xah$. In the final sections it turned out that
this method  can be ported to the manifold with corners $\FMt(n)$,
introduced by Kontsevich-Soibelman.
Just because of its corners,
the screen structure associated to a point $x$  in $\FMt(n)$
is well suited 
for modeling  a possibly degenerated Voronoi diagram for $x$. 
Recall the description of Voronoi diagrams of polynomial sites,
presented in Chapter \ref{chlimit}: we have by now two notions of
degenerated Voronoi diagram that match both combinatorially and
according to shape with the notion of classic Voronoi diagram.

From a set $S(t)$ of $n$ distinct sites 
$\{p_1(t), \dots, p_n(t)\}$ in $\R[t] \times \R[t]$ we can 
construct a point $s \in \XAH$ at $t=0$, compare Example \ref{exmoreplug}. 
A construction in the other direction can be made
as well. Suppose we are given some $x \in \XAH$. And suppose 
the $x$-screens are filled by $x$-sites $p_1(x), \dots, p_n(x)$.
Define the \bfindex{depth} of a screen $S$ as:
\beq
	\text{depth}(S) & = & 1+\# \text{~screens above} ~S.
\eeq
Let $p_j$ be the site directly above $p_i$ with respect to the
predecessor relation defined  in Definition \ref{defhookedpath}.
The polynomial site $p_i(t)$ representing $p_i(x)$ is 
defined as:
\beq
	p_1(t) & = & (0,0),\\
	p_i(t) & = & p_j(t) + t^{\text{depth}(S_{ij}) }
		~(\,\text{coordinates}(\,p_i(x)\,) ~\text{in}~S_{ij}).
\eeq
Mapping  a set $S(t)$ of polynomial sites to a point $s \in \XAH$ and back results
in this way in a normal form for $S(t)$. 

In Chapters \ref{changles} and \ref{chcont} we have studied
the compactification $\cda$. This compactification
is the closure of the graph of the directed angle map, 
cf.\ Definition \ref{dcda}.
One of the drawbacks of considering only directed angles 
is that collinear configurations can not be reconstructed
from just these angles, compare Section \ref{secreconstruct}.
This problem is solved in the compactification
$\XAH$ by adding the hooks between triples of points. 

In Chapter \ref{chcont}, Theorem \ref{hausdorffcont} we prove 
essentially that
two data sets $\gamma_n,\eta_n \in \cda$ that are Euclidean close have
Voronoi diagrams $V(\gamma_n), V(\delta_n)$ that are Hausdorff close. 
As there is an obvious continuous  map (forget all ratios) from
$\FMt(n)$ to $\cda$, this continuity result also holds for the top
screen of the Voronoi diagram $V_{FM}(x)$.
More precisely: suppose we have some $x \in \FMt(n)$. Fill the
top screen $S_{[n]}$ by $x$. The coordinates of the sites
$p_1(x), \dots, p_n(x)$  can be read off from $S_{[n]}$, 
while the directed angles between pairs of
sites  are part of $x$. These data together determine a point
$\gamma_n \in \cda$.   This establishes a map $f: \FMt \rightarrow \cda$.
We have the following situation:
\begin{displaymath}
\begin{array}{lllcl}
	\FMt(n) & \stackrel{f}{\rightarrow} &  \cda & \
\stackrel{\text{Voronoi map}}{\rightarrow} & \text{diagram space},\\
	x       & \mapsto     & \gamma_n & \mapsto & V(\gamma_n).
\end{array}
\end{displaymath}
As both $f$ and the `Voronoi map' are continuous, the composition
is continuous as well. In this short analysis, we have left out the
incorporation of the compactness condition on the Voronoi map. 
This is left as an exercise to the reader.

\newpage
\phantom{kip}
\thispagestyle{empty}
\clearpage
\addcontentsline{toc}{chapter}{\numberline{}Index.}
\printindex

\end{document}